\documentclass[12pt]{article}
\usepackage{mathtools}
\usepackage{microtype}
\usepackage[usenames,dvipsnames,svgnames,table]{xcolor}
\usepackage{amsfonts,amsmath, amssymb,amsthm,amscd}
\usepackage[height=9in,width=6.5in]{geometry}
\usepackage{verbatim}
\usepackage{tikz}
\usetikzlibrary{decorations.markings}
\tikzset{->-/.style={decoration={
			markings,
			mark=at position .6 with {\arrow{>}}},postaction={decorate}}}
\usepackage{tkz-euclide}

\usepackage{latexsym}
\usepackage[ psdextra,pdfauthor={Tang},bookmarksnumbered,hyperfigures,colorlinks=true,citecolor=BrickRed,linkcolor=BrickRed,urlcolor=BrickRed,pdfstartview=FitH]{hyperref}
\usepackage{mathrsfs}
\usepackage{graphicx}

\usepackage{enumitem}

\newtheorem{theorem}{Theorem}[section]
\newtheorem{corollary}[theorem]{Corollary}
\newtheorem{lemma}[theorem]{Lemma}
\newtheorem{proposition}[theorem]{Proposition}
\newtheorem{question}[theorem]{Question}
\newtheorem{conj}[theorem]{Conjecture}
\newtheorem{observation}[theorem]{Observation}
\newtheorem{claim}[theorem]{Claim}

\newtheorem{definition}[theorem]{Definition}
\newtheorem{remark}[theorem]{Remark}
\numberwithin{equation}{section}

\newcommand{\be}{\begin{equation}}
	\newcommand{\ee}{\end{equation}}
	
\newcommand*\circled[1]{\tikz[baseline=(char.base)]{
		\node[shape=circle,draw,inner sep=0.5pt] (char) {#1};}}	
	
\def\ba{\begin{align}}
	\def\ea{\end{align}}

\def\as {\mathrm{a.s.}}

\newcommand{\notion}[1]{{\bf  \textit{#1}}}

\DeclareMathOperator{\UST}{\mathsf{UST}}

\begin {document}
\title{Weights of uniform spanning forests\\ on nonunimodular transitive graphs}
\author{
	Pengfei Tang\thanks{Department of Mathematics, Indiana
		University. Partially supported by the National
		Science Foundation under grants DMS-1612363.
		Email: \protect\url{tangp@iu.edu}.}
}
\date{\today}
\maketitle
\begin{abstract}
Considering the wired uniform spanning forest on a nonunimodular transitive graph, we show that almost surely each tree of the wired uniform spanning forest is light. More generally we study the tilted volumes for the trees in the wired uniform spanning forest.

Regarding the free uniform spanning forest, we consider several families of nonunimodular transitive graphs. We show that the free uniform spanning forest is the same as the wired one on  Diestel--Leader graphs. 
For grandparent graphs, we show that the free uniform spanning forest is connected and has branching number bigger than one. We also show that each tree of the free uniform spanning forest is heavy and has branching number bigger than one on  a free product of a nonunimodular transitive graph with one edge when the free uniform spanning forest is not the same as the wired. 
\end{abstract}

\section{Introduction}\label{sec:intro}
 Let $G=(V(G),E(G))$ be a locally finite, connected infinite graph. The wired uniform spanning forest (WUSF) and the free uniform spanning forest (FUSF) are weak limits of the uniform spanning tree measures on an exhaustion of the graph $G$, with wired and free boundary conditions respectively. The WUSF and FUSF can be disconnected but each component is an infinite tree. Pemantle \cite{Pem1991} proved for $\mathbb{Z}^d$, the WUSF is  the same as FUSF and he also showed that the WUSF(FUSF) on $\mathbb{Z}^d$ is connected iff $d\leq 4$. Since then  WUSF and FUSF have been extensively studied. Hutchcroft studied many geometric properties of WUSF in \cite{Hutchcroft2018a}. For example Hutchcroft gave the volume growth dimension, the spectral dimension and the walk dimension for the trees in WUSF. His results are about a class of high-dimensional graphs and this class of graphs contains nonunimodular transitive graphs. In the present paper, we continue the study of the geometries of trees in WUSF and FUSF on a special family of graphs---nonunimodular transitive graphs. We are particularly interested in the geometry of trees in WUSF and FUSF with respect to the ``level structure" of the underlying nonunimodular transitive graphs. For more background on uniform spanning forests see \cite{BLPS2001} or Chapter 4 and 10 of \cite{LP2016}.
 
 Let $\textnormal{Aut}(G)$ denote the automorphism group of $G$. Suppose $\textnormal{Aut}(G)$ has a nonunimodular closed subgroup $\Gamma\subset \textnormal{Aut}(G)$ that acts transitively on $G$. In particular such $G$ must be nonamenable \cite[Proposition 8.14]{LP2016} and hence there are infinitely many components for WUSF on $G$ \cite[Corollary 10.27]{LP2016}. There is a unique left Haar measure $|\cdot|$ on $\Gamma$ (up to a multiplicative constant). For each $x\in V(G)$, let $\Gamma_x:=\{\gamma\in\Gamma:\gamma x=x\}$ denote the stabilizer of $x$ and $m(x):=|\Gamma_x|$. We call this function $m: V\rightarrow(0,\infty)$ the \notion{weight function} for $(G,\Gamma)$. For a cluster $C$, we define its weight $m(C):=\sum_{x\in C}{m(x)}$ and call $C$ a \notion{$\Gamma$-light} cluster or a \notion{$\Gamma$-heavy} cluster according to $m(C)<\infty$ or $m(C)=\infty$ respectively. For simplicity we will just say $C$ is a light (or heavy) cluster if $\Gamma$ is well-understood from the context.

 Although each component of the WUSF or FUSF on $G$ is an infinite tree, it may happen that some tree is light and has branching number bigger than one. We prove that  if there is nonunimodular subgroup $\Gamma\subset \textnormal{Aut}(G)$ that acts transitively on $G$, then each tree of the WUSF on $G$ is $\Gamma$-light $\as$

 \begin{theorem}\label{thm:light for WUSF}
 	Let $\textnormal{Aut}(G)$ denote the automorphism group of $G$. Suppose $\textnormal{Aut}(G)$ has a nonunimodular closed subgroup $\Gamma\subset \textnormal{Aut}(G)$ that acts transitively on $G$. Then each tree of the $\textnormal{WUSF}$ is $\Gamma$-light almost surely. 
 \end{theorem}
 
 From now on if $G$ is a transitive graph and $\Gamma\subset \textnormal{Aut}(G)$ is a closed nonunimodular subgroup that acts transitively on $G$, then we call $(G,\Gamma)$ is a \notion{nonunimodular transitive pair}.

 Benjamini, Lyons,  Peres and  Schramm gave several equivalent conditions for $\textnormal{WUSF}=\textnormal{FUSF}$; see \cite[Theorem 7.3]{BLPS2001}. In particular, $\textnormal{WUSF}=\textnormal{FUSF}$ for amenable transitive graphs. For the case $\textnormal{WUSF}\neq\textnormal{FUSF}$, we make the following conjecture.
 
 \begin{conj}\label{conj: 1.2}
 	Suppose $\textnormal{Aut}(G)$ has a nonunimodular closed subgroup $\Gamma\subset \textnormal{Aut}(G)$ that acts transitively on $G$ and $\textnormal{WUSF}\neq\textnormal{FUSF}$ on $G$. Then each tree in the $\textnormal{FUSF}$ is $\Gamma$-heavy and has branching number bigger than one almost surely.
 \end{conj}

However we are only able to prove this for certain examples.  

 \begin{theorem}\label{thm:heavy for FUSF}
 	Let $G$ be a free product of a nonunimodular transitive graph with one edge such that $\textnormal{WUSF}\neq\textnormal{FUSF}$ on $G$ or a grandparent graph. Then each connected tree of the $\textnormal{FUSF}$ on $G$ is heavy and has branching number bigger than one almost surely. In fact, when $G$ is a grandparent graph, the FUSF on $G$ is connected almost surely.
 \end{theorem}
\begin{remark}
	Recently Pete and Tim\'{a}r \cite{PeteTimar2020} disproved the heaviness part of   Conjecture \ref{conj: 1.2}. So it is now natural to ask the question---``which transitive nonunimodular graphs have the property that every tree in the FUSF is heavy?" Theorem \ref{thm:heavy for FUSF},  Remark \ref{rem: connected FUSF on tree cross finite graph} and Pete and Tim\'{a}r's  Theorem 1.1 and 1.3 in \cite{PeteTimar2020}  suggest this question is a quantitative issue (see Section 6 of \cite{PeteTimar2020} for some related problems). 
%

A more fundamental problem than the branching number part of Conjecture \ref{conj: 1.2} is to show that every tree in FUSF has infinitely many ends on nonunimodular transitive graphs  when $\textnormal{FUSF}\neq\textnormal{WUSF}$. This is the remaining case of Question 15.8 of \cite{BLPS2001}; the unimodular case has been answered positively in \cite{TomAsaf2017,Timar2018indis}. 
\end{remark}

 The paper is organized as follows. In Section \ref{sec: 2} we give some preliminaries---we review the definition of WUSF and FUSF, Wilson's algorithm, nonunimodular transitive graphs and the tilted mass-transport principle. In particular, we review the ``level structure" for a nonunimodular transitive pair $(G,\Gamma)$. The level structure is induced by the weight function $m:V\rightarrow (0,\infty)$. Let $u\sim v$ denote that $u$ and $v$ are neighboring vertices in $G$. Let $L_n(x)$ be the set of vertices $y$ such that $\frac{m(y)}{m(x)}\approx e^{nt_0}$, where $e^{t_0}:=\max\{\frac{m(v)}{m(u)}\colon u\sim v \}$. See Subsection \ref{sec: subsection 2.3} for the precise definition of $L_n(x)$.  We also review the so-called \notion{tilted volumes} introduced by  Hutchcroft \cite{Tom2017}. For a set  $C\subset V(G)$ and $x\in C$, the so-called tilted volumes of $C$ are defined as: $|C|_{x,\lambda}:=\sum_{y\in C}{\frac{m(y)^\lambda}{m(x)^{\lambda}}}$, where $\lambda\in \mathbb{R}$. The tilted volume can be viewed as a generalization of the size $|C|$ and weight $m(C)$: when $\lambda=0$, $|C|_{x,\lambda}$ is just the size of $C$; when $\lambda=1$, $|C|_{x,\lambda}$ is just a normalization of $m(C)$.

  In Section \ref{sec: 3} we consider WUSF on a toy model $(\mathbb{T}_{b+1},\Gamma_{\xi})$, namely regular tree together with a subgroup of automorphisms that fixes an end. 
  The WUSF on the toy model has close relations to critical percolation and Galton--Watson trees (Lemma \ref{lem: law of T_x in WUSF for the toy model}). These relations and the tree structure of the underlying graph $\mathbb{T}_{b+1}$ make the study of  WUSF on the toy model relatively easy. Moreover the results for WUSF on the toy model shed light on further studies for WUSF on general nonunimodular transitive graphs (Proposition \ref{prop: aymptotic for limsup of toy model}, Question \ref{ques: asymotopic behavior of intersection with a lower level for general case}, Proposition \ref{prop: stretched exponential decay for x-wusf intersect with L_n(x)} and Proposition \ref{prop: sharp high moments for T_x intersect L_n(x) on the toy model}).

  In Section \ref{sec: 4} we study the geometry of the trees in WUSF with respect to the level structure of the nonunimodular transitive pair. 
  We will consider quantities such as $\mathbb{E}[|T_x\cap L_n(x)|]$ and $\mathbb{P}[T_x\cap L_n(x)\neq \emptyset]$. Similar quantities have been studied for Bernoulli percolation in \cite{Tom2017}. These quantities are not only interesting in themselves but also pave the way for proving Theorem \ref{thm:light for WUSF}. We will study when the tilted volumes of the trees in WUSF are finite (Proposition \ref{prop: when the tilted volume is finite}) and what the tail behaviors of the tilted volumes are (Proposition \ref{prop: tail prob for general case}). In particular, Theorem \ref{thm:light for WUSF} is a special case of Proposition \ref{prop: when the tilted volume is finite}.
  Some other objects like the past and the future of a vertex will also be studied in Section \ref{sec: 4}.
  
  In Section \ref{sec: 5} we consider FUSF on  Diestel--Leader graphs and grandparent graphs. In Proposition \ref{prop: FUSF=WUSF on DL graphs} we show that FUSF on  a Diestel--Leader graph is the same as WUSF. For the FUSF on a grandparent graph, we show it is just one tree (Proposition \ref{prop: FUSF on grandparent graph is a tree}) and has branching number bigger than one (Proposition \ref{prop: branching number larger than one for FUSF on grandparent graph}).

  In Section \ref{sec: 6} we consider FUSF on free products of  nonunimodular transitive graphs with  $\mathbb{Z}_2$ and prove Proposition \ref{prop: free product percolation}. Theorem \ref{thm:heavy for FUSF} then follows from  Proposition \ref{prop: FUSF on grandparent graph is a tree}, \ref{prop: branching number larger than one for FUSF on grandparent graph} and \ref{prop: free product percolation}. A key ingredient of proving Proposition \ref{prop: free product percolation} is to compare the weight of a tree in the FUSF to a branching random walk and then use Biggins' theorem \cite{Lyons1997}. 
 
 In Section \ref{sec: appendix A} of the appendix we give the details of the proof the lower bounds of \eqref{eq: prob of srw intersecting a high slab with at least k vertices} and \eqref{eq: prob of srw intersecting a low slab with at least k vertices} (The lower bounds are not used in the main part).  In Section \ref{sec: appendix B} of the appendix we give the proof of Proposition \ref{prop: high moments for x-component intersecting a slab}. We also provide more quantitative results for WUSF on the toy model in Section \ref{sec: section C in the appendix} of the appendix.

 \section{Preliminaries}\label{sec: 2}
 
 \subsection{Uniform spanning forests}

 If $G$ is a finite connected graph, then it has finitely many spanning trees. The \notion{uniform spanning tree} (UST) on $G$ is the uniform measure on the set of spanning trees of $G$ and denoted by $\UST(G)$. The Aldous-Broder algorithm and Wilson's algorithm are well-known methods to generate the UST on a finite graph $G$.
 
  Suppose $G=(V,E)$ is a locally finite, connected infinite graph. An exhaustion of $G$ is a sequence of finite connected subgraphs $G_n=(V_n,E_n)$ of $G$ such that  $V_n\subset V_{n+1}$ and $G=\bigcup G_n$. The weak limit of $\UST(G_n)$ exists and is independent of the choice of the exhaustion; for example see \cite[Section 10.1]{LP2016}.  We call the weak limit  of $\UST(G_n)$ \notion{free uniform spanning forest} (FUSF). From the graph $G$, first identify the vertices outside $G_n$ to a single vertex, say $z_n$, and then remove loop-edges at $z_n$ but keep multiple edges. The graph obtained in this way is denoted by $G_n^W$. We now also assume that $G_n$ is the graph \textbf{induced}
 in $G$ by $V_n$. Then the weak limit of $\UST(G_n^W)$ exists and is independent of the choice of such exhaustion; for example see \cite[Section 10.1]{LP2016}. We call the weak limit of $\UST(G_n^W)$  \notion{wired uniform spanning forest} (WUSF). 
 
 J\'{a}rai and Redig \cite{JR2008} also introduced a $v$-WUSF on $G$, which can be roughly understood as a wired spanning forest with $v$ wired to $\infty$. Suppose $v\in V(G_n)$ for every $G_n$ in the above exhaustion and let $\widehat{G}_n$ be the graph obtained from $G_n^W$ by identifying $v$ and $z_n$. Then the $v$-WUSF on $G$ is the weak limit of $\UST(\widehat{G}_n)$. 
 
 The connected component of $v$ in the $v$-WUSF is finite almost surely if $G$ is a transient transitive graph \cite[Proposition 3.1 and Theorem 7.4]{LMS2008}.

\subsection{Wilson's algorithm}
 Next we review Wilson's algorithm for generating WUSF on a transient graph $G$. This is called Wilson's method rooted at infinity  \cite[Theorem 5.1]{BLPS2001}.  For finite graphs or recurrent graphs see section 3 of \cite{BLPS2001} for more details.
 
 For a path $w$ on $G$ that visits every vertex finitely many times, we can define a loop erasure of $w$, namely let $\textnormal{LE}(w)$ be the self-avoiding path obtained by erasing the loops chronologically as they are created.

 Suppose $G$ is a transient graph. Fix an arbitrary ordering $(v_1,v_2,\ldots)$ of the vertices of $G$. Set $\mathsf{F}_0=\emptyset$. Given $\mathsf{F}_{n-1}$ for $n\geq 1$, we construct $\mathsf{F}_n$ as follows. If $v_n\in \mathsf{F}_{n-1}$, then let $\mathsf{F}_n=\mathsf{F}_{n-1}$. Otherwise start a simple random walk from $v_n$ and  let $\tau_n$ denote the first hitting time of $\mathsf{F}_{n-1}$. In particular, $\tau_n=\infty$ if the simple random walk never hits $\mathsf{F}_{n-1}$. Let $P_n$ denote the random walk path stopped at $\tau_n$. Loop erase this random walk path and denote it by $\textnormal{LE}(P_n)$. Set $\mathsf{F}_n=\mathsf{F}_{n-1}\cup \textnormal{LE}(P_n)$.
 Finally set $\mathsf{F}:=\bigcup_{n}\mathsf{F}_n$. Then $\mathsf{F}$ has the law of WUSF on $G$, and in particular, its law does not depend on the ordering one chose.

Let $G$ be a transient network and let $\mathfrak{F}$ be a sample of WUSF on $G$ generated by using Wilson's algorithm. Then for every edge $e$ of $\mathfrak{F}$, there is a unique orientation such that $e$ is crossed by the loop-erased random walk in that direction. The resulting oriented graph is called \notion{oriented wired uniform spanning forest} and denoted by OWUSF. The OWUSF also does not depend on the ordering of the vertices. From OWUSF one can get WUSF by forgetting the orientation. It is also easy to see that in the OWUSF, every vertex has exactly one edge emanating from it.

 Suppose $G$ is a transient transitive graph. Then every tree in the WUSF on $G$ will have only one end almost surely \cite[Theorem 7.4]{LMS2008}. Let $\vec{\mathfrak{F}}$ be a sample of OWUSF on $G$ and let $\mathfrak{F}$ be the spanning forest obtained from $\vec{\mathfrak{F}}$ by forgetting its orientation. Then  $\mathfrak{F}$ has the law of WUSF. For each vertex $x\in V(G)$, let $T_x$ denote the connected component of $x$ in $\mathfrak{F}$. Since $T_x$ is one-ended almost surely, there is a unique infinite ray $\eta=(v_0,v_1,\ldots)$ starting from $v_0=x$ representing the unique end of $T_x$. Then the edge  $(v_0,v_1)$ is also the unique edge emanating from $x$ in $\vec{\mathfrak{F}}$. 
 
 Given a sample $\mathfrak{F}$ of WUSF on $G$ and $u\in V(G)$,  we define \notion{the future} of $u$ to be the unique oriented ray starting from $u$ and denote it by $\mathfrak{F}(u,\infty)$. We take the convention that $u\in\mathfrak{F}(u,\infty)$. We also define \notion{the past} of $u$ to the subgraph of $\mathfrak{F}$ spanned by those vertices  $v\in\mathfrak{F}$ such that $u\in\mathfrak{F}(v,\infty)$ and denote it by $\mathfrak{P}(u)$.

 Let $\mathfrak{F}_v$ be a sample of $v$-WUSF on $G$. Then one can generate $\mathfrak{F}_v$ by running Wilson's algorithm rooted at infinity but starting with $\mathsf{F}_0=\{v\}$, i.e. the forest with a single vertex $v$ and no edge. We can also orient $e\in\mathfrak{F}_v$ as the way it is crossed by the loop-erased random walk in the Wilson's algorithm. Then each vertex but $v$ has exactly one edge emanating from it. We can define the past and future for every vertex $u$ in $\mathfrak{F}_v$ according to this orientation and denote them by $\mathfrak{P}_v(u)$ and $\mathfrak{F}_v(u,\infty)$ respectively. In particular, the future of $v$ in $\mathfrak{F}_v$ is the single vertex $v$ itself and the past of $v$ is the connected component of $v$ in $\mathfrak{F}_v$. Let $\mathfrak{T}_v$ denote the tree containing $v$ in $\mathfrak{F}_v$. Most notation here coincides with the ones listed on page 15 of  \cite{Hutchcroft2018a} for the reader's convenience. 
 
  Given a general oriented forest $F$ of $G$, if there is an oriented path from $u$ to $v$ in $F$, then $u$ is said to be in the past of $v$ and $v$ is said to be in the future of $u$. Let $\textnormal{past}_F(v)$ denote the past of $v$ in the oriented forest $F$, namely, the set of vertices which lie in the past of $v$ in $F$. One lemma we shall need is the following stochastic domination result.
 
 \begin{lemma}[Lemma 2.1 of \cite{Hutchcroft2018a}]\label{lem:stochastic domination}
 	Let $G$ be an infinite network and $\mathfrak{F}$ be a sample of OWUSF on $G$. For each $v\in V(G)$, let $\mathfrak{F}_v$ be an oriented $v$-WUSF of $G$. Suppose $K$ is a finite set of vertices of $G$ and define $F(K):=\bigcup_{u\in K}\mathfrak{F}(u,\infty)$ and $F_v(K):=\bigcup_{u\in K}\mathfrak{F}_v(u,\infty)$. Then for every $u\in K$ and every increasing event $\mathscr{A}\subset \{0,1\}^E$ we have 
 	\be\label{eq:SD 1}
 	\mathbb{P}\left( \textnormal{past}_{\mathfrak{F}\backslash F(K)}(u)\in\mathscr{A}\mid F(K)\right)
 	\leq \mathbb{P}(\mathfrak{T}_u\in\mathscr{A})
 	\ee
 	and similarly
 	\be\label{eq:SD 2}
 	\mathbb{P}\left( \textnormal{past}_{\mathfrak{F}_v\backslash F_v(K)}(u)\in\mathscr{A}\mid F(K)\right)
 	\leq \mathbb{P}(\mathfrak{T}_u\in\mathscr{A}).
 	\ee
 \end{lemma}
 This lemma is stated for OWUSF and we will often use it for WUSF because we can first sample an OWUSF and then get WUSF by forgetting the orientation of the OWUSF.

 \subsection{Nonunimodular transitive graphs and the tilted mass-transport principle}\label{sec: subsection 2.3}
 Next we recall the tilted mass-transport principle. We will restrict to transitive graphs.
 Suppose $G$ is an infinite, locally finite connected graph and $\Gamma\subset\textnormal{Aut}(G)$ is a subgroup of automorphisms that acts transitively on $G$. For $x,y\in V(G)$, let $|\Gamma_xy|$ denote the number of vertices in the set $\{\gamma x:\gamma\in \Gamma_x \}$, where $\Gamma_x$ is the stabilizer of $x$. There is a unique (up to a multiplicative constant) nonzero left Haar measure $|\cdot|$ on $\Gamma$ and we denote by $m(x)=|\Gamma_x|$ the Haar measure of the stabilizer $\Gamma_x$. Then a simple criterion for unimodularity of $\Gamma$ is given as follows. 
\begin{proposition}[Trofimov \cite{Trofimov1985}]\label{prop:unicre}
	Suppose $\Gamma\subset \textnormal{Aut}(G)$ acts transitively on $G$. Then $\Gamma$ is unimodular if and only if for all $x,y\in V(G
	)$,  
	\[
	|\Gamma_xy|=|\Gamma_yx|.
	\]
\end{proposition}

The grandparent graph and the Diestel--Leader graph $DL(q,r)$ with $q\neq r$ are typical examples of nonunimodular transitive graphs. For more examples see section 3 of \cite{Timar2006}.

The following lemma is well known (for a proof, see for example formula (1.28) and Lemma 1.29 in \cite{Woess2000} ):
\begin{lemma}\label{lem:haar}
	Suppose $\Gamma\subset \textnormal{Aut}(G)$ acts on $G=(V,E)$ transitively, then for all $x,y\in V$
	\[
	\frac{m(x)}{m(y)}=\frac{|\Gamma_xy|}{|\Gamma_yx|}=\frac{m(\gamma x)}{m(\gamma y)},\forall \gamma\in \Gamma.
	\]
	
\end{lemma}

 Given the group $\Gamma$, define the \notion{modular function} $\Delta:V\times V\in[0,\infty]$ as follows:
 \[
 \Delta(x,y):=\frac{m(y)}{m(x)}.
 \]
 
 Then from the above lemma \ref{lem:haar} we know $\Delta$ is a $\Gamma$ diagonally invariant function, i.e. $\Delta(x,y)=\Delta(\gamma x,\gamma y), \forall \gamma\in\Gamma$. Another important property for the modular function is the \notion{cocycle identity}:
 \[
 \Delta(u,v)\Delta(v,w)=\Delta(u,w),\ \forall u,v,w\in V(G).
 \]

 For more background on modular functions defined here see Section 2.1 of \cite{Tom2017}. 
 
 For a cluster $K$ of $G$, a vertex $v\in K$ and a parameter $\lambda\in\mathbb{R}$,  Hutchcroft \cite{Tom2017} introduced the \notion{tilted volume} as follows:
 \[
 |K|_{v,\lambda}:=\sum_{y\in K}\Delta(v,y)^\lambda.
 \]
 
 The tilted mass-transport principle  is a useful technique when dealing with nonunimodular transitive graphs.  Actually the word `tilted' can be omitted, and  the mass-transport principle was defined without the word \cite{BLPS1999}. The `tilted mass-transport principle' (TMTP) was first used in \cite{Tom2017} for a different way of writing the mass-transport principle. 
 
 \begin{proposition}[Proposition 2.2  of \cite{Tom2017}]\label{prop:TMTP}
 	With the same notations as in the above definition of modular function, suppose $F:V^2\rightarrow[0,\infty]$ is invariant under the diagonal action of $\Gamma$.  Then 
 	\[
 \sum_{x\in V}{F(\rho,x)}=\sum_{x\in V}{F(x,\rho)\Delta(\rho,x)}.
 	\]
 \end{proposition}

 We will often use the following form of the above tilted mass-transport principle. Suppose $\omega$ is a $\Gamma$-invariant bond percolation process on $G=(V,E)$. Suppose $f:V^2\times\{0,1\}^E\rightarrow [0,\infty]$ is $\Gamma$-invariant, that is,
 \[
 f(x,y,\omega)=f(\gamma x,\gamma y,\gamma \omega), \, \forall\, x,y\in V, \gamma\in \Gamma. 
 \]
 Then applying the tilted mass-transport principle with $F(x,y)=\mathbb{E}[f(x,y,\omega)]$ one has
 \begin{equation}\label{eq:TMTP with percolation}
 \mathbb{E}\left[\sum_{x\in V}{f(\rho,x,\omega)}\right]=\mathbb{E}\left[\sum_{x\in V}{f(x,\rho,\omega)\Delta(\rho,x)}\right].
 \end{equation}
 When the percolation $\omega$ is clear from the context, we will often write $f(x,y)$ instead of $f(x,y,\omega)$ in the above equation.
 
 Next we recall some terminology from percolation theory. Suppose $G=(V,E)$ is a locally finite, connected graph. Let $2^E=\{0,1\}^E$ be the collection of all subsets $\eta\subset E$ and let $\mathcal{F}_E$ be the $\sigma$-field generated by sets of the form $\{\eta:e\in\eta\}$ where $e$ runs over all edges in $E$.  A bond percolation on $G$ is a pair $(\mathbf{P},\omega)$, where $\omega$ is a random element in $2^E$ and $\mathbf{P}$ is the law of $\omega$. For simplicity sometimes we will just say $\omega$ is a bond percolation. The interested reader can refer to \cite[Chapter 7 and 8]{LP2016} for more background on percolation theory. If the law $\mathbf{P}$ is invariant under a subgroup $\Gamma$ of automorphisms, then we call $(\mathbf{P},\omega)$ a $\Gamma$-invariant percolation on $G$. In particular, WUSF and FUSF can be viewed as $\textnormal{Aut}(G)$-invariant percolation processes on an infinite graph $G$.

 Suppose $(\mathbf{P},\omega)$ is an automorphism-invariant percolation on $G$ and $\mathbf{E}$ is the corresponding expectation operator. An edge $e$ is called open if $\omega(e)=1$; otherwise it is called closed. For $x\in V$, the \notion{cluster} of $x$ is the connected component of $x$ in the subgraph formed by open edges and denoted by $C_x$. In  case $\omega$ has the law of WUSF or FUSF, we will also denote the connected component of $x$ by $T_x$ since each connected component is a tree. 
 
 Note that the modular function satisfies $\Delta(y,x)=\Delta(x,y)^{-1}$. An immediate application of the tilted mass-transport principle is the following observation at the beginning of section 3 of \cite{Tom2017}: 
 \begin{lemma}\label{lem:convex and symmetry}
 	Suppose $G$ is an infinite, locally finite connected graph and $\Gamma$ is a subgroup of automorphisms that acts transitively on $G$. Suppose $(\mathbf{P},\omega)$ is a $\Gamma$-invariant percolation on $G$ and $x$ is an arbitrary vertex in $G$. Then the function $f(\lambda):=\mathbf{E}[|C_x|_{x,\lambda}]$ is symmetric about $\lambda=\frac{1}{2}$:
 	\[
 	\mathbf{E}[|C_x|_{x,\lambda}]=\mathbf{E}[|C_x|_{x,1-\lambda}].
 	\] 
 	
 	Note that $f$ is also convex, whence $f$ is decreasing on $(-\infty,\frac{1}{2}]$ and increasing on $[\frac{1}{2},\infty)$.
 \end{lemma}

 
 Finally we take a look at the level structure of nonunimodular transitive graphs; see also Section 2.2 of \cite{Tom2017}. 

 Suppose $G$ is an infinite, locally finite connected graph and $\Gamma$ is a nonunimodular subgroup of automorphisms that acts transitively on $G$. For $x,y\in V(G)$,  if $x,y$ are neighbors then we write $x\sim y $.
 
 Define \[
 t_0=t_0(G):=\max \{\log\Delta(x,y): x,y\in V(G), x\sim y \} 
 \]
 
 For each $s\leq t$ and $v\in V$ we define the \notion{slab}
 \[
 S_{s,t}(v):=\{ x\in V: s\leq \log \Delta(v,x) \leq t    \}.
 \]
 
 If $t\geq s+t_0$ and a simple path $\pi:=(v_0,\ldots,v_n)$ in $G$ crosses the slab $S_{s,t}(v)$ in the sense that 
 $\log \Delta(v,v_0)\leq s$ and $\log \Delta(v,v_n)\geq t$ then it must contain some vertex from the slab $S_{s,t}(v)$ by the definition of $t_0$. 
 
 Now we introduce a construction which will help us when applying the tilted mass-transport principle. This construction comes from \cite{Timar2006} and the setup here is borrowed from Section 2.2 of \cite{Tom2017}. 
 
 Fix an arbitrary vertex $v_0$ of $G$, let $U_{v_0}$ be a uniform $[0,1]$ random variable independent of the WUSF and $v$-WUSF we shall consider. For every other $v\in V(G)$, let 
 \[
 U_v:=U_{v_0}-\frac{1}{t_0}\log \Delta(v_0,v)\ \ \mod 1.
 \]
 Notice the law of the collection of random variables $U:=\{U_v: v\in V(G) \}$ does not depend on the choice of $v_0$.
 
 Given the collection of random variables $U$, the \notion{separating layers} are defined to be 
 \[
 L_n(v):=\{x\in V: (n+U_v-1)t_0\leq \log\Delta(v,x) \leq (n+U_v)t_0  \},n\in\mathbb{Z}.
 \]
 We also call $L_n(v)$ the $n$-th slab relative to $v$. We will use $\mathbb{E},\mathbb{P}$ to denote the expectation operator and probability measure for the joint law of $U$  and WUSF or $v$-WUSF on $G$. Note that the cocycle identity ($\Delta(x,v)\cdot\Delta(v,x)=1$) implies that  
 \be\label{eq: relation for n and -n levels}
 x\in L_n(v) \Leftrightarrow\  -n+(1-U_v)-1\leq \frac{1}{t_0}\log\Delta(x,v)\leq -n+1-U_v.
 \ee
 In particular since $1-U_v$ has the same law as $U_x$ (uniform on $[0,1]$) we have  
 \be\label{eq: probability relation for n and -n levels}
 \mathbb{P}[x\in L_n(v)]=\mathbb{P}[v\in L_{-n}(x)].
 \ee
 
 Suppose $f:V^2\times\{0,1\}^E\times[0,1]^V\rightarrow[0,\infty]$ is invariant under the diagonal action of $\Gamma$, then like \eqref{eq:TMTP with percolation}  we have the following form of tilted mass-transport principle (equation (2.2) in \cite{Tom2017})
 \be\label{eq: TMTP with lable on vertices}
 \mathbb{E}\left[\sum_{x\in V}f(\rho,x,\mathfrak{F},U)\right]=\mathbb{E}\left[\sum_{x\in V}f(x,\rho,\mathfrak{F},U)\Delta(\rho,x)\right],
 \ee
 where $\mathfrak{F}$ is a sample of WUSF and $U:=\{U_v: v\in V(G) \}$ is defined as above.
 
 If we let $\{\rho\leftrightarrow x \}$ denote the event that $\rho$ and $x$ are in the same connected component in $\mathfrak{F}$ and set  $f(\rho,x,\mathfrak{F},U)=1_{\{ \rho\leftrightarrow x, x\in L_n(\rho)   \}}$, then 
 \begin{eqnarray}\label{eq: expectation relation between n and -n level}
 \mathbb{E}[|T_\rho\cap L_n(\rho)|]&\stackrel{ \eqref{eq: TMTP with lable on vertices}}{=}&\sum_{x\in V}{\mathbb{E}[1_{\{ x\leftrightarrow \rho, \rho\in L_n(x)   \}}\Delta(\rho,x) ] }
 \stackrel{\textnormal{independence}}{=}\sum_{x\in V}{\mathbb{E}[1_{\{ x\leftrightarrow \rho\}}]\mathbb{E}[1_{\{\rho\in L_n(x)   \}}]\Delta(\rho,x)  }\nonumber\\
&\stackrel{\eqref{eq: probability relation for n and -n levels}}{=}&
\sum_{x\in V}{\mathbb{E}[1_{\{ x\leftrightarrow \rho\}}]\mathbb{E}[1_{\{x\in L_{-n}(\rho)   \}}]\Delta(\rho,x)  }\nonumber\\
&=&\mathbb{E}[\sum_{x\in V}1_{\{ x\leftrightarrow \rho, x\in L_{-n}(\rho)   \}}\Delta(\rho,x) ]
\asymp \exp(-t_0n)\mathbb{E}[|T_\rho\cap L_{-n}(\rho)|],
 \end{eqnarray}
 where the last equality holds up to a factor of $e^{\pm t_0}$.
 
 The constants $c_i$ in the present papers are always positive constants. Their values only depend on the underlying graph $G$ and the subgroup $\Gamma$. The constants $c_i$ may carry different values at different appearances.  The notations we use in this paper are similar to the ones listed on page 15 of \cite{Hutchcroft2018a}. However since we use $\Gamma$ to denote subgroup of automorphisms, we will use slightly different notations for the future of a vertex and  the paths connecting two vertices in the WUSF and $v$-WUSF. We summarize the  notations here for the reader's convenience.
 
 \begin{itemize}[leftmargin=4cm]
 	\item  [$\mathfrak{F}, \mathfrak{F}_v$] A sample of the WUSF and $v$-WUSF respectively.
 	
 	\item  [$\mathfrak{F}^f$] A sample of the FUSF.
 	
 	\item  [$x\leftrightarrow y$] The events that $x$ and $y$ are in the same connected component of $\mathfrak{F}$.
 	
 	\item [$T_v,\mathfrak{T}_v$] The connected components of $v$ in the WUSF and $v$-WUSF respectively. 
 	
 	\item [$\mathfrak{F}(x,\infty),\mathfrak{F}_{v}(x,\infty)$]  The future of $x$ in $\mathfrak{F}$ and $\mathfrak{F}_v$ respectively.
 	
 	\item  [$\mathfrak{P}(x),\mathfrak{P}_v(x)$] The past of $x$ in $\mathfrak{F}$ and $\mathfrak{F}_v$ respectively.
 	
 	
 	\item [$d_G(x,y)$] The graph distance from $x$ to $y$. When the underlying graph $G$ is clear from the context, we often write it as $d(x,y)$. 
 	
 	\item [$\{X_k^x\}_{k\geq 0}$] A simple random walk on $G$ starting from $x$. For  $y\neq z$, we take $X^y,X^z$  to be independent. 
 	
 	\item [$\sigma_y^x$] The first visit time of $y$ by a simple random walk starting from $x$. 
 	
 	\item  [$\tau_{-n}^x$] The last visit time of $L_{-n}(x)$ by a simple random walk starting from $x$.

 \item [$\asymp $] This denotes an equality that holds up to positive multiplicative constants. More precisely, for two positive functions $f,g$ on $(0,\infty)$, $f(R)\asymp g(R)$ means that there exists $R_0>0$ and $c_1,c_2>0$ such that $c_1f(R)\leq g(R)\leq c_2 f(R)$ for all $R\geq R_0$, and the implicit constants $c_1,c_2, R_0$ only depend on the graph. Sometimes we also use $f(n)=\Theta(g(n))$ to denote that $f(n)\asymp g(n)$.
 
 \item [$\asymp_{\lambda} $] Similar to the above, but the implicit constants also depend on $\lambda$.
 
 \item  [$\preceq,\preceq_\lambda $ and $\succeq, \succeq_\lambda $] are defined similarly to $\asymp,\asymp_\lambda$. 
 \end{itemize}

 \section{A toy model}\label{sec: 3} 

 Let $G$ be a regular tree $\mathbb{T}_{b+1}$ of degree $b+1$, where $b\geq2$. Given an end $\xi$ of $G$, let $\Gamma_\xi$ be the subgroup of automorphisms that fixes this end $\xi$. One can check that $\Gamma_\xi$ acts transitively on $G$ and $\Gamma_\xi$ is nonunimodular by Proposition  \ref{prop:unicre}. The \notion{toy model} $(\mathbb{T}_{b+1}, \Gamma_{\xi})$ is just this regular tree $\mathbb{T}_{b+1}$ together with the subgroup $\Gamma_\xi$.
 
 The FUSF on $\mathbb{T}_{b+1}$ is trivial, namely it equals the tree itself almost surely. The WUSF on $\mathbb{T}_{b+1}$ can be generated using  Wilson's algorithm. 
The description of $T_x$ in the following lemma is due to H\"{a}ggstr\"{o}m (1998) \cite{Hag1998}.
\begin{lemma}\label{lem: law of T_x in WUSF for the toy model}
	
	Consider WUSF on the regular tree $\mathbb{T}_{b+1}$. Let  $x$ be a fixed vertex and  $T_x$  be the tree containing $x$ in the WUSF. 
	
	Pick a ray $(x_0,x_1,x_2,\ldots)$ starting from $x_0=x$ ``uniformly" in the sense that $\mathbb{P}[(x_0,\ldots,x_n)=\gamma]=\frac{1}{b^{n-1}(b+1)}$ for any self-avoiding path $\gamma$ started from $x_0$ with length $n$. 
	Then the tree $T_x$ has the law of the connected component of $x$ in the union of this ray with an independent Bernoulli bond percolation with parameter $\frac{1}{b}$.
	
\end{lemma}
\begin{proof}
	The proof given here comes from the second paragraph of Section 11 on \cite[page 42]{BLPS2001} and we include it for readers' convenience. 
	The WUSF on $\mathbb{T}_{b+1}$ can be generated using  Wilson's algorithm:  
	For a vertex $x\in \mathbb{T}_{b+1}$, start a simple random walk from $x$ on the tree $\mathbb{T}_{b+1}$, loop erase this random walk path chronologically to get a ray $\eta$ starting from $x$. Let  $y_1,\ldots,y_b$ denote the neighbors of $x$ not on the ray $\eta$. For $i=1,\ldots,b$ start independent simple random walks from the vertex $y_i$ and let $A_i:=\{\textnormal{ the simple random walk starting from } y_i \textnormal{ hits } x \}$.  Obviously given $\eta$, the events $A_i$ are independent. If $A_i$ occurs then put the edge connecting $y_i$ and $x$ to the WUSF; otherwise $y_i,x$ will be in different components in the WUSF. On the event $A_i$, we add only the edge $(y_i,x)$ to the WUSF and repeat the process for $y_i$. Obviously $\mathbb{P}(A_i|\eta)=\frac{1}{b}$. Thus the tree of $x$ is the union of the ray $\eta$, and independent random trees attached to the vertices of $\eta$. For the root $x$, the random tree attached is a critical Galton--Watson tree with binomial progeny distribution $\textnormal{Bin}(b,1/b)$. For each other vertex on the ray $\eta$, the first generation of the random tree has binomial progeny distribution $\textnormal{Bin}(b-1,1/b)$ while the subsequent generations have progeny distribution $\textnormal{Bin}(b,1/b)$. This analysis can be extended to give the whole WUSF easily.
	
	From the above description, it is easy to see the random trees attached also have the law of clusters in an independent Bernoulli bond percolation with $p=\frac{1}{b}$.
\end{proof}

Similarly one has the following description of $\mathfrak{T}_x$ for the toy model. 
\begin{lemma}\label{lem: law of T_x in  x-wusf for the toy model}
	Let  $x$ be a fixed vertex of $\mathbb{T}_{b+1}$. Consider $x$-WUSF on the regular tree $\mathbb{T}_{b+1}$ and the tree $\mathfrak{T}_x$ that contains $x$. Then $\mathfrak{T}_x$ has the law of the cluster $C_x$ in an independent Bernoulli bond percolation with $p=\frac{1}{b}$.
\end{lemma}
  
\subsection{Two point function and first moment for the toy model}
 
  The following is a simple application of the symmetry of $\mathbb{T}_{b+1}$ and Lemma \ref{lem: law of T_x in WUSF for the toy model}.
 \begin{proposition}\label{prop: exact two point function for WUSF on regular trees}
 	Let $\mathbb{P}$ denote the law of WUSF on the regular tree $\mathbb{T}_{b+1}$, where $b\geq2$. Let $x,y\in V(G)$ be two arbitrary vertices and $\{x\leftrightarrow y\}$
 	denote the event that $x,y$ are in the same tree of WUSF. Set $n=\textnormal{dist}(x,y)$ to be the graph distance of $x,y$. Then 
 	\be\label{eq: two point function for WUSF on the toy model}
 	\mathbb{P}[x\leftrightarrow y]=\frac{1}{b^n}[1+n\cdot\frac{b-1}{b+1}].
 	\ee
 	In particular, the probability that the tree $T_x$ intersects with a high slab $L_n(x)$ decays exponentially: 
 	\[
 	\mathbb{P}[T_x\cap L_n(x)\neq\emptyset]=\frac{1}{b^n}[1+n\cdot\frac{b-1}{b+1}]\asymp ne^{-t_0n},n\geq1.
 	\]
 \end{proposition}
 \begin{proof}
 		Let $\eta_x=(x_0,x_1,x_2,\ldots)$ be the ray starting from $x$ representing the end $\xi$. 
 	Since the WUSF is invariant under the whole automorphism group of $\mathbb{T}_{b+1}$, and for $k\geq 1$ there are $(b+1)b^{k-1}$ vertices on  $\mathbb{T}_{b+1}$ with graph distance $k$ to $x$,  by symmetry  $\mathbb{P}(x_k\in \mathfrak{F}(x,\infty))=\frac{1}{(b+1)b^{k-1}}$. 
 	
 		Let $u$ denote the highest vertex in the future of $x$, i.e., the vertex  $u$ is the unique one in $\mathfrak{F}(x,\infty)$  such that $\Delta(x,u)=\sup\{\Delta(x,v): v\in\mathfrak{F}(x,\infty) \}$. 
 	Let $A_k$ denote the event that $\Delta(x,u)=b^k$ for $k\geq0$.  Let $y$ denote the highest vertex in the past of $u$ and $B_{k,k'}$ be the event that $A_k\cap \{\Delta(u,y)=b^{k'} \}$ for $k'\geq0$. See Figure \ref{fig: B_k,k'} for a typical $B_{k,k'}$.
 	
 	Note that for $k\geq0$, the event $A_k$ defined above is just the event $\{x_k\in \mathfrak{F}(x,\infty) \}\backslash \{x_{k+1}\in \mathfrak{F}(x,\infty)  \}$. Hence for $k>0$, $\mathbb{P}(A_k)=\mathbb{P}(x_k\in \mathfrak{F}(x,\infty))-\mathbb{P}(x_{k+1}\in \mathfrak{F}(x,\infty))=\frac{b-1}{b+1}\cdot\frac{1}{b^k}$. Combining with the fact that $\mathbb{P}(x_0\in\mathfrak{F}(x,\infty))=1$, one has 
 	\[
 	\mathbb{P}[A_k]=\left\{
 	\begin{array}{cc}
 	\frac{b}{b+1}, & k=0\\
 	& \\
 	\frac{b-1}{b+1}\cdot\frac{1}{b^k}, & k>0.
 	\end{array}
 	\right.
 	\]
 	In particular, we see that the future of $x$, $\mathfrak{F}(x,\infty)=(x_0,\ldots,x_n,\ldots)$ will eventually go down to the lower slabs in the sense that $\log \Delta(x_0,x_n)\rightarrow-\infty$ as $n$ tends to infinity.  In particular, this implies that for $n\geq0$, \[
 	\mathbb{P}[T_x\cap L_{-n}(x)\neq \emptyset]=	\mathbb{P}[\mathfrak{F}(x,\infty)\cap L_{-n}(x)\neq \emptyset]=1.
 	\]
 	
 By Lemma \ref{lem: law of T_x in WUSF for the toy model}  one has 
 	\[
 	\mathbb{P}[B_{k,k'}]=\mathbb{P}[A_k]\cdot\frac{b-1}{b}\cdot\frac{1}{b^{k'}}, \ \forall k, k'\geq0 .
 	\]
 	The $n=0$ case of \eqref{eq: two point function for WUSF on the toy model} is trivial. We assume $n>0$ in the following. 
 	By symmetry we can assume $y$ is the unique vertex such that $\Delta(x,y)=b^n$ and $\textnormal{dist}(x,y)=n$. Hence 
 	\[
 	\mathbb{P}[x\leftrightarrow y]=\sum_{k=n}^{\infty}\mathbb{P}[A_k ]+\sum_{k=0}^{n-1}\sum_{k'=n-k}^{\infty}\mathbb{P}[B_{k,k'}].
 	\] 
 	Then simple calculation shows that $\mathbb{P}[x\leftrightarrow y]=\frac{1}{b^n}[1+n\cdot\frac{b-1}{b+1}]$.
 \end{proof}
 
 \begin{proposition}\label{prop: first moment for T_x intersect with a slab in the toy model}
 	Consider WUSF on the toy model $(\mathbb{T}_{b+1},\Gamma_{\xi})$. Fix an arbitrary vertex $x$ of  $\mathbb{T}_{b+1}$. 
 	Then 
 	\be\label{eq: expectation of intersection of T_x with a slab for the toy model} 
 	\mathbb{E}[|T_x\cap L_n(x)|]\asymp
 	\left\{
 	\begin{array}{cc}
 	(n\vee1)e^{-t_0n}=\frac{n\vee1}{b^n} & \textnormal{ if }n\geq 0\\
 |n| & \textnormal{ if }n<0
 	\end{array}
 	\right..
 	\ee
 \end{proposition}
\begin{proof}
	Let $d(u,v)$ be the graph distance of $u,v$ in $\mathbb{T}_{b+1}$.
	Fix $n\geq0$.
	Note that for an integer $k\geq 0$, \[
	|\{y\colon d(x,y)=n+2k, \Delta(x,y)=b^n \}|=
	\left\{
	\begin{array}{cc}
1 & \textnormal{ if }k=0\\
(b-1)b^{k-1} & \textnormal{ if }k\geq1
	\end{array}
	\right..
	\]
Then by Proposition \ref{prop: exact two point function for WUSF on regular trees} one has for $n\geq0$,
\begin{multline*}
	\mathbb{E}[|T_x\cap L_n(x)|]=
	\sum_{y\in L_n(x)}\mathbb{P}[x\leftrightarrow y]\\
=1\cdot\frac{1}{b^n}[1+n\cdot\frac{b-1}{b+1}]+\sum_{k=1}^{\infty}(b-1)b^{k-1}\cdot \frac{1}{b^{n+2k}}[1+(n+2k)\cdot\frac{b-1}{b+1}]
	\asymp\frac{n\vee1}{b^n}.
\end{multline*}	
	The $n<0$ case is similar, just using the following observation instead 
	\[
		|\{y\colon d(x,y)=|n|+2k, \Delta(x,y)=b^n \}|=
	\left\{
	\begin{array}{cc}
	b^{|n|} & \textnormal{ if }k=0\\
	(b-1)b^{|n|+k-1} & \textnormal{ if }k\geq1
	\end{array}
	\right..\qedhere
	\]
\end{proof}
 
  \begin{corollary}\label{cor: expectation for WUSF on regular tree}
 	For  the toy model $(\mathbb{T}_{b+1},\Gamma_{\xi})$, $\mathbb{E}[|T_x|_{x,\lambda}]<\infty
 	$ if and only if $\lambda\in (0,1)$. 
 \end{corollary}
\begin{proof}
	If we decompose $T_x$ according to its intersection with different slabs, we get 
	\be\label{eq:decomposion w.r.t. slabs for the toy model}
	\mathbb{E}[|T_x|_{x,\lambda}]\asymp\sum_{n\in\mathbb{Z}}\mathbb{E}[|T_x\cap L_{n}(x)|]\cdot\exp(t_0\lambda n).
	\ee
	
	Then \eqref{eq:decomposion w.r.t. slabs for the toy model} and \eqref{eq: expectation of intersection of T_x with a slab for the toy model} yields the conclusion. Also for $\lambda\in(0,1)$, by \eqref{eq:decomposion w.r.t. slabs for the toy model} one has 
	\be\label{eq: asymptotic behavior of expected tilted volume when lambda is close to 0 or 1}
	\mathbb{E}[|T_x|_{x,\lambda}]\asymp 1+\sum_{n=1}^{\infty}n\big(e^{-t_0\lambda n}+e^{-t_0(1-\lambda)n}\big)\asymp\frac{1}{\lambda^2(1-\lambda)^2}.\qedhere
	\ee
\end{proof}
 
\begin{proposition}\label{prop: lightness of trees in the toy model}
	For  the toy model $(\mathbb{T}_{b+1},\Gamma_{\xi})$, if $\lambda\leq 0$, then $|T_x|_{x,\lambda}=\infty \ \as$
	If $\lambda>0$, then $|T_x|_{x,\lambda}<\infty\  \as$ 
	In particular, the $\lambda=1$ case tells us that $T_x$ is light almost surely.
\end{proposition} 
\begin{proof}
		Let $\{X_n\}_{n\geq0}$ be a simple random walk started at $x$. If we sample $\mathfrak{F}$ starting with $\{X_n\}_{n\geq0}$ and let $\mathfrak{F}(x,\infty)=(x_0,x_1,x_2,\ldots)$ be future of $x$. Then from the proof of \ref{prop: exact two point function for WUSF on regular trees} one has that $\Delta(x_0,x_n)\rightarrow0$ almost surely ($\mathbb{P}(\bigcup_{k\geq 0}A_k)=1$). 
	Thus for $\lambda\leq0 $, $|T_x|_{x,\lambda}\geq |\mathfrak{F}(x,\infty)|_{x,\lambda}=\sum_{n=0}^{\infty}\Delta(x_0,x_n)^\lambda=\infty$.
	
	For $\lambda\in(0,1)$, Corollary  \ref{cor: expectation for WUSF on regular tree}
	then implies $|T_x|_{x,\lambda}<\infty$ almost surely.
	
	For $\lambda>0$, we can write 
	\[
	|T_x|_{x,\lambda}=\sum_{y\in T_x,\Delta(x,y)\geq1}{\Delta(x,y)^\lambda}+\sum_{y\in T_x,\Delta(x,y)<1}{\Delta(x,y)^\lambda}.
	\]
	The first summation is finite almost surely since $|T_x|_{x,\frac{1}{2}}<\infty$ and then it is a summation over finitely many vertex. For $\lambda\geq1$, the second summation is bounded above by $\sum_{y\in T_x,\Delta(x,y)<1}{\Delta(x,y)^{\frac{1}{2}}}\leq |T_x|_{x,\frac{1}{2}}<\infty.$
	Thus for $\lambda\geq 1$,  one also has $|T_x|_{x,\lambda}<\infty$ almost surely.
\end{proof}

 \subsection[\texorpdfstring{On the asymptotic behavior of $T_x\cap L_{-n}(x)$.}{On the asymptotic behavior of a tree with low slabs.}]{On the asymptotic behavior of $T_x\cap L_{-n}(x)$.}\label{sec: subsection 3.2}
 
 
 We have seen that $\mathbb{E}[|T_x\cap L_{-n}(x)|]\asymp n$ for the toy model for $n>0$. A natural question is what can we say about the almost sure behavior of $|T_x\cap L_{-n}(x)|$. Since $T_x$ is light almost surely, $T_x\cap L_n(x)=\emptyset$ a.s.\ if $n$ is large. So when talking about almost sure behavior, we are only interested in the  asymptotic behavior of  $|T_x\cap L_{-n}(x)|$ when $n$ tends to infinity.
 
 
 \begin{proposition}\label{prop: asymptotic for toy model}
 	Consider WUSF on the toy model $(\mathbb{T}_{b+1},\Gamma_{\xi})$. Let $x\in \mathbb{T}_{b+1}$ be an arbitrary vertex and $T_x$ be its connected component in the WUSF. 
 	Then $	\frac{|T_x\cap L_{-n}(x)|}{n} $ converges in distribution to a random variable $Z$ and $Z$ has Gamma distribution with density function $f(z)=Cz\exp(-\frac{2z}{1-1/b})\mathbf{1}_{ \{z>0\}}$, where $C$ is a normalizing constant. Moreover
 	\[
 	\frac{\log |T_x\cap L_{-n}(x)|}{\log n} \rightarrow 1\, \, \, \,  \as
 	\]
 	
 \end{proposition}
 \begin{proof}
 	Let $v$ be the last vertex on the future of $x$ such that $\Delta(v,x)=1$. In particular, if $x$ is the highest point on its future, then $v=x$. Denote the path on the future of $x$ from $x$ to $v$ by $\pi(x,v). $Let $N$ be the largest number $k$ such that there is some vertex $x_i$ in $\pi(x,v)\backslash\{v\}$ such that the bush at $x_i$ intersects $L_{-k}(x)$. 
 	Since all the bushes are finite trees almost surely, $N$ is finite almost surely. 
 	
 	For $n>N$, $|T_x\cap L_{-n}(x)|=1+Z_n$, where $Z_n$ is the size of the $n$-th generation of a critical branching process with immigration. Indeed the $1$ on the right hand side is the contribution of the future and  $Z_n$ is given as follows:
 	
 	\[
 	Z_0=0,\,  Z_n=\sum_{j=1}^{Z_{n-1}}Y_{n,j}+I_n,\, n\geq 1,
 	\]
 	where $Y_{i,j}$'s are i.i.d.\ random variables with the progeny distribution $\textnormal{Bin}(b,\frac{1}{b})$ and the immigration $I_i$'s  are i.i.d.\ random variables with distribution $\textnormal{Bin}(b-1,\frac{1}{b})$. In fact the immigration $I_i$ is the number of children of the vertex on the future in the previous generation. 
 	
 	Since we are only interested in the normalized asymptotic behavior, it suffices to show 
 	$\frac{Z_n}{n}$ converges in distribution to a random variable with Gamma distribution. This is a classical result regarding critical branching process with immigration; for example see Theorem 3 in \cite{Pakes1971}.
 	
 	For the almost sure result, it suffices to show that $\frac{\log Z_n}{\log n}\rightarrow1\,\as$, which is due to Theorem 1.1 of \cite{Wei1991}.  
 \end{proof}

 In view of the above proposition, one might ask whether 
 $\frac{|T_x\cap L_{-n}(x)|}{n} $ converges almost surely to a  random variable with Gamma distribution. However Proposition \ref{prop: asymptotic for toy model} together with the following proposition imply that the limit does not exist. 
 
 \begin{proposition}\label{prop: aymptotic for limsup of toy model}
 	For the toy model, using the same notation as Proposition \ref{prop: asymptotic for toy model}, one has 
 	\[
 	\limsup_{n\rightarrow\infty}\frac{|T_x\cap L_{-n}(x)|}{n\log \log n}\asymp 1\, \, \as
 	\]
 	
 \end{proposition} 
 \begin{proof}
 	Similar to the proof of Proposition \ref{prop: asymptotic for toy model}, it suffices to show that 
 	\be\label{eq: 4.73upperbound}
 	\limsup_{n\rightarrow\infty}\frac{Z_n}{n\log \log n}\preceq 1\, \, \as
 	\ee
 	and
 	\be\label{eq: 4.74lowerbound}
 	\limsup_{n\rightarrow\infty}\frac{Z_n}{n\log \log n}\succeq 1\, \, \as
 	\ee
 	
 	The inequality \eqref{eq: 4.73upperbound} is a direct consequence of Remark 2.2 of \cite{Wei1991}.
 	
 	Next we show \eqref{eq: 4.74lowerbound}. 
 	
 	If $I_i=t>0$, we write $GW_{i,1},\ldots,GW_{i,t}$ to be the descendant trees of these $t$ people immigrated in generation $i$. We use  $|GW_{i,j}\cap L_{-2^k}(x)|$ to denote the contribution  of $GW_{i,j}$ to $Z_{2^k}$.
 	
 	For a constant $c>0$, define the events $A_k(c)$ for $k\geq 2$ as follows:
 	\[
 	A_k(c):= \bigcup_{i=2^{k-1}}^{2^{k-1}+2^{k-2} } \{I_i\neq 0, |GW_{i,j}\cap L_{-2^k}(x)|\geq c\cdot2^k\log k \textnormal{ for some } 1\leq j\leq I_i \}.
 	\]

 	We claim that for some small enough constant $c$, $\sum_{k=2}^{\infty}\mathbb{P}[A_k(c)]=\infty$. Since the $A_k$'s are independent events, by the Borel--Cantelli lemma one has $\mathbb{P}[ A_k \textnormal{ i.o. }]=1$. 	Notice on the event $A_k(c)$,  $Z_{2^k}\geq c\cdot2^k\log k$. Therefore \eqref{eq: 4.74lowerbound} holds.
 	
 	Now it remains to prove the claim that $\sum_{k=2}^{\infty}\mathbb{P}[A_k(c)]=\infty$ for some small constant $c>0$. 
 	
 	For simplicity, write $B_{i,k}(c)=\{I_i\neq 0, |GW_{i,j}\cap L_{-2^k}(x)|\geq c\cdot2^k\log k \textnormal{ for some } 1\leq j\leq I_i \}$. 
 	
 	For $i\in [2^{k-1},2^{k-1}+2^{k-2}]$, $n=2^k-i\in [2^{k-2},2^{k-1}]$. By  Kolmogorov's estimate (see for example Theorem 12.7 of \cite{LP2016}) one has 
 	\be\label{eq: Kolmogorow estimate}
 	\mathbb{P}[|GW_{i,j}\cap L_{-2^k}(x)|> 0]\succeq \frac{1}{n}\succeq\frac{1}{2^k}.
 	\ee
 	
 	By the inequality (2.2) on page 588 and Theorem 3.3 in \cite{PekozRollin2011} and $\textnormal{Var}(Y_{i,j})=1-\frac{1}{b}$, there exists a constant $c_{1}$ that the conditional probability $P(c):=\mathbb{P}[|GW_{i,j}\cap L_{-2^k}(x)|> c\cdot 2^k\log k\big| |GW_{i,j}\cap L_{-2^k}(x)|> 0]$ satisfies
 	\be\label{eq: Yaglom's limit}
 	P(c)
 	\geq \exp(-\frac{2}{1-1/b}\frac{c\cdot2^k\log k}{n})-1.74c_{1}\sqrt{\frac{\log n}{n}}
 	\succeq\frac{1}{k^{c\cdot c_{2}}},
 	\ee
 	where $c_{2}=\frac{8}{1-1/b}$ (since $\frac{2^k}{n}\leq 4$).
 	
 	Therefore for $i\in [2^{k-1},2^{k-1}+2^{k-2}]$
 	\[
 	\mathbb{P}[B_{i,k}(c)]\geq \mathbb{P}[I_i=1]\times\mathbb{P}[|GW_{i,j}\cap L_{-2^k}(x)|> c\cdot 2^k\log k]\stackrel{\eqref{eq: Kolmogorow estimate}, \eqref{eq: Yaglom's limit}}{\succeq} \frac{1}{2^k\cdot k^{c\cdot c_{2}}}.
 	\]
 	Using independence of $B_{i,k}(c)$ for different $i$'s, one has 
 	\[
 	\mathbb{P}[A_{k}(c)]=1-\prod_{i=2^{k-1}}^{2^{k-1}+2^{k-2}}[1-\mathbb{P}[B_{i,k}(c)]]\succeq \frac{1}{k^{c\cdot c_{2}}},
 	\]
 	whence $\sum_{k=2}^{\infty}\mathbb{P}[A_k(c)]=\infty$ for small enough constant $c>0$.
 \end{proof}
 \begin{remark}
 	From the proof of Proposition \ref{prop: aymptotic for limsup of toy model} one actually obtains  that $\limsup_{n\rightarrow\infty}\frac{Z_n}{n\log\log n}\asymp 1$ for  general critical branching process  with immigration $Z_n$ if the  conditions in Theorem 1.1 of \cite{Wei1991} are satisfied with $\delta=1$ (to use Theorem 3.3 of \cite{PekozRollin2011}).

 \end{remark}
 
 Proposition \ref{prop: asymptotic for toy model} and \ref{prop: aymptotic for limsup of toy model} implies that almost surely $\lim_{n\rightarrow\infty}\frac{|T_x\cap L_{-n}(x)|}{n}$ does not exist. Proposition \ref{prop: asymptotic for toy model} also implies that for the toy model, for every $\varepsilon>0$, $\as\,\,  |T_x\cap L_{-n}(x)|\succeq n^{1-\varepsilon}$. Can one improve this lower bound for the toy model?

 \begin{question}\label{ques: asymotopic behavior of intersection with a lower level for general case}
 	What is the asymptotic behavior of $|T_x\cap L_{-n}(x)|$ as $n\rightarrow\infty$ for WUSF on general nonunimodular transitive graphs?
 \end{question}

 \section{The geometry of trees in WUSF with respect to the level structure}\label{sec: 4}
 
 Throughout the section we will assume $(G,\Gamma)$ is a nonunimodular transitive pair.
 We are interested in understanding the geometry of the tree $T_x$ with respect to the level structure induced by the modular function. For example, we wonder whether the tree $T_x$ is light. Moreover we are also interested in detailed information such as $\mathbb{P}[T_x\cap L_n(x)\neq \emptyset]$, $\mathbb{E}[|T_x\cap L_n(x)|^k]$, $\mathbb{E}[|T_x|_{x,\lambda}]$ and $\mathbb{P}[|T_x|_{x,\lambda }>R]$. Along the way similar quantities for $\mathfrak{T}_x,\mathfrak{P}(x)$ and $\mathfrak{F}(x,\infty)$ will also be studied. 
 
\begin{figure}[h!]
 	\centering
 	\begin{tikzpicture}[scale=0.8, text height=1.5ex,text depth=.25ex]

 	\draw (0,16)--(6,16)--(6,21)--(0,21)--(0,16);
 	\fill[teal!30!white](0,16)--(6,16)--(6,21)--(0,21)--(0,16);
 	\draw (0,20)--(6,20);
 	\node[below] at(3,21-0.2) {Simple random walk};
 	
 	\draw (1,18.5)--(5,18.5)--(5,19.5)--(1,19.5)--(1,18.5); 
 	\node[below] at(3,19.5-0.2) {Lemma \ref{lem:expectation for SRW on nonunimodular transitive graph}};
 	
 	\draw (1,16.5)--(5,16.5)--(5,17.5)--(1,17.5)--(1,16.5); 
 	\node[below] at(3,17.5-0.2) {Cor.\ \ref{cor:large deviation for SRW on nonunimodular transitive graph}};

 	\draw [->, very thick,>=stealth,color=blue](3,18.4)--(3,17.6);

 	\draw (0,8)--(6,8)--(6,13)--(0,13)--(0,8);
 	\fill[purple!30!white](0,8)--(6,8)--(6,13)--(0,13)--(0,8);
 	\draw (0,12)--(6,12);
 	\node[below] at(3,13-0.2) {First moments};

 	\draw (0.1,8.1)--(5.8,8.1)--(5.8,8.9)--(0.1,8.9)--(0.1,8.1);
 	\node[below] at (3,8.8) {Prop.\ \ref{prop: slab intersection for the x-WUSF component} $\mathbb{E}|\mathfrak{T}_x\cap L_n(x)|$};
 	
 	\draw (0.1,11.1)--(5.8,11.1)--(5.8,11.9)--(0.1,11.9)--(0.1,11.1);
 	\node[below] at (3,11.8) {Prop.\ \ref{prop: slab intersection for T_x} $\mathbb{E}|T_x\cap L_n(x)|$};
 	
 	\draw (0.1,9.1)--(5.9,9.1)--(5.9,10.9)--(0.1,10.9)--(0.1,9.1); 
 	\draw (1.4,9.1)--(1.4,10.9);
 	\draw (1.4,10)--(5.9,10);
 	\node[below] at (0.8,10.6) {Prop.};
 	\node[below] at (0.8,10) {\ref{prop: slab intersection for the future and past}};
 	\node[below] at (3.7,10.8){ $\mathbb{E}|\mathfrak{P}(x)\cap L_n(x)|$};
 	\node[below] at (3.7,9.8){ $\mathbb{E}|\mathfrak{F}(x,\infty)\cap L_n(x)|$};

 	\draw (0,0)--(6,0)--(6,5)--(0,5)--(0,0);
 	\fill[blue!30!white](0,0)--(6,0)--(6,5)--(0,5)--(0,0);
 	\draw (0,4)--(6,4);
 	\node[below] at(3,5-0.2) {Expected tilted volumes};
 	
 	\draw (0.2,3.1)--(5.8,3.1)--(5.8,3.9)--(0.2,3.9)--(0.2,3.1);
 	\node[below] at (3,3.8) {Prop.\ \ref{prop: expected tilted volume for WUSF component}, $\mathbb{E}[|T_x|_{x,\lambda}]$};
 	
 	\draw (0.2,0.1)--(5.8,0.1)--(5.8,2.8)--(0.2,2.8)--(0.2,0.1); 
 	\draw (1.5,0.1)--(1.5,2.8);
 	\draw (1.5,1)--(5.8,1);
 	\draw (1.5,1.9)--(5.8,1.9);
 	\node[below] at (0.9,2.1) {Prop.};
 	\node[below] at (0.9,1.43) {\ref{prop: tilted volume bounds for x-component, past and future}};
 	\node[below] at (3.65,2.7) {$\mathbb{E}[|\mathfrak{T}_x|_{x,\lambda}]$};
 	\node[below] at (3.65,1.8) {$\mathbb{E}[|\mathfrak{P}(x)|_{x,\lambda}]$};
 	\node[below] at (3.65,0.9) {$\mathbb{E}[|\mathfrak{F}(x,\infty)|_{x,\lambda}]$};

 	\draw (12,8)--(19,8)--(19,13)--(12,13)--(12,8);
 	\fill[ green!60!white](12,8)--(19,8)--(19,13)--(12,13)--(12,8);
 	\draw (12,12)--(19,12);
 	\node[below] at(15.5,13-0.2) {Prob.\ of intersections};

 	\draw (12.15,11.1)--(18.8,11.1)--(18.8,11.9)--(12.15,11.9)--(12.15,11.1);
 	\node[below] at (15.5,11.8) {Prop.\ \ref{prop: exponential decay of the probability that T_x intersecting a high slab}, $\mathbb{P}[T_x\cap L_n(x)\neq\emptyset]$};
 	
 	\draw (12.15,10.1)--(18.8,10.1)--(18.8,10.9)--(12.15,10.9)--(12.15,10.1);
 	\node[below] at (15.5,10.8) {Prop.\ \ref{prop: decay of probability for x-WUSF component insecting a slab far away}, $\mathbb{P}[\mathfrak{T}_x\cap L_n(x)\neq\emptyset]$};

 	\draw  (12.15,8.2)--(18.8,8.2)--(18.8,9.8)--(12.15,9.8)--(12.15,8.2);
 	\draw (13.65,8.2)--(13.65,9.8);
 	\draw (13.65,9)--(18.8,9);
 	\node[below] at (12.9,9.7) {Prop.};
 	\node[below] at (12.9,9) {\ref{prop: exponentail decay for the past and future intersecting with a slab}};
 	\node[below] at (16.25,9.7) {$\mathbb{P}[\mathfrak{P}(x)\cap L_n(x)\neq\emptyset]$};
 	\node[below] at (16.25,8.9) {$\mathbb{P}[\mathfrak{F}(x,\infty)\cap L_n(x)\neq\emptyset]$};

 	\draw (12,16)--(19,16)--(19,21)--(12,21)--(12,16);
 	\fill[brown!80!white](12,16)--(19,16)--(19,21)--(12,21)--(12,16);
 	\draw (12,20)--(19,20);
 	\node[below] at(15.5,21-0.2) {High moments};
 	
 	\draw (12.1,19.1)--(18.9,19.1)--(18.9,19.9)--(12.1,19.9)--(12.1,19.1);
 	\node[below] at (15.5,19.8) {Prop.\ \ref{prop: high moments for x-component intersecting a slab}, $\mathbb{E}[|\mathfrak{T}_x\cap L_n(x)|^k]$};
 	
 	\draw (12.09,16.1)--(18.91,16.1)--(18.91,16.95)--(12.09,16.95)--(12.09,16.1);
 	\node[below] at (15.5,16.8) {Cor.\ \ref{cor: high moments for the future intersect with a slab}, $\mathbb{E}[|\mathfrak{F}(x,\infty)\cap L_n(x)|^k]$};
 	
 	\draw (12.1,17.2)--(15,17.2)--(15,18)--(12.1,18)--(12.1,17.2); 
 	\node [below] at (13.5,17.9) {Prop.\ \ref{prop: high moments for T_x intersecting a slab}};
 	\draw [->, very thick,>=stealth,color=blue](13.5,18.9)--(13.5,18.1);
 	\node[right]at (12.2,18.5) {Lem.\ \ref{lem: bounds high moments of T_x intersections in terms of the x-WUSF ones}};

 	\draw (16,17.2)--(18.9,17.2)--(18.9,18)--(16,18)--(16,17.2); 
 	\node [below] at (17.5,17.9) {Cor.\ \ref{cor: high moments for the past}};
 	\draw [->, very thick,>=stealth,color=blue](17.5,18.9)--(17.5,18.1);
 	\node[right]at (16.2,18.5) {Lem.\ \ref{lem:stochastic domination}};

 	\draw (12,0)--(19,0)--(19,5)--(12,5)--(12,0);
 	\fill[teal!60!white](12,0)--(19,0)--(19,5)--(12,5)--(12,0);
 	\draw (12,4)--(19,4);
 	\node[below] at(15.5,5-0.2) {Tail prob.\ for tilted volumes};

 	\draw (12.2,3.1)--(18.8,3.1)--(18.8,3.9)--(12.2,3.9)--(12.2,3.1);
 	\node[below] at (15.5,3.8) {Prop.\ \ref{prop: tail prob for general case}, $\mathbb{P}[|T_x|_{x,\lambda}>R]$};
 	
 	\draw (12.2,2.1)--(18.8,2.1)--(18.8,2.9)--(12.2,2.9)--(12.2,2.1);
 	\node[below] at (15.5,2.8) {Prop.\ \ref{prop: tail probability for tilted volume of x-wusf component}, $\mathbb{P}[|\mathfrak{T}_x|_{x,\lambda}>R]$};

 	\draw  (12.2,0.2)--(18.8,0.2)--(18.8,1.8)--(12.2,1.8)--(12.2,0.2);
 	\draw (13.7,0.2)--(13.7,1.8);
 	\draw (13.7,1)--(18.8,1);
 	\node[below] at (12.95,1.7) {Prop.};
 	\node[below] at (12.95,1) {\ref{prop: tail probability for tilted volume of the past and future}};
 	\node[below] at (16.25,1.7) {$\mathbb{P}[|\mathfrak{P}(x)|_{x,\lambda}>R]$};
 	\node[below] at (16.25,0.9) {$\mathbb{P}[|\mathfrak{F}(x,\infty)|_{x,\lambda}>R]$};

 	\draw (7,20)--(11,20)--(11,21)--(7,21)--(7,20);
 	\fill[brown!80!white](7,20)--(11,20)--(11,21)--(7,21)--(7,20);
 	\node[below] at (9,20.85) {Lemma \ref{lem: point function ineq with k=3}, \ref{lem: point function inequality for general k}};

 	\draw (8.5,16)--(11,16)--(11,17)--(8.5,17)--(8.5,16);
 	\fill[teal!30!white](8.5,16)--(11,16)--(11,17)--(8.5,17)--(8.5,16);
 	\node[below]at(9.75,16.85) {Prop.\ \ref{prop: intersection of a srw and a slab}};

 	\draw [->, very thick,>=stealth,color=blue](11,16.5)--(12.1,16.5);
 	
 	\draw [ very thick,color=blue](5.8,8.5)--(7,8.5)--(7,19)--(5,19);	
 	
 	\draw [->, very thick,>=stealth,color=blue](7,16.5)--(8.5,16.5);
 	
 	\draw [ very thick,color=blue](5.8,8.5)--(8,8.5)--(8,20);
 	\draw [->, very thick,>=stealth,color=blue](8,19.5)--(12.1,19.5);
 	\node[below] at (10,19.5) {tree-graph ineq.};

 	\draw [->, very thick,>=stealth,color=blue](3,15.9)--(3,13.1);
 	\node at (3,14.5) {Wilson's algorithm+TMTP};
 	
 	\draw [->, very thick,>=stealth,color=blue](3,7.9)--(3,5.1);

 	\draw [ very thick,color=blue](19.1,18.5)--(20,18.5)--(20,6.5)--(13,6.5)--(13,7.9);	
 	\draw [->, very thick,>=stealth,color=blue](16,6.5)--(16,5.1);
 	
 	
 	\draw [very thick,color=blue](6.1,10.5)--(10,10.5)--(10,14.5)--(17,14.5)--(17,15.9);
 	
 	\draw [->, very thick,>=stealth,color=blue](14.5,14.5)--(14.5,13.1);
 	
 	\end{tikzpicture}
 	\caption{The dependencies among the lemmas and propositions in Section \ref{sec: 4}.}
 	\label{fig: relations among theorems}
 \end{figure}
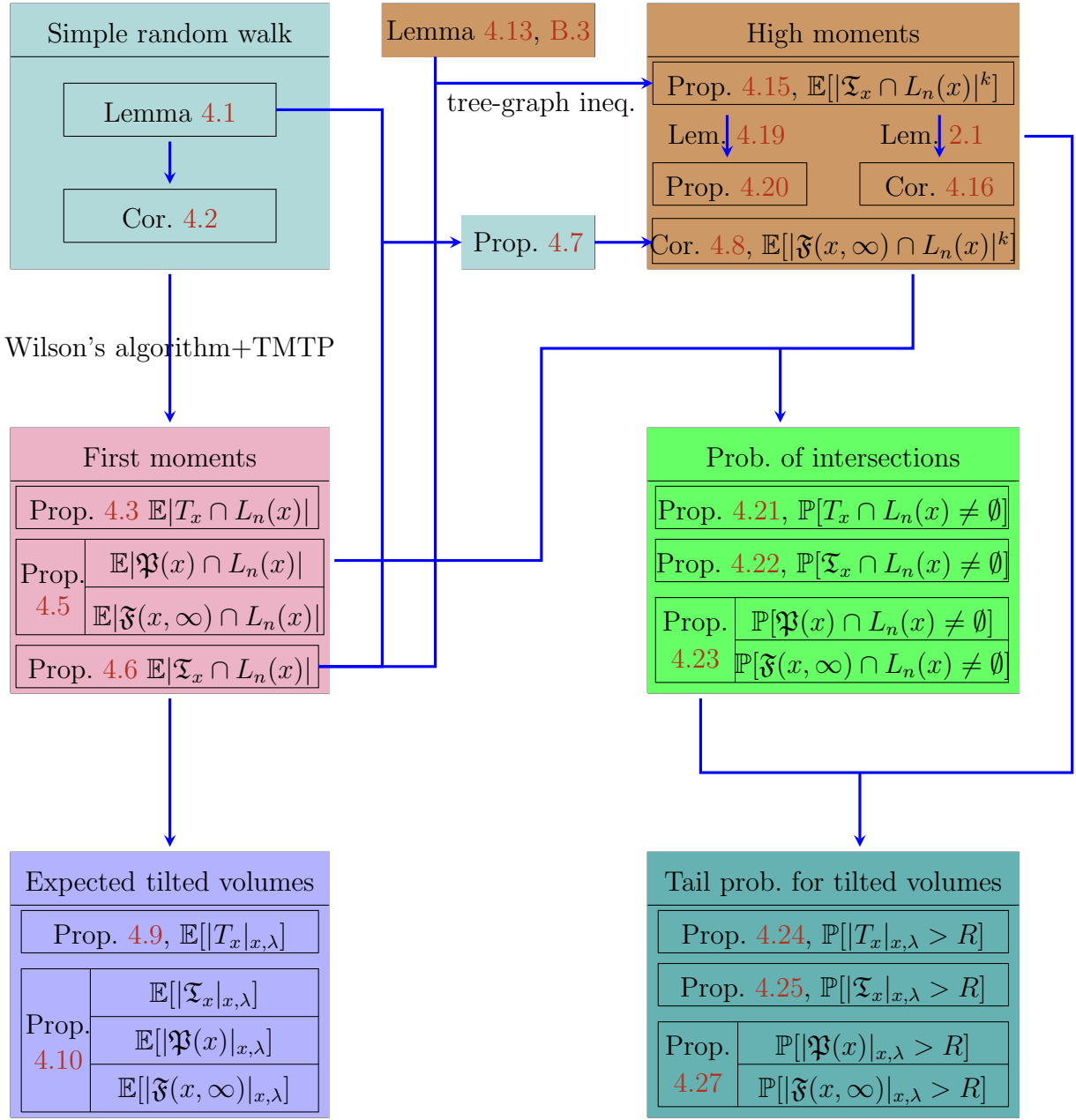

 The dependencies of the lemmas and propositions in the section can be summarized in Figure \ref{fig: relations among theorems}. 

 \subsection{Simple random walk on nonunimodular transitive graphs}
 Wilson's algorithm is a very useful tool for studying WUSF. Here we begin with a simple lemma about simple random walk on a nonunimodular transitive graph. 
 \begin{lemma}\label{lem:expectation for SRW on nonunimodular transitive graph}
 	Suppose $\{X_n\}_{n\geq0}$ is a simple random walk on $G$ and $\lambda\in\mathbb{R}$. Then one has 
 	\[
 	\mathbb{E}[\Delta(X_0,X_1)^\lambda]<1 \textnormal{ if and only if }\lambda\in(0,1).
 	\]
 	In particular, if $\lambda\in(0,1)$, then 
 	\be\label{eq: expected summation for srw}
 	\sum_{k=0}^{\infty}\mathbb{E}[\Delta(X_0,X_k)^\lambda]<\infty.
 	\ee
 	Moreover, one also has 
 	\[
 	\mathbb{E}[\log\Delta(X_0,X_1)]<0.
 	\]
 \end{lemma}
 \begin{proof}
 	Let $D$ denote the degree of the transitive graph $G$. For  each $x\in V$, let $\{X_n^x\}_{n\geq0}$ be a simple random walk starting at $X_0^x=x$.
 	
 	Define a random function $f:V^2\rightarrow[0,\infty]$ to be $f(x,y)=\mathbf{1}_{\{y=X_1^x\}}\Delta(x,y)^\lambda$, then the tilted mass-transport principle and the cocycle identity yield that 
 	\[
 	\mathbb{E}[\Delta(X_0,X_1)^\lambda\mid X_0=x]=\mathbb{E}[\sum_{y\in V}f(x,y)]=\mathbb{E}[\sum_{y\in V}f(y,x)\Delta(x,y) ]=	\mathbb{E}[\Delta(X_0,X_1)^{1-\lambda}|X_0=x].
 	\]
 	In particular, taking $\lambda=0$ one has that 
 	\[
 	1=\mathbb{E}[\sum_{y\in V}f(x,y)]=\mathbb{E}[\sum_{y\in V}f(y,x)\Delta(x,y) ]=\sum_{y\sim x}\frac{1}{D}\Delta(x,y)=\mathbb{E}[\Delta(X_0,X_1)|X_0=x].
 	\]
  Since the above two equations are true for all $x\in V(G)$, one has
  \[
  	\mathbb{E}[\Delta(X_0,X_1)^\lambda]=\mathbb{E}[\Delta(X_0,X_1)^{1-\lambda}]
  \]
 	and
 	\[
 	1=\mathbb{E}[\Delta(X_0,X_1)].
 	\]
 	For $\lambda\in (0,1)$, the function $x\mapsto x^\lambda$ is a strictly concave function on $[0,\infty]$ and  $\Delta(X_0,X_1)$ takes more than one value with positive probability by  nonunimodularity. By Jensen's inequality one has that 
 	\[
 	\mathbb{E}[\Delta(X_0,X_1)^\lambda]<\mathbb{E}[\Delta(X_0,X_1)]^\lambda=1.
 	\]
 	Using the cocycle identity and independence, one has $\mathbb{E}[\Delta(X_0,X_k)^\lambda]=\left(\mathbb{E}[\Delta(X_0,X_1)^\lambda]\right)^k$, whence \eqref{eq: expected summation for srw} follows from the fact that $\mathbb{E}[\Delta(X_0,X_1)^\lambda]<1$.
 	
 	For $\lambda>1$ or $\lambda<0$, the  function $x\mapsto x^\lambda$ is a strictly convex function on $[0,\infty]$. By Jensen's inequality one has that 
 	\[
 	\mathbb{E}[\Delta(X_0,X_1)^\lambda]>\mathbb{E}[\Delta(X_0,X_1)]^\lambda=1.
 	\]
 Since $x\mapsto \log x$ is  strictly concave on $(0,\infty)$,	$\mathbb{E}[\log\Delta(X_0,X_1)]<\log\mathbb{E}[\Delta(X_0,X_1)]=0.$\end{proof}

Lemma \ref{lem:expectation for SRW on nonunimodular transitive graph} implies that simple random walk on nonunimodular transitive graphs has a drift towards the lower slabs. 
 
 \begin{corollary}\label{cor:large deviation for SRW on nonunimodular transitive graph}
 	Fix $x\in V(G)$ and suppose $\{X_n^x\}_{n\geq0}$ is a simple random walk on $G$ starting at $x$. Let $n\vee 1=\max\{ 1,n\}$ denote the maximum of $n$ and $1$.
  	Then 
 	\be\label{eq: occupation times for srw on -n slab}
 	\mathbb{E}\left[\sum_{k=0}^{\infty}\mathbf{1}_{\{X_k^x\in L_{-n}(x) \}}\right]\asymp 1
 	\ \textnormal{ and }\ 
 	\mathbb{E}\left[\sum_{k=0}^{\infty}k\cdot\mathbf{1}_{\{X_k^x\in L_{-n}(x) \}}\right]\asymp n\vee 1.
 	\ee

 	For an integer $n\geq0$, let $\tau_{-n}^x:=\sup\{k:X_k\in L_{-n}(x)\}$ be the last visit time of the slab $L_{-n}(x)$. Then there exist constants $c_1,c_2,c_3>0$ such that 
 	\be \label{eq:tail probability for last visit time of SRW with a fixed level}
 	\mathbb{P}[\tau_{-n}^x\geq R]\leq c_1\exp(-c_2R),\ \ \forall \ R\geq c_3(n\vee 1).
 	\ee
 \end{corollary}
 \begin{proof}
 	By the cocycle identity $\log \Delta(X_0^x,X_k^x)=\sum_{j=1}^{k}\log \Delta(X_{j-1}^x,X_j^x)$. 
 	Since $\log \Delta(X_{j-1}^x,X_j^x)$ are independent, identically distributed random variables, the law of large numbers and Lemma \ref{lem:expectation for SRW on nonunimodular transitive graph}  then imply that 
 	\be\label{eq: LLN for srw}
 	\frac{1}{k}\log \Delta(X_0^x,X_k^x)=\frac{1}{k}\sum_{j=1}^{k}\log \Delta(X_{j-1}^x,X_j^x)\rightarrow
 	\mathbb{E}[\log\Delta(X_0^x,X_1^x)]<0\ \as
 	\ee
 	Also there exists a constant $c\in(0,1)$ such that 
 	\be \label{eq: positive probability for a cutpoint}
 	\mathbb{P}[\log\Delta(X_0^x,X_1^x)=-t_0,\log\Delta(X_0^x,X_k^x)\leq -2t_0,\ \forall k\geq 2]\geq c.
 	\ee
 	
 	Write $N_{-n}=\sum_{k=0}^{\infty}\mathbf{1}_{\{X_k^x\in L_{-n}(x) \}}$.  Each time a simple random walk visits $L_n(x)$, with probability bounded below by some positive constant the simple random walk will leave $L_n(x)$ and never visit it again. Then the strong Markov property implies that 
 	\be\label{eq: strong Markov property for SRW}
 	\mathbb{P}[N_{-n}\geq k]\leq (1-c)^{k-1},\,\,\forall \, k\geq1.
 	\ee
 	Hence the upper bound in the  first part of \eqref{eq: occupation times for srw on -n slab} holds.

 	Equation \eqref{eq: LLN for srw} also implies that $\tau_{-n}^x<\infty$ almost surely. The lower bound for \eqref{eq: occupation times for srw on -n slab} are straightforward since $\sum_{k=0}^{\infty}\mathbf{1}_{\{X_k^x\in L_{-n}(x) \}}\geq 1$ and $\sum_{k=0}^{\infty}k\cdot \mathbf{1}_{\{X_k^x\in L_{-n}(x) \}}\geq n-1$ almost surely.

 	For $k\geq R\geq \frac{2(n+1)t_0}{-\mathbb{E}[\log\Delta(X_0,X_1)]}$, one has $-\frac{(n+1)t_0}{k}\geq -\frac{(n+1)t_0}{R}\geq \frac{1}{2}\mathbb{E}[\log\Delta(X_0,X_1)]>\mathbb{E}[\log\Delta(X_0,X_1)]$. Therefore
 	for such $R$
 	\begin{eqnarray*}
 		\mathbb{P}[\tau_{-n}^x\geq R]&=&\mathbb{P}[\cup_{k\geq R}\{ X_{k}^x\in L_{-n}(x) \}]\leq\sum_{k\geq R}\mathbb{P}[X_{k}^x\in L_{-n}(x)]\\
 		&\leq&\sum_{k\geq R}\mathbb{P}[\log\Delta(X_0^x,X_k^x)\geq (-n-1)t_0]\\
 		&\leq&\sum_{k\geq R}\mathbb{P}\left[\frac{1}{k}\log\Delta(X_0^x,X_k^x)\geq -\frac{(n+1)t_0}{k}\right]\\
 		&\leq&\sum_{k\geq R}\mathbb{P}\left[\frac{1}{k}\log\Delta(X_0^x,X_k^x)\geq \frac{1}{2}\mathbb{E}[\log\Delta(X_0,X_1)]\right].
 	\end{eqnarray*}
 	By large deviation principle (e.g.\ \cite[Theorem 2.7.7]{Durrett2010}) we know there exists constants $c_4,c_5>0$ such that 
 	\be
 	\mathbb{P}\left[\frac{1}{k}\log\Delta(X_0^x,X_k^x)\geq \frac{1}{2}\mathbb{E}[\log\Delta(X_0^x,X_1^x)]\right]
 	\leq c_4\exp(-c_5k),\ \forall \ k\geq1.
 	\ee
 	Hence taking $c_3=\frac{4t_0}{-\mathbb{E}[\log\Delta(X_0,X_1)]}$, then for $n\geq0$ and $R\geq c_3(n\vee 1)\geq \frac{2(n+1)t_0}{-\mathbb{E}[\log\Delta(X_0,X_1)]}$ one has 
 	\be 
 	\mathbb{P}[\tau_{-n}^x\geq R]\leq \sum_{k=c_3(n\vee1)}^{\infty} c_4\exp(-c_5k) =\frac{c_4}{1-\exp(-c_5)}\exp(-c_5c_3(n\vee 1)).
 	\ee
 	Then \eqref{eq:tail probability for last visit time of SRW with a fixed level} holds if one takes $c_3=\frac{4t_0}{-\mathbb{E}[\log\Delta(X_0,X_1)]}$, $c_1=\frac{c_4}{1-\exp(-c_5)}$ and $c_2=c_5c_3$.

 	Notice that for $k\geq 1,j\geq c_3(n\vee 1)$, there exists  constants $c_6,c_7>0$ such that 
 	\be\label{eq: exponential decay for N_-n,tau-n large}
 	\mathbb{P}[N_{-n}=k,\tau_{-n}^x=j]\leq \min\{(1-c)^{k-1},c_1\exp(-c_2j) \}\leq c_6\exp(-c_7(k+j)).
 	\ee
 	Hence 
 	\begin{eqnarray}\label{eq: tau large part}
 	\sum_{k=1}^{\infty}\sum_{j=c_3(n\vee1)}^{\infty}k\cdot j\cdot \mathbb{P}[N_{-n}=k,\tau_{-n}^x=j]&\leq& c_6\sum_{k=1}^{\infty}k\exp(-c_7k)\sum_{j=1}^{\infty}j\exp(-c_7j)\nonumber\\
 	&=&\frac{c_6\exp(-2c_7)}{(1-\exp(-c_7))^4}<\infty.
 	\end{eqnarray}
 Note 
 \begin{eqnarray}\label{eq: tau small part}
 \sum_{k=1}^{\infty}\sum_{j=1}^{c_3(n\vee 1)}k\cdot j\cdot \mathbb{P}[N_{-n}=k,\tau_{-n}^x=j]&\leq&
 \sum_{k=1}^{\infty}\sum_{j=1}^{c_3(n\vee 1)}k\cdot c_3(n\vee 1)\cdot \mathbb{P}[N_{-n}=k,\tau_{-n}^x=j]\nonumber\\
 &\leq&c_3(n\vee 1)\sum_{k=1}^{\infty}k\mathbb{P}[N_{-n}=k]\nonumber\\
 &\leq &c_3(n\vee 1)\sum_{k=1}^{\infty}k(1-c)^{k-1}=\frac{c_3(n\vee 1)}{c^2}.
 \end{eqnarray}
 
 Combining \eqref{eq: tau large part} and \eqref{eq: tau small part} one has that there exists a constant $c_8>0$ such that 
 \[
 \mathbb{E}[N_{-n}\tau_{-n}^x]=\sum_{k=1}^{\infty}\sum_{j=1}^{\infty}k\cdot j\cdot \mathbb{P}[N_{-n}=k,\tau_{-n}^x=j]\leq c_8(n\vee 1).
 \]
 	
 	Note that $\sum_{k=0}^{\infty}k\cdot\mathbf{1}_{\{X_k^x\in L_{-n}(x) \}}\leq N_{-n}\tau_{-n}^x$, whence the upper bound in the second part of    \eqref{eq: occupation times for srw on -n slab} holds. 
 \end{proof}

 \subsection{First moment of intersections with a slab and the expected tilted volumes}
 
 We summarize the expectations of $\mathbb{E}[|\star \cap L_n(x)|]$ for $\star\in\{ T_x, \mathfrak{T}_x,\mathbb{P}(x),\mathfrak{F}(x,\infty) \}$ in Table \ref{table: first moments}.  Based on these first moments we can also derive estimates on the tilted volumes. 
 
 
 \begin{table}[ht]
 	\centering
 	\begin{tabular}{|c|c|c|}
 		\hline
 		Quantities $(n\in \mathbb{Z})$ & Results & Position\\
 		\hline
 		$\mathbb{E}[|T_x\cap L_n(x)|]$ & $\left\{
 		\begin{array}{cc}
 		\asymp (n\vee 1)e^{-t_0n} & \textnormal{ if }n\geq0\\
 		\asymp |n|& \textnormal{ if }n<0\\
 		\end{array}
 		\right.$ & Proposition \ref{prop: slab intersection for T_x}\\
 		\hline
 			$\mathbb{E}[|\mathfrak{P}(x)\cap L_n(x)|]$ &  $\left\{
 			\begin{array}{cc}
 			\asymp e^{-t_0n} & \textnormal{ if }n\geq0\\
 			\asymp 1& \textnormal{ if }n<0\\
 			\end{array}
 			\right.$ & Proposition \ref{prop: slab intersection for the future and past} \\
 		\hline
 		$\mathbb{E}[|\mathfrak{F}(x,\infty)\cap L_{n}(x)|]$ &  $\left\{
 		\begin{array}{cc}
 		\asymp e^{-t_0n} & \textnormal{ if }n\geq0\\
 		\asymp 1& \textnormal{ if }n<0\\
 		\end{array}
 		\right.$& Proposition \ref{prop: slab intersection for the future and past} \\
 		\hline
 		$\mathbb{E}[|\mathfrak{T}_x\cap L_{n}(x)|]$ & $\left\{
 		\begin{array}{cc}
 		\asymp e^{-t_0n} & \textnormal{ if }n\geq0\\
 		\asymp 1& \textnormal{ if }n<0\\
 		\end{array}
 		\right.$ & Proposition \ref{prop: slab intersection for the x-WUSF component}\\
 		\hline
 	\end{tabular}
 \centering
 	\caption{First moments for the intersections with a slab}
 	\label{table: first moments}
 \end{table}
 
 \begin{proposition}\label{prop: slab intersection for T_x}
 	There exist  constants $c_1>0,c_{2}>1$ such that for every $n\geq1$
 		\be\label{eq: expectation of intersection of T_x with a low slab}
 		1\leq \mathbb{E}[|T_x\cap L_{0}(x)|]\leq c_{2}\textnormal{ and }
 		c_1n\leq \mathbb{E}[|T_x\cap L_{-n}(x)|]\leq c_{2}n
 		\ee
 and
 	\be\label{eq: expectation of intersection of T_x with a high slab}
 	c_1n\exp(-t_0(n+1))\leq \mathbb{E}[|T_x\cap L_{n}(x)|]\leq c_{2}n\exp(-t_0(n-1)), \ \forall n\geq1
 	\ee	
 	
 \end{proposition}

The following lemma is a simple application of  Wilson's algorithm. The inequality \eqref{eq: upper bound on two-point function in WUSF}  comes from  the proof of Theorem 13.1 in \cite{BLPS2001} where expected quadratic growth with respect to  extrinsic graph metric was proved. Also a version of \eqref{eq: upper bound on two-point function in v-WUSF} is used in the proof of Lemma 6.6 in \cite{Hutchcroft2018a} and the Lemma 6.6 is closely related to the quadratic growth with respect to intrinsic graph metric. See Theorem 1.3, Corollary 6.4 and 6.12 in  \cite{Hutchcroft2018a} for the quadratic growth with respect to intrinsic graph metric. 
 \begin{lemma}\label{lem: upper bound on two-point functions}
 	Recall that $\{x\leftrightarrow y \}$  denotes the event that  $x,y$ are in the same connected component of the WUSF sample $\mathfrak{F}$. Then
 	\be\label{eq: upper bound on two-point function in WUSF}
 	\mathbb{P}[x\leftrightarrow y]\leq \sum_{m=0}^{\infty}(m+1)\mathbb{P}[X_m^x=y],
 	\ee
 	where $\{X_m^x\}_{m\geq0}$ is a simple random walk on $G$ started from $x$. 
 	
 	Let $y\in\mathfrak{T}_x$ denote the event that $y$ lies in the connected component of $x$ in the $x$-WUSF sample $\mathfrak{F}_x$.
 	Then 
 	\be\label{eq: upper bound on two-point function in v-WUSF}
 	\mathbb{P}[y\in\mathfrak{T}_x]=\mathbb{P}[x\in\mathfrak{T}_y]\leq \sum_{m=0}^{\infty}\mathbb{P}[X_m^x=y].
 	\ee
 	
 \end{lemma}
 \begin{proof}

 	From the reversibility of simple random walk and the regularity of $G$, one has that for any vertices $x,y$ and $k\leq m$,
 	\[
 	\sum_{z\in V}\mathbb{P}[X_k^x=z]\mathbb{P}[X_{m-k}^y=z]=
 	\sum_{z\in V}\mathbb{P}[X_k^x=z]\mathbb{P}[X_{m-k}^z=y]=\mathbb{P}[X_m^x=y].
 	\]
 	By  Wilson's algorithm, 
 	\begin{eqnarray*}
 	\mathbb{P}[x\leftrightarrow y]&\leq& \mathbb{P}[X^x \textnormal{ intersects }X^y]\\
 	&\leq& \sum_{z\in V}\sum_{m=0}^{\infty}\sum_{k=0}^{m}\mathbb{P}[X_k^x=z]\mathbb{P}[X_{m-k}^y=z]\\
 	&=& \sum_{m=0}^{\infty}\sum_{k=0}^{m}\mathbb{P}[X_m^x=y]=\sum_{m=0}^{\infty}(m+1)\mathbb{P}[X_m^x=y].
 	\end{eqnarray*}

 	 If we use Wilson's algorithm to sample $\mathfrak{F}_x$ starting with the simple random walk $X^y$, then we have that 
 	 \[
 	 \mathbb{P}[y\in\mathfrak{T}_x]=\mathbb{P}[\sigma_x^y<\infty]
 	 \]
 	 where $\sigma_x^y$ denotes the first visit time of $x$ by the simple random walk $X^y$.
 	 The reversibility of simple random walk on $G$ and the regularity of $G$ then imply that $\mathbb{P}[y\in\mathfrak{T}_x]=\mathbb{P}[\sigma_x^y<\infty]=\mathbb{P}[\sigma_y^x<\infty]=\mathbb{P}[x\in\mathfrak{T}_y]$.
 	 
 	 Therefore one has that 
 	 \[
 	 \mathbb{P}[y\in\mathfrak{T}_x]=\mathbb{P}[x\in\mathfrak{T}_y]=\mathbb{P}[\sigma_y^x<\infty] \leq\sum_{m=0}^{\infty}\mathbb{P}[X_m^x=y].\qedhere
 	 \]
 \end{proof}
 
 \begin{proof}[Proof of Proposition \ref{prop: slab intersection for T_x}]
 We first establish the upper bounds in \eqref{eq: expectation of intersection of T_x with a low slab}. 
 
 For $n\geq 0$, one has that 
 \begin{eqnarray}
 \mathbb{E}[|T_x\cap L_{-n}(x)|]
 &=&\mathbb{E}\Big[\sum_{y\in V}{\mathbf{1}_{\{y\in L_{-n}(x) \}}\cdot \mathbf{1}_{ \{y\leftrightarrow x \} } }\Big]\nonumber\\
 &\stackrel{\eqref{eq: upper bound on two-point function in WUSF} }{\leq}&
 \sum_{y\in V}\mathbb{E}[\mathbf{1}_{\{y\in L_{-n}(x) \}}]\mathbb{E}\Big[\sum_{m=0}^{\infty} (m+1)\mathbf{1}_{\{ X_m^x=y \}}\Big]\nonumber\\
 &= &\mathbb{E}\Big[\sum_{m=0}^{\infty} (m+1)\mathbf{1}_{\{ X_m^x\in L_{-n}(x) \}}\Big]\nonumber\\
 &\stackrel{\eqref{eq: occupation times for srw on -n slab}}{\leq }&c_0+c_0(n\vee 1).
 \end{eqnarray}
 Taking $c_{2}=2c_0$ one has the upper bounds in \eqref{eq: expectation of intersection of T_x with a low slab}. Since $x\in T_x\cap L_{0}(x)$, the lower bound $\mathbb{E}[T_x\cap L_0(x)]\geq 1$ is trivial.

 Next we prove the inequality $\mathbb{E}[|T_x\cap L_{-n}(x)|]\geq c_1n$ for $n\geq1$ using a similar strategy as the one used in the proof of Theorem 13.1 of \cite{BLPS2001}. 
 
 Denote by $g(u,v)=\sum_{m=0}^{\infty}\mathbb{P}[X_m^u=v]$
 the Green function for simple random walk on $G$. Since $\Gamma$ is a closed subgroup of automorphisms that acts transitively on $G$, $G$ is nonamenable (Proposition 8.14 of \cite{LP2016}) and hence the spectral radius
 $\rho(G)<1$.  By the Cauchy--Schwarz inequality $p_n(u,v)=\mathbb{P}[X_n^u=v]=\langle P^n \mathbf{1}_u,\mathbf{1}_v\rangle\leq \rho(G)^n$ and hence that 
 \[
 g(u,v)\leq \sum_{n\geq d(u,v)}\rho(G)^n\leq \frac{1}{1-\rho(G)}\rho(G)^{d(u,v)},
 \]
 where $d(u,v)$ is the graph distance of $u,v$ in $G$. 
 
 We use Wilson's algorithm to generate the WUSF sample $\mathfrak{F}$ by starting with a simple random walk. Then the future of $x$ is a subset of $X^x$. By \eqref{eq: LLN for srw} one has that $1\leq |\mathfrak{F}(x,\infty)\cap L_{-n}(x)|<\infty$  for every $n\geq0$. Also
 \[
 \mathbb{P}[y\leftrightarrow x]=\mathbb{P}[X^y \cap \mathfrak{F}(x,\infty)\neq\emptyset]
 \]

 Let $B(x,n):=\{y\in V(G)\colon d(x,y)\leq n\}$ denote the ball of radius $n$ centered at $x$. We shall show that there exists a constant $c_{1}>0$ such that for any set $S\subset B(x,n)$ that contains exactly one vertex at distance $k$ from $x$ for $1\leq k\leq n$, the following inequality holds.
 \be\label{eq: lower bound 1 of occupation time for S}
 \mathbb{E}\Big[\sum_{y\in L_{-n}(x)}{ \mathbf{1}_{ \{L_y(S)>0\} } }\Big]\geq c_{1}n,
 \ee
 where $L_y(S)$ is the total occupation time of $S$ by $X^y$, that is, $L_y(S):=\sum_{w\in S}\sum_{m\geq0}\mathbf{1}_{\{X_m^y=w\}}$. 
 
 If $S\subset \mathfrak{F}(x,\infty)\cap B(x,n)$ and $S$ contains exactly one vertex from each sphere $\partial B(x,k)$ for $ k=1,\ldots,n$, then 
 \be\label{eq: two-point fcn bigger than nonzero occupation}
 \mathbb{P}[y\leftrightarrow x]\geq \mathbb{P}[L_y(S)>0].
 \ee
 
 Since simple random walk $X^w$ on $G$ visits $S_{s,t}(w)$ for every $t\leq 0,t-s\geq t_0$ almost surely, one has that $\sum_{y\in L_{-n}(x)}{g(w,y)}\geq1$, $\forall w\in S$. Thus 
 $\sum_{y\in L_{-n}(x),w\in S}{g(w,y)}\geq |S|= n$. Hence for $n\geq1$,   
 \begin{eqnarray}\label{eq: lower bound 2 of occupation time for S}
\mathbb{E}\Big[\sum_{y\in L_{-n}(x)}L_y(S)\Big] 
&= & \mathbb{E}\Big[\sum_{y\in L_{-n}(x)}\sum_{w\in S}g(y,w)\Big] \nonumber\\
&=&\mathbb{E}\Big[\sum_{y\in L_{-n}(x)}\sum_{w\in S}g(w,y)\Big] \nonumber\\
&\geq&n.
 \end{eqnarray}
 
 For any $w,v\in S$, $d(w,v)\geq |d(w,x)-d(v,x)|$, whence for any $w\in S$ one has 
 \[
 \mathbb{E}[L_w(S)]\leq \sum_{j\geq0}2C\rho(G)^{j-1}=:c_{3}.
 \]
 The Markov property implies that $\mathbb{E}[L_y(S)|L_y(S)>0]\leq \max_{w\in S}\mathbb{E}[L_w(S)]\leq c_{3}$. Hence 
 \[
 \mathbb{E}[L_y(S)]\leq c_{3}\mathbb{P}[L_y(S)>0].
 \]
 This together with \eqref{eq: lower bound 2 of occupation time for S} implies \eqref{eq: lower bound 1 of occupation time for S} with $c_1=\frac{1}{c_{3}}$. 
 
 Thus for $n\geq 1$ by conditioning on $\mathfrak{F}(x,\infty)$ one has that 
\begin{eqnarray}
\mathbb{E}[|T_x\cap L_{-n}(x)|]&=&
\sum_{y\in V}\mathbb{P}[y\in L_{-n}(x)]\mathbb{P}[y\leftrightarrow x]\nonumber\\
&\stackrel{\eqref{eq: two-point fcn bigger than nonzero occupation}}{\geq}&
\sum_{y\in V}\mathbb{P}[y\in L_{-n}(x)]\mathbb{P}[L_y(S)>0]=\sum_{y\in V}\mathbb{E}[\mathbf{1}_{\{ y\in L_{-n}(x) \}}]\cdot \mathbb{E}[\mathbf{1}_{\{L_y(S)>0\}}]\nonumber\\
&=&
\sum_{y\in V}\mathbb{E}\left[\mathbf{1}_{\{ y\in L_{-n}(x) \}}  \mathbf{1}_{\{L_y(S)>0\}}\right] \textnormal{ by independence}\nonumber\\
&\stackrel{\eqref{eq: lower bound 1 of occupation time for S}}{\geq }&c_{1}n.
\end{eqnarray}

The inequalities \eqref{eq: expectation of intersection of T_x with a high slab} follow from \eqref{eq: expectation of intersection of T_x with a low slab} and \eqref{eq: expectation relation between n and -n level}. 
 \end{proof}

 \begin{proposition}\label{prop: slab intersection for the future and past}
 	There exists a constant $c_{1}>0$ such that for every $n\geq0$,
 	\be\label{eq: expectation of intersection of future of x with a low slab}
 	1\leq \mathbb{E}[|\mathfrak{F}(x,\infty)\cap L_{-n}(x)|]\leq c_{1}.
 	\ee
 	There exists  constants $c_{2},c_{3}>0$ such that for every $n\geq1$,
 	\be\label{eq: expectation of intersection of past of x with a low slab}
 	c_{2}\leq\mathbb{E}[|\mathfrak{P}(x)\cap L_{-n}(x)|]\leq c_{3}e^{t_0}.
 	\ee
 	Moreover one also has that
 	\be\label{eq: expectation of intersection of future of x with a high slab}
 	c_{2}\exp(-t_0(n+1))\leq \mathbb{E}[|\mathfrak{F}(x,\infty)\cap L_{n}(x)|]\leq c_{3}\exp(-t_0n)
 	\ee
 	and
 	\be\label{eq: expectation of intersection of past of x with a high slab}
 		\exp(-t_0(n+1))\leq \mathbb{E}[|\mathfrak{P}(x)\cap L_{n}(x)|]\leq c_{1}\exp(-t_0(n-1)).
    \ee
 
 \end{proposition}
 \begin{proof}
 	We will prove \eqref{eq: expectation of intersection of future of x with a low slab}, the upper bound in \eqref{eq: expectation of intersection of future of x with a high slab} and the lower bound in \eqref{eq: expectation of intersection of past of x with a low slab}. The rest will follow from these inequalities and the tilted mass-transport principle.

 	If we sample $\mathfrak{F}$ starting with $X^x$, then $\mathfrak{F}(x,\infty)\subset \{X^x_k:k\geq0\}$. Write $\mathfrak{F}(x,\infty)=(v_0=x,v_1,v_2,\ldots)$. By \eqref{eq: LLN for srw} one has that $\log\Delta(X^x_0,X^x_k)\rightarrow-\infty$ as $k\rightarrow\infty$, in particular $\log\Delta(v_0,v_k)\rightarrow-\infty$. Thus $1\leq |\mathfrak{F}(x,\infty)\cap L_{-n}(x)|$ for every $n\geq0$ almost surely, whence one has the lower bound in \eqref{eq: expectation of intersection of future of x with a low slab}: $1\leq \mathbb{E}[|\mathfrak{F}(x,\infty)\cap L_{-n}(x)|]$. 
 	
 	 The upper bound in \eqref{eq: expectation of intersection of future of x with a low slab} also follows from the above inclusion $\mathfrak{F}(x,\infty)\subset \{X^x_k:k\geq0\}$.
 	\[
 	\mathbb{E}[|\mathfrak{F}(x,\infty)\cap L_{-n}(x)|]\leq 
 	\mathbb{E}\left[\sum_{m=0}^{\infty}\mathbf{1}_{ \{X_m^x\in L_{-n}(x) \} }\right]\stackrel{\eqref{eq: occupation times for srw on -n slab}}{\leq }c_0=c_{1}.
 	\]
 	
 	Next we prove the upper bound in \eqref{eq: expectation of intersection of future of x with a high slab}. By the reversibility of simple random walk and the transitivity of $G$, one has $ \mathbb{P}[\sigma_y^x<\infty]=\mathbb{P}[\sigma_x^y<\infty]$. Hence for $n\geq1$,
 	\begin{eqnarray}\label{eq: intersection of a high slab for v-WUSF component}
 	 	\mathbb{E}\left[\sum_{y\in V}\mathbf{1}_{ \{y\in L_n(x), \sigma_y^x<\infty \} } \right]
 	 	&\stackrel{\textnormal{TMTP}}{=}&
 	 	\mathbb{E}\left[\sum_{y\in V}\mathbf{1}_{ \{x\in L_n(y), \sigma_x^y<\infty \} }\Delta(x,y) \right]\nonumber\\
 	 	&\asymp&e^{-t_0n}\mathbb{E}\left[\sum_{y\in L_{-n}(x)}\mathbf{1}_{ \{ \sigma_y^x<\infty \}}\right],
 	\end{eqnarray}
 	where the last equality holds up a multiplicative constant $e^{\pm t_0}$.
 	
 	Since $\mathbb{P}[y\in\mathfrak{F}(x,\infty)]\leq \mathbb{P}[\sigma_y^x<\infty]$, one has that 
 	\[
 	\mathbb{E}[|\mathfrak{F}(x,\infty)\cap L_{n}(x)|]\leq 
 		\mathbb{E}\left[\sum_{y\in V}\mathbf{1}_{ \{y\in L_n(x), \sigma_y^x<\infty \} } \right]
 		\stackrel{\eqref{eq: intersection of a high slab for v-WUSF component}}{\leq} e^{t_0-t_0n}\mathbb{E}\left[\sum_{m=0}^{\infty}\mathbf{1}_{ \{X_m^x\in L_{-n}(x) \} }\right]\leq c_0e^{t_0-t_0n}.
 	\]
 	Taking $c_{3}=c_0e^{t_0}$ one has the upper bound in \eqref{eq: expectation of intersection of future of x with a high slab}.

 	Finally we prove the lower bound in \eqref{eq: expectation of intersection of past of x with a low slab}.

 Let $y\in V$ and let $X^y$ be a simple random walk started at $y$. Sample  $\mathfrak{F}$ starting with the simple random walk $X^y$. Let $\mathscr{A}(y,x)$ be the event that $\sigma_x^y<\infty$, that the sets $\{X^y_m:0\leq m<\sigma_x^y\}$ and $\{X_m^y:m\geq \sigma_x^y \}$ are disjoint, so that $y\in\mathfrak{P}(x)$ on the event $\mathscr{A}(y,x)$ and hence 
 \[
 \mathbb{E}[|\mathfrak{P}(x)\cap L_{-n}(x) |]\geq \sum_{y\in V}\mathbb{P}[ y\in L_{-n}(x)]\cdot \mathbb{P}[\mathscr{A}(y,x)  ]
 \]
 	Let $Y^x$ be an independent simple random walk also started at $x$. We also use $X^x$ to denote the set of vertices $\{X_m^x:m\geq0\}$ visited by the random walk $X^x$ and write similarly $Y^x_+:=\{Y_m^x:m\geq1 \}$. 	By time-reversal one has that 
 	\[
 	\mathbb{P}[\mathscr{A}(y,x)]\geq \mathbb{P}[\sigma_y^x<\infty,X^x\cap Y^x_+=\emptyset].
 	\]
 	Since $X^x$ hits $L_{-n}(x)$ for every $n\geq0$ almost surely, one has that 
 	\[
 	 \mathbb{E}[|\mathfrak{P}(x)\cap L_{-n}(x) |]\geq \sum_{y\in V}\mathbb{P}[ y\in L_{-n}(x)]\cdot \mathbb{P}[\sigma_y^x<\infty,X^x\cap Y^x_+=\emptyset]\stackrel{\eqref{eq: occupation times for srw on -n slab}}{\geq} \mathbb{P}[X^x\cap Y^x_+=\emptyset].
 	\]
 	Since $G$ is a nonamenable transitive graph, there exists a positive constant $
 	c_{2}$ such that $\mathbb{P}[X^x\cap Y^x_+=\emptyset]\geq c_{2}$ (see \cite[Theorem 10.24]{LP2016}). Thus we have the lower bound in \eqref{eq: expectation of intersection of past of x with a low slab}.

 	Just like \eqref{eq: expectation relation between n and -n level}, the tilted mass-transport principle implies that for any $n\in\mathbb{Z}$,
 	\[
 	\mathbb{E}[|\mathfrak{P}(x)\cap L_n(x)|]=\mathbb{E}\left[\sum_{y\in V}\mathbf{1}_{\{ y\in L_{-n}(x)\}} \mathbf{1}_{\{y\in\mathfrak{F}(x,\infty) \}} \Delta(x,y) \right]\asymp
 	e^{-t_0n}\mathbb{E}[|\mathfrak{F}(x,\infty)\cap L_{-n}(x)|],
 	\]
 	where the last equality holds up to a factor $e^{\pm t_0}$.
 	This together with \eqref{eq: expectation of intersection of future of x with a low slab} implies \eqref{eq: expectation of intersection of past of x with a high slab}; this together with the upper bound in \eqref{eq: expectation of intersection of future of x with a high slab} implies the upper bound in \eqref{eq: expectation of intersection of past of x with a low slab} and this together with the lower bound in \eqref{eq: expectation of intersection of past of x with a low slab} implies the lower bound in \eqref{eq: expectation of intersection of future of x with a high slab} 
  \end{proof}

 Next we give estimates of $\mathbb{E}[|\mathfrak{T}_x\cap L_{n}(x)|]$.
 
 \begin{proposition}\label{prop: slab intersection for the x-WUSF component}
 	There exists positive constants $c_{1},c_{2}$ such that for all $n\geq0$,
 	\be\label{eq: expectation of intersection of x-WUSF component with a low slab}
 	c_{1}\leq \mathbb{E}[|\mathfrak{T}_x\cap L_{-n}(x)|]\leq c_{2}
 	\ee
 	and 
 	\be\label{eq: expectation of intersection of x-WUSF component with a high slab}
 		c_{1}\exp(-t_0(n+1))\leq \mathbb{E}[|\mathfrak{T}_x\cap L_{n}(x)|]\leq c_{2}\exp(-t_0(n-1)).
 	\ee
 \end{proposition}
 \begin{proof}
 As in the proof of Lemma \ref{lem: upper bound on two-point functions},  by  Wilson's algorithm and the reversibility of simple random walk one has that
 \[
 \mathbb{P}[y\in\mathfrak{T}_x]=\mathbb{P}[\sigma_x^y<\infty]=\mathbb{P}[\sigma_y^x<\infty],
 \]
 where $\sigma_x^y$ is the first visit time of $x$ by the simple random walk $X^y$. 
 
 Hence 
 \begin{eqnarray}\label{eq: intersection of a low slab for x-WUSF in terms of srw}
 \mathbb{E}[|\mathfrak{T}_x\cap L_{-n}(x)|]
 &=&\mathbb{E}\left[\sum_{y\in V} \mathbf{1}_{ \{ y\in L_{-n}(x) \} }\cdot\mathbf{1}_{ \{ y\in\mathfrak{T}_x   \} }\right]\nonumber\\
 &=&\mathbb{E}\left[\sum_{y\in V} \mathbf{1}_{ \{ y\in L_{-n}(x) \} }\cdot\mathbf{1}_{ \{ \sigma_y^x<\infty   \} }\right]
 \end{eqnarray}
 Since simple random walk $X^x$ visits $L_{-n}(x)$ for every $n\geq0$ almost surely, 
one has  $\sum_{y\in V} \mathbf{1}_{ \{ y\in L_{-n}(x) \} }\cdot\mathbf{1}_{ \{ \sigma_y^x<\infty   \} }\geq 1\ \as$. Then 
 \eqref{eq: intersection of a low slab for x-WUSF in terms of srw} implies that the lower bound in \eqref{eq: expectation of intersection of x-WUSF component with a low slab} holds with $c_{1}=1$.
 
 On the other hand, 
 \[
 \mathbb{E}[|\mathfrak{T}_x\cap L_{-n}(x)|]\stackrel{\eqref{eq: intersection of a low slab for x-WUSF in terms of srw}}{\leq }\mathbb{E}\left[\sum_{m=0}^{\infty}\mathbf{1}_{ \{ X_m^x\in L_{-n}(x)  \}} \right]\stackrel{\eqref{eq: occupation times for srw on -n slab}}{\leq }c_0.
 \]
 Hence the upper bound in \eqref{eq: expectation of intersection of x-WUSF component with a low slab} holds with $c_{2}=c_0$. 
 
 By Lemma \ref{lem: upper bound on two-point functions} one has $\mathbb{P}[y\in\mathfrak{T}_x]=\mathbb{P}[x\in\mathfrak{T}_y]$. Then the tilted mass-transport principle gives the relation between $\mathbb{E}[|\mathfrak{T}_x\cap L_{-n}(x)|]$ and $\mathbb{E}[|\mathfrak{T}_x\cap L_{n}(x)|]$, namely \eqref{eq: intersection of a high slab for v-WUSF component}. This together with \eqref{eq: expectation of intersection of x-WUSF component with a low slab}   yields \eqref{eq: expectation of intersection of x-WUSF component with a high slab}.
 \end{proof}

 For the intersection of a simple random walk trajectory with a slab, we have the following results.
 \begin{proposition}\label{prop: intersection of a srw and a slab}
 	Let $X^x$ denote a simple random walk on the transitive graph $G$ started from $x$. Then  for $n\geq0$ one has that
 	\be\label{eq: prob of srw intersecting a high slab with at least k vertices}
 	\mathbb{P}[|\{X_m^x:m\geq0 \}\cap L_n(x)|\geq k]=e^{-\Theta(k)-t_0n}
 	\ee
 	and 
 	\be\label{eq: prob of srw intersecting a low slab with at least k vertices}
 		\mathbb{P}[|\{X_m^x:m\geq0 \}\cap L_{-n}(x)|\geq k]=e^{-\Theta(k)}
 	\ee
 	In particular, one has that
 	\be\label{eq: prob of srw intersecting a high slab}
 	\mathbb{P}[\{X_m^x:m\geq0 \}\cap L_n(x)\neq\emptyset]\asymp \exp(-t_0n),n\geq0.
 	\ee
 	Also a simple consequence is that there exists a constant $c>0$ such that
 	\be\label{eq: k-moment of intersection of srw with a  slab}
 	\mathbb{E}[|\{X_m^x:m\geq0 \}\cap L_n(x)|^k]\leq
 	\left\{\begin{array}{cc}
 	c^kk!\exp(-t_0n)& \textnormal{ if }\,n\geq0\\
 	c^kk!& \textnormal{ if }\,n<0
 	\end{array}\right..
 	\ee
 	
 \end{proposition}
 \begin{proof}
 	Note that $\mathbb{P}[\sigma_y^x<\infty]=\mathbb{P}[\sigma_x^y<\infty]$ and
 	\[
 	|\{X_m^x:m\geq0 \}\cap L_n(x)|=\sum_{y\in L_n(x)}\mathbf{1}_{\{\sigma_y^x<\infty\}}.
 	\]
 	Hence 
 	\[
 	\mathbb{E}[|\{X_m^x:m\geq0 \}\cap L_n(x)|]=\mathbb{E}\Big[ \sum_{y\in L_n(x)}\mathbf{1}_{\{\sigma_x^y<\infty\}}\Big]=\mathbb{E}[|\mathfrak{T}_x\cap L_n(x)|]
 	\]
 	By Proposition \ref{prop: slab intersection for the x-WUSF component} one has that 
 	$	\mathbb{E}[|\mathfrak{T}_x\cap L_n(x)|]	\preceq \exp(-t_0n)$
 	and then
 	\be\label{eq: upper bound for prob of srw intersecting a high slab}
 		\mathbb{P}[\{X_m^x:m\geq0 \}\cap L_n(x)\neq\emptyset]
 		\leq 	\mathbb{E}[|\{X_m^x:m\geq0 \}\cap L_n(x)|]
 		\asymp \exp(-t_0n).
 	\ee
 	This gives the upper bound of \eqref{eq: prob of srw intersecting a high slab}.
 	

 	Write $N_n=\sum_{m=0}^{\infty}\mathbf{1}_{ \{X_m^x\in L_n(x)  \} }$ for $n\in\mathbb{Z}$ as in the proof of Corollary \ref{cor:large deviation for SRW on nonunimodular transitive graph}. Similar to \eqref{eq: strong Markov property for SRW}, by the strong Markov property of simple random walk  there exists a constant $c\in (0,1)$ such that 
 	\be\label{eq: N_n geq k conditioned on N_n big than 0}
 	\mathbb{P}[N_n\geq k|N_n>0]\leq (1-c)^{k-1},k\geq 1.
 	\ee
 	In particular,
 	\be\label{eq: k-th moments of N_n conditioned on N_n bigger than 0}
 	\mathbb{E}[N_n^k|N_n>0]\leq \int_{0}^{\infty}ky^{k-1}(1-c)^{y-1}dy=\frac{1}{1-c}\cdot \frac{k!}{(-\log(1-c) )^k}
 	\ee
 	
 	Notice that $|\{X_m^x:m\geq0 \}\cap L_n(x)|\leq N_n$ and $\mathbb{P}[N_n>0]=\mathbb{P}[\{X_m^x:m\geq0 \}\cap L_n(x)\neq\emptyset]$.
 	Therefore we have the upper bound:
 	\begin{eqnarray}\label{eq: second moment for SRW intersect with a high slab}
 		\mathbb{E}[|\{X_m^x:m\geq0 \}\cap L_n(x)|^2]&\leq & \mathbb{P}[N_n>0]\mathbb{E}[N_n^2|N_n>0]\stackrel{\eqref{eq: k-th moments of N_n conditioned on N_n bigger than 0}}{\preceq}\mathbb{P}[N_n>0]\nonumber\\
 		&=& \mathbb{P}[\{X_m^x:m\geq0 \}\cap L_n(x)\neq\emptyset]\stackrel{\eqref{eq: upper bound for prob of srw intersecting a high slab}}{\preceq} \exp(-t_0n).
 	\end{eqnarray}
 	
 	By second moment method one has the lower bound of \eqref{eq: prob of srw intersecting a high slab}, namely 
 	\[
 	\mathbb{P}[\{X_m^x:m\geq0 \}\cap L_n(x)\neq\emptyset]
 	\geq \frac{\big(\mathbb{E}[|\{X_m^x:m\geq0 \}\cap L_n(x)|]\big)^2}{\mathbb{E}[|\{X_m^x:m\geq0 \}\cap L_n(x)|^2]}
 	\stackrel{\eqref{eq: upper bound for prob of srw intersecting a high slab}, \eqref{eq: second moment for SRW intersect with a high slab}}{\succeq} \exp(-t_0n).
 	\]
 	
 Hence for $n\geq 0$, by \eqref{eq: N_n geq k conditioned on N_n big than 0} and \eqref{eq: prob of srw intersecting a high slab} there exists a constant $c_1>0$ such that 
 \be\label{eq: tail probability for srw intersection with a high slab}
 \mathbb{P}[|\{X_m^x:m\geq0 \}\cap L_{n}(x)|\geq k]\leq \mathbb{P}[N_n\geq k]\leq e^{-t_0n-c_1k}.
 \ee
 	For $n\geq0$, the simple random walk $X^x$ hits every $L_{-n}(x)$ almost surely. Thus $\mathbb{P}[N_{-n}>0]=1$. Hence  by \eqref{eq: N_n geq k conditioned on N_n big than 0}
 	\be\label{eq: tail probability for srw intersection with a low slab}
 	\mathbb{P}[|\{X_m^x:m\geq0 \}\cap L_{-n}(x)|\geq k]\leq \mathbb{P}[N_{-n}\geq k]\leq e^{-c_1k}.
 	\ee
 	Recall for a nonnegative random variable $Y$ and $p>0$, $\mathbb{E}[Y^p]=\int_{0}^{\infty}py^{p-1}\mathbb{P}[Y>y]dy$ (for example see  Lemma 2.2.13 in \cite{Durrett2010}). This fact and estimates (\ref{eq: tail probability for srw intersection with a high slab}, \ref{eq: tail probability for srw intersection with a low slab}) imply \eqref{eq: k-moment of intersection of srw with a  slab}.
 	
 It remains to show that  there exists $c_2>0$ such that for all $n\in\mathbb{Z},k\geq1$
 	\be\label{eq: exponential decay lower bound in k for srw intersect with a slab}
 	\mathbb{P}[|\{X_m^x:m\geq0 \}\cap L_{n}(x)|\geq k]\geq 
 	\left\{
 	\begin{array}{cc}
 	e^{-t_0n-c_2k} & \textnormal{ if }n\geq0\\
 	e^{-c_2k} & \textnormal{ if }n<0\\
 	\end{array}
 	\right..
 	\ee
 	and we defer its proof to the appendix.
  \end{proof}

  \begin{corollary}\label{cor: high moments for the future intersect with a slab}
  	For the future $\mathfrak{F}(x,\infty)$, there exists a constant $c>0$ such that 
  	\[
  	\mathbb{E}[|\mathfrak{F}(x,\infty)\cap L_{n}(x)|^k]\leq 
  	\left\{
  	\begin{array}{cc}
  	c^kk!e^{-t_0n} & \textnormal{ if }n\geq0\\
  	c^kk!& \textnormal{ if }n<0
  	\end{array}
  	\right.
  	\]
  \end{corollary}

 Next we extend Corollary \ref{cor: expectation for WUSF on regular tree} to all nonunimodular transitive graphs. 
 \begin{proposition}\label{prop: expected tilted volume for WUSF component}
 	The expected tilted volume $\mathbb{E}[|T_x|_{x,\lambda}]$ is finite if and only if $\lambda\in (0,1)$. 
 \end{proposition}
 \begin{proof}
 	The proof is the same as Corollary \ref{cor: expectation for WUSF on regular tree}. One just need to replace  \eqref{eq: expectation of intersection of T_x with a slab for the toy model}  with the corresponding estimates \eqref{eq: expectation of intersection of T_x with a low slab} and \eqref{eq: expectation of intersection of T_x with a high slab} in the general case. 
 \end{proof}

 Similarly Proposition \ref{prop: slab intersection for the x-WUSF component} and \ref{prop: slab intersection for the future and past} yield the following proposition. We omit its proof. 
 
 \begin{proposition}\label{prop: tilted volume bounds for x-component, past and future}
 	The expected tilted volume $\mathbb{E}[|\mathfrak{T}_x|_{x,\lambda}]$ is finite if and only if $\lambda\in (0,1)$.
 Similarly,
 	\begin{itemize}
 		\item $\mathbb{E}[|\mathfrak{F}(x,\infty)|_{x,\lambda}]<\infty$ is finite if and only if $\lambda\in (0,1)$.

 		\item $\mathbb{E}[|\mathfrak{P}(x)|_{x,\lambda}]<\infty$ is finite if and only if $\lambda\in (0,1)$.

 	\end{itemize}
 	
 \end{proposition}

 Now we are ready to prove Theorem \ref{thm:light for WUSF}. In fact it is just the  $\lambda=1$  case in the following proposition. 
 \begin{proposition}\label{prop: when the tilted volume is finite}
 	If $\lambda\leq 0$, then $|T_x|_{x,\lambda}=\infty \ \as$
 	If $\lambda>0$, then $|T_x|_{x,\lambda}<\infty\  \as$  
 \end{proposition}
 \begin{proof}
 	
 	Let $\{X_n\}_{n\geq0}$ be a simple random walk started at $x$. If we sample $\mathfrak{F}$ starting with $\{X_n\}_{n\geq0}$ and let $\mathfrak{F}(x,\infty)=(x_0,x_1,x_2,\ldots)$ be future of $x$. Then Corollary \ref{cor:large deviation for SRW on nonunimodular transitive graph} implies that $\Delta(x_0,x_n)\rightarrow0$ almost surely. 
 	Thus for $\lambda\leq0 $, $|T_x|_{x,\lambda}\geq |\mathfrak{F}(x,\infty)|_{x,\lambda}=\infty$.
 	
 	For $\lambda>0$, use the same proof as Proposition \ref{prop: lightness of trees in the toy model} with Proposition  \ref{prop: expected tilted volume for WUSF component} replacing the role of Corollary \ref{cor: expectation for WUSF on regular tree}.	
 \end{proof}
  Moreover \eqref{eq:bounds for tail prob of general case} in Proposition \ref{prop: tail prob for general case} gives a quantitative tail bound of the weight of the tree $T_x$ in the WUSF. 

 \begin{remark}\label{rem: update-tolerance proof of the lightness}
 	Given Proposition \ref{prop: slab intersection for T_x} there is also another way of proving Theorem \ref{thm:light for WUSF}. We just sketch the idea here.  One can use the update-tolerance \cite{Tom2016wired} and TMTP to show that if $T_x$ is heavy with positive probability, then $T_x$ will intersect $L_n(x)$ with infinitely many vertices with positive probability, which contradicts the fact that $\mathbb{E}[|T_x\cap L_n(x)|]<\infty$.
 \end{remark}

 \subsection{High moments}

 We shall establish some upper bounds  on high moments of $|\mathfrak{T}_x\cap L_n(x)|$ and $|T_x\cap L_n(x)|$. The high moments for $|\mathfrak{F}(x,\infty)\cap L_n(x)|$ have already been given in Corollary \ref{cor: high moments for the future intersect with a slab}. The high moments for the past $|\mathfrak{P}(x)\cap L_n(x)|$ can be obtained from the corresponding result for $\mathfrak{T}_x$ using Lemma \ref{lem:stochastic domination}.  See Table \ref{table: high moments} for a summary of these results.

  \begin{table}[ht]
 	\centering
 	\begin{tabular}{|c|c|c|}
 		\hline
 		Quantities $(n\in \mathbb{Z},k\geq2)$ & Results & Position\\
 			\hline
 		$\mathbb{E}[|\mathfrak{F}(x,\infty)\cap L_{n}(x)|^k]$ &  $\left\{
 		\begin{array}{cc}
 		\leq c^kk! e^{-t_0n} & \textnormal{ if }n\geq0\\
 		\leq c^kk! & \textnormal{ if }n<0\\
 		\end{array}
 		\right.$& Corollary \ref{cor: high moments for the future intersect with a slab} \\
 		\hline
 		$\mathbb{E}[|\mathfrak{T}_x\cap L_n(x)|^k]$ & $\left\{
 		\begin{array}{cc}
 		\leq c^k(k!)^2e^{-t_0n} & \textnormal{ if }n\geq0\\
 		\leq c^k(k!)^2|n|^{k-1}& \textnormal{ if }n<0\\
 		\end{array}
 		\right.$ & Proposition \ref{prop: high moments for x-component intersecting a slab}\\
 		\hline
 		$\mathbb{E}[|\mathfrak{P}(x)\cap L_n(x)|^k]$ &  $\left\{
 		\begin{array}{cc}
 		\leq c^k(k!)^2e^{-t_0n} & \textnormal{ if }n\geq0\\
 		\leq c^k(k!)^2|n|^{k-1}& \textnormal{ if }n<0\\
 		\end{array}
 		\right.$ & Corollary \ref{cor: high moments for the past} \\
 		\hline
 		$\mathbb{E}[|T_x\cap L_{n}(x)|^k]$ & $\left\{
 		\begin{array}{cc}
 		\leq c^k(k!)^5 (n\vee1)e^{-t_0n} & \textnormal{ if }n\geq0\\
 		\leq c^k(k!)^5|n|^{2k}& \textnormal{ if }n<0\\
 		\end{array}
 		\right.$ & Proposition \ref{prop: high moments for T_x intersecting a slab}\\
 		\hline
 	\end{tabular}
 \caption{High moments for the intersections with a slab}
 \label{table: high moments}
 \end{table}

 Write $\tau(x,y):=\mathbb{P}[y\in\mathfrak{T}_x]$. By transitivity one has $\tau(x,y)=\mathbb{P}[\sigma_x^y<\infty]=\mathbb{P}[\sigma_y^x<\infty]=\tau(y,x)$, where $\sigma_x^y$ is hitting time of $x$ by a simple random walk started from $y$. In the following we will use the convention $\tau(x,x)=1$. Write $\tau(x_0,x_1,\ldots,x_k):=\mathbb{P}\big[\bigcap_{i=1}^{k}\{x_i\in \mathfrak{T}_{x_0}\}\big]$. The following lemma is an analogue of Lemma 6.89 in \cite{Grimmett1999}. 
 \begin{lemma}\label{lem: point function ineq with k=3}
 	For all vertices $x_0,x_1,x_2$ we have 
 	\be\label{eq: point function ineq with k=3}
 	\tau(x_0,x_1,x_2)\leq \sum_{u\in V} \tau(x_0,u)\tau(u,x_1)\tau(u,x_2)
 	\ee
 \end{lemma}
 \begin{proof}
 	Rooted the tree $\mathfrak{T}_{x_0}$ at $x_0$. Let $\gamma_i$ be the path on the tree from $x_0$ to $x_i$ for $i=1,2$. Let $s(x_0,x_1,x_2)$ be the last vertex on the path $\gamma_1\cap \gamma_2$. If we sample $\mathfrak{T}_{x_0}$ from independent simple random walks started from $u,x_1,x_2,\ldots$, then we know that 
 	\be
 	\mathbb{P}\big[x_1,x_2\in\mathfrak{T}_{x_0}\textnormal{ and }s(x_0,x_1,x_2)=u \big]\leq
 	\tau(x_0,u)\tau(u,x_1)\tau(u,x_2)
 	\ee
 	Notice  if $x_1=x_2$, then $s(x_0.x_1,x_2)=x_1=x_2$ the above inequality is still true. If  $x_1=x_0\neq x_2$, then $s(x_0,x_1,x_2)=x_0=x_1$ and so the above inequality is still true. 
 	Summing this inequality over all the possible choice of $s(x_0,x_1,x_2)$ one obtains the desired inequality \eqref{eq: point function ineq with k=3}.
 \end{proof}
\begin{observation}\label{obse: relative level relation}
If $u\in L_m(x)$ and $x_1\in L_n(x)$, then \[
	x_1\in\bigcup_{j=n-m-1}^{n-m+1}L_m(u).
	\]
	In other words,
	\[
	L_n(x)\subset \bigcup_{j=n-m-1}^{n-m+1}L_m(u).
	\]
	Also by Proposition \ref{prop: slab intersection for the x-WUSF component} one has 
	\be\label{eq: relative expectation when x_1 runs over L_n(x)}
	\mathbb{E}\Big[\sum_{x_1\in L_n(x)}{\tau(u,x_1)}\bigm| U\Big]\asymp \exp(-t_0(n-m)\vee 0)
	\ee
\end{observation}
 \begin{proposition}\label{prop: high moments for x-component intersecting a slab}
 	There exists $c_{1}>0$ such that for all $n\geq0,k\geq 2$, one has 
 	\be\label{eq: second moment of x-component intersecting a high slab}
 	\mathbb{E}[|\mathfrak{T}_x\cap L_n(x)|^k]\leq c_{1}^k(k!)^{2}\exp(-t_0n)
 	\ee
 	and 
 	\be\label{eq: second moment of x-component intersecting a low slab}
 		\mathbb{E}[|\mathfrak{T}_x\cap L_{-n}(x)|^k]\leq c_{1}^k(k!)^{2} (n\vee 1)^{k-1}.
 	\ee
 \end{proposition}

 \begin{proof}	
Here we only prove the case $k=2$. The proof of the general case will be given in the appendix.

Note $|\mathfrak{T}_x\cap L_n(x)|^2=\sum_{x_1,x_2\in L_n(x)}{\tau(x,x_1,x_2)}$. Using Lemma \ref{lem: point function ineq with k=3} one has 
\begin{eqnarray}
\mathbb{E}[|\mathfrak{T}_x\cap L_n(x)|^2]&=&\mathbb{E}\bigg[\sum_{x_1,x_2\in L_n(x)}{\tau(x,x_1,x_2)}\bigg]\nonumber\\
&\leq&\mathbb{E}\bigg[\sum_{u\in V,x_1,x_2\in L_n(x)}{\tau(x,u)\tau(u,x_1)\tau(u,x_2)}\bigg]\nonumber\\
&\leq&\sum_{j=-\infty}^{\infty}\mathbb{E}\bigg[\sum_{u\in L_j(x),x_1,x_2\in L_n(x)}{\tau(x,u)\tau(u,x_1)\tau(u,x_2)}\bigg]\nonumber\\
&\stackrel{\eqref{eq: relative expectation when x_1 runs over L_n(x)} }{\asymp }&
\sum_{j=-\infty}^{\infty}\exp(-t_0(j\vee0))\cdot \exp(-2t_0(n-j)\vee0)
\end{eqnarray}

Hence for $n\geq0$, 
one has 
\[
\mathbb{E}[|\mathfrak{T}_x\cap L_n(x)|^2]\preceq \sum_{j=-\infty}^{0}e^{-2t_0n}e^{2t_0j}+
\sum_{j=0}^{n}e^{-2t_0n}e^{t_0j}+\sum_{j=n+1}^{\infty}e^{-t_0j}\asymp e^{-t_0n}.
\]
Similarly for $n\leq-1$, one can compute that 
\[
\mathbb{E}[|\mathfrak{T}_x\cap L_{-n}(x)|^2]\preceq |n|.\qedhere
\]
 \end{proof}
 
 An immediate corollary of Proposition \ref{prop: high moments for x-component intersecting a slab} is the upper bound on $|\mathfrak{P}(x)\cap L_n(x)|$ using Lemma \ref{lem:stochastic domination}. 
 \begin{corollary}\label{cor: high moments for the past}
 	There exists a constant $c>0$ such that for all $n\geq0,k\geq2$, 
 	\[
 	\mathbb{E}[|\mathfrak{P}(x)\cap L_n(x)|^k]\leq \left\{
 	\begin{array}{cc}
 	c^k(k!)^2e^{-t_0n} &\textnormal{ if }n\geq0\\
 		c^k(k!)^2|n|^{k-1} & \textnormal{ if }n<0
 	\end{array}
 	\right..
 	\]
 \end{corollary}
 
 \begin{corollary}\label{cor: tail probability for x-wusf intersect with L_n with many vertices}
 	There exist  constants $c_1,c_2,c_3,c_4>0$ such that for all $n\geq0,k\geq1$,
 	\be\label{eq: stretched exponential tail for x-wusf intersect with L_n with many vertices}
 	\mathbb{P}[|\mathfrak{T}_x\cap L_n(x)|\geq k]\leq c_2e^{-t_0n}\cdot e^{-c_1\sqrt{k}}
 	\ee
 	and for all $n>0,k\geq1$
 	\be\label{eq: stretched exponential tail for x-wusf intersect with L_(-n) with many vertices}
 	\mathbb{P}[|\mathfrak{T}_x\cap L_{-n}(x)|\geq k]\leq \frac{c_4}{n}e^{-c_3\sqrt{\frac{k}{n}}}
 	\ee
 	
 \end{corollary}
\begin{proof}
	Taking $c_1=\frac{1}{2\sqrt{c}}$, then
	\begin{eqnarray}
	\mathbb{E}\big[\exp\big(c_1\sqrt{|\mathfrak{T}_x\cap L_n(x)|}\big)\big]&\leq 
	&\mathbb{E}\big[\exp\big(c_1\sqrt{|\mathfrak{T}_x\cap L_n(x)|}\big)\big]+\mathbb{E}\big[\exp\big(-c_1\sqrt{|\mathfrak{T}_x\cap L_n(x)|}\big)\big]\nonumber\\
	&=&2\sum_{k=0}^{\infty}\frac{c_1^{2k}}{(2k)!}\mathbb{E}\big[|\mathfrak{T}_x\cap L_n(x)|^k\big]\nonumber\\
	&\stackrel{\eqref{eq: second moment of x-component intersecting a high slab}}{\leq}&
	2\sum_{k=0}^{\infty}\frac{c_1^{2k}}{(2k)!} c^k \cdot (k!)^2\cdot e^{-t_0n}\nonumber\\
	&\leq&2\sum_{k=0}^{\infty}(c_1^2c)^k\cdot e^{-t_0n}<  3e^{-t_0n}.
	\end{eqnarray}
	Hence by Markov's inequality  for all $n\geq0,k\geq1$ one has that 
	\begin{eqnarray*}
	\mathbb{P}[|\mathfrak{T}_x\cap L_n(x)|\geq k]&=&\mathbb{P}[\exp\big(c_1\sqrt{|\mathfrak{T}_x\cap L_n(x)|}\big)\geq \exp(c_1\sqrt{k})]\\
	&\leq& \frac{3e^{-t_0n}}{\exp(c_1\sqrt{k})}\leq c_2e^{-t_0n}\cdot e^{-c_1\sqrt{k}}.\qedhere
	\end{eqnarray*}
The proof of \eqref{eq: stretched exponential tail for x-wusf intersect with L_(-n) with many vertices} is similar. For small $c_3>0$, one has that  $\mathbb{E}[\exp(c_3\sqrt{\frac{|\mathfrak{T}_x\cap L_{-n}(x)|}{n}})]\leq \frac{c_4}{n}$ and then use Markov's inequality. 
\end{proof}
\begin{question}
	Is it the case that for all $n\geq0,k\geq1$,
	\[
		\mathbb{P}[|\mathfrak{T}_x\cap L_n(x)|\geq k]=e^{-t_0n-\Theta(\sqrt{k})}?
	\]
	The above Corollary \ref{cor: tail probability for x-wusf intersect with L_n with many vertices} establishes the upper bound. As we will see in the appendix, the corresponding lower bound is also true for the toy model. 
\end{question}
 
 Using the estimates for $\mathbb{E}[|\mathfrak{T}_x\cap L_n(x)|^k]$ we can derive upper bounds for $\mathbb{E}[|T_x\cap L_n(x)|^k]$. Recall $U=\{U_v:v\in V\}$ are the labels we used to define the slabs. 
 
 \begin{lemma}\label{lem: bounds high moments of T_x intersections in terms of the x-WUSF ones}
 	Suppose $\mathfrak{F}(x,\infty)=(x_0,x_1,x_2,\ldots)$ is the future of $x$ in the WUSF sample $\mathfrak{F}$. Let $\mathfrak{T}_{x_i}$ be the tree of $x_i$ in  $x_i$-WUSF, sampled independently of  $\mathfrak{F}$, $U$ and each other. Then for every $n\in\mathbb{Z},k\geq 2$, one has
 	\be\label{eq: k-moments bounds between T_x and x-WUSF}
 	\mathbb{E}\big[|T_x\cap L_n(x)|^k\bigm|\mathfrak{F}(x,\infty),U\big]
 	\leq k! \sum_{\substack{k_0,k_1,\ldots \geq 0:\\ k_0+ k_1+\cdots  =k}} \, \prod_{i:k_i\neq 0} \frac{1}{k_i!}\mathbb{E}\big[|\mathfrak{T}_{x_i}\cap L_n(x)|^{k_i}\bigm|U\big]. 
 	\ee
 \end{lemma}
 \begin{proof}
 	The proof is similar to the one of Lemma 6.8 in \cite{Hutchcroft2018a}. We present the details for reader's convenience. In the following we fix $n\in\mathbb{Z}$.
 	
 	Given the future $\mathfrak{F}(x,\infty)=(x_0,x_1,x_2,\ldots)$ of $x$ and $i\geq0$, we call the connected component of $x_i$ in $T_x\backslash\{x_0,x_1,\ldots,x_{i-1},x_{i+1},x_{i+2},\ldots\}$ the $i$-th \notion{bush} of $T_x$ and denote it by $\textnormal{Bush}_i(x)$. 
 	Denote by $N_i$ the number of vertices in $L_n(x)\cap \textnormal{Bush}_i(x)$. In particular, $\textnormal{Bush}_0(x)=\mathfrak{P}(x)$.  By the lightness of $T_x$, almost surely only finitely many $N_i$'s are nonzero. 
 	
 	Notice that 
 	\be\label{eq: lem4.15-1}
 	|T_x\cap L_n(x)|^k=(N_0+N_1+\cdots)^k=\sum_{\substack{k_0,k_1,\ldots \geq 0:\\ k_0+ k_1+\cdots  =k}} \,k!\prod_{i:k_i\neq 0} \frac{N_i^{k_i}}{k_i!}
 	\ee
 	
 	Since $N_i=\sum_{y\in V}\mathbf{1}_{ \{ y\in L_n(x)\cap \textnormal{Bush}_i(x) \} }$, 
 	$N_i^{k_i}=\sum_{y_j: 1\leq j\leq k_i}\mathbf{1}_{ \{y_j\in L_n(x)\cap \textnormal{Bush}_i(x) \textnormal{ for }j=1,\ldots,k_i  \} }$ for $k_i>0$ and then for $k_i>0$
 	\[
 	\mathbb{E}[N_i^{k_i}\mid \mathfrak{F}(x,\infty),U]=\sum_{y_j: 1\leq j\leq k_i}\mathbb{P}[\bigcap_{1\leq j\leq k_i} y_j\in L_n(x)\cap \textnormal{Bush}_i(x) \mid \mathfrak{F}(x,\infty),U ]. 
 	\]
 
 Similarly for any  sequence $(k_i)_{i\geq0}$ such that $k_i\geq 0,\sum_{i=0}^{\infty}k_i=k$ one has that 
 	\be\label{eq: lem4.15-2}
 	\mathbb{E}\Big[\prod_{i:k_i\neq 0}N_i^{k_i}\bigm|\mathfrak{F}(x,\infty),U\Big]
 	=\prod_{i:k_i\neq 0} \sum_{\substack{y_{i,j}:\\ 1\leq j\leq k_i}}\mathbb{P}\Big[\bigcap_{i:k_i\neq0}\bigcap_{1\leq j\leq k_i} y_{i,j}\in L_n(x)\cap \textnormal{Bush}_i(x)\bigm|\mathfrak{F}(x,\infty),U\Big]
 	\ee
 	
 	For each $i\geq0$, let $Y_i=\{y_{i,1},\ldots,y_{i,k_i}\}$ be a finite (possibly with multiplicity) collection of vertices of $G$ and $W_i=\{w_{i,1},\ldots,w_{i,m_i} \}$ the corresponding set of vertices of $Y_i$ without multiplicity. In particular if $k_i=0$ then $m_i=0$ and $W_i$ is an empty set. Let $\mathscr{A}_i$ be the event that for every vertex $w\in W_i$, $w\in \textnormal{Bush}_i(x)$. 
 	
 	Let $\{X^{i,j}:i\geq0,m_i\neq 0,1\leq j\leq m_i \}$ be a collection of independent simple random walks, independent of $\mathfrak{F}(x,\infty)$, such that $X_0^{i,j}=w_{i,j}$ for each 
 	$i\geq0$ with $m_i\neq0$ and $1\leq j\leq m_i$. 
 	
 	For each $i\geq0$ such that $m_i\neq0$, let $\mathscr{B}_i$ be the event that, if we sample $x_i$-WUSF using Wilson's algorithm, starting with the random walks $X^{i,1},\ldots,X^{i,m_i}$, then for every vertex $w\in W_i$, $w$ is connected to $x_i$ in $x_i$-WUSF.
 	
 	It is easy to see that if we sample $\mathfrak{F}$ conditional on $\mathfrak{F}(x,\infty)$ using Wilson's algorithm starting with $X^{0,1},\ldots,X^{0,m_0}$, then $X^{1,1},\ldots,X^{1,m_1}$, and so on, then we have $\mathscr{A}_i\subset \mathscr{B}_i$. 
 	Therefore
 	\be\label{eq: lem4.15-3}
 	\mathbb{P}\big[\bigcap_{i:k_i\neq0}\mathscr{A}_i \bigm|\mathfrak{F}(x,\infty),U\big]
 	\leq \prod_{i:k_i\neq0}\mathbb{P}\big[\mathscr{B}_i \bigm|\mathfrak{F}(x,\infty),U\big]
 	\ee
 	
 	Summing over all the possible choices of the sets $Y_i$ such that $Y_i\subset L_n(x),i\geq0$ and $\sum_{i=0}^{\infty}k_i=k$, by \eqref{eq: lem4.15-1}, \eqref{eq: lem4.15-2} and \eqref{eq: lem4.15-3} we obtain \eqref{eq: k-moments bounds between T_x and x-WUSF}.
 	 \end{proof}
 
 	Considering the high moments for the intersection $T_x\cap L_n(x)$, we have the following upper bounds.
 \begin{proposition}\label{prop: high moments for T_x intersecting a slab}
  For all $n\geq0,k\geq2$ there exists a constant $c_{1}>0$ such that 
 	\be\label{eq: high moment for T_x intersecting a high slab }
 	\mathbb{E}[|T_x\cap L_n(x)|^k]\leq  c_{1}^k(k!)^5(n\vee 1)e^{-t_0n}
 	\ee
 	and
 	\be\label{eq: high moment for T_x intersecting a low slab }
 	\mathbb{E}[|T_x\cap L_{-n}(x)|^k]\leq c_{1}^k(k!)^5(n\vee 1)^{2k}.
 	\ee
 	
 \end{proposition}

\begin{proof}
		By Proposition \ref{prop: slab intersection for the x-WUSF component} and \ref{prop: high moments for x-component intersecting a slab} one has that there exists $c>0$ such that for all $k\geq1$,
	\be\label{eq: all moments for x-wusf intersect with L_n(x)}
	\mathbb{E}[|\mathfrak{T}_{x}\cap L_n(x)|^k]\preceq 
	c^k(k!)^2(|n\wedge0|\vee1)^{k-1}\exp(-t_0(n\vee0))=
	\left\{
	\begin{array}{cc}
		c^k(k!)^2|n|^{k_i-1} & \textnormal{if }n<0\\
		c^k(k!)^2e^{-t_0n} & \textnormal{if }n\geq 0
	\end{array}
	\right.
	\ee
	Suppose  $\mathfrak{F}(x,\infty)=(x_0,x_1,x_2,\ldots)$ is the future of $x$.  From the Observation \ref{obse: relative level relation} we know that if $u\in L_m(x)$ for some integer $m$, then $L_n(x)\subset  \bigcup_{j=n-m-1}^{n-m+1}L_j(u)$. Hence for any positive integer $t$
    \begin{eqnarray*}
	\mathbb{E}\big[|\mathfrak{T}_{u}\cap L_n(x)|^t\bigm|U\big]
	&\leq& 3^t \sum_{j=n-m-1}^{n-m+1}\mathbb{E}[|\mathfrak{T}_{u}\cap L_j(u)|^t]\\
	&\leq& 3^tc^t(t!)^{2} \sum_{j=n-m-1}^{n-m+1}(|j\wedge0|\vee 1)^{t-1}\exp(-t_0(j\vee 0))\nonumber\\
	&\preceq& (6c)^t(t!)^{2}(|(n-m)\wedge0|\vee 1)^{t-1}\exp(-t_0(0\vee(n-m)))
	\end{eqnarray*}

Let $m(x_i)$ denote the integer $m$ such that $x_i\in L_{m}(x)$ (if there are two such $m$'s, taking the smaller one). Then using Lemma \ref{lem: bounds high moments of T_x intersections in terms of the x-WUSF ones}
one has that for any $n\in\mathbb{Z},k\geq 2$
\begin{eqnarray}\label{eq: 4.49 in terms of S_t}
\mathbb{E}[|T_x\cap L_n(x)|^k]
&\preceq&
(6c)^kk! \sum_{\substack{k_0,k_1,\ldots \geq 0:\\ k_0+ k_1+\cdots  =k}} \,\mathbb{E}\left[ \prod_{i:k_i\neq 0} \frac{(k_i!)^{2}}{k_i!}
(|(n-m(x_i))\wedge0|\vee 1)^{k_i-1}\left(\frac{e^{t_0m(x_i)}}{e^{t_0n}}\wedge1\right)\right]\nonumber\\
&\preceq&c_{2}^k(k!)^2\sum_{\substack{k_0,k_1,\ldots \geq 0:\\ k_0+ k_1+\cdots  =k}} \,\mathbb{E}\left[ \prod_{i:k_i\neq 0} 
(|(n-m(x_i))\wedge0|\vee 1)^{k_i-1}\left(\frac{e^{t_0m(x_i)}}{e^{t_0n}}\wedge1\right)\right]\nonumber\\
&\preceq&c_{2}^k(k!)^2\sum_{\substack{k_0,k_1,\ldots \geq 0:\\ k_0+ k_1+\cdots  =k}} \,\mathbb{E}\left[ \prod_{i:k_i\neq 0} 
(|(n-m(x_i))\wedge0|\vee 1)^{k_i}\left(\frac{e^{t_0m(x_i)}}{e^{t_0n}}\wedge1\right)\right]
\end{eqnarray} 
	where in the second inequality we use $\prod_{i:k_i\neq0}k_i!\leq k!$.
	
Write $I:=\{i\colon k_i\neq 0 \}$. Then $1\leq |I|\leq k$ and for a fixed set $I$ of indices, the number of positive solutions $(k_i\colon i\in I)$ for the equation $\sum_{i\in I}{k_i}=k$ is $\binom{k}{|I|-1}\leq 2^k$. Set $f(k_i,m,n):=\left\{\begin{array}{cc}
1 & \textnormal{ if }m\leq n\\
k_i & \textnormal{ if }m>n
\end{array}\right.$.Then for each  sequence of integers $(k_0,k_1,\ldots)$ such that $k_i\geq0$ and $\sum_{i=0}^{\infty}k_i=k$, the term $ \prod_{i:k_i\neq 0} 
(|(n-m(x_i))\wedge0|\vee 1)^{k_i}\big(\frac{e^{t_0m(x_i)}}{e^{t_0n}}\wedge1\big)$ appears in the expansion of $\left(
\sum_{i=0}^{\infty}(|(n-m(x_i))\wedge0|\vee 1)\big(\frac{e^{t_0m(x_i)}}{e^{t_0n}}\wedge1\big)
\right)^t$, where $t=\sum_{i\in I}{f(k_i,m(x_i),n)}\leq \sum_{i\in I}{k_i}=k$. Since  this term appears at most $\binom{k}{|I|-1}\leq 2^k$ times in $$\sum_{\substack{k_0,k_1,\ldots \geq 0:\\ k_0+ k_1+\cdots  =k}} \,\mathbb{E}\left[ \prod_{i:k_i\neq 0} 
(|(n-m(x_i))\wedge0|\vee 1)^{k_i}\left(\frac{e^{t_0m(x_i)}}{e^{t_0n}}\wedge1\right)\right],$$ one has that 
\begin{eqnarray*}\label{eq: 4.49' in terms of S_t}
\mathbb{E}[|T_x\cap L_n(x)|^k]
&\preceq&2^kc_{2}^k(k!)^2\sum_{t=1}^{k}\mathbb{E}\left(
\sum_{i=0}^{\infty}(|(n-m(x_i))\wedge0|\vee 1)\left(\frac{e^{t_0m(x_i)}}{e^{t_0n}}\wedge1\right)
\right)^t\nonumber\\
&=&2^kc_{2}^k(k!)^2\sum_{t=1}^{k}\mathbb{E}\bigg(\sum_{j=-\infty}^{\infty} (|(n-j)\wedge0|\vee1)\big(\frac{e^{t_0j}}{e^{t_0n}}\wedge1\big)|\mathfrak{F}(x,\infty)\cap L_j(x)| \bigg)^t.
\end{eqnarray*}

Set $S_t:=\mathbb{E}\left(\sum_{j=-\infty}^{\infty} (|(n-j)\wedge0|\vee1)\left(\frac{e^{t_0j}}{e^{t_0n}}\wedge1\right)|\mathfrak{F}(x,\infty)\cap L_j(x)| \right)^t$ for $1\leq t\leq k$. 
The above inequality becomes 
\be\label{eq: bound k-th moment of T_x intersect with L_n(x) in terms of S_t}
\mathbb{E}[|T_x\cap L_n(x)|^k]\leq c_{3}^k(k!)^2(S_1+\cdots+S_k)
\ee

By Proposition \ref{prop: slab intersection for the future and past} one has that for $n\geq 0$
\be\label{eq: bound S_1 when n geq 0}
S_1\preceq \sum_{j=-\infty}^{-1}\frac{e^{t_0j}}{e^{t_0n}}+
\sum_{j=0}^{n}\frac{1}{e^{t_0n}}+\sum_{j=n+1}^{\infty}(j-n)e^{-t_0j} 
\preceq (n\vee 1)e^{-t_0n}.
\ee
Similarly if $n<0$, then 
\be\label{eq: bound S_1 when n leq 0}
S_1\preceq \sum_{j=-\infty}^{n-1}\frac{e^{t_0j}}{e^{t_0n}}+\sum_{j=n}^{0}((j-n)\vee 1)+
\sum_{j=1}^{\infty}(j-n)e^{-t_0j}\preceq|n|^2.
\ee	
	
	By H\"{o}lder's inequality one has that for any $t\geq2$ and nonnegative sequences $(a_j)_{j\in\mathbb{Z}},(b_j)_{j\in\mathbb{Z}}$
	\[
	\left(\sum_{j=-\infty}^{\infty}a_jb_j\right)^t
	\leq\left(\sum_{j=-\infty}^{\infty}a_j^{\frac{t}{t-1}}\right)^{t-1}\left(\sum_{j=-\infty}^{\infty}b_j^t\right),
	\]
	 with equality when $b_j=ca_j^{\frac{1}{t-1}}$.
	 
For $2\leq t\leq k$, applying the above H\"{o}lder's inequality with $a_j>0$ and 
$b_j=\frac{1}{a_j}\times (|(n-j)\wedge0|\vee 1)\left(\frac{e^{t_0j}}{e^{t_0n}}\wedge1\right)|\mathfrak{F}(x,\infty)\cap L_j(x)|$ and then taking expectations one has that 
\begin{eqnarray}\label{eq: bound S_t using holder's inequality}
S_t&\leq& \left(\sum_{j=-\infty}^{\infty}a_j^{\frac{t}{t-1}}\right)^{t-1}\cdot \sum_{j=-\infty}^{\infty}
\frac{1}{a_j^t}(|(n-j)\wedge0|\vee1)^t\big(\frac{e^{t_0j}}{e^{t_0n}}\wedge1\big)^t
\mathbb{E}|\mathfrak{F}(x,\infty)\cap L_j(x)|^t\nonumber\\
&\leq&\left(\sum_{j=-\infty}^{\infty}a_j^{\frac{t}{t-1}}\right)^{t-1}\cdot \sum_{j=-\infty}^{\infty}
\frac{1}{a_j^t}(|(n-j)\wedge0|\vee1)^t\big(\frac{e^{t_0j}}{e^{t_0n}}\wedge1\big)^t(ct)^te^{-t_0(j\vee0)},
\end{eqnarray}
	 where in the second inequality we use \eqref{eq: k-moment of intersection of srw with a  slab} and $t!\leq t^t$.
Pick $a_j$ such that 
\be\label{eq: best choice of values of a_j in compact form}
a_j^{\frac{1}{t-1}}=\frac{1}{a_j}(|(n-j)\wedge0|\vee1)\big(\frac{e^{t_0j}}{e^{t_0n}}\wedge1\big)(ct)e^{-\frac{t_0(j\vee0)}{t}}=:\widetilde{b}_j,
\ee
i.e., for $n\geq0$,
\be\label{eq: best choice of values of a_j when n geq 0}
	a_j^{\frac{t}{t-1}}=\left\{\begin{array}{cc}
		(j-n)\cdot 1\cdot (ct)e^{-\frac{t_0j}{t}} & \textnormal{ if }j>n\\
		& \\
		1\cdot \frac{e^{t_0j}}{e^{t_0n}}\cdot (ct)e^{-\frac{t_0j}{t}} & \textnormal{ if }0\leq j\leq n\\
		& \\
		1\cdot \frac{e^{t_0j}}{e^{t_0n}}\cdot (ct) & \textnormal{ if }j<0
	\end{array}\right.
\ee	
and for $n<0$,
\be\label{eq: best choice of values of a_j when n less than 0}
a_j^{\frac{t}{t-1}}=\left\{\begin{array}{cc}
	(j-n)\cdot 1\cdot (ct)e^{-\frac{t_0j}{t}} & \textnormal{ if }j>0\\
	&\\
	(j-n)\cdot 1\cdot (ct) & \textnormal{ if }n<j\leq 0\\
	&\\
	1\cdot \frac{e^{t_0j}}{e^{t_0n}}\cdot (ct) & \textnormal{ if }j\leq n
\end{array}\right.
\ee	
With $a_j$ given above and $\widetilde{b}_j$ as defined in \eqref{eq: best choice of values of a_j in compact form},  $\widetilde{b}_j=a_j^{\frac{1}{t-1}}$ and  \eqref{eq: bound S_t using holder's inequality} becomes
\begin{eqnarray}\label{eq: bound S_t with best choice of a_j}
S_t&\leq&\left(\sum_{j=-\infty}^{\infty}a_j^{\frac{t}{t-1}}\right)^{t-1}\cdot \sum_{j=-\infty}^{\infty} \widetilde{b}_j^t\nonumber\\
&=&\left(\sum_{j=-\infty}^{\infty}a_j^{\frac{t}{t-1}}\right)^{t-1}\cdot \sum_{j=-\infty}^{\infty} a_j^{\frac{t}{t-1}}=\left(\sum_{j=-\infty}^{\infty}a_j^{\frac{t}{t-1}}\right)^{t}
\end{eqnarray}

Hence for $n\geq0,2\leq t\leq k$, by \eqref{eq: best choice of values of a_j when n geq 0} and \eqref{eq: bound S_t with best choice of a_j} one has 
\begin{eqnarray}\label{eq: bound of S_t when n geq 0 and t geq 2}
S_t&\leq&\left(\sum_{j=-\infty}^{-1}ct\frac{e^{t_0j}}{e^{t_0n}}+\sum_{j=1}^{n}ct\frac{e^{t_0j}}{e^{t_0n}}e^{-\frac{t_0j}{t}} +\sum_{j=n+1}^{\infty}(j-n)cte^{-\frac{t_0j}{t}} \right)^t\nonumber\\
&\leq&(ct)^t\left(\frac{c_4}{e^{t_0n}}+\frac{c_4e^{t_0(1-\frac{1}{t})n}}{e^{t_0n}(e^{t_0(1-\frac{1}{t})}-1)}+c_4e^{-\frac{t_0n}{t}}\int_{0}^{\infty}xe^{-\frac{t_0x}{t}}dx\right)^t\nonumber\\
&\leq&(ct)^t\left(\frac{c_4}{e^{t_0n}}+\frac{c_5}{e^{\frac{t_0n}{t}}}+c_4e^{-\frac{t_0n}{t}}\big(\frac{t}{t_0}\big)^2\right)^t\nonumber\\
&\leq&c_6^tt^{3t}e^{-t_0n}\leq c_7^t(t!)^3e^{-t_0n},
\end{eqnarray}
where we use  $t^t\leq c_8^tt!$ for some $c_8>1$ (by Stirling's formula for example) in the last inequality.

For $n<0,2\leq t\leq k$, by \eqref{eq: best choice of values of a_j when n less than 0} and \eqref{eq: bound S_t with best choice of a_j} one has 
\begin{eqnarray}\label{eq: bound of S_t when n leq 0 and t geq 2}
S_t&\leq&\left(\sum_{j=-\infty}^{n}ct\frac{e^{t_0j}}{e^{t_0n}}+\sum_{j=n+1}^{0}(j-n)ct +\sum_{j=1}^{\infty}(j-n)cte^{-\frac{t_0j}{t}} \right)^t\nonumber\\
&\leq&(ct)^t\left(c_9+|n|^2+c_9|n|\int_{0}^{\infty}e^{-\frac{t_0x}{t}}dx+c_9\int_{0}^{\infty}xe^{-\frac{t_0x}{t}}dx\right)^t\nonumber\\
&=&(ct)^t\left(c_9+|n|^2+c_9|n|\frac{t}{t_0}+c_9\big(\frac{t}{t_0}\big)^2\right)^t\nonumber\\
&\leq&(ct)^t\big(c_{10}(|n|+t)^2\big)^t\leq c_{11}^tt^t(|n|t)^{2t}\nonumber\\
&\leq& c_{12}^t(t!)^3|n|^{2t},
\end{eqnarray}

Now we are ready to get the final conclusion. For $n\geq0$, by \eqref{eq: bound k-th moment of T_x intersect with L_n(x) in terms of S_t}, \eqref{eq: bound S_1 when n geq 0} and \eqref{eq: bound of S_t when n geq 0 and t geq 2} we get 
\[
\mathbb{E}[|T_x\cap L_n(x)|^k]\leq c_{1}^k(k!)^5(n\vee 1)e^{-t_0n}.
\]
For $n<0$, by \eqref{eq: bound k-th moment of T_x intersect with L_n(x) in terms of S_t}, \eqref{eq: bound S_1 when n leq 0} and \eqref{eq: bound of S_t when n leq 0 and t geq 2} we get 
\[
\mathbb{E}[|T_x\cap L_n(x)|^k]\leq c_{1}^k(k!)^5|n|^{2k}.\qedhere
\]
\end{proof}

 \subsection{Probabilities for intersections with a slab}

 Another natural question to consider is the decay of the probability of a vertex connecting to a certain slab. Similar question was analyzed for Bernoulli percolation \cite[Lemma 5.2]{Tom2017}.
 
  \begin{table}[ht]
 	\centering
 	\begin{tabular}{|c|c|c|}
 		\hline
 		Quantities $(n\in\mathbb{Z})$ & Results & Position\\
 		\hline
 		$\mathbb{P}[T_x\cap L_n(x)\neq \emptyset]$ & $\left\{
 		\begin{array}{cc}
 		\asymp ne^{-t_0n} & \textnormal{ if }n>0\\
 		\stackrel{\eqref{eq: almost surely nonempty intersection with lower slabs}}{=}1& \textnormal{ if }n\leq 0\\
 		\end{array}
 		\right.$ & Proposition \ref{prop: exponential decay of the probability that T_x intersecting a high slab}\\
 			\hline
 		$\mathbb{P}[\mathfrak{T}_x\cap L_{n}(x)\neq\emptyset]$ & $\left\{
 		\begin{array}{cc}
 		\asymp e^{-t_0n} & \textnormal{ if }n\geq0\\
 		\asymp \frac{1}{|n|}& \textnormal{ if }n<0\\
 		\end{array}
 		\right.$ & Proposition \ref{prop: decay of probability for x-WUSF component insecting a slab far away}\\
 		\hline
 		$\mathbb{P}[\mathfrak{P}(x)\cap L_n(x)\neq\emptyset]$ &  $\left\{
 		\begin{array}{cc}
 		\asymp e^{-t_0n} & \textnormal{ if }n\geq0\\
 		\asymp \frac{1}{n}& \textnormal{ if }n<0\\
 		\end{array}
 		\right.$ & Proposition \ref{prop: exponentail decay for the past and future intersecting with a slab} \\
 		\hline
 		$\mathbb{P}[\mathfrak{F}(x,\infty)\cap L_{n}(x)\neq\emptyset]$ &  $\left\{
 		\begin{array}{cc}
 		\asymp e^{-t_0n} & \textnormal{ if }n>0\\
 		\stackrel{\eqref{eq: almost surely nonempty intersection with lower slabs}}{=}1& \textnormal{ if }n\leq 0\\
 		\end{array}
 		\right.$& Proposition \ref{prop: exponentail decay for the past and future intersecting with a slab}\\
 		\hline
 	\end{tabular}
 \centering
 	\caption{Probabilities of nonempty intersections with a slab}
 	\label{table: probability of nonempty intersections}
 \end{table}
 
 \begin{proposition}\label{prop: exponential decay of the probability that T_x intersecting a high slab}
 	Let $\{T_x\cap L_n(x)\neq \emptyset\}$ denote the event that $T_x$ has a nonempty intersection with the slab $L_n(x)$. Then for $n\geq1$
 	\be\label{eq: exponential decay for T_x intersecting a high slab}
  \mathbb{P}[T_x\cap L_n(x)\neq \emptyset]\asymp n\exp(-t_0n).
 	\ee
 	In particular
 	\be
 	\lim_{n\rightarrow\infty}-\frac{1}{t_0n}\log\mathbb{P}[T_x\cap L_n(x)\neq \emptyset]=1.
 	\ee
 	
 \end{proposition}
\begin{proof}
	By Lemma \ref{lem:expectation for SRW on nonunimodular transitive graph} one has that $\mathfrak{F}(x,\infty)$ intersects every $L_n(x)$ almost surely. Since $\mathfrak{F}(x,\infty)\subset T_x$, one has 
	\be\label{eq: almost surely nonempty intersection with lower slabs}
	\mathbb{P}[T_x\cap L_n(x)\neq\emptyset]=\mathbb{P}[\mathfrak{F}(x,\infty)\cap L_n(x)\neq\emptyset]=1,\,\,\forall \ n\leq 0.
	\ee

	On the one hand, for $n\geq1$ using second moment method and  \eqref{eq: expectation of intersection of T_x with a high slab} and \eqref{eq: high moment for T_x intersecting a high slab }  one has 
	\[
	\mathbb{P}[T_x\cap L_n(x)\neq \emptyset]
	\geq 
	\frac{\big(\mathbb{E}[|T_x\cap L_n(x)|]\big)^2}{\mathbb{E}[|T_x\cap L_n(x)|^2]}
	\succeq ne^{-t_0n}.
	\]
	On the other hand, note when $T_x\cap L_n(x)\neq \emptyset$, $|T_x\cap L_n(x)|\geq1$. Thus by Markov's inequality one has that 
	\[
		\mathbb{P}[T_x\cap L_n(x)\neq \emptyset]=\mathbb{P}[|T_x\cap L_n(x)|\geq 1]
		\leq \mathbb{E}[|T_x\cap L_{n}(x)|]\stackrel{\eqref{eq: expectation of intersection of T_x with a high slab}}{\preceq } n\exp(-t_0n).
	\]
	Combining the above two inequalities one has the conclusion. 
\end{proof}

Since the tree $\mathfrak{T}_x$ in the $x$-WUSF is almost surely finite, $\mathbb{P}[\mathfrak{T}_x\cap L_{n}(x)\neq \emptyset]$ and $\mathbb{P}[\mathfrak{T}_x\cap L_{-n}(x)\neq \emptyset]$  both decay to zero as $n\rightarrow\infty$. 
\begin{proposition}\label{prop: decay of probability for x-WUSF component insecting a slab far away}
	For the tree $\mathfrak{T}_x$ in the $x$-WUSF and $n\geq 0$, one  has
	\be\label{eq: exponential decay for x-WUSF component intersecting a high slab}
	\mathbb{P}[\mathfrak{T}_x\cap L_{n}(x)\neq \emptyset]\asymp \exp(-t_0n)
	\ee
	
	In particular,  
	\be
	\lim_{n\rightarrow\infty}-\frac{1}{t_0n}\log\mathbb{P}[\mathfrak{T}_x\cap L_{n}(x)\neq \emptyset]=1.
	\ee
	Moreover for $n\geq 1$, the decay of $\mathbb{P}[\mathfrak{T}_x\cap L_{-n}(x)\neq \emptyset]$ is much slower.
	\be\label{eq: polynomial decay for x-WUSF component intersecting a low slab}
	\mathbb{P}[\mathfrak{T}_x\cap L_{-n}(x)\neq \emptyset]\asymp\frac{1}{n}.
	\ee
\end{proposition}
\begin{proof}
	The proof of \eqref{eq: exponential decay for x-WUSF component intersecting a high slab} is similar to \eqref{eq: exponential decay for T_x intersecting a high slab} and thus we omit it.
	
	By second moment method one has the lower bound in \eqref{eq: polynomial decay for x-WUSF component intersecting a low slab}, namely 
	\[
	\mathbb{P}[\mathfrak{T}_x\cap L_{-n}(x)\neq \emptyset]\geq \frac{\left(\mathbb{E}[|\mathfrak{T}_x\cap L_{-n}(x)|]\right)^2}{\mathbb{E}[|\mathfrak{T}_x\cap L_{-n}(x)|^2]}
	\stackrel{\eqref{eq: expectation of intersection of x-WUSF component with a low slab},\eqref{eq: second moment of x-component intersecting a low slab}}{\succeq } \frac{1}{n}.
	\]
	
	Let $\textnormal{diam}_{\textnormal{int}}(\mathfrak{T}_x)$ denote the intrinsic diameter of the finite tree $\mathfrak{T}_x$. If $\mathfrak{T}_x\cap L_{-n}(x)\neq \emptyset$, then $\textnormal{diam}_{\textnormal{int}}(\mathfrak{T}_x)\geq n-1$. Hence for $n\geq 1$  one has 
	\[
	\mathbb{P}[\mathfrak{T}_x\cap L_{-n}(x)\neq \emptyset]\leq\mathbb{P}[\textnormal{diam}_{\textnormal{int}}(\mathfrak{T}_x)\geq n-1]\preceq\frac{1}{n},
	\]
	where the last inequality is due to Theorem 7.1 of \cite{Hutchcroft2018a}.
\end{proof}

For  the future $\mathfrak{F}(x,\infty)$, it intersects every $L_{-n}(x),n\geq0$ almost surely, and the probability $\mathbb{P}[\mathfrak{F}(x,\infty)\cap  L_n(x)\neq\emptyset]$ decays exponentially. For the past $\mathfrak{P}(x)$, it is a finite tree and hence $\mathbb{P}[\mathfrak{P}(x)\cap  L_n(x)\neq\emptyset]$ and $\mathbb{P}[\mathfrak{P}(x)\cap  L_{-n}(x)\neq\emptyset]$ both decay to zero as $n$ tending to infinity. We summarize these in the following Proposition.
\begin{proposition}\label{prop: exponentail decay for the past and future intersecting with a slab}
	Suppose $n\geq0$. For the future one has the following asymptotic behavior:
	\be\label{eq: exponential decay for future intersect a high slab}
	\mathbb{P}[\mathfrak{F}(x,\infty)\cap  L_n(x)\neq\emptyset]\asymp \exp(-t_0n)
	\ee
	For the past one has 
	\be\label{eq: exponential decay for past intersect a high slab}
	   \mathbb{P}[\mathfrak{P}(x)\cap  L_n(x)\neq\emptyset]\asymp \exp(-t_0n)
	\ee
	and
	\be\label{eq: polynomial decay for past intersect a low slab}
	\mathbb{P}[\mathfrak{P}(x)\cap  L_{-n}(x)\neq\emptyset]\asymp\frac{1}{n\vee 1}.
	\ee

\end{proposition}
\begin{proof}
By \eqref{eq: expectation of intersection of future of x with a high slab} one has $\mathbb{E}[|\mathfrak{F}(x,\infty)\cap L_n(x)|]\asymp e^{-t_0n}$. From Corollary \ref{cor: high moments for the future intersect with a slab} one has $\mathbb{E}[|\mathfrak{F}(x,\infty)\cap L_n(x)|^2]\preceq e^{-t_0n}$. Using the first-moment and second-moment method as in the proof of Proposition \ref{prop: exponential decay of the probability that T_x intersecting a high slab} again one obtains \eqref{eq: exponential decay for future intersect a high slab}.

	Similarly for the past, by \eqref{eq: expectation of intersection of past of x with a high slab} one has $\mathbb{E}[|\mathfrak{P}(x)\cap L_n(x)|]\asymp e^{-t_0n}$. From Corollary \ref{cor: high moments for the past} one has $\mathbb{E}[|\mathfrak{P}(x)\cap L_n(x)|^2]\preceq e^{-t_0n}$. Using the first-moment and second-moment method  again one obtains \eqref{eq: exponential decay for past intersect a high slab}.
	
	By Lemma \ref{lem:stochastic domination} and Proposition \ref{prop: high moments for x-component intersecting a slab} one has \[
	\mathbb{E}[|\mathfrak{P}(x)\cap L_n(x)|^2]\preceq n\exp(-t_0n).
	\]
	From Proposition \ref{prop: slab intersection for the future and past} one has 
	$\mathbb{E}[|\mathfrak{P}(x)\cap L_n(x)|]\succeq \exp(-t_0n)$.
	Using Markov's inequality one then has  the lower bound of \eqref{eq: exponential decay for past intersect a high slab}. 
	
	By Lemma \ref{lem:stochastic domination} and Proposition \ref{prop: decay of probability for x-WUSF component insecting a slab far away} one has the upper bound in \eqref{eq: polynomial decay for past intersect a low slab}:
	\[
	\mathbb{P}[\mathfrak{P}(x)\cap  L_{-n}(x)\neq\emptyset]\leq	\mathbb{P}[\mathfrak{T}_x\cap  L_{-n}(x)\neq\emptyset] \preceq \frac{1}{n\vee 1}.
	\]
	
	By \eqref{eq: expectation of intersection of past of x with a low slab} one has  $\mathbb{E}[|\mathfrak{P}(x)\cap L_{-n}(x)|]\asymp 1$.  From Corollary \ref{cor: high moments for the past} one has $\mathbb{E}[|\mathfrak{P}(x)\cap L_{-n}(x)|^2]\preceq n\vee1$. Then the proof of the lower bound of \eqref{eq: polynomial decay for past intersect a low slab} is an application of second-moment method again.
\end{proof}

 \subsection{Tail bounds for tilted volumes}
The results for this subsection are summarized in Table \ref{table: tail probability for tilted volumes}.
 \begin{table}[!htbp]
	\centering
	\begin{tabular}{|c|c|c|}
		\hline
		Quantities $(R\geq 2)$ & Results & Position\\
		\hline
		$\mathbb{P}[|T_x|_{x,\lambda}\geq R]$ & $\left\{
		\begin{array}{cc}
		\asymp_\lambda R^{-\frac{1}{\lambda}}\log R & \textnormal{ if }\lambda>0\\
		\stackrel{\eqref{eq: tilted volume is infinite for T_x and the future with nonpositive lambda}}{=}1& \textnormal{ if }\lambda\leq 0\\
		\end{array}
		\right.$ & Proposition \ref{prop: tail prob for general case}\\
		\hline
		$\mathbb{P}[|\mathfrak{T}_x|_{x,\lambda}\geq R]$ & $\left\{
		\begin{array}{cc}
		\asymp_\lambda R^{-\frac{1}{\lambda}} & \textnormal{ if }\lambda>0\\
		\asymp \frac{1}{\sqrt{R}} & \textnormal{ if }\lambda=0\\
		\asymp_\lambda \frac{1}{\log R}& \textnormal{ if }\lambda <0\\
		\end{array}
		\right.$ & Proposition \ref{prop: tail probability for tilted volume of x-wusf component}\\
		\hline
		$\mathbb{P}[|\mathfrak{P}(x)|_{x,\lambda}\geq R]$ &  $\left\{
		\begin{array}{cc}
		\asymp_\lambda R^{-\frac{1}{\lambda}} & \textnormal{ if }\lambda>0\\
		\asymp \frac{1}{\sqrt{R}} & \textnormal{ if }\lambda=0\\
		\asymp_\lambda \frac{1}{\log R}& \textnormal{ if }\lambda <0\\
		\end{array}
		\right.$  & Proposition \ref{prop: tail probability for tilted volume of the past and future} \\
		\hline
		$\mathbb{P}[|\mathfrak{F}(x,\infty)|_{x,\lambda}\geq R]$ &  $\left\{
		\begin{array}{cc}
		\asymp_\lambda R^{-\frac{1}{\lambda}}\log R & \textnormal{ if }\lambda>0\\
		\stackrel{\eqref{eq: tilted volume is infinite for T_x and the future with nonpositive lambda}}{=}1& \textnormal{ if }\lambda\leq 0\\
		\end{array}
		\right.$& Proposition \ref{prop: tail probability for tilted volume of the past and future}\\
		\hline
	\end{tabular}
\centering
	\caption{Probabilities of nonempty intersections with a slab}
	\label{table: tail probability for tilted volumes}
\end{table}

\begin{proposition}\label{prop: tail prob for general case}
	For $\lambda>0$, the tilted volume $|T_x|_{x,\lambda}$ has the following property.
	\be\label{eq:bounds for tail prob of general case}
		 \mathbb{P}[|T_x|_{x,\lambda}\geq R]\asymp_{\lambda} R^{-\frac{1}{\lambda}}\log R , \ \forall\  R\geq 2.
		\ee

\end{proposition}

 \begin{proof}
   By Lemma \ref{lem:expectation for SRW on nonunimodular transitive graph} the future $\mathfrak{F}(x,\infty)$ intersects every $L_n(x),n\leq 0$ almost surely. In particular this implies that for every $\lambda\leq0$,  
   \be\label{eq: tilted volume is infinite for T_x and the future with nonpositive lambda}
   \mathbb{P}[|T_x|_{x,\lambda}\geq R]=\mathbb{P}[|\mathfrak{F}(x,\infty)|_{x,\lambda}\geq R]=1,\,\,R\geq 2.
   \ee

 		Taking $n_0=\lceil1+\frac{\log R}{t_0\lambda} \rceil$, then $1+\frac{\log R}{t_0\lambda} \leq n_0\leq 2+\frac{\log R}{t_0\lambda}$ and 
 		\be\label{eq: choice of n_0 in terms of R}
 		e^{-2t_0}R^{-\frac{1}{\lambda}}\leq \exp(-t_0n_0)\leq 	e^{-t_0}R^{-\frac{1}{\lambda}}
 		\ee
 		On the event that $\{T_x\cap L_{n_0}(x)\neq \emptyset \}$, one has $|T_x|_{x,\lambda}\geq   \exp(t_0\lambda (n_0-1))\geq R$. Hence by Proposition \ref{prop: exponential decay of the probability that T_x intersecting a high slab}
 		\[
 		\mathbb{P}[|T_x|_{x,\lambda}\geq R]\geq \mathbb{P}[T_x\cap L_{n_0}(x)\neq \emptyset ]
 		\succeq n_0\exp(-t_0n_0)\succeq\frac{1}{\lambda} R^{-\frac{1}{\lambda}}\log R.
 		\]

 		Now we derive the upper bound. If $T_x\cap L_{n_0}(x)= \emptyset$, then $T_x\cap L_{j}(x)= \emptyset$ for all $j\geq n_0$. If $y\in L_j(x)$, then $\Delta(x,y)^\lambda\leq e^{t_0(j+1)\lambda}$. Thus on the event $T_x\cap L_{n_0}(x)= \emptyset$, $|T_x|_{x,\lambda}\leq \sum_{j=\infty}^{n_0-1}e^{t_0(j+1)\lambda}|T_x\cap L_j(x)|$. 
 		Using the union bounds we have 
 		\begin{eqnarray}
 		\mathbb{P}[|T_x|_{x,\lambda}\geq R]
 		&\leq& \mathbb{P}[T_x\cap L_{n_0}(x)\neq \emptyset]+\mathbb{P}\left[
 		\sum_{j=-\infty}^{n_0}e^{t_0(j+1)\lambda}|T_x\cap L_j(x)|\geq R\right]\nonumber\\
 		&=:&J_1+J_2
 		\end{eqnarray}
 		
 		By  Proposition \ref{prop: exponential decay of the probability that T_x intersecting a high slab}
 		we have 
 		\be
 		J_1\preceq  n_0\exp(-t_0n_0)\stackrel{\eqref{eq: choice of n_0 in terms of R} }{\preceq} \frac{2t_0\lambda+\log R}{\lambda} R^{-\frac{1}{\lambda}}\preceq_\lambda R^{-\frac{1}{\lambda}}\log R.
 		\ee
 		Taking a small constant $c>0$ and $h_j=\frac{c}{(|j-n_0|\vee 1)^2}$ such that $\sum_{j=-\infty}^{n_0}h_j<\sum_{j=-\infty}^{\infty}h_j\leq 1$, using the union bound we get
 		\be
 		J_2\leq \sum_{j=-\infty}^{n_0}\mathbb{P}[e^{t_0(j+1)\lambda}|T_x\cap L_j(x)|\geq Rh_j].
 		\ee
 		Taking a positive integer $k=\lceil\frac{1}{\lambda}\rceil+1>\frac{1}{\lambda}$, by Markov's inequality one has 
 		\be
 		J_2\leq \sum_{j=-\infty}^{n_0}\frac{1}{R^kh_j^k}\mathbb{E}[e^{t_0(j+1)\lambda k}|T_x\cap L_j(x)|^k]
 		\ee
 		By Proposition \ref{prop: high moments for T_x intersecting a slab} we have 
 		\begin{eqnarray}
 		J_2&\preceq_\lambda&\sum_{j=-\infty}^{0}\frac{e^{t_0(j+1)\lambda k}}{R^kh_j^k}(|j|\vee1)^{2k}+\sum_{j=1}^{n_0}\frac{e^{t_0(j+1)\lambda k}}{R^kh_j^k}j\exp(-t_0j)\nonumber\\
 		&\preceq_\lambda&\frac{1}{R^k}\sum_{j=-\infty}^{0}e^{t_0\lambda kj}(|j|\vee1)^{2k}(|j-n_0|\vee1)^{2k}+\frac{1}{R^k}\sum_{j=1}^{n_0}j(|j-n_0|\vee1)^{2k}e^{t_0(\lambda k-1)j}\nonumber\\
 		&\preceq_\lambda&\frac{1}{R^k}\sum_{j=-\infty}^{0}e^{t_0\lambda kj}(|j|\vee1)^{4k}n_0^{2k}+\frac{n_0e^{t_0(\lambda k-1)n_0}}{R^k}+\frac{n_0}{R^k}\sum_{j=1}^{n_0-1}(n_0-j)^{2k}e^{t_0(\lambda k-1)j}\nonumber\\
 		&\preceq_\lambda&\frac{n_0^{2k}}{R^k}+n_0e^{-t_0n_0}+\frac{n_0e^{t_0(\lambda k-1)n_0}}{R^k}\int_{0}^{n_0}(n_0-x)^{2k}e^{t_0(\lambda k-1)(x-n_0)}dx\nonumber\\
 		&\preceq_\lambda&\frac{\big(\log R\big)^{2+2\lceil \frac{1}{\lambda}\rceil}}{R^{\frac{1}{\lambda}+1}}+n_0e^{-t_0n_0}+n_0e^{-t_0n_0}\cdot
 		\int_{0}^{\infty}y^{2k}e^{-t_0(\lambda k-1)y}dx
 		 \nonumber\\
 		&\preceq_\lambda&\frac{\big(\log R\big)^{2+2\lceil \frac{1}{\lambda}\rceil}}{R^{\frac{1}{\lambda}+1}}+n_0e^{-t_0n_0}+n_0e^{-t_0n_0}\nonumber\\
 			&\preceq_\lambda& R^{-\frac{1}{\lambda}}\log R,
 		\end{eqnarray}	
 		where in the fourth to the sixth inequalities we use $R^k\asymp_\lambda e^{t_0\lambda kn_0}$ by the choice of $n_0$.

 		Therefore 
 		\[
 		\mathbb{P}[|T_x|_{x,\lambda}\geq R]\leq J_1+J_2\preceq_\lambda R^{-\frac{1}{\lambda}}\log R.\qedhere
 		\]
	
 \end{proof}

Since the tree $\mathfrak{T}_x$ in the $x$-WUSF is almost surely finite, $|\mathfrak{T}_x|_{x,\lambda}<\infty$ for all $\lambda\in\mathbb{R}$. 
\begin{proposition}\label{prop: tail probability for tilted volume of x-wusf component}
	The tail probability  $\mathbb{P}[|\mathfrak{T}_x|_{x,\lambda}\geq R]$ satisfies the following inequalities. 
	\begin{itemize}
		\item If $\lambda<0$, then 
		\be\label{eq: tail probability of tilted volume of x-WUSF component in case lambda less than zero}
		\mathbb{P}[|\mathfrak{T}_x|_{x,\lambda}\geq R]\asymp_{\lambda} \frac{1}{\log R},\, R\geq 2.
		\ee
		
		\item  If $\lambda=0$, then $|\mathfrak{T}_x|_{x,\lambda}=|\mathfrak{T}_x|$ is the size of the tree $\mathfrak{T}_x$ and 
		\be\label{eq: tail probability of tilted volum of x-WUSF component large in case lambda equal to zero}
			\mathbb{P}[|\mathfrak{T}_x|\geq R]\asymp \frac{1}{\sqrt{R}}, \,R\geq 2.
		\ee
		
		\item  If $\lambda>0$, then 
		\be\label{eq: tail probaiblity of x-WUSF component in case lambda larger than zero}
			\mathbb{P}[|\mathfrak{T}_x|_{x,\lambda}\geq R]\asymp_{\lambda} R^{-\frac{1}{\lambda}}, \,R\geq 2.
		\ee

	\end{itemize}
\end{proposition}
\begin{proof}
	We only prove \eqref{eq: tail probability of tilted volume of x-WUSF component in case lambda less than zero} since \eqref{eq: tail probability of tilted volum of x-WUSF component large in case lambda equal to zero} was proved to hold for more general networks in Theorem 7.1 of \cite{Hutchcroft2018a} and the proof of  \eqref{eq: tail probaiblity of x-WUSF component in case lambda larger than zero} is similar to the ones in Proposition \ref{prop: tail prob for general case}. 
	
	Now we assume $\lambda<0$ and $R\geq 2$. Taking  $n_0=\lceil \frac{\log R}{-t_0\lambda}+1\rceil$, then $\exp(-t_0(n_0-1)\lambda)\geq R$. Hence if $\mathfrak{T}_x\cap L_{-n_0}(x)\neq \emptyset$, then 
$|\mathfrak{T}_x|_{x,\lambda}\geq \exp(-t_0(n_0-1)\lambda)\geq R$. Therefore by Proposition \ref{prop: decay of probability for x-WUSF component insecting a slab far away} one has that 
\[
	\mathbb{P}[|\mathfrak{T}_x|_{x,\lambda}\geq R]\geq 
	\mathbb{P}[\mathfrak{T}_x\cap L_{-n_0}(x)\neq \emptyset]\succeq\frac{1}{n_0}\succeq_\lambda
 \frac{1}{\log R}.
\]
	On the other hand, there exists a constant $C>0$ small enough such that for all $R\geq 2D$
	\[
	\sum_{y\in B(x,C\log R)}\Delta(x,y)^\lambda
	< D^{C\log R+1}\exp(-t_0\lambda C\log R)
	\leq  R.
	\]
	Therefore, for $R\geq 2D$
	\[
		\mathbb{P}[|\mathfrak{T}_x|_{x,\lambda}\geq R]
		\leq
		\mathbb{P}[\textnormal{diam}_{\textnormal{int}}(\mathfrak{T}_x)\geq C\log R] \preceq \frac{1}{\log R},
	\] 
	where the last inequality is again due to Theorem 7.1 of \cite{Hutchcroft2018a}.
\end{proof}

\begin{remark}
	As $\lambda\rightarrow0$, the tilted volume $|\mathfrak{T}_x|_{x,\lambda}\rightarrow|\mathfrak{T}_x|$. So for properly related $\lambda$ and $R$, the probabilities $ \mathbb{P}[|\mathfrak{T}_x|_{x,\lambda}\geq R]$  and $ \mathbb{P}[|\mathfrak{T}_x|\geq R]$   should be close. This is due to the dependence of $\lambda$ implicitly in \eqref{eq: tail probability of tilted volume of x-WUSF component in case lambda less than zero} and \eqref{eq: tail probaiblity of x-WUSF component in case lambda larger than zero}. Indeed if $\lambda=-\frac{\log R}{\sqrt{R}}$, then 
	$\mathbb{P}[|\mathfrak{T}_x|_{x,\lambda}\geq R]\succeq\frac{1}{n_0}\succeq\frac{1}{\sqrt{R}}$, where $n_0=\lceil \frac{\log R}{-t_0\lambda}+1\rceil$. Moreover for $\lambda=-\frac{\log R}{\sqrt{R}}$, let $j_0=\frac{\log 2}{t_0\lambda}$, if $j\geq j_0$, then $\exp(t_0j\lambda)\leq 2$. Thus
	$\mathbb{P}[|\mathfrak{T}_x|_{x,\lambda}\geq R]\leq \mathbb{P}[\mathfrak{T}_x\cap L_{j_0}(x)\neq\emptyset]+\mathbb{P}[2|\mathfrak{T}_x|\geq R]\preceq \frac{1}{|j_0|}+\frac{1}{\sqrt{R}}\preceq \frac{\log R}{\sqrt{R}}$. 
\end{remark}

For the future and past of $x$ in the WUSF sample $\mathfrak{F}$, one can also consider the tail probability for the titled volumes. We summarize the results in the following proposition and omit the proofs since they are similar to the ones for $T_x$ and $\mathfrak{T}_x$. 
\begin{proposition}\label{prop: tail probability for tilted volume of the past and future}
	If $\lambda>0$, then 
	\[
\mathbb{P}[|\mathfrak{F}(x,\infty)|_{x,\lambda}\geq R]\asymp_\lambda  R^{-\frac{1}{\lambda}}
	\]
	and 
	\[
		\mathbb{P}[|\mathfrak{P}(x)|_{x,\lambda}\geq R]
	\asymp_{\lambda} R^{-\frac{1}{\lambda}}.
	\]
	
	If $\lambda=0$, then 
	\[
	\mathbb{P}[|\mathfrak{P}(x)|\geq R]\asymp\frac{1}{\sqrt{R}}.
	\]
	
	If $\lambda<0$, then 
	\[
	\mathbb{P}[|\mathfrak{P}(x)|_{x,\lambda}\geq R]\asymp_\lambda\frac{1}{\log R}.
	\]
\end{proposition}

 \section{FUSF on Diestel--Leader graphs and grandparent graphs}\label{sec: 5}

 \subsection{FUSF=WUSF on Diestel--Leader graphs}
 Diestel--Leader graphs are a family of transitive graphs constructed by Diestel and Leader in \cite{DL2001} as possible examples of transitive graphs that are not roughly isometric to any Cayley graph. Later Eskin, Fisher, and Whyte  \cite{EFW2012} showed that the Diestel--Leader graph $DL(q,r)$ with $q\neq r$ is indeed not roughly isometric to any Cayley graph. 
 
 Next we give a precise definition of the Diestel--Leader graph $DL(q,r)$. For more details see \cite{Woess2005}. 
	 Given a regular tree $\mathbb{T}_{b+1}$ with $b\geq2$, fix a root $o\in \mathbb{T}_{b+1}$ and an end $\xi$. For each vertex $v\in \mathbb{T}_{b+1}$, there is a unique ray $\eta_v$ starting from $v$ and representing $\xi$. We call the unique neighbor of $v$ on the ray $\eta_v$ the \notion{parent} of $v$ with respect to $\xi$, and denote it by $v^-$. For $x,y\in \mathbb{T}_{b+1}$, we define $x\curlywedge y$ (w.r.t. $\xi$) to be the first intersection vertex of $\eta_x$ and $\eta_y$.  We define the \notion{Busemann function} $\mathfrak{h}:\mathbb{T}_{b+1}\rightarrow\mathbb{Z}$ with respect to $o,\xi$ as follows:
\[
\mathfrak{h}(x)=d(o,x\curlywedge o)-d(x,x\curlywedge o), 
\]
where $d(u,v)$ is the graph distance between $u$ and $v$ on  $\mathbb{T}_{b+1}$. We also define the \notion{horocycles} $H_k$ with respect to $o,\xi$ as $H_k=H_k(o,\xi):=\{x\in \mathbb{T}_{b+1}: \mathfrak{h}(x)=k \}$.  Note that changing the base $o$ will only change the Busemann function by adding a constant.

 Suppose $q\geq r\geq2$ are two positive integers and $\mathbb{T}_{q+1},\mathbb{T}_{r+1}$ are two regular trees with degree $q+1,r+1$ respectively. Fix roots $o_1,o_2$ and reference ends $\omega_1,\omega_2$ for $\mathbb{T}_{q+1},\mathbb{T}_{r+1}$ respectively. 
 \begin{definition}
 	The set of vertices of Diestel--Leader graph $DL(q,r)$ is given by 
 	\[
   DL(q,r)=\{ x_1x_2\in \mathbb{T}_{q+1}\times\mathbb{T}_{r+1}: \mathfrak{h}(x_1)+\mathfrak{h}(x_2)=0 \},
 	\]
 	where we use $\mathfrak{h}$ as the Busemann functions on $\mathbb{T}_{q+1},\mathbb{T}_{r+1}$ w.r.t. to $o_1,\omega_1$ and $o_2,\omega_2$.
 	
 	The neighborhood relation is given by 
 	\[
 	x_1x_2\sim y_1y_2 \textnormal{ if and only if } x_1\sim y_1 \textnormal{ and }y_1\sim y_2.
 	\]
 \end{definition}
 A way of visualizing of $DL(q,r)$ is described on page 418 of \cite{Woess2005}; see also Figure 2 there for an example $DL(2,2)$.

 Let $P$ denote the transition operator corresponding to simple random walk on $DL(q,r)$. Namely for $x_1x_2\in DL(q,r)$, $p(x_1x_2,y_1y_2)=\frac{\mathbf{1}_{\{x_1x_2\sim y_1y_2\}}}{q+r}$and for functions $h:DL(q,r)\rightarrow \mathbb{R}$, $Ph(x_1x_2)=\sum_{y_1y_2}p(x_1x_2,y_1y_2)h(y_1y_2)$.
 
 A ($P$-)\notion{harmonic} function $h$ is the one that satisfies $Ph=h$.
 
 Let $P_1,P_2$ denote the projection of $P$ on $\mathbb{T}_{q+1}$ and $\mathbb{T}_{r+1}$ respectively: \\
 $p_1(x_1,y_1)=\left\{
 \begin{array}{cc}
 \frac{1}{q+r} & \textnormal{ if } y_1^-=x_1\\
  & \\
 \frac{r}{q+r} & \textnormal{ if }y_1=x_1^-\\
  & \\
 0 & \textnormal{ otherwise,}
 \end{array}
 \right. , \,\, 
 p_2(x_2,y_2)=\left\{
 \begin{array}{cc}
 \frac{1}{q+r} & \textnormal{ if } y_2^-=x_2\\
 & \\
 \frac{q}{q+r} & \textnormal{ if }y_2=x_2^-\\
 & \\
 0 & \textnormal{ otherwise.}
 \end{array}
 \right.
 $\\
 Woess proved the following decomposition theorem about positive harmonic functions.
 \begin{theorem}[Theorem 2.3 of \cite{Woess2005}]\label{thm: decom of PH}
 	If $h$ is a non-negative $P$-harmonic function on $DL(q,r)$, then there are non-negative $P_i$-harmonic functions $h_i$, $i=1,2$ on $\mathbb{T}_{q+1}$ and $\mathbb{T}_{r+1}$ respectively, such that \[
 	h(x_1x_2)=h_1(x_1)+h_2(x_2),\, \, \forall x_1x_2\in DL(q,r)
 	\]
 \end{theorem}
 
 \begin{proposition}\label{prop: FUSF=WUSF on DL graphs}
 	FUSF is the same as WUSF on Diestel--Leader graphs. 
 \end{proposition}
 \begin{proof}
 	By Theorem 7.3 of \cite{BLPS2001}, it suffices to show that every harmonic Dirichlet fuction on $DL(q,r)$ is a constant function and this is an easy consequence of Theorem \ref{thm: decom of PH}. 
 	
 	In fact, suppose there are non-constant harmonic Dirichlet functions on $DL(q,r)$. Then there are non-constant bounded harmonic Dirichlet functions on $DL(q,r)$ (Theorem 3.73 of \cite{Soardi1994}), whence there are non-constant non-negative harmonic Dirichlet  functions on $DL(q,r)$. 
 	
 	Let $h$ be a  non-constant non-negative harmonic Dirichlet  function on $DL(q,r)$.  By Theorem \ref{thm: decom of PH}, there exist non-negative functions $h_1$ and $h_2$ on 
 	$\mathbb{T}_{q+1}$ and $\mathbb{T}_{r+1}$ such that $	h(x_1x_2)=h_1(x_1)+h_2(x_2),\, \, \forall x_1x_2\in DL(q,r)$.
 	
 	Since $h$ is not a constant function, at least one of $h_1,h_2$ is also not a constant. Without loss of generality, we assume that $h_1$ is not a constant. Suppose $x_1,y_1\in\mathbb{T}_{q+1}$ are two neighboring vertices such that $y_1^-=x_1$ and $h_1(x_1)\neq h_1(y_1)$. We first show that for any $z_1\in\mathbb{T}_{q+1}$ such that $z_1^-=x_1$, one has $h_1(z_1)=h_1(y_1)$. Suppose $z_1\neq y_1$. 
 	
 	Since $h$ is a harmonic Dirichlet function, 
 	\begin{eqnarray}
 	\infty&>&\sum_{x_2\in \mathbb{T}_{r+1}:\mathfrak{h}(x_2)=-\mathfrak{h}(x_1)}|h(x_1x_2)-h(y_1x_2^-)|^2+
 	|h(x_1x_2)-h(z_1x_2^-)|^2\nonumber\\
 	&\geq&\sum_{x_2\in \mathbb{T}_{r+1}:\mathfrak{h}(x_2)=-\mathfrak{h}(x_1)}\frac{1}{2}|h(y_1x_2^-)-h(z_1x_2^-)|^2
 	\nonumber\\
    &=&	\sum_{x_2\in \mathbb{T}_{r+1}:\mathfrak{h}(x_2)=-\mathfrak{h}(x_1)}\frac{1}{2}|h_1(y_1)-h_1(z_1)|^2.
 	\end{eqnarray}
Since there are infinitely many $x_2\in \mathbb{T}_{r+1}$ such that $\mathfrak{h}(x_2)=-\mathfrak{h}(x_1)$, we have $h_1(y_1)=h_1(z_1)$. From this and the fact that $h_1$ is $P_1$ harmonic, one has that $h_1$ is constant on each horocycle of $\mathbb{T}_{q+1}$.  Similarly $h_2$ is also a constant on each horocycle of $\mathbb{T}_{r+1}$, whence $h(x_1x_2)$ only depends on which horocycle $x_1$ lies in. Then $h$ must be a constant function on $DL(q,r)$ to have finite Dirichlet energy. This contradicts with the choice that $h$ is a non-constant function. 
 \end{proof}

Now by Proposition \ref{prop: FUSF=WUSF on DL graphs}, the study of FUSF on $DL(qr)$ reduces to WUSF. In particular we know that  each component of FUSF on $DL(q,r)$ with $q>r\geq2$ is one-ended and light almost surely.

 \subsection{FUSF on grandparent graphs}
 
 We first recall the definition of grandparent graphs. For more details see Section 7.1 of \cite{LP2016}. 
 
 Consider a regular tree $\mathbb{T}_{b+1}$ with degree $b+1\geq3$. Let $\xi$ be a fixed end of $\mathbb{T}_{b+1}$. As in the previous subsection, for each $v$ in $\mathbb{T}_{b+1}$, there is a unique ray $\eta_v=(v_0,v_1,v_2,\ldots)$ that represents the end $\xi$ starting at $v_0=v$.  
 We call $v_2$ the $\xi$-\notion{grandparent} of $v$. Throughout this subsection we let $G$ be the graph obtained from $\mathbb{T}_{b+1}$ by adding the edges $(v,v_2)$ between $v$ and its $\xi$-grandparent for all $v\in\mathbb{T}_{b+1}$. It is well known that $G$ is a nonunimodular transitive graph. For two vertices $x,y$ in $G$, we denote by $d_G(x,y),d_{\mathbb{T}}(x,y)$ the graph distance of $x,y$ in $G$ and $\mathbb{T}_{b+1}$ respectively. 
 
Fix a base point $v$ and let  $\eta_v=(v_0,v_1,v_2,\ldots)$ be the unique ray that represents the end $\xi$ starting at $v_0=v$.  We consider the following exhaustion of $G$.  For $n\geq 1$, let $G_n$ be the subgraph of $G$ induced by vertices $\{x:d_{\mathbb{T}}(x,v)\leq n \}$.
 
 For $k,n\geq 1$, let $\mathscr{P}_{k,n}$ denote the set of self-avoiding paths that connect $v_0,v_1$ in $G_n$ with length $k$.
 
 We start with an estimate on the size of $\mathscr{P}_{k,n}$. 
 
\begin{lemma}\label{lem: number of paths connecting two neighbors}
	If $1\leq k\leq n+1$, then  $	|\mathscr{P}_{k,n}|\asymp  b^k$.
	 If $k\geq n+2$, then $\mathscr{P}_{k,n}=\emptyset$.
\end{lemma}
 \begin{proof}
 	For a self-avoiding path $\pi=(w_0,w_1,\ldots,w_k)$ in $G$, if $\Delta(w_i,w_{i+1})=b^{-2}$, then we say the step from $w_i$ to $w_{i+1}$ is \notion{downward} 2 levels, and denote it by $w_i\stackrel{-2}{\rightarrow}w_{i+1}$. Similarly we define downward 1 level, upward 1 level and upward 2 levels. If $w_i\stackrel{-1}{\rightarrow}w_{i+1}$ or $w_i\stackrel{+1}{\rightarrow}w_{i+1}$, then the edge $e=(w_i,w_{i+1})$ is an edge in $\mathbb{T}_{b+1}$, we call it a \notion{tree} edge; if $w_i\stackrel{-2}{\rightarrow}w_{i+1}$ or $w_i\stackrel{+2}{\rightarrow}w_{i+1}$, then the edge $e=(w_i,w_{i+1})$ is an edge connecting a vertex to its grandparent, we call it a \notion{grandparent} edge.
 	
 	The following relations are obvious:
 	\begin{itemize}
 		\item  if $w_{i-1}\stackrel{+1}{\rightarrow}w_{i}$, then $w_i$ is the parent of $w_{i-1}$;
 		\item  if $w_{i-1}\stackrel{+2}{\rightarrow}w_{i}$, then $w_i$ is the grandparent of $w_{i-1}$;
 		\item  if $w_{i-1}\stackrel{-1}{\rightarrow}w_{i}$, then $w_{i-1}$ is the parent of $w_i$;
 		\item  if $w_{i-1}\stackrel{-2}{\rightarrow}w_{i}$, then $w_{i-1}$ is the grandparent of $w_{i}$.
 	\end{itemize}

 	\begin{observation}\label{obs: 5.5 self-avoiding path after x and its parent}
 		If $\pi=(w_0,w_1,\ldots,w_k)\in\mathscr{P}_{k,n}$ has visited $x$ and its parent $x^-$ by time $i$, then the remaining part of the path $(w_{i+1},\ldots,w_k)$ must lie in the connected component of $v_1$ in the induced subgraph of $G\backslash\{x,x^-\}$. This is because $\{x,x^-\}$ is a cutset for $G$ and $\pi$ is self-avoiding. 
 	\end{observation}
 	
 	\begin{claim}\label{claim: forms of paths}
 			A self-avoiding path $\pi=(w_0,w_1,\ldots,w_k)$ in $G$ connecting $v_0$ and $v_1$ has one of the following forms. 
 		If $k$ is even,  then exactly one of the following statements holds. See Figure \ref{fig: five classes adbee' of paths} and \ref{fig: four classes of paths cc'ff'} for an illustration of typical paths of each form.
 		
 		\begin{itemize}
 			\item[(a)] $\pi$ has the form $w_0\stackrel{-2}{\rightarrow}\cdots \stackrel{-2}{\rightarrow}w_{t-1}\stackrel{-1}{\rightarrow}w_{t}\stackrel{+2}{\rightarrow}w_{t+1}\stackrel{+2}{\rightarrow}\cdots \stackrel{+2}{\rightarrow}w_{2t+2}$. Here $w_0=v_0$ and $w_{2t+2}=v_1$.
 			
 			\item[(b)] $\pi$ has the form $w_0\stackrel{+2}{\rightarrow}\cdots \stackrel{+2}{\rightarrow}w_{t}\stackrel{-1}{\rightarrow}w_{t+1}\stackrel{-2}{\rightarrow}\cdots \stackrel{-2}{\rightarrow}w_{2t}$.
 			
 			\item[(b')] $\pi$ has the form $w_0\stackrel{+2}{\rightarrow}\cdots \stackrel{+2}{\rightarrow}w_{t}\stackrel{-1}{\rightarrow}w_{t+1}\stackrel{+2}{\rightarrow}w_{t+2}\stackrel{-2}{\rightarrow}\cdots \stackrel{-2}{\rightarrow}w_{2t+2}$, where $w_{t+1}$ and $w_{t+3}$ are different children of $w_{t+2}$.
 			
 			\item[(c)]  For some positive integers $\alpha,\beta$ such that $k=2+2\alpha+2\beta$, $\pi$ has the form 
 			$w_0\stackrel{+2}{\rightarrow}\cdots \stackrel{+2}{\rightarrow}w_{\alpha}\stackrel{-2}{\rightarrow}w_{\alpha+1}\stackrel{-2}{\rightarrow}\cdots \stackrel{-2}{\rightarrow} w_{\alpha+\beta}\stackrel{-1}{\rightarrow}w_{\alpha+\beta+1}\stackrel{+2}{\rightarrow}\cdots \stackrel{+2}{\rightarrow}w_{\alpha+2\beta+2}\stackrel{-2}{\rightarrow}\cdots\stackrel{-2}{\rightarrow}w_{2+2\alpha+2\beta}$, where $w_{\alpha-1}$ and $w_{\alpha+1}$ have different parents.

 			\item[(c')] For some positive integers $\alpha,\beta$ such that $k=2\alpha+2\beta$, $\pi$ has the form 
 			$w_0\stackrel{+2}{\rightarrow}\cdots \stackrel{+2}{\rightarrow}w_{\alpha}\stackrel{-2}{\rightarrow}w_{\alpha+1}\stackrel{-2}{\rightarrow}\cdots \stackrel{-2}{\rightarrow} w_{\alpha+\beta}\stackrel{-1}{\rightarrow}w_{\alpha+\beta+1}\stackrel{+2}{\rightarrow}\cdots \stackrel{+2}{\rightarrow}w_{\alpha+2\beta+1}\stackrel{-2}{\rightarrow}\cdots\stackrel{-2}{\rightarrow}w_{2\alpha+2\beta}$, where $w_{\alpha-1}$ and $w_{\alpha+1}$ have the same parent (e.g. $w_{\alpha+2\beta+1}$).
 		\end{itemize}

 		If $k$ is odd,  then exactly one of the following statements holds.
 		\begin{itemize}
 			\item[(d)] $\pi$ has the form $w_0\stackrel{-2}{\rightarrow}\cdots \stackrel{-2}{\rightarrow}w_{t}\stackrel{+1}{\rightarrow}w_{t+1}\stackrel{+2}{\rightarrow}w_{t+2}\stackrel{+2}{\rightarrow}\cdots \stackrel{+2}{\rightarrow}w_{2t+1}$. 
 			
 			\item[(e)] $\pi$ has the form $w_0\stackrel{+2}{\rightarrow}\cdots \stackrel{+2}{\rightarrow}w_{t}\stackrel{+1}{\rightarrow}w_{t+1}\stackrel{-2}{\rightarrow}\cdots \stackrel{-2}{\rightarrow}w_{2t+1}$.
 			
 			\item[(f)]  For some positive integers $\alpha,\beta$ such that $k=1+2\alpha+2\beta$, $\pi$ has the form 
 			$w_0\stackrel{+2}{\rightarrow}\cdots \stackrel{+2}{\rightarrow}w_{\alpha}\stackrel{-2}{\rightarrow}w_{\alpha+1}\stackrel{-2}{\rightarrow}\cdots \stackrel{-2}{\rightarrow} w_{\alpha+\beta}\stackrel{+1}{\rightarrow}w_{\alpha+\beta+1}\stackrel{+2}{\rightarrow}\cdots \stackrel{+2}{\rightarrow}w_{\alpha+2\beta+1}\stackrel{-2}{\rightarrow}\cdots\stackrel{-2}{\rightarrow}w_{1+2\alpha+2\beta}$, where $w_{\alpha-1}$ and $w_{\alpha+1}$ have different parents.
 			
 			\item[(f')] For some positive integers $\alpha,\beta$ such that $k=2\alpha+2\beta-1$, $\pi$ has the form 
 			$w_0\stackrel{+2}{\rightarrow}\cdots \stackrel{+2}{\rightarrow}w_{\alpha}\stackrel{-2}{\rightarrow}w_{\alpha+1}\stackrel{-2}{\rightarrow}\cdots \stackrel{-2}{\rightarrow} w_{\alpha+\beta}\stackrel{+1}{\rightarrow}w_{\alpha+\beta+1}\stackrel{+2}{\rightarrow}\cdots \stackrel{+2}{\rightarrow}w_{\alpha+2\beta}\stackrel{-2}{\rightarrow}\cdots\stackrel{-2}{\rightarrow}w_{2\alpha+2\beta-1}$, where $w_{\alpha-1}$ and $w_{\alpha+1}$ have the same parent (e.g. $w_{\alpha+2\beta}$).

 		\end{itemize}

 	\end{claim}

 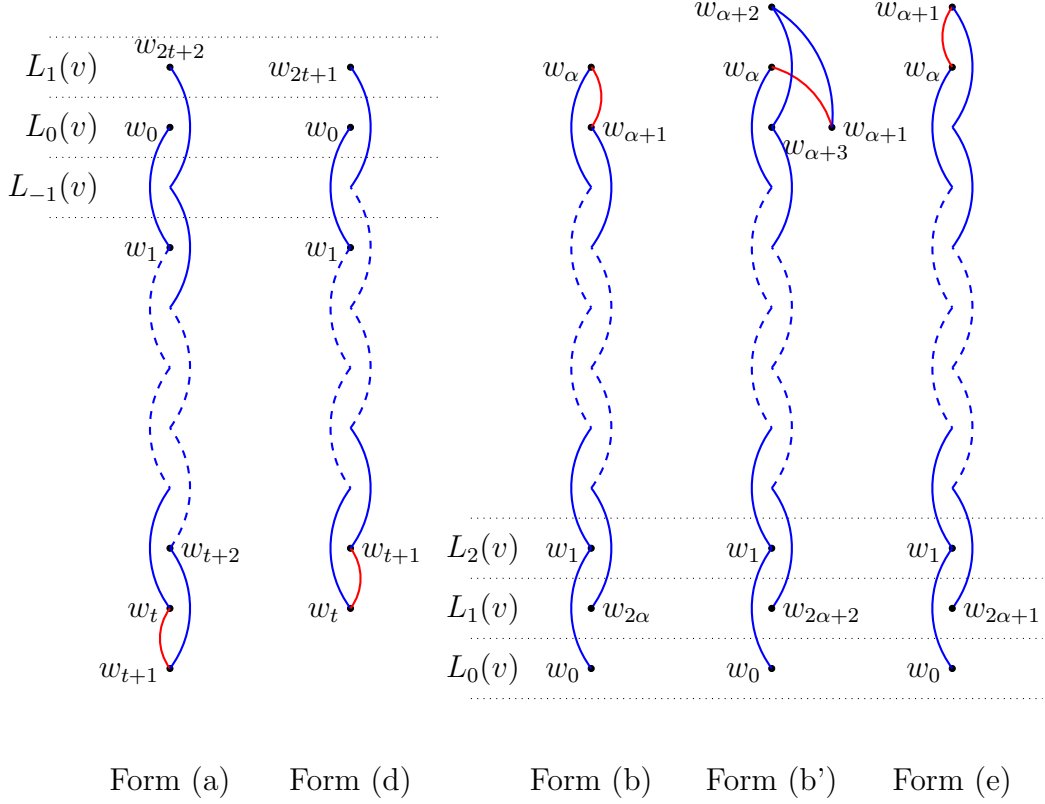
\begin{figure}[h!]
 	\centering
 	\begin{tikzpicture}[scale=0.8, text height=1.5ex,text depth=.25ex] 
 	\draw [dotted,thin] (-2,9.5-1)--(4.5,9.5-1);
 	\draw [dotted,thin] (-2,9.5)--(4.5,9.5);
 	\draw [dotted,thin] (-2,9.5+1)--(4.5,9.5+1);
 	\draw [dotted,thin] (-2,9.5+2)--(4.5,9.5+2);
 	\node[left] at (-1,10) {$L_0(v)$};
 	\node[left] at (-1,11) {$L_1(v)$};
 	\node[left] at (-1,9) {$L_{-1}(v)$};

 	\draw [dotted,thin] (5,0.5)--(14.5,0.5);
 	\draw [dotted,thin] (5,1.5)--(14.5,1.5);
 	\draw [dotted,thin] (5,2.5)--(14.5,2.5);
 	\draw [dotted,thin] (5,3.5)--(14.5,3.5);
 	\node[left] at (6,1) {$L_0(v)$};
 	\node[left] at (6,2) {$L_1(v)$};
 	\node[left] at (6,3) {$L_{2}(v)$};

 	
 	\draw[fill=black] (0,10) circle [radius=0.05];
 	\node[left] at (0,10) {$w_0$};
 	\draw[fill=black] (0,8) circle [radius=0.05];
 	\node[left] at (0,8) {$w_1$};
 	
 	\draw[fill=black] (0,2) circle [radius=0.05];
 	\node[left] at (0,2) {$w_{t}$};
 	
 	\draw[fill=black] (0,1) circle [radius=0.05];
 	\node[left] at (0,1) {$w_{t+1}$};

 	\draw[fill=black] (0,3) circle [radius=0.05];
 	\node[right] at (0,3) {$w_{t+2}$};
 	
 	\draw[fill=black] (0,11) circle [radius=0.05];
 	\node[above] at (0,11) {$w_{2t+2}$};
 	
 	\begin{scope}
 	\tkzDefPoint(0,10){a}\tkzDefPoint(0,8){b}\tkzDefPoint(3,9){C}
 	\tkzCircumCenter(a,b,C)\tkzGetPoint{O}
 	\tkzDrawArc[color=blue,thick](O,a)(b)
 	\end{scope}

 	\begin{scope}
 	\tkzDefPoint(0,8){a}\tkzDefPoint(0,6){b}\tkzDefPoint(3,7){C}
 	\tkzCircumCenter(a,b,C)\tkzGetPoint{O}
 	\tkzDrawArc[color=blue,dashed,thick](O,a)(b)
 	\end{scope}

 	\begin{scope}
 	\tkzDefPoint(0,6){a}\tkzDefPoint(0,4){b}\tkzDefPoint(3,5){C}
 	\tkzCircumCenter(a,b,C)\tkzGetPoint{O}
 	\tkzDrawArc[color=blue,,dashed,thick](O,a)(b)
 	\end{scope}

 	\begin{scope}
 	\tkzDefPoint(0,4){a}\tkzDefPoint(0,2){b}\tkzDefPoint(3,3){C}
 	\tkzCircumCenter(a,b,C)\tkzGetPoint{O}
 	\tkzDrawArc[color=blue,thick](O,a)(b)
 	\end{scope}

 	\begin{scope}
 	\tkzDefPoint(0,2){a}\tkzDefPoint(0,1){b}\tkzDefPoint(1.5,1.5){C}
 	\tkzCircumCenter(a,b,C)\tkzGetPoint{O}
 	\tkzDrawArc[color=red,thick](O,a)(b)
 	\end{scope}

 	\begin{scope}
 	\tkzDefPoint(0,1){a}\tkzDefPoint(0,3){b}\tkzDefPoint(-3,2){C}
 	\tkzCircumCenter(a,b,C)\tkzGetPoint{O}
 	\tkzDrawArc[color=blue,thick](O,a)(b)
 	\end{scope}
 	
 	\begin{scope}
 	\tkzDefPoint(0,3){a}\tkzDefPoint(0,5){b}\tkzDefPoint(-3,4){C}
 	\tkzCircumCenter(a,b,C)\tkzGetPoint{O}
 	\tkzDrawArc[color=blue,dashed,thick](O,a)(b)
 	\end{scope}
 	
 	\begin{scope}
 	\tkzDefPoint(0,5){a}\tkzDefPoint(0,7){b}\tkzDefPoint(-3,6){C}
 	\tkzCircumCenter(a,b,C)\tkzGetPoint{O}
 	\tkzDrawArc[color=blue,dashed,thick](O,a)(b)
 	\end{scope}
 	
 	\begin{scope}
 	\tkzDefPoint(0,7){a}\tkzDefPoint(0,9){b}\tkzDefPoint(-3,8){C}
 	\tkzCircumCenter(a,b,C)\tkzGetPoint{O}
 	\tkzDrawArc[color=blue,thick](O,a)(b)
 	\end{scope}

 	\begin{scope}
 	\tkzDefPoint(0,9){a}\tkzDefPoint(0,11){b}\tkzDefPoint(-3,10){C}
 	\tkzCircumCenter(a,b,C)\tkzGetPoint{O}
 	\tkzDrawArc[color=blue,thick](O,a)(b)
 	\end{scope}

 	\node[below] at (0,-0.5) {Form (a)};

 	
 	\draw[fill=black] (0+3,10) circle [radius=0.05];
 	\node[left] at (0+3,10) {$w_0$};
 	\draw[fill=black] (0+3,8) circle [radius=0.05];
 	\node[left] at (0+3,8) {$w_1$};

 	\draw[fill=black] (0+3,2) circle [radius=0.05];
 	
 	\node[left] at (0+3,2) {$w_{t}$};

 	\draw[fill=black] (0+3,3) circle [radius=0.05];
 	\node[right] at (0+3,3) {$w_{t+1}$};
 	
 	\draw[fill=black] (0+3,11) circle [radius=0.05];
 	\node[left] at (0+3,11) {$w_{2t+1}$};
 	
 	\begin{scope}
 	\tkzDefPoint(0+3,10){a}\tkzDefPoint(0+3,8){b}\tkzDefPoint(3+3,9){C}
 	\tkzCircumCenter(a,b,C)\tkzGetPoint{O}
 	\tkzDrawArc[color=blue,thick](O,a)(b)
 	\end{scope}

 	\begin{scope}
 	\tkzDefPoint(0+3,8){a}\tkzDefPoint(0+3,6){b}\tkzDefPoint(3+3,7){C}
 	\tkzCircumCenter(a,b,C)\tkzGetPoint{O}
 	\tkzDrawArc[color=blue,dashed,thick](O,a)(b)
 	\end{scope}

 	\begin{scope}
 	\tkzDefPoint(0+3,6){a}\tkzDefPoint(0+3,4){b}\tkzDefPoint(3+3,5){C}
 	\tkzCircumCenter(a,b,C)\tkzGetPoint{O}
 	\tkzDrawArc[color=blue,,dashed,thick](O,a)(b)
 	\end{scope}

 	\begin{scope}
 	\tkzDefPoint(0+3,4){a}\tkzDefPoint(0+3,2){b}\tkzDefPoint(3+3,3){C}
 	\tkzCircumCenter(a,b,C)\tkzGetPoint{O}
 	\tkzDrawArc[color=blue,thick](O,a)(b)
 	\end{scope}

 	\begin{scope}
 	\tkzDefPoint(0+3,2){a}\tkzDefPoint(0+3,3){b}\tkzDefPoint(-1.5+3,2.5){C}
 	\tkzCircumCenter(a,b,C)\tkzGetPoint{O}
 	\tkzDrawArc[color=red,thick](O,a)(b)
 	\end{scope}

 	\begin{scope}
 	\tkzDefPoint(0+3,3){a}\tkzDefPoint(0+3,5){b}\tkzDefPoint(-3+3,4){C}
 	\tkzCircumCenter(a,b,C)\tkzGetPoint{O}
 	\tkzDrawArc[color=blue,thick](O,a)(b)
 	\end{scope}
 	
 	\begin{scope}
 	\tkzDefPoint(0+3,5){a}\tkzDefPoint(0+3,7){b}\tkzDefPoint(-3+3,6){C}
 	\tkzCircumCenter(a,b,C)\tkzGetPoint{O}
 	\tkzDrawArc[color=blue,dashed,thick](O,a)(b)
 	\end{scope}
 	
 	\begin{scope}
 	\tkzDefPoint(0+3,7){a}\tkzDefPoint(0+3,9){b}\tkzDefPoint(-3+3,8){C}
 	\tkzCircumCenter(a,b,C)\tkzGetPoint{O}
 	\tkzDrawArc[color=blue,dashed,thick](O,a)(b)
 	\end{scope}

 	\begin{scope}
 	\tkzDefPoint(0+3,9){a}\tkzDefPoint(0+3,11){b}\tkzDefPoint(-3+3,10){C}
 	\tkzCircumCenter(a,b,C)\tkzGetPoint{O}
 	\tkzDrawArc[color=blue,thick](O,a)(b)
 	\end{scope}

 	\node[below] at (0+3,-0.5) {Form (d)};

 	
 	\begin{scope}[shift={(-1,0)}]

 	\draw[fill=black] (8,1) circle [radius=0.05];
 	\node[left] at (8,1) {$w_0$};
 	
 	\draw[fill=black] (8,3) circle [radius=0.05];
 	\node[left] at (8,3) {$w_1$};
 	
 	\begin{scope}
 	\tkzDefPoint(8,1){a}\tkzDefPoint(8,3){b}\tkzDefPoint(8+3,2){C}
 	\tkzCircumCenter(a,b,C)\tkzGetPoint{O}
 	\tkzDrawArc[color=blue,thick](O,b)(a)
 	\end{scope}

 	\begin{scope}
 	\tkzDefPoint(8,3){a}\tkzDefPoint(8,5){b}\tkzDefPoint(8+3,4){C}
 	\tkzCircumCenter(a,b,C)\tkzGetPoint{O}
 	\tkzDrawArc[color=blue,thick](O,b)(a)
 	\end{scope}

 	\begin{scope}
 	\tkzDefPoint(8,5){a}\tkzDefPoint(8,7){b}\tkzDefPoint(8+3,6){C}
 	\tkzCircumCenter(a,b,C)\tkzGetPoint{O}
 	\tkzDrawArc[color=blue,dashed,thick](O,b)(a)
 	\end{scope}
 	
 	\begin{scope}
 	\tkzDefPoint(8,7){a}\tkzDefPoint(8,9){b}\tkzDefPoint(8+3,8){C}
 	\tkzCircumCenter(a,b,C)\tkzGetPoint{O}
 	\tkzDrawArc[color=blue,dashed,thick](O,b)(a)
 	\end{scope}
 	
 	\begin{scope}
 	\tkzDefPoint(8,9){a}\tkzDefPoint(8,11){b}\tkzDefPoint(8+3,10){C}
 	\tkzCircumCenter(a,b,C)\tkzGetPoint{O}
 	\tkzDrawArc[color=blue,thick](O,b)(a)
 	\end{scope}
 	
 	\draw[fill=black] (8,11) circle [radius=0.05];
 	\node[left] at (8,11) {$w_\alpha$};
 	\draw[fill=black] (8,10) circle [radius=0.05];
 	\node[right] at (8,10) {$w_{\alpha+1}$};
 	
 	\begin{scope}
 	\tkzDefPoint(8,11){a}\tkzDefPoint(8,10){b}\tkzDefPoint(8-1.5,10.5){C}
 	\tkzCircumCenter(a,b,C)\tkzGetPoint{O}
 	\tkzDrawArc[color=red,thick](O,b)(a)
 	\end{scope}
 	
 	\begin{scope}
 	\tkzDefPoint(8,10){a}\tkzDefPoint(8,8){b}\tkzDefPoint(8-3,9){C}
 	\tkzCircumCenter(a,b,C)\tkzGetPoint{O}
 	\tkzDrawArc[color=blue,thick](O,b)(a)
 	\end{scope}
 	
 	\begin{scope}
 	\tkzDefPoint(8,8){a}\tkzDefPoint(8,6){b}\tkzDefPoint(8-3,7){C}
 	\tkzCircumCenter(a,b,C)\tkzGetPoint{O}
 	\tkzDrawArc[color=blue,dashed,thick](O,b)(a)
 	\end{scope}
 	
 	\begin{scope}
 	\tkzDefPoint(8,6){a}\tkzDefPoint(8,4){b}\tkzDefPoint(8-3,5){C}
 	\tkzCircumCenter(a,b,C)\tkzGetPoint{O}
 	\tkzDrawArc[color=blue,dashed,thick](O,b)(a)
 	\end{scope}

 	\begin{scope}
 	\tkzDefPoint(8,4){a}\tkzDefPoint(8,2){b}\tkzDefPoint(8-3,3){C}
 	\tkzCircumCenter(a,b,C)\tkzGetPoint{O}
 	\tkzDrawArc[color=blue,thick](O,b)(a)
 	\end{scope}
 	
 	\draw[fill=black] (8,2) circle [radius=0.05];
 	\node[right] at (8,2) {$w_{2\alpha}$};

 	\node[below] at (8,-0.5) {Form (b)};
 	
 	\end{scope}

 	
 	\begin{scope}[shift={(2,0)}]

 	\draw[fill=black] (8,1) circle [radius=0.05];
 	\node[left] at (8,1) {$w_0$};
 	
 	\draw[fill=black] (8,3) circle [radius=0.05];
 	\node[left] at (8,3) {$w_1$};
 	
 	\begin{scope}
 	\tkzDefPoint(8,1){a}\tkzDefPoint(8,3){b}\tkzDefPoint(8+3,2){C}
 	\tkzCircumCenter(a,b,C)\tkzGetPoint{O}
 	\tkzDrawArc[color=blue,thick](O,b)(a)
 	\end{scope}

 	\begin{scope}
 	\tkzDefPoint(8,3){a}\tkzDefPoint(8,5){b}\tkzDefPoint(8+3,4){C}
 	\tkzCircumCenter(a,b,C)\tkzGetPoint{O}
 	\tkzDrawArc[color=blue,thick](O,b)(a)
 	\end{scope}

 	\begin{scope}
 	\tkzDefPoint(8,5){a}\tkzDefPoint(8,7){b}\tkzDefPoint(8+3,6){C}
 	\tkzCircumCenter(a,b,C)\tkzGetPoint{O}
 	\tkzDrawArc[color=blue,dashed,thick](O,b)(a)
 	\end{scope}
 	
 	\begin{scope}
 	\tkzDefPoint(8,7){a}\tkzDefPoint(8,9){b}\tkzDefPoint(8+3,8){C}
 	\tkzCircumCenter(a,b,C)\tkzGetPoint{O}
 	\tkzDrawArc[color=blue,dashed,thick](O,b)(a)
 	\end{scope}
 	
 	\begin{scope}
 	\tkzDefPoint(8,9){a}\tkzDefPoint(8,11){b}\tkzDefPoint(8+3,10){C}
 	\tkzCircumCenter(a,b,C)\tkzGetPoint{O}
 	\tkzDrawArc[color=blue,thick](O,b)(a)
 	\end{scope}
 	
 	\draw[fill=black] (8,11) circle [radius=0.05];
 	\node[left] at (8,11) {$w_\alpha$};
 	
 	\draw[fill=black] (9,10) circle [radius=0.05];
 	\node[right] at (9,10) {$w_{\alpha+1}$};
 	
 	\draw[fill=black] (8,12) circle [radius=0.05];
 	\node[left ] at (8,12) {$w_{\alpha+2}$};

 	\draw[fill=black] (8,10) circle [radius=0.05];
 	\node[right] at (8,10-0.3) {$w_{\alpha+3}$};
 	
 	\begin{scope}
 	\tkzDefPoint(8,11){a}\tkzDefPoint(9,10){b}\tkzDefPoint(8-1.5,10.5){C}
 	\tkzCircumCenter(a,b,C)\tkzGetPoint{O}
 	\tkzDrawArc[color=red,thick](O,b)(a)
 	\end{scope}
 	
 	\begin{scope}
 	\tkzDefPoint(8,12){a}\tkzDefPoint(9,10){b}\tkzDefPoint(8-3,11){C}
 	\tkzCircumCenter(a,b,C)\tkzGetPoint{O}
 	\tkzDrawArc[color=blue,thick](O,b)(a)
 	\end{scope}

 	\begin{scope}
 	\tkzDefPoint(8,12){a}\tkzDefPoint(8,10){b}\tkzDefPoint(8-3,11){C}
 	\tkzCircumCenter(a,b,C)\tkzGetPoint{O}
 	\tkzDrawArc[color=blue,thick](O,b)(a)
 	\end{scope}
 	
 	\begin{scope}
 	\tkzDefPoint(8,10){a}\tkzDefPoint(8,8){b}\tkzDefPoint(8-3,9){C}
 	\tkzCircumCenter(a,b,C)\tkzGetPoint{O}
 	\tkzDrawArc[color=blue,thick](O,b)(a)
 	\end{scope}
 	
 	\begin{scope}
 	\tkzDefPoint(8,8){a}\tkzDefPoint(8,6){b}\tkzDefPoint(8-3,7){C}
 	\tkzCircumCenter(a,b,C)\tkzGetPoint{O}
 	\tkzDrawArc[color=blue,dashed,thick](O,b)(a)
 	\end{scope}
 	
 	\begin{scope}
 	\tkzDefPoint(8,6){a}\tkzDefPoint(8,4){b}\tkzDefPoint(8-3,5){C}
 	\tkzCircumCenter(a,b,C)\tkzGetPoint{O}
 	\tkzDrawArc[color=blue,dashed,thick](O,b)(a)
 	\end{scope}

 	\begin{scope}
 	\tkzDefPoint(8,4){a}\tkzDefPoint(8,2){b}\tkzDefPoint(8-3,3){C}
 	\tkzCircumCenter(a,b,C)\tkzGetPoint{O}
 	\tkzDrawArc[color=blue,thick](O,b)(a)
 	\end{scope}
 	
 	\draw[fill=black] (8,2) circle [radius=0.05];
 	\node[right] at (8,2) {$w_{2\alpha+2}$};

 	\node[below] at (8,-0.5) {Form (b')};
 	
 	\end{scope}

 	
 	\begin{scope}[shift={(2,0)}]

 	\draw[fill=black] (8+3,1) circle [radius=0.05];
 	\node[left] at (8+3,1) {$w_0$};
 	
 	\draw[fill=black] (8+3,3) circle [radius=0.05];
 	\node[left] at (8+3,3) {$w_1$};
 	
 	\begin{scope}
 	\tkzDefPoint(8+3,1){a}\tkzDefPoint(8+3,3){b}\tkzDefPoint(8+3+3,2){C}
 	\tkzCircumCenter(a,b,C)\tkzGetPoint{O}
 	\tkzDrawArc[color=blue,thick](O,b)(a)
 	\end{scope}

 	\begin{scope}
 	\tkzDefPoint(8+3,3){a}\tkzDefPoint(8+3,5){b}\tkzDefPoint(8+3+3,4){C}
 	\tkzCircumCenter(a,b,C)\tkzGetPoint{O}
 	\tkzDrawArc[color=blue,thick](O,b)(a)
 	\end{scope}

 	\begin{scope}
 	\tkzDefPoint(8+3,5){a}\tkzDefPoint(8+3,7){b}\tkzDefPoint(8+3+3,6){C}
 	\tkzCircumCenter(a,b,C)\tkzGetPoint{O}
 	\tkzDrawArc[color=blue,dashed,thick](O,b)(a)
 	\end{scope}
 	
 	\begin{scope}
 	\tkzDefPoint(8+3,7){a}\tkzDefPoint(8+3,9){b}\tkzDefPoint(8+3+3,8){C}
 	\tkzCircumCenter(a,b,C)\tkzGetPoint{O}
 	\tkzDrawArc[color=blue,dashed,thick](O,b)(a)
 	\end{scope}
 	
 	\begin{scope}
 	\tkzDefPoint(8+3,9){a}\tkzDefPoint(8+3,11){b}\tkzDefPoint(8+3+3,10){C}
 	\tkzCircumCenter(a,b,C)\tkzGetPoint{O}
 	\tkzDrawArc[color=blue,thick](O,b)(a)
 	\end{scope}
 	
 	\draw[fill=black] (8+3,11) circle [radius=0.05];
 	\node[left] at (8+3,11) {$w_\alpha$};
 	\draw[fill=black] (8+3,12) circle [radius=0.05];
 	\node[left] at (8+3,12) {$w_{\alpha+1}$};
 	
 	\begin{scope}
 	\tkzDefPoint(8+3,11){a}\tkzDefPoint(8+3,12){b}\tkzDefPoint(8+1.5+3,11.5){C}
 	\tkzCircumCenter(a,b,C)\tkzGetPoint{O}
 	\tkzDrawArc[color=red,thick](O,b)(a)
 	\end{scope}
 	
 	\begin{scope}
 	\tkzDefPoint(8+3,12){a}\tkzDefPoint(8+3,10){b}\tkzDefPoint(8-3+3,11){C}
 	\tkzCircumCenter(a,b,C)\tkzGetPoint{O}
 	\tkzDrawArc[color=blue,thick](O,b)(a)
 	\end{scope}
 	
 	\begin{scope}
 	\tkzDefPoint(8+3,10){a}\tkzDefPoint(8+3,8){b}\tkzDefPoint(8-3+3,9){C}
 	\tkzCircumCenter(a,b,C)\tkzGetPoint{O}
 	\tkzDrawArc[color=blue,thick](O,b)(a)
 	\end{scope}
 	
 	\begin{scope}
 	\tkzDefPoint(8+3,8){a}\tkzDefPoint(8+3,6){b}\tkzDefPoint(8-3+3,7){C}
 	\tkzCircumCenter(a,b,C)\tkzGetPoint{O}
 	\tkzDrawArc[color=blue,dashed,thick](O,b)(a)
 	\end{scope}
 	
 	\begin{scope}
 	\tkzDefPoint(8+3,6){a}\tkzDefPoint(8+3,4){b}\tkzDefPoint(8-3+3,5){C}
 	\tkzCircumCenter(a,b,C)\tkzGetPoint{O}
 	\tkzDrawArc[color=blue,dashed,thick](O,b)(a)
 	\end{scope}

 	\begin{scope}
 	\tkzDefPoint(8+3,4){a}\tkzDefPoint(8+3,2){b}\tkzDefPoint(8-3+3,3){C}
 	\tkzCircumCenter(a,b,C)\tkzGetPoint{O}
 	\tkzDrawArc[color=blue,thick](O,b)(a)
 	\end{scope}
 	
 	\draw[fill=black] (8+3,2) circle [radius=0.05];
 	\node[right] at (8+3,2) {$w_{2\alpha+1}$};

 	\node[below] at (8+3,-0.5) {Form (e)};

 	\end{scope}

 	\end{tikzpicture}
 	\caption{A systematic drawing of the five classes (a),(d),(b),(b'),(e) of paths for $\mathscr{P}_{k,n}$, where $w_0=v$.}
 	\label{fig: five classes adbee' of paths}
 	
 \end{figure}

 \begin{proof}[Proof of Claim \ref{claim: forms of paths}]
The proof of this claim is a case-by-case study and repeated applications of Observation \ref{obs: 5.5 self-avoiding path after x and its parent} with appropriate choices of $x,x^-$.

 Suppose $\pi=(w_0,\ldots,w_k)$ is a self-avoiding path connecting $v_0$ and $v_1$ with length $k$. By parity $\pi$ must use at least one tree edge. 
 
 \textbf{Case one:} (First step is downward 2 levels)\\
 If the first step of $\pi$ is downward 2 levels, i.e. $w_0\stackrel{-2}{\rightarrow}w_1$, then the next step cannot be upward 2 levels, otherwise it will not be self-avoiding. By the same reason the path can only go downward 2 levels  each step before encountering the tree edge.
Let $t$ be the number of steps before encountering a tree edge, i.e., $w_0\stackrel{-2}{\rightarrow}\cdots \stackrel{-2}{\rightarrow}w_{t}$. 

Now look at the next step $w_t\rightarrow w_{t+1}$, it must be a tree edge by the choice of $t$. 
If $w_t\stackrel{-1}{\rightarrow}w_{t+1}$, then apply Observation \ref{obs: 5.5 self-avoiding path after x and its parent} with $x=w_{t+1}$ and $x^-=w_t$ since $w_t\stackrel{-1}{\rightarrow}w_{t+1}$. Note that $v_{1}$ and the descendants of $w_{t+1}$ are in different connected component of $G\backslash\{w_t,w_{t+1} \}$. Thus the next step must be $w_{t+1}\stackrel{+2}{\rightarrow}w_{t+2}$. 
Apply Observation \ref{obs: 5.5 self-avoiding path after x and its parent} with $x=w_{t+2}$ and $x^-=w_{t-1}$, the next step must be $w_{t+2}\stackrel{+2}{\rightarrow}w_{t+3}$. Continue in this way and then $\pi$ must have form (a).

Using similar reasoning,  if $w_t\stackrel{+1}{\rightarrow}w_{t+1}$, to avoid cycle on $\pi$ the remaining steps after time $t+1$ can only be upward 2 levels until reaching $v_1$.

In sum if the first step of $\pi$ is downward 2 levels, then $\pi$ must be of the form (a) or (d).

 	\textbf{Case two:} (First step is upward 2 levels)\\
 	Suppose the first step of $\pi$ is upward 2 levels, i.e. $w_0\stackrel{+2}{\rightarrow}w_1$. Let $\alpha$ be the number of steps of $\pi$ which is upward 2 levels until another type of step is encountering, i.e. $w_0\stackrel{+2}{\rightarrow}\cdots \stackrel{+2}{\rightarrow}w_{\alpha}$ but $w_\alpha\rightarrow w_{\alpha+1}$ is not upward 2 levels.

 	\textbf{Sub-case 2(a)}
 	Suppose the step $w_\alpha\rightarrow w_{\alpha+1}$ is a tree edge. 
 	
 	If $w_{\alpha}\stackrel{+1}{\rightarrow}w_{\alpha+1}$, then applying Observation \ref{obs: 5.5 self-avoiding path after x and its parent} with $x=w_\alpha,x^-=w_{\alpha+1}$ it must be the case that  $w_{\alpha+1}\stackrel{-2}{\rightarrow} w_{\alpha+2}$ and $w_{\alpha+2}$ is an ancestor of $v_1$. So on and so forth, the remaining steps can only be downward 2 levels until hitting $v_1$. Thus in this sub-case the path $\pi$  must be of the form  (e).  
 	
 	 If $w_{\alpha}\stackrel{-1}{\rightarrow}w_{\alpha+1}$ and $w_{\alpha+1}$ is an ancestor of $v_1$ (i.e. $w_{\alpha+1}=w_{\alpha-1}^-$),
 	 then applying Observation \ref{obs: 5.5 self-avoiding path after x and its parent} with $x=w_{\alpha-1},x^-=w_{\alpha+1}$, it must be the case that  $w_{\alpha+1}\stackrel{-2}{\rightarrow} w_{\alpha+2}$ and $w_{\alpha+2}$ is an ancestor of $v_1$. So on and so forth, the remaining steps can only be downward 2 levels until hitting $v_1$. Thus in this sub-case the path $\pi$  must be of the form  (b).

 	 	If $w_{\alpha}\stackrel{-1}{\rightarrow}w_{\alpha+1}$ and $w_{\alpha+1}$ is not an ancestor of $v_1$ (i.e. $w_{\alpha+1}\neq w_{\alpha-1}^-$),
 	 then applying Observation \ref{obs: 5.5 self-avoiding path after x and its parent} with $x=w_{\alpha+1},x^-=w_{\alpha}$, it must be the case that  $w_{\alpha+1}\stackrel{+2}{\rightarrow} w_{\alpha+2}$. Next applying Observation \ref{obs: 5.5 self-avoiding path after x and its parent} with $x=w_\alpha,x^-=w_{\alpha+2}$, it must be the case that  $w_{\alpha+2}\stackrel{-2}{\rightarrow} w_{\alpha+3}$ and $w_{\alpha+3}$ is an ancestor of $v_1$ (i.e. $w_{\alpha-1}^-=w_{\alpha+3}$). Applying Observation \ref{obs: 5.5 self-avoiding path after x and its parent} with $x=w_{\alpha-1},x^-=w_{\alpha+3}$, it must be the case that  $w_{\alpha+3}\stackrel{-2}{\rightarrow} w_{\alpha+4}$ and $w_{\alpha+4}$ is an ancestor of $v_1$.  So on and so forth, the remaining steps can only be downward 2 levels until hitting $v_1$. Thus in this sub-case the path $\pi$  must be of the form  (b').

 	\textbf{Sub-case 2(b)}
 	Suppose the step $w_\alpha\rightarrow w_{\alpha+1}$ is downward 2 levels. To avoid cycle, $w_{\alpha+1}$ must be a different grandchildren of $w_\alpha$ other than $w_{\alpha-1}$. 
 	Similar as in the first case, from the time $\alpha$ the path can only go downward 2 levels each step before encountering a tree edge, say there are $\beta$ such steps.
 	Now the path looks like $w_0\stackrel{+2}{\rightarrow}\cdots \stackrel{+2}{\rightarrow}w_{\alpha}\stackrel{-2}{\rightarrow}w_{\alpha+1}\stackrel{-2}{\rightarrow}\cdots \stackrel{-2}{\rightarrow} w_{\alpha+\beta}\rightarrow\cdots$
 	By the choice of $\beta$, the next step $w_{\alpha+\beta}\rightarrow w_{\alpha+\beta+1}$ is a tree edge. Using the same reasoning as in the first case, after time $\alpha+\beta+1$ the path can only be upward 2 levels each step until it hits the ray  $\eta_v$. Then from that hitting vertex, to avoid cycles the path can only go downward 2 levels each step until hitting $v_1$. 
 	Thus in this sub-case the path $\pi$  must be of the form (c) or (c') or (f) or (f') depending on the parity of the length of $\pi$ and whether $w_{\alpha-1}$ and $w_{\alpha+1}$ has the same parent.
 	We leave the details of this \textbf{Sub-case 2(b)} to the reader.

 	\begin{figure}[h!]
 		\centering
 		\begin{tikzpicture}[scale=0.8, text height=1.5ex,text depth=.25ex] 

 		\draw [dotted, thin] (-2,-1.5)--(18,-1.5);
 		
 		\draw [dotted, thin] (-2,-0.5)--(18,-0.5);
 		
 		\draw [dotted, thin] (-2,0.5)--(18,0.5);
 		
 		\draw [dotted, thin] (-2,1.5)--(18,1.5);
 		
 		\node[left] at (-1,-1) {$L_0(v)$};
 		\node[left] at (-1,0) {$L_1(v)$};
 		\node[left] at (-1,1) {$L_2(v)$};

 		
 		\draw[fill=black] (0,-1) circle [radius=0.05];
 		\node[left] at (0,-1) {$w_0$};
 		\draw[fill=black] (0,1) circle [radius=0.05];
 		\node[left] at (0,1) {$w_1$};
 		\draw[fill=black] (0,3) circle [radius=0.05];
 		\node[left] at (0,3) {$w_2$};
 		
 		\begin{scope}
 		\tkzDefPoint(0,-1){a}\tkzDefPoint(0,1){b}\tkzDefPoint(3,0){C}
 		\tkzCircumCenter(a,b,C)\tkzGetPoint{O}
 		\tkzDrawArc[color=blue,thick](O,b)(a)
 		\end{scope}
 		
 		\begin{scope}
 		\tkzDefPoint(0,1){a}\tkzDefPoint(0,3){b}\tkzDefPoint(3,2){C}
 		\tkzCircumCenter(a,b,C)\tkzGetPoint{O}
 		\tkzDrawArc[color=blue,thick](O,b)(a)
 		\end{scope}
 		
 		\begin{scope}
 		\tkzDefPoint(0,3){a}\tkzDefPoint(0,5){b}\tkzDefPoint(3,4){C}
 		\tkzCircumCenter(a,b,C)\tkzGetPoint{O}
 		\tkzDrawArc[color=blue,dashed,thick](O,b)(a)
 		\end{scope}
 		
 		\begin{scope}
 		\tkzDefPoint(0,5){a}\tkzDefPoint(0,7){b}\tkzDefPoint(3,6){C}
 		\tkzCircumCenter(a,b,C)\tkzGetPoint{O}
 		\tkzDrawArc[color=blue,dashed,thick](O,b)(a)
 		\end{scope}

 		\begin{scope}
 		\tkzDefPoint(0,7){a}\tkzDefPoint(0,9){b}\tkzDefPoint(3,8){C}
 		\tkzCircumCenter(a,b,C)\tkzGetPoint{O}
 		\tkzDrawArc[color=blue,dashed,thick](O,b)(a)
 		\end{scope}
 		
 		\begin{scope}
 		\tkzDefPoint(0,9){a}\tkzDefPoint(0,11){b}\tkzDefPoint(3,10){C}
 		\tkzCircumCenter(a,b,C)\tkzGetPoint{O}
 		\tkzDrawArc[color=blue,thick](O,b)(a)
 		\end{scope}
 		
 		\draw[fill=black] (0,9) circle [radius=0.05];
 		\node[left] at (0,9) {$w_{\alpha-1}$};
 		\draw[fill=black] (0,11) circle [radius=0.05];
 		\node[left] at (0,11) {$w_{\alpha}$};
 		
 		\draw[fill=black] (2,9) circle [radius=0.05];
 		\node[right] at (2,9) {$w_{\alpha+1}$};
 		
 		\begin{scope}
 		\tkzDefPoint(0,11){a}\tkzDefPoint(2,9){b}\tkzDefPoint(1.4,10.4){C}
 		\tkzCircumCenter(a,b,C)\tkzGetPoint{O}
 		\tkzDrawArc[color=blue,thick](O,b)(a)
 		\end{scope}
 		
 		\begin{scope}
 		\tkzDefPoint(2,9){a}\tkzDefPoint(2,7){b}\tkzDefPoint(-1,8){C}
 		\tkzCircumCenter(a,b,C)\tkzGetPoint{O}
 		\tkzDrawArc[color=blue,dashed,thick](O,b)(a)
 		\end{scope}
 		
 		\begin{scope}
 		\tkzDefPoint(2,7){a}\tkzDefPoint(2,5){b}\tkzDefPoint(-1,6){C}
 		\tkzCircumCenter(a,b,C)\tkzGetPoint{O}
 		\tkzDrawArc[color=blue,dashed,thick](O,b)(a)
 		\end{scope}
 		
 		\begin{scope}
 		\tkzDefPoint(2,5){a}\tkzDefPoint(2,3){b}\tkzDefPoint(-1,4){C}
 		\tkzCircumCenter(a,b,C)\tkzGetPoint{O}
 		\tkzDrawArc[color=blue,thick](O,b)(a)
 		\end{scope}
 		
 		\draw[fill=black] (2,3) circle [radius=0.05];
 		\node[right] at (2,3) {$w_{\alpha+\beta}$};
 		\draw[fill=black] (2,2) circle [radius=0.05];
 		\node[right] at (2,2) {$w_{\alpha+\beta+1}$};
 		
 		\begin{scope}
 		\tkzDefPoint(2,3){a}\tkzDefPoint(2,2){b}\tkzDefPoint(0.5,2.5){C}
 		\tkzCircumCenter(a,b,C)\tkzGetPoint{O}
 		\tkzDrawArc[color=red,thick](O,b)(a)
 		\end{scope}

 		\begin{scope}
 		\tkzDefPoint(2,2){a}\tkzDefPoint(2,4){b}\tkzDefPoint(5,3){C}
 		\tkzCircumCenter(a,b,C)\tkzGetPoint{O}
 		\tkzDrawArc[color=blue,thick](O,b)(a)
 		\end{scope}
 		
 		\begin{scope}
 		\tkzDefPoint(2,4){a}\tkzDefPoint(2,6){b}\tkzDefPoint(5,5){C}
 		\tkzCircumCenter(a,b,C)\tkzGetPoint{O}
 		\tkzDrawArc[color=blue,dashed,thick](O,b)(a)
 		\end{scope}

 		\begin{scope}
 		\tkzDefPoint(2,6){a}\tkzDefPoint(2,8){b}\tkzDefPoint(5,7){C}
 		\tkzCircumCenter(a,b,C)\tkzGetPoint{O}
 		\tkzDrawArc[color=blue,dashed,thick](O,b)(a)
 		\end{scope}
 		
 		\begin{scope}
 		\tkzDefPoint(2,8){a}\tkzDefPoint(1,10){b}\tkzDefPoint(4,9){C}
 		\tkzCircumCenter(a,b,C)\tkzGetPoint{O}
 		\tkzDrawArc[color=blue,thick](O,b)(a)
 		\end{scope}
 		
 		\draw[fill=black] (1,10) circle [radius=0.05];
 		\node[right] at (1,10) {$w_{\alpha+2\beta+1}$};
 		
 		\begin{scope}
 		\tkzDefPoint(1,10){a}\tkzDefPoint(0,12){b}\tkzDefPoint(0.7,11){C}
 		\tkzCircumCenter(a,b,C)\tkzGetPoint{O}
 		\tkzDrawArc[color=blue,thick](O,a)(b)
 		\end{scope}
 		
 		\draw[fill=black] (0,12) circle [radius=0.05];
 		\node[left] at (0,12) {$w_{\alpha+2\beta+2}$};

 		\begin{scope}
 		\tkzDefPoint(0,12){a}\tkzDefPoint(0,10){b}\tkzDefPoint(-3,11){C}
 		\tkzCircumCenter(a,b,C)\tkzGetPoint{O}
 		\tkzDrawArc[color=blue,thick](O,b)(a)
 		\end{scope}
 		
 		\begin{scope}
 		\tkzDefPoint(0,10){a}\tkzDefPoint(0,8){b}\tkzDefPoint(-3,9){C}
 		\tkzCircumCenter(a,b,C)\tkzGetPoint{O}
 		\tkzDrawArc[color=blue,thick](O,b)(a)
 		\end{scope}

 		\begin{scope}
 		\tkzDefPoint(0,8){a}\tkzDefPoint(0,6){b}\tkzDefPoint(-3,7){C}
 		\tkzCircumCenter(a,b,C)\tkzGetPoint{O}
 		\tkzDrawArc[color=blue,dashed,thick](O,b)(a)
 		\end{scope}
 		
 		\begin{scope}
 		\tkzDefPoint(0,6){a}\tkzDefPoint(0,4){b}\tkzDefPoint(-3,5){C}
 		\tkzCircumCenter(a,b,C)\tkzGetPoint{O}
 		\tkzDrawArc[color=blue,dashed,thick](O,b)(a)
 		\end{scope}
 		
 		\begin{scope}
 		\tkzDefPoint(0,4){a}\tkzDefPoint(0,2){b}\tkzDefPoint(-3,3){C}
 		\tkzCircumCenter(a,b,C)\tkzGetPoint{O}
 		\tkzDrawArc[color=blue,dashed,thick](O,b)(a)
 		\end{scope}

 		\begin{scope}
 		\tkzDefPoint(0,2){a}\tkzDefPoint(0,0){b}\tkzDefPoint(-3,1){C}
 		\tkzCircumCenter(a,b,C)\tkzGetPoint{O}
 		\tkzDrawArc[color=blue,thick](O,b)(a)
 		\end{scope}
 		\draw[fill=black] (0,0) circle [radius=0.05];
 		\node[right] at (0,0) {$w_{2\alpha+2\beta+2}$};
 		
 		\node[below] at (1,-2) {Form (c)};

 		\draw[fill=red] (0,10) circle [radius=0.07];

 		
 		\begin{scope}[shift={(5,0)}]

 		\draw[fill=black] (0,-1) circle [radius=0.05];
 		\node[left] at (0,-1) {$w_0$};
 		\draw[fill=black] (0,1) circle [radius=0.05];
 		\node[left] at (0,1) {$w_1$};
 		\draw[fill=black] (0,3) circle [radius=0.05];
 		\node[left] at (0,3) {$w_2$};
 		
 		\begin{scope}
 		\tkzDefPoint(0,-1){a}\tkzDefPoint(0,1){b}\tkzDefPoint(3,0){C}
 		\tkzCircumCenter(a,b,C)\tkzGetPoint{O}
 		\tkzDrawArc[color=blue,thick](O,b)(a)
 		\end{scope}
 		
 		\begin{scope}
 		\tkzDefPoint(0,1){a}\tkzDefPoint(0,3){b}\tkzDefPoint(3,2){C}
 		\tkzCircumCenter(a,b,C)\tkzGetPoint{O}
 		\tkzDrawArc[color=blue,thick](O,b)(a)
 		\end{scope}
 		
 		\begin{scope}
 		\tkzDefPoint(0,3){a}\tkzDefPoint(0,5){b}\tkzDefPoint(3,4){C}
 		\tkzCircumCenter(a,b,C)\tkzGetPoint{O}
 		\tkzDrawArc[color=blue,dashed,thick](O,b)(a)
 		\end{scope}
 		
 		\begin{scope}
 		\tkzDefPoint(0,5){a}\tkzDefPoint(0,7){b}\tkzDefPoint(3,6){C}
 		\tkzCircumCenter(a,b,C)\tkzGetPoint{O}
 		\tkzDrawArc[color=blue,dashed,thick](O,b)(a)
 		\end{scope}

 		\begin{scope}
 		\tkzDefPoint(0,7){a}\tkzDefPoint(0,9){b}\tkzDefPoint(3,8){C}
 		\tkzCircumCenter(a,b,C)\tkzGetPoint{O}
 		\tkzDrawArc[color=blue,dashed,thick](O,b)(a)
 		\end{scope}
 		
 		\begin{scope}
 		\tkzDefPoint(0,9){a}\tkzDefPoint(1,11){b}\tkzDefPoint(3,10){C}
 		\tkzCircumCenter(a,b,C)\tkzGetPoint{O}
 		\tkzDrawArc[color=blue,thick](O,b)(a)
 		\end{scope}
 		
 		\draw[fill=black] (0,9) circle [radius=0.05];
 		\node[left] at (0,9) {$w_{\alpha-1}$};
 		\draw[fill=black] (1,11) circle [radius=0.05];
 		\node[above] at (1,11) {$w_{\alpha}$};
 		
 		\draw[fill=black] (2,9) circle [radius=0.05];
 		\node[right] at (2,9) {$w_{\alpha+1}$};
 		
 		\begin{scope}
 		\tkzDefPoint(1,11){a}\tkzDefPoint(2,9){b}\tkzDefPoint(-1,10){C}
 		\tkzCircumCenter(a,b,C)\tkzGetPoint{O}
 		\tkzDrawArc[color=blue,thick](O,b)(a)
 		\end{scope}
 		
 		\begin{scope}
 		\tkzDefPoint(2,9){a}\tkzDefPoint(2,7){b}\tkzDefPoint(-1,8){C}
 		\tkzCircumCenter(a,b,C)\tkzGetPoint{O}
 		\tkzDrawArc[color=blue,dashed,thick](O,b)(a)
 		\end{scope}
 		
 		\begin{scope}
 		\tkzDefPoint(2,7){a}\tkzDefPoint(2,5){b}\tkzDefPoint(-1,6){C}
 		\tkzCircumCenter(a,b,C)\tkzGetPoint{O}
 		\tkzDrawArc[color=blue,dashed,thick](O,b)(a)
 		\end{scope}
 		
 		\begin{scope}
 		\tkzDefPoint(2,5){a}\tkzDefPoint(2,3){b}\tkzDefPoint(-1,4){C}
 		\tkzCircumCenter(a,b,C)\tkzGetPoint{O}
 		\tkzDrawArc[color=blue,thick](O,b)(a)
 		\end{scope}
 		
 		\draw[fill=black] (2,3) circle [radius=0.05];
 		\node[right] at (2,3) {$w_{\alpha+\beta}$};
 		\draw[fill=black] (2,2) circle [radius=0.05];
 		\node[right] at (2,2) {$w_{\alpha+\beta+1}$};
 		
 		\begin{scope}
 		\tkzDefPoint(2,3){a}\tkzDefPoint(2,2){b}\tkzDefPoint(0.5,2.5){C}
 		\tkzCircumCenter(a,b,C)\tkzGetPoint{O}
 		\tkzDrawArc[color=red,thick](O,b)(a)
 		\end{scope}

 		\begin{scope}
 		\tkzDefPoint(2,2){a}\tkzDefPoint(2,4){b}\tkzDefPoint(5,3){C}
 		\tkzCircumCenter(a,b,C)\tkzGetPoint{O}
 		\tkzDrawArc[color=blue,thick](O,b)(a)
 		\end{scope}
 		
 		\begin{scope}
 		\tkzDefPoint(2,4){a}\tkzDefPoint(2,6){b}\tkzDefPoint(5,5){C}
 		\tkzCircumCenter(a,b,C)\tkzGetPoint{O}
 		\tkzDrawArc[color=blue,dashed,thick](O,b)(a)
 		\end{scope}

 		\begin{scope}
 		\tkzDefPoint(2,6){a}\tkzDefPoint(2,8){b}\tkzDefPoint(5,7){C}
 		\tkzCircumCenter(a,b,C)\tkzGetPoint{O}
 		\tkzDrawArc[color=blue,dashed,thick](O,b)(a)
 		\end{scope}
 		
 		\begin{scope}
 		\tkzDefPoint(2,8){a}\tkzDefPoint(1,10){b}\tkzDefPoint(1.2,9){C}
 		\tkzCircumCenter(a,b,C)\tkzGetPoint{O}
 		\tkzDrawArc[color=blue,thick](O,b)(a)
 		\end{scope}
 		
 		\draw[fill=black] (1,10) circle [radius=0.05];
 		\node[right] at (1,10) {$w_{\alpha+2\beta+1}$};
 		
 		\begin{scope}
 		\tkzDefPoint(1,10){a}\tkzDefPoint(0,8){b}\tkzDefPoint(0.7,9){C}
 		\tkzCircumCenter(a,b,C)\tkzGetPoint{O}
 		\tkzDrawArc[color=blue,thick](O,b)(a)
 		\end{scope}


 		\begin{scope}
 		\tkzDefPoint(0,8){a}\tkzDefPoint(0,6){b}\tkzDefPoint(-3,7){C}
 		\tkzCircumCenter(a,b,C)\tkzGetPoint{O}
 		\tkzDrawArc[color=blue,dashed,thick](O,b)(a)
 		\end{scope}
 		
 		\begin{scope}
 		\tkzDefPoint(0,6){a}\tkzDefPoint(0,4){b}\tkzDefPoint(-3,5){C}
 		\tkzCircumCenter(a,b,C)\tkzGetPoint{O}
 		\tkzDrawArc[color=blue,dashed,thick](O,b)(a)
 		\end{scope}
 		
 		\begin{scope}
 		\tkzDefPoint(0,4){a}\tkzDefPoint(0,2){b}\tkzDefPoint(-3,3){C}
 		\tkzCircumCenter(a,b,C)\tkzGetPoint{O}
 		\tkzDrawArc[color=blue,dashed,thick](O,b)(a)
 		\end{scope}

 		\begin{scope}
 		\tkzDefPoint(0,2){a}\tkzDefPoint(0,0){b}\tkzDefPoint(-3,1){C}
 		\tkzCircumCenter(a,b,C)\tkzGetPoint{O}
 		\tkzDrawArc[color=blue,thick](O,b)(a)
 		\end{scope}
 		\draw[fill=black] (0,0) circle [radius=0.05];
 		\node[right] at (0,0) {$w_{2\alpha+2\beta}$};
 		
 		\node[below] at (1,-2) {Form (c')};

 		\end{scope}


 		\begin{scope}[shift={(10,0)}]

 		\draw[fill=black] (0,-1) circle [radius=0.05];
 		\node[left] at (0,-1) {$w_0$};
 		\draw[fill=black] (0,1) circle [radius=0.05];
 		\node[left] at (0,1) {$w_1$};
 		\draw[fill=black] (0,3) circle [radius=0.05];
 		\node[left] at (0,3) {$w_2$};
 		
 		\begin{scope}
 		\tkzDefPoint(0,-1){a}\tkzDefPoint(0,1){b}\tkzDefPoint(3,0){C}
 		\tkzCircumCenter(a,b,C)\tkzGetPoint{O}
 		\tkzDrawArc[color=blue,thick](O,b)(a)
 		\end{scope}
 		
 		\begin{scope}
 		\tkzDefPoint(0,1){a}\tkzDefPoint(0,3){b}\tkzDefPoint(3,2){C}
 		\tkzCircumCenter(a,b,C)\tkzGetPoint{O}
 		\tkzDrawArc[color=blue,thick](O,b)(a)
 		\end{scope}
 		
 		\begin{scope}
 		\tkzDefPoint(0,3){a}\tkzDefPoint(0,5){b}\tkzDefPoint(3,4){C}
 		\tkzCircumCenter(a,b,C)\tkzGetPoint{O}
 		\tkzDrawArc[color=blue,dashed,thick](O,b)(a)
 		\end{scope}
 		
 		\begin{scope}
 		\tkzDefPoint(0,5){a}\tkzDefPoint(0,7){b}\tkzDefPoint(3,6){C}
 		\tkzCircumCenter(a,b,C)\tkzGetPoint{O}
 		\tkzDrawArc[color=blue,dashed,thick](O,b)(a)
 		\end{scope}

 		\begin{scope}
 		\tkzDefPoint(0,7){a}\tkzDefPoint(0,9){b}\tkzDefPoint(3,8){C}
 		\tkzCircumCenter(a,b,C)\tkzGetPoint{O}
 		\tkzDrawArc[color=blue,dashed,thick](O,b)(a)
 		\end{scope}
 		
 		\begin{scope}
 		\tkzDefPoint(0,9){a}\tkzDefPoint(0,11){b}\tkzDefPoint(3,10){C}
 		\tkzCircumCenter(a,b,C)\tkzGetPoint{O}
 		\tkzDrawArc[color=blue,thick](O,b)(a)
 		\end{scope}
 		
 		\draw[fill=black] (0,9) circle [radius=0.05];
 		\node[left] at (0,9) {$w_{\alpha-1}$};
 		\draw[fill=black] (0,11) circle [radius=0.05];
 		\node[left] at (0,11) {$w_{\alpha}$};
 		
 		\draw[fill=black] (2,9) circle [radius=0.05];
 		\node[right] at (2,9) {$w_{\alpha+1}$};
 		
 		\begin{scope}
 		\tkzDefPoint(0,11){a}\tkzDefPoint(2,9){b}\tkzDefPoint(1.4,10.4){C}
 		\tkzCircumCenter(a,b,C)\tkzGetPoint{O}
 		\tkzDrawArc[color=blue,thick](O,b)(a)
 		\end{scope}
 		
 		\begin{scope}
 		\tkzDefPoint(2,9){a}\tkzDefPoint(2,7){b}\tkzDefPoint(-1,8){C}
 		\tkzCircumCenter(a,b,C)\tkzGetPoint{O}
 		\tkzDrawArc[color=blue,dashed,thick](O,b)(a)
 		\end{scope}
 		
 		\begin{scope}
 		\tkzDefPoint(2,7){a}\tkzDefPoint(2,5){b}\tkzDefPoint(-1,6){C}
 		\tkzCircumCenter(a,b,C)\tkzGetPoint{O}
 		\tkzDrawArc[color=blue,dashed,thick](O,b)(a)
 		\end{scope}
 		
 		\begin{scope}
 		\tkzDefPoint(2,5){a}\tkzDefPoint(2,3){b}\tkzDefPoint(-1,4){C}
 		\tkzCircumCenter(a,b,C)\tkzGetPoint{O}
 		\tkzDrawArc[color=blue,thick](O,b)(a)
 		\end{scope}
 		
 		\draw[fill=black] (2,3) circle [radius=0.05];
 		\node[below] at (2,3) {$w_{\alpha+\beta}$};
 		\draw[fill=black] (2,4) circle [radius=0.05];
 		\node[left] at (2,4) {$w_{\alpha+\beta+1}$};
 		
 		\begin{scope}
 		\tkzDefPoint(2,3){a}\tkzDefPoint(2,4){b}\tkzDefPoint(3.5,3.5){C}
 		\tkzCircumCenter(a,b,C)\tkzGetPoint{O}
 		\tkzDrawArc[color=red,thick](O,b)(a)
 		\end{scope}

 		\begin{scope}
 		\tkzDefPoint(2,4){a}\tkzDefPoint(2,6){b}\tkzDefPoint(5,5){C}
 		\tkzCircumCenter(a,b,C)\tkzGetPoint{O}
 		\tkzDrawArc[color=blue,dashed,thick](O,b)(a)
 		\end{scope}

 		\begin{scope}
 		\tkzDefPoint(2,6){a}\tkzDefPoint(2,8){b}\tkzDefPoint(5,7){C}
 		\tkzCircumCenter(a,b,C)\tkzGetPoint{O}
 		\tkzDrawArc[color=blue,dashed,thick](O,b)(a)
 		\end{scope}
 		
 		\begin{scope}
 		\tkzDefPoint(2,8){a}\tkzDefPoint(1,10){b}\tkzDefPoint(4,9){C}
 		\tkzCircumCenter(a,b,C)\tkzGetPoint{O}
 		\tkzDrawArc[color=blue,thick](O,b)(a)
 		\end{scope}
 		
 		\draw[fill=black] (1,10) circle [radius=0.05];
 		\node[right] at (1,10) {$w_{\alpha+2\beta}$};
 		
 		\begin{scope}
 		\tkzDefPoint(1,10){a}\tkzDefPoint(0,12){b}\tkzDefPoint(0.7,11){C}
 		\tkzCircumCenter(a,b,C)\tkzGetPoint{O}
 		\tkzDrawArc[color=blue,thick](O,a)(b)
 		\end{scope}
 		
 		\draw[fill=black] (0,12) circle [radius=0.05];
 		\node[left] at (0,12) {$w_{\alpha+2\beta+1}$};

 		\begin{scope}
 		\tkzDefPoint(0,12){a}\tkzDefPoint(0,10){b}\tkzDefPoint(-3,11){C}
 		\tkzCircumCenter(a,b,C)\tkzGetPoint{O}
 		\tkzDrawArc[color=blue,thick](O,b)(a)
 		\end{scope}
 		
 		\begin{scope}
 		\tkzDefPoint(0,10){a}\tkzDefPoint(0,8){b}\tkzDefPoint(-3,9){C}
 		\tkzCircumCenter(a,b,C)\tkzGetPoint{O}
 		\tkzDrawArc[color=blue,thick](O,b)(a)
 		\end{scope}

 		\begin{scope}
 		\tkzDefPoint(0,8){a}\tkzDefPoint(0,6){b}\tkzDefPoint(-3,7){C}
 		\tkzCircumCenter(a,b,C)\tkzGetPoint{O}
 		\tkzDrawArc[color=blue,dashed,thick](O,b)(a)
 		\end{scope}
 		
 		\begin{scope}
 		\tkzDefPoint(0,6){a}\tkzDefPoint(0,4){b}\tkzDefPoint(-3,5){C}
 		\tkzCircumCenter(a,b,C)\tkzGetPoint{O}
 		\tkzDrawArc[color=blue,dashed,thick](O,b)(a)
 		\end{scope}
 		
 		\begin{scope}
 		\tkzDefPoint(0,4){a}\tkzDefPoint(0,2){b}\tkzDefPoint(-3,3){C}
 		\tkzCircumCenter(a,b,C)\tkzGetPoint{O}
 		\tkzDrawArc[color=blue,dashed,thick](O,b)(a)
 		\end{scope}

 		\begin{scope}
 		\tkzDefPoint(0,2){a}\tkzDefPoint(0,0){b}\tkzDefPoint(-3,1){C}
 		\tkzCircumCenter(a,b,C)\tkzGetPoint{O}
 		\tkzDrawArc[color=blue,thick](O,b)(a)
 		\end{scope}
 		\draw[fill=black] (0,0) circle [radius=0.05];
 		\node[right] at (0,0) {$w_{2\alpha+2\beta+1}$};

 		\draw[fill=red] (0,10) circle [radius=0.07];
 		
 		\node[below] at (1,-2) {Form (f)};

 		\end{scope}


 		\begin{scope}[shift={(15,0)}]

 		\draw[fill=black] (0,-1) circle [radius=0.05];
 		\node[left] at (0,-1) {$w_0$};
 		\draw[fill=black] (0,1) circle [radius=0.05];
 		\node[left] at (0,1) {$w_1$};
 		\draw[fill=black] (0,3) circle [radius=0.05];
 		\node[left] at (0,3) {$w_2$};
 		
 		\begin{scope}
 		\tkzDefPoint(0,-1){a}\tkzDefPoint(0,1){b}\tkzDefPoint(3,0){C}
 		\tkzCircumCenter(a,b,C)\tkzGetPoint{O}
 		\tkzDrawArc[color=blue,thick](O,b)(a)
 		\end{scope}
 		
 		\begin{scope}
 		\tkzDefPoint(0,1){a}\tkzDefPoint(0,3){b}\tkzDefPoint(3,2){C}
 		\tkzCircumCenter(a,b,C)\tkzGetPoint{O}
 		\tkzDrawArc[color=blue,thick](O,b)(a)
 		\end{scope}
 		
 		\begin{scope}
 		\tkzDefPoint(0,3){a}\tkzDefPoint(0,5){b}\tkzDefPoint(3,4){C}
 		\tkzCircumCenter(a,b,C)\tkzGetPoint{O}
 		\tkzDrawArc[color=blue,dashed,thick](O,b)(a)
 		\end{scope}
 		
 		\begin{scope}
 		\tkzDefPoint(0,5){a}\tkzDefPoint(0,7){b}\tkzDefPoint(3,6){C}
 		\tkzCircumCenter(a,b,C)\tkzGetPoint{O}
 		\tkzDrawArc[color=blue,dashed,thick](O,b)(a)
 		\end{scope}

 		\begin{scope}
 		\tkzDefPoint(0,7){a}\tkzDefPoint(0,9){b}\tkzDefPoint(3,8){C}
 		\tkzCircumCenter(a,b,C)\tkzGetPoint{O}
 		\tkzDrawArc[color=blue,dashed,thick](O,b)(a)
 		\end{scope}
 		
 		\begin{scope}
 		\tkzDefPoint(0,9){a}\tkzDefPoint(1,11){b}\tkzDefPoint(3,10){C}
 		\tkzCircumCenter(a,b,C)\tkzGetPoint{O}
 		\tkzDrawArc[color=blue,thick](O,b)(a)
 		\end{scope}
 		
 		\draw[fill=black] (0,9) circle [radius=0.05];
 		\node[left] at (0,9) {$w_{\alpha-1}$};
 		\draw[fill=black] (1,11) circle [radius=0.05];
 		\node[above] at (1,11) {$w_{\alpha}$};
 		
 		\draw[fill=black] (2,9) circle [radius=0.05];
 		\node[right] at (2,9) {$w_{\alpha+1}$};
 		
 		\begin{scope}
 		\tkzDefPoint(1,11){a}\tkzDefPoint(2,9){b}\tkzDefPoint(-1,10){C}
 		\tkzCircumCenter(a,b,C)\tkzGetPoint{O}
 		\tkzDrawArc[color=blue,thick](O,b)(a)
 		\end{scope}
 		
 		\begin{scope}
 		\tkzDefPoint(2,9){a}\tkzDefPoint(2,7){b}\tkzDefPoint(-1,8){C}
 		\tkzCircumCenter(a,b,C)\tkzGetPoint{O}
 		\tkzDrawArc[color=blue,dashed,thick](O,b)(a)
 		\end{scope}
 		
 		\begin{scope}
 		\tkzDefPoint(2,7){a}\tkzDefPoint(2,5){b}\tkzDefPoint(-1,6){C}
 		\tkzCircumCenter(a,b,C)\tkzGetPoint{O}
 		\tkzDrawArc[color=blue,dashed,thick](O,b)(a)
 		\end{scope}
 		
 		\begin{scope}
 		\tkzDefPoint(2,5){a}\tkzDefPoint(2,3){b}\tkzDefPoint(-1,4){C}
 		\tkzCircumCenter(a,b,C)\tkzGetPoint{O}
 		\tkzDrawArc[color=blue,thick](O,b)(a)
 		\end{scope}
 		
 		\draw[fill=black] (2,3) circle [radius=0.05];
 		\node[below] at (2,3) {$w_{\alpha+\beta}$};
 		\draw[fill=black] (2,4) circle [radius=0.05];
 		\node[left] at (2,4) {$w_{\alpha+\beta+1}$};
 		
 		\begin{scope}
 		\tkzDefPoint(2,3){a}\tkzDefPoint(2,4){b}\tkzDefPoint(3.5,3.5){C}
 		\tkzCircumCenter(a,b,C)\tkzGetPoint{O}
 		\tkzDrawArc[color=red,thick](O,b)(a)
 		\end{scope}

 		\begin{scope}
 		\tkzDefPoint(2,4){a}\tkzDefPoint(2,6){b}\tkzDefPoint(5,5){C}
 		\tkzCircumCenter(a,b,C)\tkzGetPoint{O}
 		\tkzDrawArc[color=blue,dashed,thick](O,b)(a)
 		\end{scope}

 		\begin{scope}
 		\tkzDefPoint(2,6){a}\tkzDefPoint(2,8){b}\tkzDefPoint(5,7){C}
 		\tkzCircumCenter(a,b,C)\tkzGetPoint{O}
 		\tkzDrawArc[color=blue,dashed,thick](O,b)(a)
 		\end{scope}
 		
 		\begin{scope}
 		\tkzDefPoint(2,8){a}\tkzDefPoint(1,10){b}\tkzDefPoint(1.2,9){C}
 		\tkzCircumCenter(a,b,C)\tkzGetPoint{O}
 		\tkzDrawArc[color=blue,thick](O,b)(a)
 		\end{scope}
 		
 		\draw[fill=black] (1,10) circle [radius=0.05];
 		\node[right] at (1,10) {$w_{\alpha+2\beta}$};
 		
 		\begin{scope}
 		\tkzDefPoint(1,10){a}\tkzDefPoint(0,8){b}\tkzDefPoint(0.7,9){C}
 		\tkzCircumCenter(a,b,C)\tkzGetPoint{O}
 		\tkzDrawArc[color=blue,thick](O,b)(a)
 		\end{scope}


 		\begin{scope}
 		\tkzDefPoint(0,8){a}\tkzDefPoint(0,6){b}\tkzDefPoint(-3,7){C}
 		\tkzCircumCenter(a,b,C)\tkzGetPoint{O}
 		\tkzDrawArc[color=blue,dashed,thick](O,b)(a)
 		\end{scope}
 		
 		\begin{scope}
 		\tkzDefPoint(0,6){a}\tkzDefPoint(0,4){b}\tkzDefPoint(-3,5){C}
 		\tkzCircumCenter(a,b,C)\tkzGetPoint{O}
 		\tkzDrawArc[color=blue,dashed,thick](O,b)(a)
 		\end{scope}
 		
 		\begin{scope}
 		\tkzDefPoint(0,4){a}\tkzDefPoint(0,2){b}\tkzDefPoint(-3,3){C}
 		\tkzCircumCenter(a,b,C)\tkzGetPoint{O}
 		\tkzDrawArc[color=blue,dashed,thick](O,b)(a)
 		\end{scope}

 		\begin{scope}
 		\tkzDefPoint(0,2){a}\tkzDefPoint(0,0){b}\tkzDefPoint(-3,1){C}
 		\tkzCircumCenter(a,b,C)\tkzGetPoint{O}
 		\tkzDrawArc[color=blue,thick](O,b)(a)
 		\end{scope}
 		\draw[fill=black] (0,0) circle [radius=0.05];
 		\node[right] at (0,0) {$w_{2\alpha+2\beta-1}$};
 		
 		\node[below] at (1,-2) {Form (f')};

 		\end{scope}

 		\end{tikzpicture}
 		\caption{A systematic drawing of the four classes (c),(c'),(f),(f') of paths for $\mathscr{P}_{k,n}$, where $w_0=v$. In form (c) and (f), the parent of $w_{\alpha-1}$ is denoted by relatively large and  red node.}
 		\label{fig: four classes of paths cc'ff'}
 		
 	\end{figure}
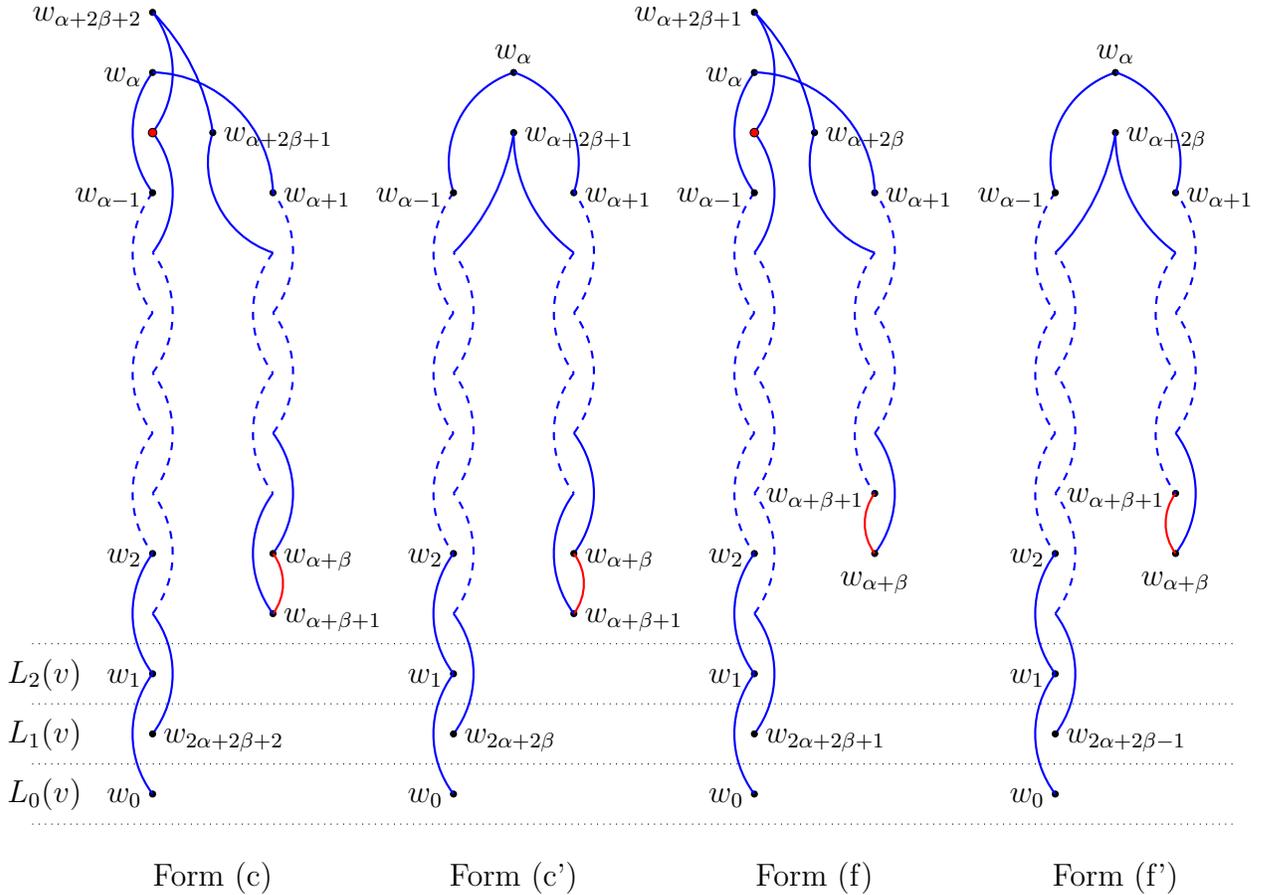
 	
 	In particular from the above analysis, we see that each self-avoiding path $\pi$ that connecting $v_0,v_1$ uses exactly one tree edge. 
 \end{proof}

 	Also from the forms of the path in Claim \ref{claim: forms of paths}, we see that if a self-avoiding path $\pi=(w_0,\ldots,w_k)$ connecting $v_0$ and $v_1$ with length $k$ is one of the forms (a), (b'),(c), (d) and (f), then 
 	\[
 	\max\{d_{\mathbb{T}}(w_0,w_i):1\leq i\leq k \}=(k-1)\vee 1;
 	\] 
 	If $\pi$ is of the form (c') or (f'), then 
 	\[
 	\max\{d_{\mathbb{T}}(w_0,w_i):1\leq i\leq k \}=k+1.
 	\] 
 	If $\pi$ is of the form (b) or (e), then 
 	\[
 	\max\{d_{\mathbb{T}}(w_0,w_i):1\leq i\leq k \}=k.
 	\] 
 	
 	Therefore, if $k\geq n+2$ then $\mathscr{P}_{k,n}=\emptyset$.

 Suppose $k$ is even and $k=2t+2$. The number of paths of the form (a) is $b^{2t+1}=b^{k-1}$ because for the first $t$ steps one has $b^2$ choices for the grandchildren and $b$ choices for $w_{t+1}$. Once $w_0,\ldots,w_{t+1}$ are fixed, the remaining vertices $w_{t+2},\ldots,w_{2t+2}$ are fixed. The number of paths of the form (b) is just $1$. The number of paths of the form (b') is  $b-1$ since one has $b-1$ choices for $w_{\alpha+1}$. The number of paths of the form (c) is $\sum_{\beta=1}^{t-1}(b-1)b^{2\beta}$ (Given $\alpha,\beta$, $w_0,\ldots,w_\alpha$ are fixed, one has $b(b-1)$ choices for $w_{\alpha+1}$, and $b^2$ choices for each of $w_{\alpha+2},\ldots,w_{\alpha+\beta}$ and $b$ choices for  $w_{\alpha+\beta+1}$. Once $w_0,\ldots,w_{\alpha+\beta+1}$ are fixed, the remaining path  is fixed). Using similar counting argument, we have that the number of paths of the form (c') is $\sum_{\beta=1}^{t}(b-1)b^{2\beta-1}$.

 	In sum if $k\in[1,n)$ is even, then $|\mathscr{P}_{k,n}|=b^{k-1}+1+\sum_{\beta=1}^{(k-4)/2}(b-1)b^{2\beta}+\sum_{\beta=1}^{(k-2)/2}(b-1)b^{2\beta-1}\asymp b^k$. 
 	
 	Similarly if $k\in[1,n)$ is odd, one also has 
 	$|\mathscr{P}_{k,n}|\asymp b^k$. 
 	
 	If $k=n\textnormal{ or }n+1$, the estimate 	$|\mathscr{P}_{k,n}|\asymp b^k$ is also true because we only need to deduct the contributions of cases (b), (c'), (e) and (f') from the previous expression for $|\mathscr{P}_{k,n}|$.
 \end{proof}

 Next we recall some notation and results regarding loop-erased random walk on a graph from \cite{Lawler1999}.
 
 Let $S(t)$ be a discrete time Markov chain on a countable state space $X$ with transition probabilities $p(x,y)$. For a subset $A$ of $X$, define the hitting time 
 \[\tau_A:=\inf\{t\geq0: S(t)\in A \} \]
 and the Green function 
 \[
 G(x,A):=\sum_{j=0}^{\infty}\mathbf{P}^x[S(j)=x;S(t)\notin A,t=0,\ldots,j],
 \]
 where $\mathbf{P}^x$ denote the measure of the Markov chain $S(t)$ started from $x$. 
 
 Fix a base point $o\in X$ and we assume that $S(t)$ is irreducible. Let $B\subset X$ be a subset with $\mathbf{P}^o[\tau_B<\infty]=1$. Suppose $o\notin B$, start the Markov chain at $o$, let it run until it hits $B$, and then erase the loops chronologically. Then we get a probability measure $\mu$ on the self-avoiding paths from $o$ to $B$. In particular, if $w=[w_0,\ldots,w_k]$ is such a self-avoiding path from $o$ to $B$ with $w_0,\ldots,w_{k-1}\notin B$, then 
 \be\label{eq: probability of a loop-erased path}
 \mu(w)=\left[\prod_{j=1}^{k}p(w_{j-1},w_j)\right]\cdot
        \left[\prod_{j=0}^{k-1}G(w_{j},B\cup A_{j-1}) \right],
 \ee
 where $A_{-1}=\emptyset,A_j=\{w_0,\ldots,w_j \}$ for $j\geq 1$ (see Proposition 3.2 of \cite{Lawler2018}). 
 
 \begin{proposition}[Proposition 3.3 of  \cite{Lawler2018}]\label{prop: symmetry about ordering}
 	Suppose $S$ is irreducible, $B\neq\emptyset$ and $w_0,\ldots,w_{k-1}\notin B$, define
 	\[
 	f(w_0,\ldots,w_{k-1};B)=\prod_{j=0}^{k-1}G(w_j,B\cup A_{j-1}),
 	\]
 	where $A_{-1}=\emptyset$ and $A_j=\{w_0,\ldots,w_j \}$.
 	Then $f(w_0,\ldots,w_{k-1};B)$ is a symmetric function of $w_0,\ldots,w_{k-1}$.
 \end{proposition}

\begin{proposition}\label{prop: FUSF on grandparent graph is a tree}
	The FUSF on the grandparent graph $G$ is connected almost surely.
\end{proposition}
\begin{proof}
	Recall $G_n$ is the sub-graph of $G$  induced by vertices $\{x:d_{\mathbb{T}}(x,v_0)\leq n \}$. Let $T_n$ be a uniform spanning tree on $G_n$.

	We will show that there exists positive constants $c_{23},c_{24}$ such that
	\be\label{eq: exponential decay of large distance in FUSF on grandparent graph}
	\mathbb{P}[d_{{T}_n}(v_0,v_1)\geq k]\leq  c_{24}\exp{(-c_{23}k)}.
	\ee
	
	Start a simple random walk on $G_n$ from $v_0$ and stop at the first hit of $v_1$, then loop-erase this random walk path. Then the self-avoiding path we get has the law of the path from $v_0$ to $v_1$ in $T_n$. Let $S(t)$ be the simple random walk on $G_n$, let $o$ be $v_0$ and $B=\{v_1\}$. By \eqref{eq: probability of a loop-erased path}
	\be\label{eq: 5.3-tail prob for large distance in FUSF on grandparent graph}
	\mathbb{P}[d_{{T}_n}(v_0,v_1)\geq k]= \sum_{j=k}^{\infty}\sum_{w\in \mathscr{P}_{j,n}}{\mu (w)}. 
	\ee
	By Lemma \ref{lem: number of paths connecting two neighbors}, \eqref{eq: exponential decay of large distance in FUSF on grandparent graph} holds trivially for $k\geq n+2$. In the following we assume $k\leq n+1$.

	Notice that in our case the first product in \eqref{eq: probability of a loop-erased path} equals $\prod_{j=0}^{k-1}\frac{1}{\textnormal{deg}_{G_n}(w_j)}$, which is again symmetric for $w_0,\ldots,w_{k-1}$. Hence by Proposition \ref{prop: symmetry about ordering}, we have $\mu(w)$ is a symmetric function of $w_0,\ldots,w_{k-1}$. Using the connection between effective resistance and Green function (see Proposition 2.1 of \cite{LP2016}) one has that 
	for any reordering $w_0',\ldots,w_{k-1}'$ of $w_0,\ldots,w_{k-1}$ 
	\be\label{eq: probability of a self-avoiding path in terms of resistance}
	\mu(w)=\prod_{j=0}^{k-1}\mathscr{R}(w_j'\leftrightarrow B\cup A_{j-1}'),
	\ee
	where $A_{-1}'=\emptyset$, $A_j'=\{w_0',\ldots,w_j'\}$ and $\mathscr{R}(w_j'\leftrightarrow B\cup A_{j-1}')$ is the effective resistance from $w_j'$ to $B\cup A_{j-1}'$ in $G_n$.

For $w=(w_0,\ldots,w_k)\in\mathscr{P}_{k,n} $, since $w$ has one of the forms listed in the proof of Lemma \ref{lem: number of paths connecting two neighbors}, using a case by case analysis it is easy to see that we can reorder $w_0,\ldots,w_{k-1}$ as $w_0',\ldots,w_{k-1}'$ such that except at most $8$ $j'$s in $\{0,\ldots,k-1\}$, one has $w_{j-2}',w_{j-1}'$ are the grandparent and parent of $w_j'$ and the grandchildren of $w_j'$ are also in $G_n$. For such an ordering and an index $j$ such that $w_{j-2}',w_{j-1}'$ are the grandparent and parent of $w_j'$ and the grandchildren of $w_j'$ are in $G_n$, one has that 
\[
\mathscr{R}(w_j'\leftrightarrow B\cup A_{j-1}')\leq \frac{b+4}{b^2+4b+8}<\frac{1}{b}.
\]
Indeed, from the local structure one has that 
\[
\mathscr{R}(w_j'\leftrightarrow B\cup A_{j-1}')\leq 
\mathscr{R}(w_j'\leftrightarrow \{w_{j-2}',w_{j-1}'\})\stackrel{\textnormal{Fig.}\ref{fig:local estimate for resistance}}{\leq} 
\frac{1}{1+1+b\cdot\frac{1}{1+\frac{1}{1+b/2}} }=\frac{b+4}{b^2+4b+8}<\frac{1}{b}.
\]
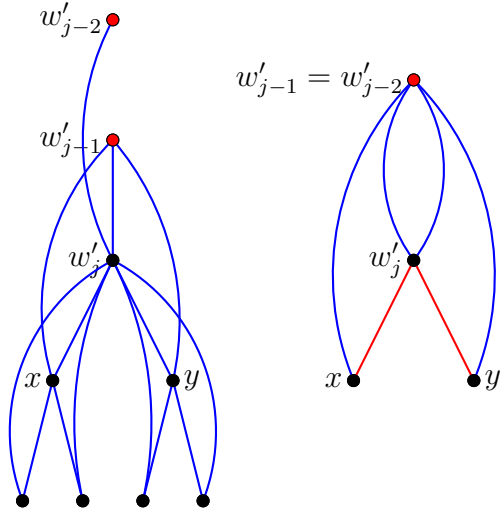
\begin{figure}[h!]
	
	\centering
	\begin{tikzpicture}[scale=0.8]
	
	\draw[fill=black] (5,5) circle [radius=0.10];
	\node [left] at (5,5) { $w_j'$};
	\draw[fill=red] (5,7) circle [radius=0.10];
	\node [left] at (5,7) { $w_{j-1}'$};
	\draw[fill=red] (5,9) circle [radius=0.10];
	\node [left] at (5,9) { $w_{j-2}'$};
	
	\node [left] at (4,3) {$x$};
	\node [right] at (6,3) {$y$};

	\draw[fill=black] (4,3) circle [radius=0.10];
	\draw[fill=black] (6,3) circle [radius=0.10];
	\draw[fill=black] (3.5,1) circle [radius=0.10];
	\draw[fill=black] (4.5,1) circle [radius=0.10];
	\draw[fill=black] (5.5,1) circle [radius=0.10];
	\draw[fill=black] (6.5,1) circle [radius=0.10];

	\begin{scope}
	\tkzDefPoint(5,9){a}\tkzDefPoint(5,5){b}\tkzDefPoint(4.5,7){C}
	\tkzCircumCenter(a,b,C)\tkzGetPoint{O}
	\tkzDrawArc[color=blue,thick](O,a)(b)
	\end{scope}
	\draw [color=blue,thick](5,5)--(5,7);
	
	\draw [color=blue,thick](5,5)--(4,3);
	\draw [color=blue,thick](5,5)--(6,3);
	
	\draw [color=blue,thick](4,3)--(3.5,1);
	\draw [color=blue,thick](4,3)--(4.5,1);
	\draw [color=blue,thick](6,3)--(5.5,1);
	\draw [color=blue,thick](6,3)--(6.5,1);
	
	\begin{scope}
	\tkzDefPoint(5,7){a}\tkzDefPoint(4,3){b}\tkzDefPoint(3.9,5){C}
	\tkzCircumCenter(a,b,C)\tkzGetPoint{O}
	\tkzDrawArc[color=blue,thick](O,a)(b)
	\end{scope}

	\begin{scope}
	\tkzDefPoint(5,7){a}\tkzDefPoint(6,3){b}\tkzDefPoint(6,5){C}
	\tkzCircumCenter(a,b,C)\tkzGetPoint{O}
	\tkzDrawArc[color=blue,thick](O,b)(a)
	\end{scope}

	\begin{scope}
	\tkzDefPoint(5,5){a}\tkzDefPoint(3.5,1){b}\tkzDefPoint(3.4,3){C}
	\tkzCircumCenter(a,b,C)\tkzGetPoint{O}
	\tkzDrawArc[color=blue,thick](O,a)(b)
	\end{scope}

	\begin{scope}
	\tkzDefPoint(5,5){a}\tkzDefPoint(4.5,1){b}\tkzDefPoint(4.4,3){C}
	\tkzCircumCenter(a,b,C)\tkzGetPoint{O}
	\tkzDrawArc[color=blue,thick](O,a)(b)
	\end{scope}

	\begin{scope}
	\tkzDefPoint(5,5){a}\tkzDefPoint(5.5,1){b}\tkzDefPoint(5.6,3){C}
	\tkzCircumCenter(a,b,C)\tkzGetPoint{O}
	\tkzDrawArc[color=blue,thick](O,b)(a)
	\end{scope}

	\begin{scope}
	\tkzDefPoint(5,5){a}\tkzDefPoint(6.5,1){b}\tkzDefPoint(6.6,3){C}
	\tkzCircumCenter(a,b,C)\tkzGetPoint{O}
	\tkzDrawArc[color=blue,thick](O,b)(a)
	\end{scope}

	\draw[fill=black] (10,5) circle [radius=0.10];
	\node [left] at (10,5) { $w_j'$};
	\draw[fill=red] (10,8) circle [radius=0.10];
	\node [left] at (10,8) { $w_{j-1}'=w_{j-2}'$};
	\begin{scope}
	\tkzDefPoint(10,8){a}\tkzDefPoint(10,5){b}\tkzDefPoint(9.5,6.5){C}
	\tkzCircumCenter(a,b,C)\tkzGetPoint{O}
	\tkzDrawArc[color=blue,thick](O,a)(b)
	\end{scope}

	\begin{scope}
	\tkzDefPoint(10,8){a}\tkzDefPoint(10,5){b}\tkzDefPoint(10.5,6.5){C}
	\tkzCircumCenter(a,b,C)\tkzGetPoint{O}
	\tkzDrawArc[color=blue,thick](O,b)(a)
	\end{scope}

	\draw[fill=black] (9,3) circle [radius=0.10];
	\node [left] at (9,3) { $x$};
	\draw[fill=black] (11,3) circle [radius=0.10];
	\node [right] at (11,3) { $y$};
	\draw [color=red,thick] (10,5)--(9,3);
	\draw [color=red,thick] (10,5)--(11,3);
	
	\begin{scope}
	\tkzDefPoint(10,8){a}\tkzDefPoint(9,3){b}\tkzDefPoint(9,6.5){C}
	\tkzCircumCenter(a,b,C)\tkzGetPoint{O}
	\tkzDrawArc[color=blue,thick](O,a)(b)
	\end{scope}

	\begin{scope}
	\tkzDefPoint(10,8){a}\tkzDefPoint(11,3){b}\tkzDefPoint(11,6.5){C}
	\tkzCircumCenter(a,b,C)\tkzGetPoint{O}
	\tkzDrawArc[color=blue,thick](O,b)(a)
	\end{scope}

	
	\draw[fill=black] (5,5) circle [radius=0.10];
	\draw[fill=red] (5,7) circle [radius=0.10];
	\draw[fill=red] (5,9) circle [radius=0.10];
	
	
	\draw[fill=black] (4,3) circle [radius=0.10];
	\draw[fill=black] (6,3) circle [radius=0.10];
	\draw[fill=black] (3.5,1) circle [radius=0.10];
	\draw[fill=black] (4.5,1) circle [radius=0.10];
	\draw[fill=black] (5.5,1) circle [radius=0.10];
	\draw[fill=black] (6.5,1) circle [radius=0.10];

	\draw[fill=black] (10,5) circle [radius=0.10];
	\draw[fill=red] (10,8) circle [radius=0.10];
	\draw[fill=black] (9,3) circle [radius=0.10];
	\draw[fill=black] (11,3) circle [radius=0.10];

	\end{tikzpicture}

	\caption{Local structure for estimating $\mathscr{R}(w_j'\leftrightarrow \{w_{j-2}',w_{j-1}'\})$, the right half is a network reduction with red edges with conductance  ${1+b/2}$. All blue edges has conductance $1$.}
	\label{fig:local estimate for resistance}
	
\end{figure}

	 For the other $j$'s we use trivial estimates 
	 \[
	 \mathscr{R}(w_j'\leftrightarrow B\cup A_{j-1}')\leq k+1.
	 \]
	
	Now by \eqref{eq: probability of a self-avoiding path in terms of resistance} there exists a constant 
  $c_{25}>0$ such that
	\[
	\max_{w\in \mathscr{P}_{k,n}}\mu(w)\leq (k+1)^8\left(\frac{b+4}{b^2+4b+8}\right)^{k-8}
	\leq c_{25} \left(\frac{b+4}{b^2+4b+7}\right)^{k}.
	\]
	
	This combined with \eqref{eq: 5.3-tail prob for large distance in FUSF on grandparent graph} and Lemma \ref{lem: number of paths connecting two neighbors} yields \eqref{eq: exponential decay of large distance in FUSF on grandparent graph}.

	Since $G_n$ is an exhaustion of $G$, $T_n$ converges weakly to the FUSF on $G$. Thus \eqref{eq: exponential decay of large distance in FUSF on grandparent graph} implies that in the FUSF sample $\mathfrak{F}^f$ on $G$, the probability that $d_{\mathfrak{F}^f}(v_0,v_1)\geq k$ decays exponentially in $k$. In particular, $v_0$ and $v_1$ are connected almost surely in $\mathfrak{F}^f$. By transitivity, almost surely any vertex of $G$ is in the same connected component as its parent in $\mathfrak{F}^f$. Therefore the FUSF $\mathfrak{F}^f$ is connected almost surely. 
\end{proof}

Now we know that the FUSF sample $\mathfrak{F}^f$ on a grandparent graph $G$ is just a tree. Next we consider the branching number of $\mathfrak{F}^f$.

\begin{proposition}\label{prop: branching number larger than one for FUSF on grandparent graph}
	The FUSF on a grandparent graph $G$ has branching number strictly larger than one.
\end{proposition}
\begin{proof}
	From Proposition \ref{prop: FUSF on grandparent graph is a tree} we know that there is at least one tree edge in $\mathfrak{F}^f$, otherwise there would be at least two trees in $\mathfrak{F}^f$.
	
	For $x\in V(G)$, let $y$ be the parent of $x$ and $z$ be the grandparent of $x$. Let $x_1,\ldots,x_b$ be the children of $x$. Note that $\{x,y\}$ is a cutset for $G$, and $G\backslash\{x,y\}$ has $b+1$ connected components, one containing $z$ and other $b$ ones each containing a unique child of $x$. We denote the connected component containing $z$ by $K_0(x)$ and the  connected components containing $x_i$ by  $K_i(x)$ for $i=1,\ldots,b$.
	Let 
	$\widehat{K}_i(x)$ be the subgraph of $G$ induced by $K_i(x)\cup\{x,y\}$ for $i=0,1,\ldots,b$.

	Conditioned on the event that the tree edge $e=(x,y)\in \mathfrak{F}^f$, one has the following observation: 
	\be\label{eq: prop 5.7-conditional distribution given a tree edge}
	\textnormal{The conditional distributions of } \mathfrak{F}^f\cap \widehat{K}_i(x),i=0,1,\ldots,b \textnormal{ are independent.}
	\ee
	To see \eqref{eq: prop 5.7-conditional distribution given a tree edge}, consider an exhaustion of $G$ by the balls $B(x,n)$. Sample UST on $B(x,n)$ using Wilson's algorithm  with $\mathsf{F}_0=\{x\}$ and the first simple random walk $X^y$ starting from $y$. Then the edge $e=(y,x)$ belongs to the uniform spanning tree $\textnormal{UST}(B(x,n))$ if and only if  the first time $X^y$ hits $x$ via the edge $(y,x)$. Suppose the event  $e=(y,x)\in \textnormal{UST}(B(x,n))$ occurs and continue Wilson's algorithm. Suppose we already obtain $F_{k-1}$ and the next simple random walk starts from $w$. Since $\{x,y\}$ is a cutset for $B(x,n)$, the simple random walk $X^w$ will hit $F_{k-1}$ before visiting other connected component of $B(x,n)\backslash\{x,y\}$. Therefore the $b+1$ trees $\textnormal{UST}(B(x,n))\cap \widehat{K}_i(x),i=0,1,\ldots,b$ are independent conditioned on $e=(y,x)\in \textnormal{UST}(B(x,n))$. By the definition of FUSF, one has \eqref{eq: prop 5.7-conditional distribution given a tree edge}.

	Next we show that conditioned on the event $e=(x,y)\in\mathfrak{F}^f$, for each $i=1,\ldots,b$, almost surely there is another tree edge in $\mathfrak{F}^f\cap \widehat{K}_i(x)$.
	 If not, we define a mass transport as follows: \[f(u,v):=\mathbf{1}_{ \{ v \textnormal{ is the nearest ancestor of } u \textnormal{ such that } (v,v^-)\in\mathfrak{F}^f \} },\]
		where $v^-$ denotes the parent of $v$.
		
		Then the mass sent out from a vertex is at most one. But if conditioned on the event $e=(x,y)\in \mathfrak{F}^f$, with positive probability there is no other tree edge in $\mathfrak{F}^f\cap \widehat{K}_i(x)$, then $x$ will receive infinite mass with positive probability. This contradicts with the tilted mass-transport principle. 
	
Thus there is a large constant $M>0$ such that the following inequality holds
	\be\label{eq: prop 5.7-2}
	\mathbb{P}\big[ \exists e'=(v,v^-)\in\mathfrak{F}^f\cap \widehat{K}_i(x) \textnormal{ s.t. }d_{\mathfrak{F}^f}(v,x)\in [1,M]  \bigm| e=(x,y)\in \mathfrak{F}^f\big]>\frac{1}{b}.
	\ee


Also observe that conditioned on $e=(x,y)\in \mathfrak{F}^f$ and $e'=(v,v^-)\in\mathfrak{F}^f\cap \widehat{K}_i(x)$, the conditional distributions of $\mathfrak{F}^f\cap \widehat{K}_i(v)$ are independent and are the same as the distribution $\mathfrak{F}^f\cap \widehat{K}_i(v)$ conditioned only on $e'=(v,v^-)\in\mathfrak{F}^f$.  

Conditioned on the tree edge $e=(x,y)\in\mathfrak{F}^f$, we call a component $\widehat{K}_i(x)$ \textbf{good} if there exists some edge $e'=(v,v^-)\in\mathfrak{F}^f\cap \widehat{K}_i(x) \textnormal{ s.t. }d_{\mathfrak{F}^f}(v,x)\in [1,M] $ and we also call this edge $e'=(v,v^-)$ a \textbf{good} edge for $e=(x,y)$. In a good component we pick an arbitrary good edge. Note that the size of good edges form  a supercritical Galton--Watson tree in light of \eqref{eq: prop 5.7-2}.
With positive probability the supercritical Galton--Watson tree has branching number bigger than one \cite[Cor. 5.10]{LP2016}. Since in $\mathfrak{F}^f$ two good edges in neighboring generations of the Galton--Watson tree have distance at most $M$,
   with positive probability $\mathfrak{F}^f$ also has branching number strictly larger than one. Since the branching number of  $\mathfrak{F}^f$ is a constant almost surely (Theorem 10.18 of \cite{LP2016}), it is strictly larger than one almost surely. 
	\end{proof}

\begin{remark}
Proposition \ref{prop: branching number larger than one for FUSF on grandparent graph} is non-trivial in the sense that there exist spanning trees of $G$ with branching number equals to one. In fact one can even find recurrent spanning trees of $G$. It is also of interest to find the exact value of this branching number. 
\end{remark}
One natural further question is the following:
\begin{question}\label{ques: invariant recurrent spanning tree for grandparent graphs}
	Are there invariant spanning trees of the grandparent graphs with branching number equals to 1 or even being recurrent? 
\end{question}


\begin{remark}\label{rem: connected FUSF on tree cross finite graph}
	The conclusion of Proposition \ref{prop: FUSF on grandparent graph is a tree} and \ref{prop: branching number larger than one for FUSF on grandparent graph} also hold for the Cartesian product $\mathbb{T}_{b+1}\square \mathbb{Z}_2$, where the $\mathbb{Z}_2$-edges in  $\mathbb{T}_{b+1}\square \mathbb{Z}_2$ will play the role of the tree edges in the above proofs. We conjecture that for any finite connected graph $H$, almost surely the FUSF on the Cartesian product $\mathbb{T}_{b+1}\square H$ is connected. Note that every tree in the FUSF on the Cartesian product $\mathbb{T}_{b+1}\square H$  has branching number larger than one since $\mathbb{T}_{b+1}\square H$ is unimodular \cite{TomAsaf2017,Timar2018indis}. 
\end{remark}

Quite recently Pete and Tim\'{a}r  \cite{PeteTimar2020} disproved the conjecture in the Remark \ref{rem: connected FUSF on tree cross finite graph}.

 \section[ \texorpdfstring{FUSF on free products of  nonunimodular transitive graphs with $\mathbb{Z}_2$}{FUSF on free products of nonunimodular transitive graphs with one edge}]{FUSF on  free products of nonunimodular transitive graphs with $\mathbb{Z}_2$}\label{sec: 6}

The free product of two Cayley graphs is well known. More generally one can define the free product of two transitive graphs $G_1,G_2$. For more details, see the description on page 2349 of \cite{Timar2006}. 

Suppose $G_0=(V_0,E_0)$ is a nonunimodular transitive graph and $\mathbb{Z}_2$ is the graph of two vertices connected by one edge. We now give the detailed definition of the free product $G_0*\mathbb{Z}_2$ just like \cite{Timar2006}. First take a copy of $G_0$ and countably many copies of $\mathbb{Z}_2$. Fix a bijection from the vertices of this copy of $G_0$ to the copies of $\mathbb{Z}_2$. Identify each vertex of this copy of $G_0$ with an arbitrary vertex in its image under the bijection. Call the resulting graph $H_1$, it is  formed by attaching an edge to each vertex of the copy of $G_1$. Let $I_1$ denote the set of  vertices on the  $\mathbb{Z}_2$ edges that are not identified with a vertex of $G_0$. Fix a bijection between $I_1$ and countably many new copies of $G_0$. Identify every vertex of $I_1$ with an arbitrary vertex in its image by the bijection to obtain a graph $H_2$. So $H_2$ is formed by attaching a $G_0$ copy to each vertex in $I_1$. 
Continue this process similarly, given $H_i$, and if $I_i$ is the set of vertices in $H_i$ that were not born by identification in some previous steps, then fix a bijection between $I_i$ to a set of infinitely many copies of $G_0$ if $i$ is odd  (or infinitely many copies of $\mathbb{Z}_2$
if $i$ is even). Identify every vertex in $I_i$ with an arbitrary vertex in  its image by the bijection to obtain $H_{i+1}$. If we view $H_i$ as a subgraph of $H_{i+1}$, then finally the free product $G_0*\mathbb{Z}_2:=\bigcup H_i$. It is easy to see $G$ is still a transitive graph. 
Also we call an edge in the free product $G_0*\mathbb{Z}_2$ a $G_0$-edge if its two endpoints lies in the same copy of $G_0$ in the above construction. Similarly we call the other edges  $\mathbb{Z}_2$-edges. The free product $G=G_0*\mathbb{Z}_2$ can be viewed as countably many disjoint copies of $G_0$, say $G_0^{(i)}\,(i=1,2,\ldots)$, connected by countably many $\mathbb{Z}_2$ edges in a certain way. So we can fix a labeling of the vertices of $G$: $V(G)=\{(v,i)\colon v\in V_0,i\in\mathbb{Z}^+\}$, where $(v,i)$ indicates  the vertex $v$ in the copy $G_0^{(i)}$. Let $\Phi_i: G_0^{(i)}\rightarrow G_0$ denote the projection map, namely for each vertex $(u,i)\in V(G_0^{(i)})$ and edge $e=((v,i),(w,i))\in E(G_0^{(i)})$, $\Phi_i((u,i))=u\in V_0$ and $\Phi_i(e)=(v,w)\in E_0$.


\begin{definition}
	Suppose $G_0$ is a transitive graph with a closed subgroup $\Gamma$ of automorphisms that acts transitively on $G_0$. Let $G=G_0*\mathbb{Z}_2$. Suppose $\omega_0$ is a $\Gamma$-invariant percolation process on $G_0$. We view $\omega_0$ as a random subgraph of $G_0$. For each   $G_0$ copy $G_0^{(i)}$ in $G$, we take an independent percolation $\omega_i$ on $G_0^{(i)}$ such that $\Phi_i(\omega_i)$ has the same law as $\omega_0$. 
	Let $\omega$ be the union of these independent percolation subgraphs and all the $\mathbb{Z}_2$-edges. We call $\omega$ the \notion{free product percolation} of $\omega_0$ on $G=G_0*\mathbb{Z}_2$. 
\end{definition}

Suppose $(G_0,\Gamma)$ is a nonunimodular transitive pair. For the free product $G=G_0*\mathbb{Z}_2$ defined above, let $\widetilde{\Gamma}$ be the set of automorphisms $\gamma_G$ of $G$ with the following two properties:
	\begin{itemize}
		\item $\gamma_G$ maps $\mathbb{Z}_2$ edges to $\mathbb{Z}_2$ edges and $G_0$ copies to $G_0$ copies;
		\item for each $i\in\mathbb{Z}^+$, the bijection $\varphi_i=\varphi_i(\gamma_G)$ is an element of $\Gamma$, where $\varphi_i$ is given as follows. The automorphism $\gamma_G$ maps the copy $G_0^{(i)}$ to another $G_0$ copy, say $G_0^{(j)}$. Then restricted to $G_0^{(i)}$ and its image $G_0^{(j)}$, $\gamma_G$ is a bijection  from $G_0^{(i)}$ to $G_0^{(j)}$. Projecting $\gamma_G\big|_{G_0^{(i)}}$ to the first coordinate one gets the bijection  $\varphi_i:G_0\rightarrow G_0$, namely $\varphi_i=\Phi_j\circ \gamma_G\big|_{G_0^{(i)}} \circ \Phi_i^{-1}$. 
	\end{itemize}

\begin{observation}\label{obse: action on the free product}
	Suppose $(G_0,\Gamma)$ is a nonunimodular transitive pair. 	Let $\widetilde{\Gamma}$ be the set of automorphisms of $G=G_0*\mathbb{Z}_2$ as defined above.  Then $\widetilde{\Gamma}$ has the following properties.
	\begin{enumerate}
		\item The set $\widetilde{\Gamma}$ is a closed subgroup of $\textnormal{Aut}(G)$ and $\widetilde{\Gamma}$ acts on $G$ transitively .
		\item For each copy $G_0^{(i)}$ and any two vertices $(u,i),(v,i)$ in that copy, one has $|\widetilde{\Gamma}_{(u,i)}(v,i)|=|\Gamma_uv|$. In particular by Lemma \ref{lem:haar}  this implies that $\widetilde{\Gamma}$ is also nonunimodular.
		\item If $(u,i)$ and $(v,j)$ are connected in $G$ by a $\mathbb{Z}_2$-edge, then $|\widetilde{\Gamma}_{(u,i)}(v,j)|=1$.
	\end{enumerate}
\end{observation}

\begin{proposition}\label{prop: free product percolation}
	Suppose $(G_0,\Gamma)$ is a nonunimodular transitive pair.  Suppose $\omega_0$ is a $\Gamma$-invariant percolation process on $G_0$ and   $\omega$ is the {free product percolation} of $\omega_0$ on $G=G_0*\mathbb{Z}_2$. If almost surely every connected component of $\omega_0$ is infinite, then each connected component of $\omega$ is $\widetilde{\Gamma}$-heavy and has branching number bigger than one. 
\end{proposition}

Given Proposition \ref{prop: free product percolation} it is easy  to prove Theorem \ref{thm:heavy for FUSF}.
\begin{proof}[Proof of Theorem \ref{thm:heavy for FUSF}]
	Notice that the FUSF on $G=G_0*\mathbb{Z}_2$ is an example of {free product percolation} of $\mathfrak{F}^f(G_0)$ on $G$, where $\mathfrak{F}^f(G_0)$ is a FUSF sample on $G_0$.
Theorem \ref{thm:heavy for FUSF} follows from the  combination of Proposition \ref{prop: FUSF on grandparent graph is a tree}, \ref{prop: branching number larger than one for FUSF on grandparent graph} and Proposition \ref{prop: free product percolation}. 
\end{proof}

\begin{proof}[Proof of Proposition \ref{prop: free product percolation}]
	The part that each connected component of $\omega$ has branching number bigger than one is obvious and we omit the details. 
	
	The part that each connected component of $\omega$ is heavy can be proved using comparison to branching random walks. 
	
	
	For a fix vertex $u\in V(G_0)$, since each connected component of $\omega_0$ is infinite, writing $K_{G_0}(u)$ of the connected component of $u$ in $\omega_0$, one has 
	\[
	\mathbb{E}\Big[\sum_{v\in V(G_0)}\mathbf{1}_{  \{v\in K_{G_0}(u)\} }\Big]=\infty.
	\] 
	Using TMTP and noting $v\in K_{G_0}(u)$ if and only if $u\in K_{G_0}(v)$, one has 
	\be\label{eq: prop 6.2-1}
	\mathbb{E}\Big[\sum_{v\in V(G_0)} \mathbf{1}_{  \{v\in K_{G_0}(u)\}}\frac{m(v)}{m(u)} \Big]=\infty. 
	\ee

	By monotone convergence theorem there exists a large constant $M$ such that 
	\be\label{eq: prop 6.2-2}
		\mathbb{E}\Big[\sum_{v\in V(G_0)} \mathbf{1}_{  \{v\in K_{B_{G_0}(u,M)}(u)\}}\frac{m(v)}{m(u)} \Big]>e,
	\ee
	where $B_{G_0}(u,M)$ denotes the ball in $G_0 $ with center $u$ and radius $M$ and $K_{B_{G_0}(u,M)}(u)$ denote the connected component of $u$ when one consider the percolation of $\omega_0$ restricted in $B_{G_0}(u,M)$. 
	
	For $x=(u,i)\in V(G)$, there is a unique $\mathbb{Z}_2$-edge incident to $x$ and we suppose the other vertex of the $\mathbb{Z}_2$-edge is $x'=(u',j)$. Let $K(x)$ denote the connected component of $x$ in the free product percolation $\omega$. Let $K_h(x)$ denote the connected component of $x$ if we delete the $\mathbb{Z}_2$-edge $e=(x,x')$ from $K(x)$. 
	
	For a fixed constant $M>0$, we truncate $K_h(x)$ as follows. First we truncate all the edges $e=((a,i),(b,i))$ in the same copy $G_0^{(i)}$ as $x=(u,i)$ if $\max\{d_{G_0}(a,u),d_{G_0}(b,u) \}\geq M$. For all the vertex $y$ in the same copy $G_0^{(i)}$  as $x$ that can be connected to $x$ by an $\omega$-open path staying in $B_{G_0^{(i)}}(x,M)$ (the ball in the copy $G_0^{(i)}$  with center $x$ and radius $M$), we keep the $\mathbb{Z}_2$-edge $(y,y')$ and do the same truncation procedure for $y'$ as previous for $x$. Keep doing this procedure and in the end we get an infinite random graph $K_h^M(x)$. This random graph $K_h^M(x)$ is just the sub-graph of $K_h(x)$ induced by those vertices that have an $\omega$-open path to $x$  such  that between any two consecutive $\mathbb{Z}_2$-edges on this path there is at most $M$ other ($G_0$)-edges .

	Now we show that $K_h^M(x)$ is heavy with positive probability.
	
	We first briefly recall the definition of a branching random walk and a related result; see \cite{Lyons1997} for more details. Let $\mathcal{L}:=\{X_i\}_{i=1}^{N}$ be a random $N$-tuple of real numbers, where $N$ is also random. We can view the branching random walk as an ordered point process on the real line. An initial point $x$ is located at the origin. It gives birth to $N$ children $x_1,\ldots,x_N$ with random displacements $X_1,\ldots,X_N$. Then each $x_i$ gives birth to a random number of particles with random displacement relative to the position of $x_i$ according to the same law as $\mathcal{L}$ and independently of one another and of the initial displacements. This procedure continues forever or until no more particles are born. 
	
	For a particle $u$, let $|u|$ be the generation of $u$ and $X(u)$ its displacement from its parent, and $S(u)$  its position (relative to the origin). Denote the initial particle at origin by $0$. If $u$ is an ancestor of $v$, write $u<v$. Then $S(v)=\sum_{0<u\leq v}{X(u)}$. 
	
	For $\alpha\in\mathbb{R}$, write $\langle \alpha,\mathcal{L}\rangle:=\sum_{i=1}^{N}e^{-\alpha X_i}$ and $\lambda(\alpha):=\mathbb{E}[\langle \alpha,\mathcal{L}\rangle]\in (0,\infty]$. Assume $\lambda(0)>1$ so that the extinction probability $q<1$. If $\lambda(\alpha)<\infty$ for some $\alpha$, then $W_n(\alpha):=\frac{\sum_{|u|=n}{e^{-\alpha S(u)}}}{\lambda(\alpha)^n}$ is a martingale with a.s. limit $W(\alpha)$. 

	Set 
	\[
	\lambda'(\alpha):=\mathbb{E}\Big[\sum_{i=1}^{N}X_ie^{-\alpha X_i}\Big]
	\]
	when the integral exists as a Lebesgue integral. Biggins' theorem tells us that if there exists a real number $\alpha$ such that the following three conditions holds, then the limit $W(\alpha)$ satisfies $\mathbb{P}[W(\alpha)=0]=q<1$.
	
	\begin{itemize}
		\item [$(1)$] $\lambda(\alpha)<\infty$ and $\lambda'(\alpha)$ exists and is finite;
		\item[$(2)$] $\mathbb{E}[\langle \alpha,\mathcal{L}\rangle\log ^+\langle \alpha,\mathcal{L}\rangle]<\infty$ and 
		\item[$(3)$] $\alpha\lambda'(\alpha)/\lambda(\alpha)<\log \lambda(\alpha)$.
	\end{itemize} 

Now given $K_h^M(x)$, we construct a corresponding branching random walk in the following manner. Let $N+1$ be the number of  vertices in the connected component of $x=(u,i)$ in $K_h^M(x)$ intersecting with the  copy  $G_0^{(i)}$, which we denote by $K_h^M(x,B_{G_0^{(i)}})$. Note that $N\in [M, |B_{G_0}(x,M)|]$ is a finite random number. Let $x_1,\ldots,x_N$ be an arbitrary ordering of the vertices in $K_h^M(x,B_{G_0^{(i)}})\backslash\{x \}$. Now for the corresponding branching random walk, we let the initial particle give birth to $N$ children, each with displacement $X_i=\log\frac{m(x_i)}{m(x)}$. For each $x_i$, let $x_i'=(u_i,j_i)$ be the vertex incident to $x_i$ by a $\mathbb{Z}_2$-edge. Let $N_{x_i}$ denote the number of vertices in $K_h^M(x_i',B_{G_0^{(j_i)}})\backslash\{x_i' \}$. In the corresponding branching random walk, we let the particle corresponding to $x_i$ give birth to $N_{x_i}$ new particles each with a relative displacement $\log\frac{m(y)}{m(x_i')},\, \forall \, y\in K_h^M(x_i',B_{G_0^{(j_i)}})\backslash\{x_i' \}$.  By the third item in Observation \ref{obse: action on the free product} and Lemma \ref{lem:haar} one has $m(x_i)=m(x_i')$, thus $\log \frac{m(y)}{m(x_i')}=\log \frac{m(y)}{m(x_i)}$.

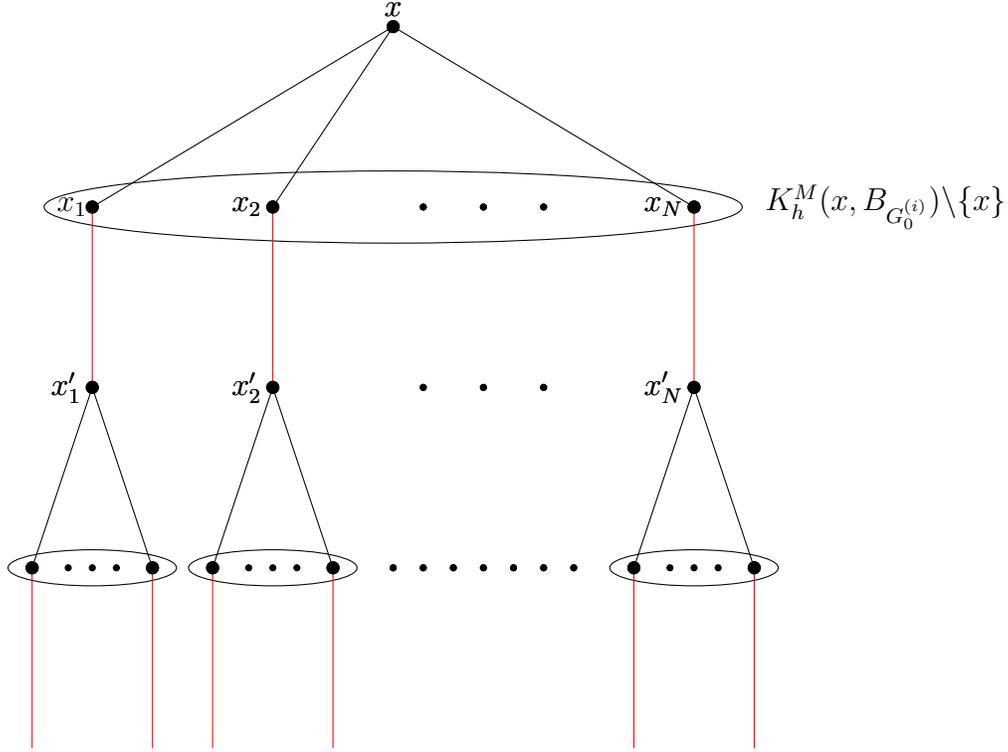
\begin{figure}[h!]
	\centering
	\begin{tikzpicture}[scale=0.8]\label{pic:comparison to branching random walk}

	\draw[fill=black] (8,15) circle [radius=0.10];
	\node [above] at (8,15) {$x$};

	\draw[fill=black] (3,12) circle [radius=0.10];
	\draw[fill=black] (6,12) circle [radius=0.10];
	\node [left] at (6,12) {$x_2$};
	\draw[fill=black] (13,12) circle [radius=0.10];
	\node [left] at (13,12) {$x_N$};
	
	\draw[fill=black] (8.5,12) circle [radius=0.05];
	\draw[fill=black] (9.5,12) circle [radius=0.05];
	\draw[fill=black] (10.5,12) circle [radius=0.05];
	\draw[fill=black] (8.5,9) circle [radius=0.05];
	\draw[fill=black] (9.5,9) circle [radius=0.05];
	\draw[fill=black] (10.5,9) circle [radius=0.05];
	\draw[fill=black] (8.5,6) circle [radius=0.05];
	\draw[fill=black] (9.5,6) circle [radius=0.05];
	\draw[fill=black] (10.5,6) circle [radius=0.05];
	\draw[fill=black] (8,6) circle [radius=0.05];
	\draw[fill=black] (9,6) circle [radius=0.05];
	\draw[fill=black] (10,6) circle [radius=0.05];
	\draw[fill=black] (11,6) circle [radius=0.05];
	
	\node [right] at (14,12) {$K_h^M(x,B_{G_0^{(i)}})\backslash \{ x\}$};
	
	\draw (8,12)  ellipse (5.8 and 0.6);

	\draw(8,15)--(3,12);
	\draw (8,15)--(6,12);
	\draw (8,15)--(13,12); 
	
	\draw[fill=black] (3,9) circle [radius=0.10];
	\node [left] at (3,9) {$x_1'$};
	\draw [color=red](3,12)--(3,9);
	\draw[fill=black] (6,9) circle [radius=0.10];
	\node [left] at (6,9) {$x_2'$};
	\draw [color=red](6,12)--(6,9);
	\draw[fill=black] (13,9) circle [radius=0.10];
	\node [left] at (13,9) {$x_N'$};
	\draw [color=red](13,12)--(13,9);

	\draw[fill=black] (2,6) circle [radius=0.10];
	\draw[fill=black] (4,6) circle [radius=0.10];
	\draw (3,9)--(2,6);
	\draw (3,9)--(4,6);
	\draw[fill=black] (2.6,6) circle [radius=0.05];
	\draw[fill=black] (3,6) circle [radius=0.05];
	\draw[fill=black] (3.4,6) circle [radius=0.05];
	\draw (3,6) ellipse (1.4 and 0.3);
	\draw [color=red](2,6)--(2,3);
	\draw [color=red](4,6)--(4,3);

	\draw[fill=black] (5,6) circle [radius=0.10];
	\draw[fill=black] (7,6) circle [radius=0.10];
	\draw (6,9)--(5,6);
	\draw (6,9)--(7,6);
	\draw[fill=black] (5.6,6) circle [radius=0.05];
	\draw[fill=black] (6,6) circle [radius=0.05];
	\draw[fill=black] (6.4,6) circle [radius=0.05];
	\draw (6,6) ellipse (1.4 and 0.3);
	\draw[color=red](5,6)--(5,3);
	\draw[color=red](7,6)--(7,3);

	\draw[fill=black] (12,6) circle [radius=0.10];
	\draw[fill=black] (14,6) circle [radius=0.10];
	\draw (13,9)--(12,6);
	\draw (13,9)--(14,6);
	\draw[fill=black] (12.6,6) circle [radius=0.05];
	\draw[fill=black] (13,6) circle [radius=0.05];
	\draw[fill=black] (13.4,6) circle [radius=0.05];
	\draw (13,6) ellipse (1.4 and 0.3);
	\draw [color=red](12,6)--(12,3);
	\draw [color=red](14,6)--(14,3);

	\draw[fill=black] (8,15) circle [radius=0.10];
	\node [above] at (8,15) {$x$};

	\draw[fill=black] (3,12) circle [radius=0.10];
	\node [left] at (3.1,12) {$x_1$};
	\draw[fill=black] (6,12) circle [radius=0.10];
	\node [left] at (6,12) {$x_2$};
	\draw[fill=black] (13,12) circle [radius=0.10];
	\node [left] at (13,12) {$x_N$};
	
	\draw[fill=black] (8.5,12) circle [radius=0.05];
	\draw[fill=black] (9.5,12) circle [radius=0.05];
	\draw[fill=black] (10.5,12) circle [radius=0.05];
	\draw[fill=black] (8.5,9) circle [radius=0.05];
	\draw[fill=black] (9.5,9) circle [radius=0.05];
	\draw[fill=black] (10.5,9) circle [radius=0.05];
	\draw[fill=black] (8.5,6) circle [radius=0.05];
	\draw[fill=black] (9.5,6) circle [radius=0.05];
	\draw[fill=black] (10.5,6) circle [radius=0.05];
	\draw[fill=black] (8,6) circle [radius=0.05];
	\draw[fill=black] (9,6) circle [radius=0.05];
	\draw[fill=black] (10,6) circle [radius=0.05];
	\draw[fill=black] (11,6) circle [radius=0.05];

	\draw[fill=black] (3,9) circle [radius=0.10];
	\node [left] at (3,9) {$x_1'$};
	\draw[fill=black] (6,9) circle [radius=0.10];
	\node [left] at (6,9) {$x_2'$};
	\draw[fill=black] (13,9) circle [radius=0.10];
	\node [left] at (13,9) {$x_N'$};

	\draw[fill=black] (2,6) circle [radius=0.10];
	\draw[fill=black] (4,6) circle [radius=0.10];
	\draw[fill=black] (5,6) circle [radius=0.10];
	\draw[fill=black] (7,6) circle [radius=0.10];
	\draw[fill=black] (12,6) circle [radius=0.10];
	\draw[fill=black] (14,6) circle [radius=0.10];
	\end{tikzpicture}

	\caption{A systematic drawing of the family tree of the corresponding branching random walk for $K_h^M(x)$ if one contracts all the red edges}
	\label{fig: branching comparison family tree}
\end{figure}

We use  $y$ to denote a vertex  in $K_h^M(x)$ or its corresponding vertex in the family tree of the corresponding branching random walk.  
Now we take $\alpha=-1$, then $e^{-\alpha X_i}=\frac{m(x_i)}{m(x)}$, $e^{-\alpha S(y)}=\frac{m(y)}{m(x)}$. Notice that as a subgraph of $G_0^{(i)}$, $K_h^M(x,B_{G_0^{(i)}})$ has the same law as $K_{B_{G_0}(u,M)}(u)$ as a subgraph of $G_0$. By Lemma \ref{lem:haar} and the second item in Observation \ref{obse: action on the free product}
\begin{eqnarray}\label{eq: prop 6.2-3}
\lambda(-1)&=&
\mathbb{E}\Big[\sum_{y\in K_h^M(x,B_{G_0^{(i)}})\backslash\{x\} } \frac{m(y)}{m(x)} \Big]\nonumber\\
&=&	\mathbb{E}\Big[\sum_{y\in V(G_0),y\neq x} \mathbf{1}_{  \{y\in K_{B_{G_0}(x,M)}(x)\}}\frac{m(y)}{m(x)} \Big]\stackrel{\eqref{eq: prop 6.2-2}}{>}e. 
\end{eqnarray}

Now we verify that the conditions listed above to use Biggins' theorem hold. 
Conditions $(1)$ and $(2)$ are trivial since in our case $X_i$ and $N$ are both bounded.
Since $\alpha=-1$ and $\lambda(-1)>0$, condition $(3)$ is just $\lambda'(-1)+\lambda(-1)\log\lambda(-1)> 0$. Since $\lambda(-1)>e$, it suffices to show that 
\be\label{eq: prop 6.2-4}
\lambda'(-1)+\lambda(-1)\geq  0.
\ee

By definition and the equivalence of $y\in K_{B_{G_0}(x,M)}(x)$ and $x\in K_{B_{G_0}(y,M)}(y)$ one has
\begin{eqnarray}\label{eq: prop 6.2-5}
\lambda'(-1)&=& \mathbb{E}\Big[\sum_{y\in V(G_0),y\neq x}  \mathbf{1}_{  \{y\in K_{B_{G_0}(x,M)}(x)\}}\frac{m(y)}{m(x)} \log \frac{m(y)}{m(x)} \Big] \nonumber\\
&\stackrel{\textnormal{TMTP}}{=}&
\mathbb{E}\Big[\sum_{y\in V(G_0),y\neq x}  \mathbf{1}_{  \{y\in K_{B_{G_0}(x,M)}(x)\}}\log \frac{m(x)}{m(y)} \Big] 
\end{eqnarray}
Note that $f(t)=t+\log\frac{1}{t}\geq f(1)=1$ on $(0,\infty)$. 
Hence by \eqref{eq: prop 6.2-3} and \eqref{eq: prop 6.2-5} one has \eqref{eq: prop 6.2-4}:
\[
\lambda'(-1)+\lambda(-1)=\mathbb{E}\Big[\sum_{y\in V(G_0),y\neq x}  \mathbf{1}_{  \{y\in K_{B_{G_0}(x,M)}(x)\}} f(\frac{m(y)}{m(x)})\Big] \geq \mathbb{E}[N]\geq M>0. 
\]

Therefore Biggins' theorem yields that with positive probability $W(\alpha)>0$. On the event that  $W(\alpha)>0$, the corresponding $m(K_h^M(x)):=\sum_{y\in K_h^M(x)}{m(y)}=\infty$ since 
$\sum_{|u|=n}{e^{-\alpha S(u)}}=\sum_{|u|=n}\frac{m(u)}{m(x)}$ tends to infinity because $\lambda(-1)>1$ and $W(\alpha)>0$. Using the standard trick (see e.g.\ Proposition 5.6 in \cite{LP2016}) one has  $m(K_h^M(x))=\infty$ almost surely. 

Therefore $K(x)$ is heavy almost surely.
\end{proof}


\begin{remark}
	The free minimal spanning forests on $G=G_0*\mathbb{Z}_2$ is also a free product percolation that satisfies the condition of Proposition \ref{prop: free product percolation}. Hence each tree of the free minimal spanning forests on $G=G_0*\mathbb{Z}_2$ is also heavy and has branching number bigger than one. Interested readers can refer to Chapter 11 of \cite{LP2016} for more background on the free minimal spanning forests.
\end{remark}

We conclude with two further open questions on FUSF on nonunimodular transitive graphs. The first question is about the number of trees in the FUSF. Benjamini et al asked whether the number of trees of the FUSF is $1$ or $\infty$ almost surely (Question 15.6 in \cite{BLPS2001}). Hutchcroft and Nachmias \cite{TomAsaf2017} answered this question positively for unimodular transitive graphs and the nonunimodular case remains open.  Another question one can consider is about the indistinguishability of the trees in the FUSF on a nonunimodular transitive graph. Since the trees in the WUSF on a nonunimodular transitive graph is light almost surely, they are distinguishable by automorphism-invariant properties, e.g.\ the sum of degrees of the highest points in the components. So the  interesting case is for the trees in the FUSF on nonunimodular transitive graphs with the property that $\textnormal{FUSF}\neq \textnormal{WUSF}$. The techniques from \cite{TomAsaf2017, Timar2018indis,Tang2018} might be useful for this question.

\subsection*{Acknowledgments}
We thank Russell Lyons for many helpful discussions and suggestions. We thank \'{A}d\'{a}m Tim\'{a}r and G\'{a}bor Pete for sharing their recent work \cite{PeteTimar2020} and  pointing out the branching number part of Remark \ref{rem: connected FUSF on tree cross finite graph}.  We thank the referees for their helpful comments and suggestions. Especially one referee pointed out a serious mistake in the previous version and proposed a useful rescue, namely the tree-graph inequality method.  This also leads to sharper bounds on certain estimates in Section \ref{sec: 4}. The study of $\mathbb{P}[|\mathfrak{T}_x\cap L_n(x)|>k]$ in Proposition \ref{prop: stretched exponential decay for x-wusf intersect with L_n(x)} was also suggested by the referee. Question \ref{ques: invariant recurrent spanning tree for grandparent graphs} was suggested by another referee. 

\appendix
\section{Proof of \eqref{eq: exponential decay lower bound in k for srw intersect with a slab} in Proposition \ref{prop: intersection of a srw and a slab}}\label{sec: appendix A}

We have already seen the following property of simple random walk from Lemma \ref{lem:expectation for SRW on nonunimodular transitive graph} and \eqref{eq: prob of srw intersecting a high slab}:
	\be\label{eq: prob of intersection with a slab for SRW}
	\mathbb{P}[\{X_m^x:m\geq 0\}\cap L_n(x)\neq \emptyset]
	\left\{
	\begin{array}{cc}
	\asymp e^{-t_0n} & \textnormal{ if }n>0\\
	=1 & \textnormal{ if }n\leq0
	\end{array}
	\right.
	\ee

Using \eqref{eq: prob of intersection with a slab for SRW} and strong Markov property for simple random walk, to show \eqref{eq: exponential decay lower bound in k for srw intersect with a slab} it suffices to show the case $n=0$, i.e., there exists $c>0$ such that 
\be\label{eq: lower bound of exponential decay for srw intersect with L_0(x)}
	\mathbb{P}[|\{X_m^x:m\geq 0\}\cap L_0(x)|\geq k]\geq e^{-ck},\,\,\forall\, k\geq 1
\ee
The basic idea is to show that there is a large constant $c>0$ such that for every $k\geq 1$, one can construct a simple  random walk trajectory $\{X_m^x\colon m\geq0\}$ starting from $x$ such that  $|\{X_m^x\colon m=0,\ldots,ck\}\cap L_0(x)|\geq k$.

Recall $t_0=\max\{\log\Delta(x,y)\colon x\sim y\}>0$. For $k\in \mathbb{Z}$, let $N_{k}(x):=\{y\colon d_G(x,y)=|k|, \textnormal{ and }\log \Delta(x,y)=kt_0 \}$. In particular, $N_0(x)=\{x\}$. For a set $A\subset V(G)$, write $N_k(A):=\bigcup_{x\in A} N_k(x)$. 

Obviously for $k\geq 1$,  $N_k(x)\subset N_1(N_{k-1}(x))$. On the other hand, if $y\in N_1(N_{k-1}(x))$, then there exists $z\in N_{k-1}(x)$ such that $y\in N_1(z)$. Hence 
$d_G(x,z)=k-1,d_G(z,y)=1$ and $\Delta(x,z)=e^{(k-1)t_0}, \Delta(z,y)=e^{t_0}$. Thus 
$d_G(x,y)\leq d_G(x,z)+d_G(z,y)=k$ and $\Delta(x,y)=\Delta(x,z)\cdot \Delta(z,y)=e^{kt_0}$. By the choice of $t_0$, one must have $d_G(x,y)=k$ and thus $y\in N_k(x)$. Therefore one has the reverse direction $N_k(x)\supset N_1(N_{k-1}(x))$. Hence  for all $k\geq1$, $N_k(x)=N_1(N_{k-1}(x))$. 

Using similar idea, one has $N_{k+m}(x)=N_k(N_m(x))$ and $N_{-k-m}(x)=N_{-k}(N_{-m}(x))$ for all $k,m\geq0$. In particular, $|N_{k+m}(x)|\leq |N_k(x)|\cdot |N_{m}(x)|$ and $N_{-k-m}(x)=|N_{-k}(x)|\cdot |N_{-m}(x)|$ for all $k,m\geq0$. Also the fact $N_{k+m}(x)=N_k(N_m(x)),k,m\geq0$ implies that 
\be\label{eq: increasing for N_k}
|N_{k}(x)|\geq |N_{j}(x)|,\,\,\forall\, k\geq j\geq 0.
\ee


Consider a mass-transport function as $F_k(x,y)=\frac{\mathbf{1}_{\{y\in N_{k}(x)\} }}{|N_{k}(x)|}$. Notice that $x\in N_{k}(y)$ if and only if $y\in N_{-k}(x)$. By the tilted mass-transport principle we have 
\[
1=\sum_{y\in V}{F_k(x,y)}=\sum_{y\in V}{F_k(y,x)\Delta(x,y)}=\sum_{y\in N_{-k}(x)}{\frac{1}{|N_{k}(x)|}e^{-kt_0}}.
\]
Thus 
\be\label{eq: N_(-k) equal exp(kt_0) N_(k)}
|N_{-k}(x)|=e^{kt_0}|N_{k}(x)|,\,\,\forall\, k\geq 1. 
\ee

By \eqref{eq: increasing for N_k} and \eqref{eq: N_(-k) equal exp(kt_0) N_(k)} one has that there exists a positive constant $c_1$ such that for all $k\geq1$,
\be\label{eq: bounded total population by k-th level}
\sum_{j=0}^{k}|N_{-j}(x)|=|N_{-k}(x)|\sum_{j=0}^{k}e^{(j-k)t_0}|\frac{|N_j(x)|}{|N_k(x)|}
\leq |N_{-k}(x)|\sum_{j=0}^{k}e^{(j-k)t_0}\leq c_1|N_{-k}(x)|.
\ee

Consider the following subgraph $T(x)$ of $G$. Let $V(T(x))$ be the union $\bigcup_{k\geq0}N_{-k}(x)$. For each $k\geq 0$ and $y\in N_{-k-1}(x)$, there is at least one edge linking from $N_{-k}(x)$ to  $y$. If there are more one such edges, just pick one arbitrarily. Let $E(T(x))$ be the collection of edges picked. Then obviously $T(x)$ is a tree. Write $T_{\leq n}(x):=\{ y\in V(T)\colon d_{T(x)}(y,x)\leq n \}$. 

By \eqref{eq: N_(-k) equal exp(kt_0) N_(k)}, $|N_{-j}(x)|$ grows exponentially fast in $j$. Now for a fixed integer $k\geq 1$, let $j$ be the smallest nonnegative integer such that $|N_{-j}(x)|> k$. 
Then \[
|N_{-j+1}(x)|\leq k< |N_{-j}(x)|\leq |N_{-1}(x)|\cdot |N_{-j+1}(x)|\leq |N_{-1}(x)|\cdot k.
\]

Pick $v\in N_j(x)$. Let the simple random walk $\{X_m^x\colon m\geq0 \}$ travels according the depth-first search of the tree $T_{\leq j}(v)$ until it visits all vertices in $N_{-j}(v)$. Under this event one has that $|\{X_m^x\colon m\geq0 \}\cap L_0(x)|\geq |N_{-j}(v)|>k$.  Notice that there are $\sum_{i=1}^{j}|N_{-i}(v)|\stackrel{\eqref{eq: bounded total population by k-th level}}{\leq }c_1|N_{-j}(v)|\leq c_1|N_{-1}(x)|k$ edges in  $T_{\leq j}(v)$. Also notice each edge in the tree $T_{\leq j}(v)$ is been crossed by the depth-first search at most twice.
Denoting by $D$  the degree of $G$, one obtains \eqref{eq: lower bound of exponential decay for srw intersect with L_0(x)}:
\[
\mathbb{P}[|\{X_m^x\colon m\geq0 \}\cap L_0(x)|>k]\geq \frac{1}{D^{2c_1|N_{-1}(x)|k}}\geq \exp(-c_2k).
\]

\section{Proof of Proposition \ref{prop: high moments for x-component intersecting a slab} for general $k$}\label{sec: appendix B}
The proof of Proposition \ref{prop: high moments for x-component intersecting a slab} uses similar ideas and notations from the proof of \cite[Theorem 6.75]{Grimmett1999}, where  high moments of a cluster were estimated in the Bernoulli percolation case. 

Here we recall some definitions from \cite[Page 135]{Grimmett1999}. 
\begin{definition}
	A tree is called a \notion{skeleton} if each vertex has degree $1$ or $3$. Suppose $S$ is a skeleton, let $I(S)$ denote the set of vertices with degree $3$ and call these vertices \notion{interior vertices}. Call the other vertices of $S$ \notion{exterior vertices} (they are just leaves of the tree $S$). It is easy to see that if a skeleton with $k$ exterior vertices, then it must have $k-2$ interior vertices. Write $V(S), E(S)$ as the vertex set and edge set of $S$ respectively as usual (in particular $E(S)$ is not the set of exterior vertices).

	A skeleton with $k$ exterior vertices is called \notion{labelled} if there is an assignment of the numbers $0,1,\ldots,k-1$ to the exterior vertices. Two labelled skeleton $S_1$ and $S_2$ are called \notion{isomorphic} if there is a graph isomorphism from $S_1$ to $S_2$ which also preserves the labels of the exterior vertices. 
\end{definition}

\begin{figure}[h!]
	\centering
	\begin{tikzpicture}[scale=0.8, text height=1.5ex,text depth=.25ex] 
	\draw[fill=black] (4,5) circle [radius=0.05]; 
	\node[above] at (4,5) {$0$};
	\draw[fill=black] (4,4) circle [radius=0.05];
	\draw[fill=black] (2,3) circle [radius=0.05];
	\draw[fill=black] (6,3) circle [radius=0.05];
	\draw[fill=black] (1,2) circle [radius=0.05];
	\draw[fill=black] (3,2) circle [radius=0.05];
	\node[left] at (3,2) {$5$};
	\draw[fill=black] (5,2) circle [radius=0.05];
	\draw[fill=black] (7,2) circle [radius=0.05];
	\node[left] at (7,2) {$2$};
	\draw[fill=black] (0,1) circle [radius=0.05];
	\node[left] at (0,1) {$1$};
	\draw[fill=black] (2,1) circle [radius=0.05];
	\node[left] at (2,1) {$4$};
	\draw[fill=black] (4,1) circle [radius=0.05];
	\node[left] at (4,1) {$6$};
	\draw[fill=black] (6,1) circle [radius=0.05];
	\node[left] at (6,1) {$3$};

	\draw (4,5)--(4,4); 
	\draw (4,4)--(2,3); 
	\draw (4,4)--(6,3); 
	\draw (2,3)--(1,2); 
	\draw (2,3)--(3,2); 
	\draw (6,3)--(5,2); 
	\draw (6,3)--(7,2); 
	\draw (1,2)--(0,1); 
	\draw (1,2)--(2,1); 
	\draw (5,2)--(4,1); 
	\draw (5,2)--(6,1); 
	
	\end{tikzpicture}
	\caption{A labelled skeleton with $6+1$ exterior vertices and $6-1$ interior vertices.}
\end{figure}
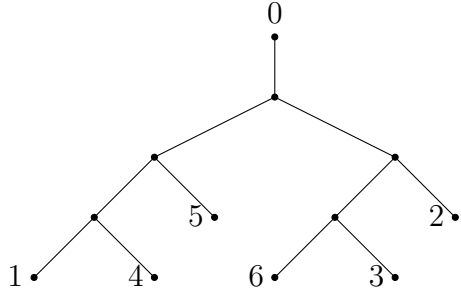 

The following lemma is an adaption of the claim from \cite[Page 135-136, Section 6.3]{Grimmett1999}.
\begin{lemma}\label{lem: skeleton map}
	Suppose $x_0,x_1,\ldots,x_k$ are (possibly non-distinct) vertices of $G$ such that $x_i\in\mathfrak{T}_{x_0}$ for $i=1,\ldots,k$. Write $\mathbf{x}=(x_0,x_1,\ldots,x_k)$. Then there is a labelled skeleton $S$ with $k+1$ exterior vertices together with  a map $\varphi_{\mathbf{x}}:V(S)\rightarrow V(G)$ such that 
	\begin{itemize}
		\item the exterior vertex with label $i$ maps to $x_i$ by $\varphi_{\mathbf{x}}$ for $i=0,1,\ldots,k$ and 
		\item the edges of $S$ correspond to edge-disjoint paths in $\mathfrak{T}_{x_0}$.
	\end{itemize}  
	Here if $(u,v)$ is an edge of $S$ such that $\varphi_{\mathbf{x}}(u)=\varphi_{\mathbf{x}}(v)$, then its corresponding path is an empty path and we use the convention that the empty path is edge disjoint with any other path. Such a map is called an \notion{admissible map}. 
\end{lemma}

\begin{lemma}\label{lem: point function inequality for general k}
	Suppose $x_0,x_1,\ldots,x_k$ are (possibly non-distinct) vertices of $G$. Then 
	\be\label{eq:  point function inequality for general k}
	\tau(x_0,x_1,\ldots,x_k)\leq \sum_S\sum_{\varphi_{\mathbf{x}}}\prod_{(u,v)\in E(S)}\tau(\varphi_{\mathbf{x}}(u),\varphi_{\mathbf{x}}(v))
	\ee 
	where the first sum is over all labelled skeleton with $k+1$ exterior vertices and the second sum is over all admissible maps $\varphi_{\mathbf{x}}$ from $V(S)$ to $V(G)$.
\end{lemma}
\begin{proof}
	This is the generalization of Lemma \ref{lem: point function ineq with k=3}.
	By Lemma \ref{lem: skeleton map}, one has 
	\begin{align*}
	\tau(x_0,x_1,\ldots,x_k)\leq \sum_S\sum_{\varphi_{\mathbf{x}}}\mathbb{P}\big[
	\varphi_{\mathbf{x}}(u)\in \mathfrak{T}_{x_0} \textnormal{ for all }u\in V(S)\textnormal{ and the paths in } \mathfrak{T}_{x_0} \textnormal{ connecting } \\
	\varphi_{\mathbf{x}}(u),\varphi_{\mathbf{x}}(v)  \textnormal{ are edge-disjoint for all }(u,v)\in E(S)
	\big]
	\end{align*}
	Place the skeleton $S$ in the lower half plane in some arbitrary way. Order the vertices of $S$ first according to the graph distance to the exterior vertex with label $0$ (for simply we will denote this vertex $0$ hereafter) and if two vertices are at the same distance to $0$, then order them according left to right. Call this order $Or$.
	Sample  $\mathfrak{T}_{x_0}$ using independent simple random walks started from the order $\varphi_{\mathbf{x}}(Or)$. Then one can see that \eqref{eq:  point function inequality for general k} holds. 
\end{proof}

\begin{figure}[h!]
	\centering
	\begin{tikzpicture}[scale=0.8, text height=1.5ex,text depth=.25ex] 
	\draw[fill=black] (4,5) circle [radius=0.05]; 
	\node[above] at (4,5) {$x_0$};
	\draw[fill=black] (4,4) circle [radius=0.05];
	\node[left] at (4,4) {$u_1$};
	\draw[fill=black] (2,3) circle [radius=0.05];
	\node[left] at (2,3) {$u_2$};
	\draw[fill=black] (6,3) circle [radius=0.05];
	\node[left] at (6,3) {$u_3$};
	\draw[fill=black] (1,2) circle [radius=0.05];
	\node[left] at (1,2) {$u_4$};
	\draw[fill=black] (3,2) circle [radius=0.05];
	\node[left] at (3,2) {$x_5$};
	\draw[fill=black] (5,2) circle [radius=0.05];
	\node[left] at (5,2) {$u_5$};
	\draw[fill=black] (7,2) circle [radius=0.05];
	\node[left] at (7,2) {$x_2$};
	\draw[fill=black] (0,1) circle [radius=0.05];
	\node[left] at (0,1) {$x_1$};
	\draw[fill=black] (2,1) circle [radius=0.05];
	\node[left] at (2,1) {$x_4$};
	\draw[fill=black] (4,1) circle [radius=0.05];
	\node[left] at (4,1) {$x_6$};
	\draw[fill=black] (6,1) circle [radius=0.05];
	\node[left] at (6,1) {$x_3$};

	\draw [->-](4,5)--(4,4); 
	\draw [->-](4,4)--(2,3); 
	\draw [->-](4,4)--(6,3); 
	\draw [->-](2,3)--(1,2); 
	\draw [->-,blue](2,3)--(3,2); 
	\draw [->-](6,3)--(5,2); 
	\draw [->-,blue](6,3)--(7,2); 
	\draw [->-,blue](1,2)--(0,1); 
	\draw [->-,blue](1,2)--(2,1); 
	\draw [->-,blue](5,2)--(4,1); 
	\draw [->-,blue](5,2)--(6,1); 
	
	\end{tikzpicture}
	\caption{A labelled skeleton with edges oriented pointing out the root $0$ and vertices named under the image of $\varphi_{\mathbf{x}}$.}
	\label{fig: oriented skeleton}
\end{figure}
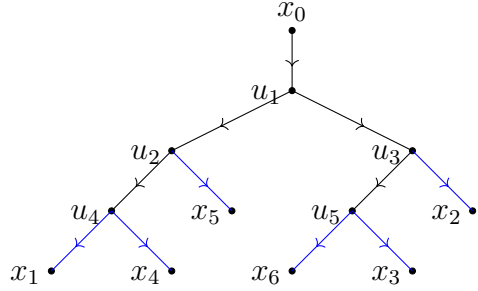

\begin{definition}\label{def: oriented skeleton sum}
	Suppose $x_0,x_1,\ldots,x_k$ are (possibly non-distinct) vertices of $G$.
	Given a labelled skeleton $S$ and an admissible map  $\varphi_{\mathbf{x}}$ from $V(S)$ to $V(G)$. Name the  vertices of $S$ as the image under $\varphi_{\mathbf{x}}$. In particular, write $I(S)=\{u_1,\ldots,u_{k-1}\}$. Let $l_v$ be the integer such that $\varphi_{\mathbf{x}}(v)\in L_{l_v}(x_0)$. Write $\vec{l}=(l_v\colon v\in V(S))$. Oriented the edges of $S$ such that they all lead out the root $x_0$ and write the set of oriented edges as $E(\vec{S})$; see Figure \ref{fig: oriented skeleton} for an example. Definite $f(S,\vec{l})$ to be
	\[
	f(S,\vec{l}):=\sum_{\langle u,v\rangle\in E(\vec{S})}{(l_v-l_u)\vee0}.
	\]
\end{definition}

\begin{lemma}\label{lem: A.5 in appendix}
	Use same notations as in Definition \ref{def: oriented skeleton sum}. Fix an arbitrary integer $n$.  Let
	$m(S,n,j)$ be the number of $k-1$ tuples $(l_{1},\ldots,l_{k-1})$ such that the corresponding $\vec{l}$ with $l_{u_i}=l_i,i=1,\ldots,k-1$, $l_{x_0}=0$ and $l_{x_i}=n,i=1,\ldots,k$ satisfying 
	\[
	f(S,\vec{l})=j,
	\] where $j\geq n\vee 0$.
	For $n\geq0$, one has $m(S,n,n)=1$.
	Also there exists a constant $c>0$ such that 
	\be
	m(S,n,j)\leq c^k(j-n)^{k-1}\textnormal{ for all } n\geq0,\,  j>n
	\ee
	and 
	\be
	m(S,n,j)\leq c^k|n|^{k-1}(j\vee 1)^{k-1}\textnormal{ for all } n<0,\, j\geq0
	\ee
\end{lemma}
\begin{proof}
	For each oriented path $\gamma=\langle v_1,\ldots,v_t\rangle$ in $S$, one has that \[
	\sum_{\langle u,v\rangle \in E(\gamma)}{(l_v-l_u)\vee0}\geq 
	\sum_{\langle u,v\rangle \in E(\gamma)}{(l_v-l_u)}=l_{v_t}-l_{v_1}.
	\]
	First consider the case $n\geq0,j\geq n$. For each inner vertex $u_i$,  considering the oriented path $\eta_{x_0,u_i}$ from $x_0$ to $u_i$ in $S$, one has that 
	\[
	l_i=l_{u_i}-l_{x_0}\leq \sum_{\langle u,v\rangle \in E(\eta(x_0,u_i))}{(l_v-l_u)\vee0}
	\leq f(S,\vec{l})=j.
	\]
	On the other hand, there are at least two disjoint oriented paths starting from $u_i$ to the leaves in $S$, say $\eta_1,\eta_2$. These two paths are also disjoint with the oriented path from $x_0$ to $u_i$.
	Hence 
	\[
	2(n-l_i)+l_i\leq \sum_{\langle u,v\rangle \in E(\eta_1)\cup E(\eta_2)}{(l_v-l_u)\vee0}+\sum_{\langle u,v\rangle \in E(\eta(x_0,u_i))}{(l_v-l_u)\vee0}\leq f(S,\vec{l})=j
	\]
	Hence the above two inequalities imply that $l_i\in[2n-j,j],i=1,\ldots,k-1$. In particular, for $j=n\geq0$, the only tuple $(l_{1},\ldots,l_{k-1})$ such that $f(S,\vec{l})=j$ is $(n,\ldots,n)$ and thus $m(S,n,n)=1$. For $j>n\geq0$, the number of tuples such that $f(S,\vec{l})=j$ is at most $(2j-2n+1)^{k-1}\leq c^k(j-n)^{k-1}$. 
	
	Next for the case $n<0,j\geq0$, one just need to observe that there exists a large constant $c>0$ such that $2j-2n+1\leq c|n|(j\vee 1)$. 
\end{proof}

\begin{proof}[Proof of Proposition \ref{prop: high moments for x-component intersecting a slab}]
	Since for $k\geq2$, $|\mathfrak{T}_x\cap L_n(x)|^k=\sum_{x_1,\ldots,x_k\in L_n(x)}{\tau(x,x_1,\ldots,x_k)}$, using Lemma \ref{lem: point function inequality for general k} one has 
	\begin{eqnarray}\label{eq: A.4 in appendix}
	\mathbb{E}[|\mathfrak{T}_x\cap L_n(x)|^k]&\leq&\mathbb{E}\bigg[
	\sum_{x_1,\ldots,x_k\in L_n(x)}\sum_{S}\sum_{\varphi_{\mathbf{x}}} \prod_{(u,v)\in E(S)}\tau(\varphi_{\mathbf{x}}(u),\varphi_{\mathbf{x}}(v))\bigg]\nonumber\\
	&\leq& \mathbb{E}\bigg[\sum_{S}\sum_{l_1,\ldots,l_{k-1}}
	\sum_{x_1,\ldots,x_k\in L_n(x)}\sum_{\varphi_{\mathbf{x}}:\varphi_{\mathbf{x}}(u_i)\in L_{l_i}(x)} \prod_{(u,v)\in E(S)}\tau(\varphi_{\mathbf{x}}(u),\varphi_{\mathbf{x}}(v))\bigg]\nonumber\\
	&\preceq& c^k\sum_{S}\sum_{l_1,\ldots,l_{k-1}}\exp(-t_0f(S,\vec{l})),
	\end{eqnarray}
	where in the last step we sum over  $x_i\in L_n(x),u_i\in L_{l_i}(x)$ in a reverse order to the one given in the proof of Lemma \ref{lem: point function inequality for general k} and use the equation \eqref{eq: relative expectation when x_1 runs over L_n(x)} in  Observation \ref{obse: relative level relation} $(2k-1)$ times (there are $2k-1$ edges in $E(S)$).
	
From Formula (6.96) (\cite[Page 138, Section 6.3]{Grimmett1999}) we know that the number of labelled skeleton with $k+1$ exterior vertices is $\frac{(2k-2)!}{2^{k-1}(k-1)!}\preceq 2^kk!$. 
Hence by \eqref{eq: A.4 in appendix}  one has that 
\begin{equation}\label{eq: A.5 in appendix}
\mathbb{E}[|\mathfrak{T}_x\cap L_n(x)|^k]\leq
(2c)^kk!\sum_{j\geq n\vee 0}{m(S,n,j)\exp(-t_0j)}
\end{equation}

For $n\geq0$, by Lemma \ref{lem: A.5 in appendix}  there exists a constant $c_1>0$ such that 
\begin{eqnarray*}
\mathbb{E}[|\mathfrak{T}_x\cap L_n(x)|^k]&\leq& c_1^kk! \sum_{j=n}^{\infty}((j-n)\vee 1)^{k-1}\exp(-t_0j)\nonumber\\
&\preceq &c_1^kk!  e^{-t_0n}\int_{0}^{\infty}x^{k-1}e^{-t_0x}dx\nonumber\\
&=&\big(\frac{c_1}{t_0}\big)^k(k!)^2e^{-t_0n}.
\end{eqnarray*}

For $n<0$, by Lemma \ref{lem: A.5 in appendix}  there exists a constant $c_2>0$ such that 
\begin{eqnarray*}
	\mathbb{E}[|\mathfrak{T}_x\cap L_n(x)|^k]&\leq& c_2^kk!|n|^{k-1} \sum_{j=0}^{\infty}(j\vee1)^{k-1}\exp(-t_0j)\nonumber\\
	&\preceq &c_2^kk!|n|^{k-1}\int_{0}^{\infty}x^{k-1}e^{-t_0x}dx\nonumber\\
	&=&\big(\frac{c_2}{t_0}\big)^k(k!)^2|n|^{k-1}.\qedhere
\end{eqnarray*}
\end{proof}

\section{More quantitative results for the toy model}\label{sec: section C in the appendix}

 \subsection{The intersection  $\mathfrak{T}_x\cap L_n(x)$ when $n\geq0$}

By Lemma \ref{lem: law of T_x in  x-wusf for the toy model} we know the tree $\mathfrak{T}_x$ in the $x$-WUSF has the same law of $C_x$, where $C_x$ denotes  the cluster of $x$ in an independent critical Bernoulli bond percolation on the regular tree $\mathbb{T}_{b+1}$. Many things are known about the cluster of $C_x$ at criticality; for example see \cite[Theorem 1.6]{Tom2017}. 

Here we are particularly interested in the tail probability of $|C_x\cap L_n(x)|$ when $n\geq0$. 

\begin{proposition}\label{prop: stretched exponential decay for x-wusf intersect with L_n(x)}
	Consider $x$-WUSF on the toy model $(\mathbb{T}_{b+1},\Gamma_{\xi})$. Then for $n\geq0,k\geq1$,
	\be\label{eq: stretched exponential decay for C_x intersecting with L_n(x)}
	\mathbb{P}[|\mathfrak{T}_x\cap L_n(x)|\geq k]=e^{-t_0n-\Theta(\sqrt{k})}
	\ee
\end{proposition}

\begin{question}
	Can one generalize Proposition \ref{prop: stretched exponential decay for x-wusf intersect with L_n(x)} to WUSF on all nonunimodular transitive graphs? 
\end{question}
In light of Lemma \ref{lem: law of T_x in  x-wusf for the toy model}, the following question is also natural to ask.
\begin{question}
	Consider critical Bernoulli percolation on a nonunimodular transitive graph. Do we have the following estimate: 
	\[
		\mathbb{P}[|C_x\cap L_n(x)|\geq k]=e^{-t_0n-\Theta(\sqrt{k})},n\geq0,k\geq1?
	\]
\end{question}

We start with some lemmas on critical  Galton--Watson trees.

\begin{lemma}\label{lem: upper bound on the k-th moment of Z_n for a critical GW tree}
	Consider a critical Galton--Watson tree with binomial progeny distribution   $\textnormal{Bin}(b,1/b)$, where $b\geq2$. Let $Z_n$ be the size of the $n$-th generation of this critical Galton--Watson tree and $Z_0=1$.  Then there exists a constant $C>1$ such that for all $k\geq 1,n\geq 1$,
	\[
	\mathbb{E}[Z_n^k]\leq C^k n^{k-1} k!
	\]
\end{lemma}
The intuition behind this lemma is to think about the limit of $Z_n$. By classical results of branching process (for example Theorem 1 and 2 in \cite[page 19-20]{AthreyaNey1972}),  
\[
\lim_{n\rightarrow\infty}\frac{\mathbb{P}[Z_n>0]}{\frac{2}{(1-\frac{1}{b})n}}=1
\]
and 
\[
\lim_{n\rightarrow\infty}\mathbb{P}\bigg[\frac{Z_n}{n}>x\mid Z_n>0\bigg]=\exp\big(-\frac{2x}{(1-\frac{1}{b})}\big),\,x\geq 0.
\]
If a nonnegative random variable $Z$ has the following distribution
\begin{itemize}
	\item $\mathbb{P}[Z>0]=1-\mathbb{P}[Z=0]=\frac{2}{(1-\frac{1}{b})n}$ and
	\item $\mathbb{P}\bigg[\frac{Z}{n}>x\mid Z>0\bigg]=\exp\big(-\frac{2x}{(1-\frac{1}{b})}\big),\,x\geq 0.$.
\end{itemize}
then it's easy to compute that $\mathbb{E}[Z^k]= \frac{n^{k-1}\cdot k!}{\lambda^{k-1}},\,\forall\, k\geq 1$, where $\lambda=\frac{2}{1-\frac{1}{b}}$.

\begin{proof}
	For a random variable $W$ with binomial distribution $\textnormal{Bin}(m,p)$, its moment generating function is 
	\[
	M_W(t):=\mathbb{E}[e^{tW}]=[1-p+pe^{t}]^m.
	\] 

	Notice that conditioned on $Z_{m-1}$, $Z_m$ has binomial distribution  $\textnormal{Bin}(bZ_{m-1},\frac{1}{b})$. Hence for $m\geq2$,
	\begin{eqnarray}
	\mathbb{E}[e^{tZ_m}]&=&\mathbb{E}\big[\mathbb{E}[e^{tZ_m}\mid Z_{m-1}]\big]\nonumber\\
	&=&\mathbb{E}\big[\big(1-\frac{1}{b}+\frac{1}{b}e^t\big)^{bZ_{m-1}}\big]\nonumber\\
	&=&\mathbb{E}\big[e^{\big( b\log(1-\frac{1}{b}+\frac{1}{b}e^t)\big)\cdot Z_{m-1} }\big]
	\end{eqnarray}
	
	Using Taylor expansion at $0$ one has
	\[
	b\log(1-\frac{1}{b}+\frac{1}{b}e^t)=t+\big(\frac{1}{2}-\frac{1}{2b}\big)t^2+o(t^2).
	\]	
	Hence there exists a large constant $C_0>0$ such that for all $t\in\big[0,\frac{1}{C_0}\big]$, 
	\[
	t\leq b\log(1-\frac{1}{b}+\frac{1}{b}e^t)\leq t+\frac{1}{2}t^2.
	\]	
	Hence for $t\in \big[0,\frac{1}{C_0}\big]$,  one  has 
	\be\label{eq: generation induction step}
	\mathbb{E}[e^{tZ_m}]=\mathbb{E}\big[e^{\big( b\log(1-\frac{1}{b}+\frac{1}{b}e^t)\big)\cdot Z_{m-1} }\big]
	\leq \mathbb{E}\big[e^{t\big(1+\frac{t}{2}\big)Z_{m-1}}\big]
	\ee
	If $t$ is such that $t\cdot e<\frac{2}{n}\wedge \frac{1}{C_0}$, then one has $$t\big(1+\frac{t}{2}\big)<t(1+\frac{1}{n}).$$
	Take $C\geq C_0$ such that  $t=\frac{1}{nC}$ satisfying $t\cdot e<\frac{2}{n}\wedge \frac{1}{C_0}$ for all $n\geq1$, then by \eqref{eq: generation induction step} one has  
	\[
	\mathbb{E}[e^{tZ_n}]\leq \mathbb{E}\big[e^{t\big(1+\frac{1}{n}\big)Z_{n-1}}\big]
	\]
	By our choice of $C$ and $t$, one has $t\big(1+\frac{1}{n}\big)<t\cdot e<\frac{2}{n}\wedge \frac{1}{C_0}$. Hence by \eqref{eq: generation induction step} again one has 
	\[
	\mathbb{E}[e^{tZ_n}]\leq \mathbb{E}\big[e^{t\big(1+\frac{1}{n}\big)Z_{n-1}}\big]\leq 
	\mathbb{E}\big[e^{t\big(1+\frac{1}{n}\big)^2Z_{n-2}}\big]
	\]
	Repeating this one obtains for $t=\frac{1}{nC}$
	\[
	\mathbb{E}[e^{tZ_n}]\leq\mathbb{E}\big[e^{t\big(1+\frac{1}{n}\big)^{n-1}Z_{1}}\big]\leq \mathbb{E}\big[e^{teZ_1}\big]\leq \mathbb{E}\big[e^{\frac{2Z_1}{n}}\big]
	\]
	
	Since $Z_1$ has binomial distribution $\textnormal{Bin}(b,\frac{1}{b})$,
	one has for $t=\frac{1}{nC}$,
	\be\label{eq: moment generating function of Z_n  with t small}
	\mathbb{E}[e^{tZ_n}]\leq \big[1-\frac{1}{b}+\frac{1}{b}e^{\frac{2}{n}}\big]^b\leq 1+\frac{C_1}{n}
	\ee
	
	Hence for $k\geq1$,
	\[
	\frac{t^k\mathbb{E}[Z_n^k]}{k!}\leq 	\mathbb{E}[e^{tZ_n}]-1\leq \frac{C_1}{n}.
	\]
	Rewriting this one obtains 
	\[
	\mathbb{E}[Z_n^k]\leq C_1\cdot C^kn^{k-1} \cdot k!.\qedhere
	\]
\end{proof}
\begin{remark}
	Notice that $Z_n\leq b^n$ almost surely. So for fixed $n$, $\mathbb{E}[Z_n^k]$ grows at most exponentially in $k$. 
	Hence there is no constant $C>0$ such that for all $n\geq1,k\geq1$
	\[
	\mathbb{E}[Z_n^k]\geq \frac{n^{k-1} \cdot k!}{C^k}.
	\]
\end{remark}

\begin{lemma}\label{lem: bounds on the n-th generation size of a critical GW tree bigger than square of n}
	Consider a a critical Galton--Watson tree with binomial progeny distribution   $\textnormal{Bin}(b,1/b)$. Let $Z_n$ be the size of the $n$-th generation of this critical Galton--Watson tree.  Then there exist constants $c_1,c_2,c_3,c_4>0$ such that for all $n\geq 5$,
	\[
	c_1\exp(-c_2n)\leq \mathbb{P}[Z_n\geq n^2]\leq c_3\exp(-c_4n).
	\]
	Here the requirement of $n\geq 5$ is just to ensure $\mathbb{P}[Z_n\geq n^2]>0$.
\end{lemma}
\begin{proof}
	 By \eqref{eq: moment generating function of Z_n  with t small} we have seen that for a large enough constant $C$,
	\[
	\mathbb{E}[e^{\frac{Z_n}{nC}}]\leq 1+\frac{C_1}{n}.
	\]
	Using Markov's inequality one has that 
	\[
	\mathbb{P}[Z_n\geq n^2]=\mathbb{P}\big[e^{\frac{Z_n}{nC}}\geq e^{\frac{n}{C}}\big]
	\leq \frac{1+\frac{C_1}{n}}{e^{\frac{n}{C}}}\leq (1+C_1)e^{-\frac{n}{C}}.
	\]

	It is easy to compute that 
		\[
		\mathbb{E}[Z_n^2]=1+\frac{(b-1)n}{b}.
		\]
		and \[
		\textnormal{Var}(Z_n)=\frac{(b-1)n}{b}.
		\]
	\begin{claim}\label{claim: 2}
	There exist constants $m_0\geq 5,c_5,c_6>0$ such that for all $m\ge m_0$,
	\be\label{eq: induction step for X_n geq square of n with even n}
	\mathbb{P}\big[Z_{2m}\geq 4m^2 \mid Z_{m}\geq m^2\big]\geq c_5e^{-c_6m}.
	\ee
	and
	\be\label{eq: induction step for X_n geq square of n with odd n}
	\mathbb{P}\big[Z_{2m+1}\geq (2m+1)^2 \mid Z_{m}\geq m^2\big]\geq c_5e^{-c_6m}.
	\ee
\end{claim} 
	Now for $n=2^k$, multiplying \eqref{eq: induction step for X_n geq square of n with even n} for $m=2^3,\ldots,2^{k-1}$ one obtains that there exists $c_7,c_8>0$ such that 
	\[
	\mathbb{P}[Z_n\geq n^2]\geq \mathbb{P}[Z_8\geq 8^2]\cdot c_5^k\exp(-c_6(2^3+\cdots+2^k))\geq c_7\exp(-c_8n).
	\]
	For general $n\geq5$, multiplying \eqref{eq: induction step for X_n geq square of n with even n} or \eqref{eq: induction step for X_n geq square of n with odd n} for $m=f(n):=\lfloor \frac{n}{2}\rfloor, f^{\circ 2}(n),\cdots$ one obtains
	\[
	\mathbb{P}[Z_n\geq n^2]\geq c\cdot c_5^{\log n}\exp(-c_6(\lfloor \frac{n}{2}\rfloor+\cdots +))
	\geq c_7\exp(-c_8n)
	\]
	So it suffices to show Claim \ref{claim: 2}. We will just show \eqref{eq: induction step for X_n geq square of n with even n}, the case \eqref{eq: induction step for X_n geq square of n with odd n} being similar.
	
	On the event $\{Z_{m}\geq m^2\}$, let  $Y_1,Y_2,\ldots,Y_{m^2}$ denote the number of children at generation $2m$ respectively of the first $m^2$ children of $Z_0$ at generation $m$. 
	The random variables $Y_i$'s are independent and satisfy $\mathbb{P}[Y_i>0]=\frac{2}{(1-1/b)m}(1+o(1))$
	and $\textnormal{Var}(Y_i)=\frac{(b-1)m}{b}\asymp m$. Thus we have 
	\be\label{eq: Z_2m at least 4 times square of m bounded below in terms of Y_i}
	\mathbb{P}\big[Z_{2m}\geq 4m^2 \mid Z_{m}\geq m^2\big]\geq\mathbb{P}[Y_1+\cdots+Y_{m^2}\geq 4m^2]
	\ee
	%

	Write $J:=\{j\in[1,m^2]\colon Y_j>0\}$. Then there exists a constant $c_9>0$
	such that $\mathbb{P}[|J|>c_9m]\geq c_9$. In fact $\mathbb{E}[|J|]=\sum_{j=1}^{m^2}\mathbb{P}[Y_j>0]\asymp m$. By independence of $Y_j$'s, $\mathbb{E}[|J|^2]=\mathbb{E}[|J|]^2+\textnormal{Var}(|J|)\asymp m^2+m\cdot\frac{\lambda}{m}(1-\frac{\lambda}{m})\asymp m^2$, where $\lambda=\frac{2}{1-1/b}$. By Paley-Zygmund inequality one has $\mathbb{P}[|J|>c_9m]\geq c_9$ for some small $c_9>0$.
	
Suppose  $W_j=W_j(m)$'s are i.i.d.\ random variables with  the law of $\frac{Y_j}{m}$ conditioned on $ Y_j>0$. Then by Theorem 1 and 2 in \cite[page 19-20]{AthreyaNey1972} $W_j$ converges in distribution to an exponential distribution $\exp(\lambda)$ with $\lambda=\frac{2}{1-1/b}$ as $m\rightarrow\infty$. Hence there exists $m_0\geq 5$ and $c_{10}>0$ such that $\mathbb{P}\big(W_j\geq \frac{4}{c_9}\big)\geq e^{-c_{10}}$ for all $m\geq m_0$. 


	Therefore for $m\geq m_0$,
	\begin{eqnarray}\label{eq: sum of Y_i from 1 to squre of m is at least 4 times square of m}
		\mathbb{P}[Y_1+\cdots+Y_{m^2}\geq 4m^2]&=& \mathbb{P}\bigg[\sum_{j\in J}{Y_j}\geq 4m^2\bigg]\nonumber\\
		&\geq&\mathbb{P}\bigg[|J|>c_9m,\sum_{j\in J}{\frac{Y_j}{m}}\geq 4m\bigg]\nonumber\\
		&\geq& c_9\cdot \mathbb{P}\bigg[\sum_{j=1}^{c_9m}W_j\geq 4m\bigg]\nonumber\\
		&\geq &c_9\cdot \big[\mathbb{P}\big(W_j\geq \frac{4}{c_9}\big)\Big]^{c_9m}\geq  c_9\exp(-c_9c_{10}m).
	\end{eqnarray}

By \eqref{eq: Z_2m at least 4 times square of m bounded below in terms of Y_i} and \eqref{eq: sum of Y_i from 1 to squre of m is at least 4 times square of m} one obtains \eqref{eq: induction step for X_n geq square of n with even n}. 
\end{proof}

\begin{proof}[Proof of Propsition \ref{prop: stretched exponential decay for x-wusf intersect with L_n(x)}]
The upper bound has already been proved for general nonunimodular transitive graphs in Corollary \ref{cor: tail probability for x-wusf intersect with L_n with many vertices}. Next we focus on the lower bound.  
	
		Let $\eta_x:=(x,x_1,x_2,\ldots)$ be the ray starting from $x$ that represents $\xi$. Let $E_m:=\{ x\leftrightarrow x_m, x\not\leftrightarrow x_{m+1}\}$ denote the event that the path from $x$ to $x_m$ is open but the edge $(x_{m},x_{m+1})$ is not open in Bernoulli$(\frac{1}{b})$ percolation on $\mathbb{T}_{b+1}$. Let $C_x$ denote the cluster of $x$ in this Bernoulli$(\frac{1}{b})$ percolation. 
	Then 
	\[
	\mathbb{P}[E_m]=\frac{1}{b^m}\big(1-\frac{1}{b}\big).
	\]
	Conditioned on $E_m$, one has $|C_x\cap L_0(x)|=X_0'+X_1'+\cdots+X_m'$ where $X_0'\equiv1$, $X_1',\ldots, X_m'$ are independent and $X_i'$ has the law of the size of the $i$-th generation of a ``critical" Galton--Watson tree. The first generation of this Galton--Watson tree has binomial distribution $\textnormal{Bin}(b-1,\frac{1}{b})$ and the other generations have binomial progeny distribution $\textnormal{Bin}(b,\frac{1}{b})$. 
	
	Similar to the proof of Lemma \ref{lem: bounds on the n-th generation size of a critical GW tree bigger than square of n}, one can show that  there exist constants $m_0>0,c>0$ such that  
	\[
	\mathbb{P}[X_m'\geq m^2]\geq e^{-cm},\,\,\forall \, m\geq m_0.
	\]
	Now for $k\geq m_0^2$ and taking $m=\lceil \sqrt{k}\rceil$,
\[
\mathbb{P}[|C_x\cap L_0(x)|\geq k]\geq \mathbb{P}[E_m]\cdot \mathbb{P}[X_m'\geq m^2]\geq 
c_1\exp(-c_2\sqrt{k}).
\]
Taking sufficiently small $c_1,c_2$ the above inequality would also hold for $k\in[1,m_0^2]$. 
Since for the toy model the distribution of $|C_x\cap L_n(x)|$ given $|C_x\cap L_n(x)|>0$ is the same as $|C_x\cap L_0(x)|$ and $\mathbb{P}[|C_x\cap L_n(x)|>0]=e^{-t_0n}=\frac{1}{b^n}$, one has the desired lower bound for all $n\geq0,k\geq 1$:
\[
\mathbb{P}[|C_x\cap L_n(x)|\geq k]\geq 
c_1\exp(-t_0n-c_2\sqrt{k}).
\] 
By Lemma \ref{lem: law of T_x in  x-wusf for the toy model} we are done.	
\end{proof}

\begin{corollary}\label{cor: sharp high moments for x-wusf intersect L_n(x) on the toy model}
	For $x$-WUSF on the toy model $(\mathbb{T}_{b+1},\Gamma_{\xi})$, there exists a constant $C>1$ such that for all $n\geq0,k\geq 1$,
	\be\label{eq: sharp high moments for x-wusf intersect L_n(x) on the toy model}
	\frac{(k!)^2e^{-t_0n}}{C^k}\leq \mathbb{E}[|\mathfrak{T}_x\cap L_n(x)|^k]\leq C^k(k!)^2e^{-t_0n}
	\ee
\end{corollary}
\begin{proof}
	The upper bound has already been proved in \eqref{eq: second moment of x-component intersecting a high slab} for general nonunimodular transitive graphs. In fact we use these $k$th moments to prove the upper bound in Proposition \ref{prop: stretched exponential decay for x-wusf intersect with L_n(x)}.
	
	Now we look at the lower bound. 
	
	By  Proposition \ref{prop: stretched exponential decay for x-wusf intersect with L_n(x)}, there exists a constant $c>0$ such that for all $n\geq0,k\geq1$,
	\[
	\mathbb{P}[|\mathfrak{T}_x\cap L_n(x)|\geq y]\geq e^{-t_0n-c\sqrt{y}}.
	\]
   Hence 
   \begin{eqnarray}\label{eq: c_10 in appendix}
   \mathbb{E}[|\mathfrak{T}_x\cap L_n(x)|^k]&=&
   \int_{0}^{\infty}ky^{k-1}	\mathbb{P}[|\mathfrak{T}_x\cap L_n(x)|\geq y]dy\nonumber\\
   &\geq&\int_{0}^{\infty}ky^{k-1}e^{-t_0n-c\sqrt{y}}dy\nonumber\\
   &\geq&	\frac{(k!)^2e^{-t_0n}}{C^k}
   \end{eqnarray}
   for some large constant $C>0$.
\end{proof}

 \subsection{The intersection  $T_x\cap L_n(x)$ when $n\geq0$}
 
 \begin{proposition}\label{prop: sharp high moments for T_x intersect L_n(x) on the toy model}
 	For WUSF on the toy model $(\mathbb{T}_{b+1},\Gamma_{\xi})$, there exists a constant $C>1$ such that for all $n\geq0,k\geq 1$,
 	\be\label{eq: sharp high moments for T_x intersect L_n(x) on the toy model}
 	\frac{(k!)^2(n\vee 1)e^{-t_0n}}{C^k}\leq \mathbb{E}[|T_x\cap L_n(x)|^k]\leq C^k(k!)^2(n\vee 1)e^{-t_0n}
 	\ee
 and
	\be\label{eq: stretched exponential decay for T_x intersecting with L_n(x)}
	\mathbb{P}[|T_x\cap L_n(x)|\geq k]=(n\vee 1)e^{-t_0n-\Theta(\sqrt{k})}
	\ee
 \end{proposition}

\begin{lemma}\label{lem: k-th moment of Z_n plus dot to Z_n+m}
	Let $Z_n$ denote the size of $n$-th generation of a critical Galton--Watson tree with progeny distribution $\textnormal{Bin}(b,\frac{1}{b})$. Suppose $X_n$ has the same law as $Z_n$ and $X_i$'s are independent. Then there is a constant $c>0$ such that for all $k\geq 1,n\geq1,m\geq0$,
	\be\label{eq: k-th moment of independent Z_n plus dot to Z_n+m}
	\mathbb{E}\big[\big(X_n+\cdots+X_{n+m}\big)^k\big]\leq c^k(k!)(m+1)(n+m)^{k-1}.
	\ee
	Similarly there is a constant $c>0$ such that for all $k\geq 1,m\geq1$,
	\be\label{eq: k-th moment of independent Z_0 plus dot to Z_m}
	\mathbb{E}\big[\big(X_0+\cdots+X_{m}\big)^k\big]\leq c^k(k!)m^k.
	\ee
\end{lemma}
\begin{proof}
	We will only show \eqref{eq: k-th moment of independent Z_n plus dot to Z_n+m}. Since $X_0\equiv1$, the inequality \eqref{eq: k-th moment of independent Z_0 plus dot to Z_m} can be obtained from \eqref{eq: k-th moment of independent Z_n plus dot to Z_n+m} easily. 
 Observe that
	\begin{eqnarray}\label{eq: decomposition of the k-th power of the sum of X_n, X_n+m}
	\mathbb{E}\big[\big(X_n+X_{n+1}+\cdots+X_{n+m}\big)^k\big]
	&=&\sum_{\substack{k_0,\ldots,k_m\geq0\\k_0+\cdots+k_m=k} }{\frac{k!}{k_0!\cdots k_m!}}\prod_{i=0}^{m}\mathbb{E}[X_{n+i}^{k_i}]\nonumber\\
	&\stackrel{\textnormal{Lem.}\ref{lem: upper bound on the k-th moment of Z_n for a critical GW tree}}{\leq } & \sum_{\substack{k_0,\ldots,k_m\geq0\\k_0+\cdots+k_m=k} }\frac{k!}{k_0!\cdots k_m!}\prod_{i=0}^{m}\big(C^{k_i}k_i!(n+i)^{k_i-1}\vee 1\big)\nonumber\\
	&=&C^kk!\sum_{\substack{k_0,\ldots,k_m\geq0\\k_1+\cdots+k_m=k} }\prod_{i=0}^{m}\big((n+i)^{k_i-1}\vee 1\big).
	\end{eqnarray}
	\begin{claim}\label{claim: C9}
		For all $n\geq1,m\geq0,k\geq1$,
		\be\label{eq: claim 2}
		\sum_{\substack{k_0,\ldots,k_m\geq0\\k_0+\cdots+k_m=k} }\prod_{i=0}^{m}\big((n+i)^{k_i-1}\vee 1\big)\leq 3^k(m+1)(n+m)^{k-1}
		\ee
	\end{claim}	

		We prove  Claim \ref{claim: C9} by induction on $m$. When $m=0$, the left hand side of \eqref{eq: claim 2} equals $n^{k-1}\vee1=n^{k-1}$ and the right hand side equals $3^kn^{k-1}$, so \eqref{eq: claim 2} holds for all $m=0$.
		
		Now suppose $m\geq1$ and \eqref{eq: claim 2} holds for $m-1$. Using the convention that $\sum_{k_m=1}^{k-1}\cdot =0$ when $k=1$  one has the left hand side of  \eqref{eq: claim 2} satisfies
		\begin{eqnarray}
		\textnormal{LHS}&=&\sum_{k_m=0}^{k}\big((n+m)^{k_m-1}\vee 1\big)\cdot \sum_{\substack{k_0,\ldots,k_{m-1}\geq0\\k_0+\cdots+k_{m-1}=k-k_m} }\prod_{i=0}^{m-1}\big((n+i)^{k_i-1}\vee 1\big)\nonumber\\
		&=&(n+m)^{k-1}+\sum_{\substack{k_0,\ldots,k_{m-1}\geq0\\k_0+\cdots+k_{m-1}=k} }\prod_{i=0}^{m-1}\big((n+i)^{k_i-1}\vee 1\big)\nonumber\\
		&&+
		\sum_{k_m=1}^{k-1}\big((n+m)^{k_m-1}\vee 1\big)\cdot \sum_{\substack{k_0,\ldots,k_{m-1}\geq0\\k_0+\cdots+k_{m-1}=k-k_m} }\prod_{i=0}^{m-1}\big((n+i)^{k_i-1}\vee 1\big)\nonumber\\
		&\leq&(n+m)^{k-1}+3^km(n+m-1)^{k-1}+
		\sum_{k_m=1}^{k-1}(n+m)^{k_m-1}\cdot 3^{k-k_m}m(n+m-1)^{k-k_m-1}\nonumber\\
		&\leq&(n+m)^{k-1}+3^km(n+m-1)^{k-1}+
		\sum_{k_m=1}^{k-1} 3^{k-k_m}m(n+m)^{k-2}\nonumber\\
		&=&3^k(m+1)(n+m)^{k-1}\bigg[\frac{1}{3^k(m+1)}+\frac{m}{m+1}\big(1-\frac{1}{n+m}\big)^{k-1}+\frac{m}{(m+1)(n+m)}	\sum_{k_m=1}^{k-1} 3^{-k_m}\bigg]\nonumber\\
		&\leq &3^k(m+1)(n+m)^{k-1}\bigg[\frac{1}{3(m+1)}+\frac{m}{m+1}+\frac{m}{(m+1)(n+m)}\cdot \frac{1}{2}\bigg]\nonumber\\
		&<&3^k(m+1)(n+m)^{k-1}.
		\end{eqnarray}
		Hence by induction the claim holds.

Now the conclusion \eqref{eq: k-th moment of independent Z_n plus dot to Z_n+m} follows from \eqref{eq: decomposition of the k-th power of the sum of X_n, X_n+m} and \eqref{eq: claim 2}.
\end{proof}

\begin{proof}[Proof of Proposition \ref{prop: sharp high moments for T_x intersect L_n(x) on the toy model}]
Given the upper bound in \eqref{eq: sharp high moments for T_x intersect L_n(x) on the toy model}, the proof of the upper bound of \eqref{eq: stretched exponential decay for T_x intersecting with L_n(x)} is the same as Corollary \ref{cor: tail probability for x-wusf intersect with L_n with many vertices} and thus we omit the details.

Given the lower bound of \eqref{eq: stretched exponential decay for T_x intersecting with L_n(x)}, the lower bound of \eqref{eq: sharp high moments for T_x intersect L_n(x) on the toy model} can be obtained in the same way as \eqref{eq: c_10 in appendix}.

So it suffices to show the upper bound of \eqref{eq: sharp high moments for T_x intersect L_n(x) on the toy model} and the lower bound of \eqref{eq: stretched exponential decay for T_x intersecting with L_n(x)}.

	Let $\eta_x=(x_0,x_1,x_2,\ldots)$ be the ray starting from $x$ representing the end $\xi$. Recall the events $A_k$ and $B_{k,k'}$ we defined in the proof of Proposition \ref{prop: exact two point function for WUSF on regular trees}; see Figure \ref{fig: B_k,k'} for an example of $B_{k,k'}$.
	
	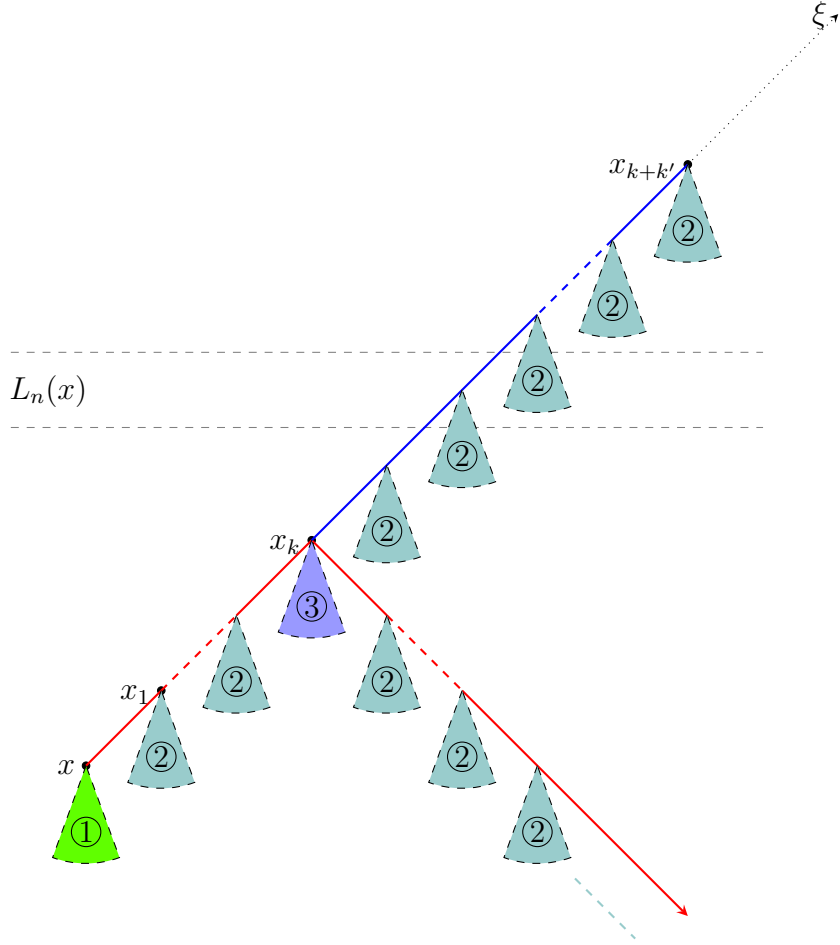
\begin{figure}[h!]
		\centering
		\begin{tikzpicture}[scale=1, text height=1.5ex,text depth=.25ex] 
		\draw [dashed, color=gray](1,7.5)--(11,7.5);
		\draw [dashed, color=gray](1,6.5)--(11,6.5);
		\node[below] at (1.5,7.3) {$L_n(x)$};
		
		\draw[fill=black] (2,2) circle [radius=0.05];
		\node[left] at (2,2) {$x$};
		
		\draw[fill=black] (3,3) circle [radius=0.05];
		\node[left] at (3,3) {$x_1$};

		\draw[fill=black] (5,5) circle [radius=0.05];
		\node[left] at (5,5) {$x_k$};
		
		\draw[fill=black] (10,10) circle [radius=0.05];
		\node[left] at (10,10) {$x_{k+k'}$};
		
		\draw [thick,color=red](2,2)--(3,3);
		\draw [thick,dashed,color=red](3,3)--(4,4);
		\draw [thick,color=red](4,4)--(5,5);
		\draw [thick,color=red](6,4)--(5,5);
		\draw [thick,dashed,color=red](6,4)--(7,3);
		\draw [->,  thick,>=stealth,color=red] (7,3)--(10,0);
		
		\draw [ thick, color=blue] (5,5)--(8,8);
		\draw [ thick, color=blue] (9,9)--(10,10);
		\draw [ thick,dashed,color=blue] (9,9)--(8,8);
		
		\draw [->, >=stealth,dotted,  color=black] (10,10)--(12,12);
		\node[left] at (12,12) {$\xi$};

		\begin{scope}[shift={(2,2)}]
		\draw[dashed,fill=green!50!lime] (0,0) -- (250:1.3) arc (250:290:1.3)--(0,0);
		\end{scope}
		\node[below] at (2,1.4) {$\circled{1}$};

		\begin{scope}[shift={(3,3)}]
		\draw[dashed,fill=teal!40!white] (0,0) -- (250:1.3) arc (250:290:1.3)--(0,0);
		\end{scope}
		\node[below] at (3,2.4) {$\circled{2}$};

		\begin{scope}[shift={(5,5)}]
		\draw[dashed,fill=blue!40!white] (0,0) -- (250:1.3) arc (250:290:1.3)--(0,0);
		\end{scope}
		\node[below] at (5,4.4) {$\circled{3}$};

		\begin{scope}[shift={(4,4)}]
		\draw[dashed,fill=teal!40!white] (0,0) -- (250:1.3) arc (250:290:1.3)--(0,0);
		\end{scope}
		\node[below] at (4,3.4) {$\circled{2}$};

		\begin{scope}[shift={(6,6)}]
		\draw[dashed,fill=teal!40!white] (0,0) -- (250:1.3) arc (250:290:1.3)--(0,0);
		\end{scope}
		\node[below] at (6,5.4) {$\circled{2}$};
		
		\begin{scope}[shift={(10,10)}]
		\draw[dashed,fill=teal!40!white] (0,0) -- (250:1.3) arc (250:290:1.3)--(0,0);
		\end{scope}
		\node[below] at (10,9.4) {$\circled{2}$};
		\begin{scope}[shift={(9,9)}]
		\draw[dashed,fill=teal!40!white] (0,0) -- (250:1.3) arc (250:290:1.3)--(0,0);
		\end{scope}
		\node[below] at (9,8.4) {$\circled{2}$};

		\begin{scope}[shift={(8,8)}]
		\draw[dashed,fill=teal!40!white] (0,0) -- (250:1.3) arc (250:290:1.3)--(0,0);
		\end{scope}
		\node[below] at (8,7.4) {$\circled{2}$};
		
		\begin{scope}[shift={(7,7)}]
		\draw[dashed,fill=teal!40!white] (0,0) -- (250:1.3) arc (250:290:1.3)--(0,0);
		\end{scope}
		\node[below] at (7,6.4) {$\circled{2}$};

		\begin{scope}[shift={(6,4)}]
		\draw[dashed,fill=teal!40!white] (0,0) -- (250:1.3) arc (250:290:1.3)--(0,0);
		\end{scope}
		\node[below] at (6,3.4) {$\circled{2}$};
		
		\begin{scope}[shift={(7,3)}]
		\draw[dashed,fill=teal!40!white] (0,0) -- (250:1.3) arc (250:290:1.3)--(0,0);
		\end{scope}
		\node[below] at (7,2.4) {$\circled{2}$};
		
		\begin{scope}[shift={(8,2)}]
		\draw[dashed,fill=teal!40!white] (0,0) -- (250:1.3) arc (250:290:1.3)--(0,0);
		\end{scope}
		\node[below] at (8,1.4) {$\circled{2}$};
		
		\draw [color=teal!40!white,dashed,thick] (8.5,0.5)--(9.3,-0.3);
		
		\end{tikzpicture}
		\caption{A sketch of the tree $T_x$ conditioned on a typical $B_{k,k'}$. The red line denotes the future of $x$. There are three types of independent ``critical" Galton--Watson trees attached to the future of $x$ and vertices $x_{k+1},\ldots,x_{k+k'}$, and the only difference among them is the progeny distribution of the first generation. The distribution of the first generation distribution for the three types are $\textnormal{Bin}(b,\frac{1}{b})$, $\textnormal{Bin}(b-1,\frac{1}{b})$ and $\textnormal{Bin}(b-2,\frac{1}{b})$ respectively.}
		\label{fig: B_k,k'}
	\end{figure}
	
	 Also recall that 
	\[
	\mathbb{P}[A_k]= \left\{
	\begin{array}{cc}
	\frac{b}{b+1}, & k=0\\
	& \\
	\frac{b-1}{b+1}\cdot\frac{1}{b^k}, & k>0.
	\end{array}
	\right.
	\]
and 
\[
\mathbb{P}[B_{k,k'}]=\mathbb{P}[A_k]\cdot \frac{b-1}{b}\cdot \frac{1}{b^{k'}}\asymp \frac{1}{b^{k+k'}},\,\,\forall \, k,k'\geq0.
\]

\begin{observation}\label{obs: T_x under B_k,k'}
	Suppose  $X_n,X_n',X_n''$ have the law of the size of the $n$-th generation of ``critical"  Galton--Watson trees with first generation progeny distributions  $\textnormal{Bin}(b,\frac{1}{b})$, $\textnormal{Bin}(b-1,\frac{1}{b})$ and $\textnormal{Bin}(b-2,\frac{1}{b})$ respectively (for generation at least $2$, the progeny distribution is always $\textnormal{Bin}(b,\frac{1}{b})$.)  In particular, when $b=2$, $X_0''=1$ and $X_n''=0$ for $n\geq1$. Moreover assume all the random variables $\{X_i,X_i',X_i''\colon i=0,1,2,\ldots \}$ are independent. 
	
	It is easy to see from Figure \ref{fig: B_k,k'} that conditioned on the event $B_{k,k'}$, one has 
	\begin{enumerate}
		\item if $n>k+k'$, then $|T_x\cap L_n(x)|=0$;
		\item if $0\leq k<n\leq k+k'$, then $|T_x\cap L_n(x)|$ has the same distribution as $\sum_{j=n}^{k+k'}X_{j-n}'$;
		\item if $0<n\leq k\leq k+k'$, then $|T_x\cap L_n(x)|$ has the same distribution as $\sum_{j=n}^{k-1}X_{j-n}'+X_{k-n}''+\sum_{j=k+1}^{k+k'}X_{j-n}'+\sum_{j=n}^{k-1}\widetilde{X}_{j-n}'$, where all these random variables  are independent and $\widetilde{X}_j'$ has the same distribution as $X_j'$;
		\item if $n\leq0$, then $|T_x\cap L_n(x)|$ has the same distribution as $X_{-n}+\sum_{j=1}^{k-1}X_{j-n}'+\sum_{j=k+1}^{k+k'}X_{j-n}'+X_{k-n}''+\sum_{j=n}^{k-1}\widetilde{X}_{j-n}'$, where all these random variables  are independent and $\widetilde{X}_j'$ has the same distribution as $X_j'$.
	\end{enumerate}
\end{observation}

Note that  the random variables $X_j',X_j''$ and $\widetilde{X}_j'$ can be stochastically dominated by $X_j$, which has the same law as $Z_j$ in Lemma \ref{lem: upper bound on the k-th moment of Z_n for a critical GW tree}. Also note that $B_{k,k'}$ are disjoint and $\mathbb{P}\big[\bigcup_{k,k'\geq0} B_{k,k'}\big]=1$. 

Therefore for $n>0,t\geq2$, one has that 
\begin{eqnarray}\label{eq: c18 for high moments of T_x cap L_n}
\mathbb{E}[|T_x\cap L_n(x)|^t]&=&
\sum_{k=0}^{n-1}\sum_{k'=n-k}^{\infty}\mathbb{P}[B_{k,k'}]\cdot \mathbb{E}\bigg[\Big(\sum_{j=n}^{k+k'}X_{j-n}'\Big)^t\bigg] \nonumber\\
&&+\sum_{k=n}^{\infty}\sum_{k'=0}^{\infty}\mathbb{P}[B_{k,k'}]\cdot \mathbb{E}\bigg[\Big(\sum_{j=n}^{k-1}X_{j-n}'+\sum_{j=k+1}^{k+k'}X_{j-n}'+X_{k-n}''+\sum_{j=n}^{k-1}\widetilde{X}_{j-n}'\Big)^t\bigg] \nonumber\\
&\leq &C^t\sum_{k=0}^{n-1}\sum_{k'=n-k}^{\infty}\frac{1}{b^{k+k'}}\cdot \mathbb{E}\bigg[\Big(\sum_{j=n}^{k+k'}X_{j-n}\Big)^t\bigg]\nonumber\\
&& +C^t\sum_{k=n}^{\infty}\sum_{k'=0}^{\infty}\frac{1}{b^{k+k'}}\cdot \mathbb{E}\bigg[\Big(\sum_{j=n}^{k+k'}X_{j-n}\Big)^t\bigg] \nonumber\\
&\stackrel{\eqref{eq: k-th moment of independent Z_0 plus dot to Z_m}}{\leq}& 
C^t\sum_{k=0}^{n-1}\sum_{k'=n-k}^{\infty}\frac{1}{b^{k+k'}}\cdot c^t(t!)[(k+k'-n)\vee 1]^{t}\nonumber\\
&&+C^t\sum_{k=n}^{\infty}\sum_{k'=0}^{\infty}\frac{1}{b^{k+k'}}\cdot c^t(t!)[(k+k'-n)\vee 1]^{t}\nonumber\\&\leq& c_2^t(t!)\sum_{j=0}^{\infty}(n+1)\cdot \frac{1}{b^{n+j}}\cdot (j\vee 1)^t\nonumber\\
&\leq & c_3^t(t!)^2\frac{n}{b^n}.
\end{eqnarray}
This proves the upper bound of \eqref{eq: sharp high moments for T_x intersect L_n(x) on the toy model} for $n\geq1$. The case of $n=0$ can be dealt in a similar way and we omit the details. 

Next we turn to the proof of the lower bound of \eqref{eq: stretched exponential decay for T_x intersecting with L_n(x)}. 

As mentioned earlier using a similar proof of Lemma \ref{lem: bounds on the n-th generation size of a critical GW tree bigger than square of n},  there exist constants $m_0>0,c>0$ such that  
\[
\mathbb{P}[X_m'\geq m^2]\geq e^{-cm},\,\,\forall \, m\geq m_0.
\]
Notice for different pairs of $(k,k')$, the event $B_{k,k'}$ are disjoint. Conditioned on $B_{k,k'}$ with $k+k'>n$, one has that $|T_x\cap L_n(x)|$ stochastically dominates $X_{k+k'-n}'$. Hence for $t\geq m_0^2$, one has the desired lower bound of \eqref{eq: stretched exponential decay for T_x intersecting with L_n(x)}:
\begin{eqnarray}
\mathbb{E}[|T_x\cap L_n(x)|\geq t]
&\geq&\sum_{k+k'=n+\lceil\sqrt{t}\rceil}\mathbb{P}[B_{k,k'}]\cdot \mathbb{P}\big[X_{\lceil\sqrt{t}\rceil}'\geq \lceil\sqrt{t}\rceil^2\big] \nonumber\\
&\succeq&\sum_{k+k'=n+\lceil\sqrt{t}\rceil}\frac{1}{b^{n+\lceil\sqrt{t}\rceil}}\cdot e^{-c\lceil\sqrt{t}\rceil}\nonumber\\
&\geq&\frac{(n\vee 1)}{b^{n}}\cdot e^{-c_1\sqrt{t}}
\end{eqnarray} 
For $t\in[1,m_0^2]$, taking a small constant $c_1$ would be enough. 
\end{proof}

\begin{remark}
	One can use \eqref{eq: k-th moment of independent Z_n plus dot to Z_n+m} and a similar calculation to \eqref{eq: c18 for high moments of T_x cap L_n} in the above proof to give another proof of Proposition \ref{prop: high moments for x-component intersecting a slab} for the toy model. Just use the events $E_m$ in the proof of Proposition \ref{prop: stretched exponential decay for x-wusf intersect with L_n(x)} instead of $B_{k,k'}$. 
\end{remark}

\subsection{Intersections with $L_{-n}(x)$ when $n>0$}

\begin{proposition}\label{prop: sharp high moments for T_x intersect L_-n(x) on the toy model}
		For WUSF on the toy model $(\mathbb{T}_{b+1},\Gamma_{\xi})$, there exist constants $c_1,c_2,c_3$ such that for all $n>0,k\geq 1$,
	\be\label{eq: sharp high moments for T_x intersect L_-n(x) on the toy model}
\mathbb{E}[|T_x\cap L_{-n}(x)|^k]\leq c_1^k(k!)^2n^k
	\ee
	and
	\be\label{eq: stretched exponential decay for T_x intersecting with L_-n(x)}
	\mathbb{P}[|T_x\cap L_{-n}(x)|\geq k]\leq c_2e^{-c_3\sqrt{\frac{k}{n}}}
	\ee
\end{proposition}
\begin{remark}
	Comparing the upper bound of  \eqref{eq: sharp high moments for T_x intersect L_n(x) on the toy model}  for the toy model with \eqref{eq: high moment for T_x intersecting a high slab } in  Proposition \ref{prop: high moments for T_x intersecting a slab}, the power of $k!$ in \eqref{eq: high moment for T_x intersecting a high slab } for general nonunimodular transitive graphs  in  Proposition \ref{prop: high moments for T_x intersecting a slab} might be not optimal.  Similarly the power of $k!$  on the right hand of  \eqref{eq: sharp high moments for T_x intersect L_-n(x) on the toy model} suggests that the power of $k!$ in \eqref{eq: high moment for T_x intersecting a low slab } might be not optimal. Besides, the power of $n$ in   \eqref{eq: sharp high moments for T_x intersect L_-n(x) on the toy model} suggests that the power of $n$ in \eqref{eq: high moment for T_x intersecting a low slab }  might also be not optimal.
\end{remark}
\begin{remark}
	The corresponding upper bounds for $\mathfrak{T}_x\cap L_{-n}(x)$ have been proved for general nonunimodular graphs; see  \eqref{eq: second moment of x-component intersecting a low slab} in Proposition \ref{prop: high moments for x-component intersecting a slab} and \eqref{eq: stretched exponential tail for x-wusf intersect with L_(-n) with many vertices} in  Corollary \ref{cor: tail probability for x-wusf intersect with L_n with many vertices}. 
\end{remark}

\begin{remark}
For the toy model, one can get a lower bound with the form of the upper bound in \eqref{eq: stretched exponential tail for x-wusf intersect with L_(-n) with many vertices}  when $1\leq k\leq c_4 n$.  To see this, note that for $1\leq k\leq c_4 n$, $\mathbb{P}[|\mathfrak{T}_x\cap L_{-n}(x)|\geq k]\geq \mathbb{P}[X_n\geq c_4n]\geq \frac{c_5}{n}\geq \frac{c_5}{n} e^{-c_6\sqrt{\frac{k}{n}}}$ by Theorem 1 on page 19 of \cite{AthreyaNey1972}, where $c_5$ depends on $c_4$.

	Similarly for the toy model,  when $1\leq k\leq c_4 n$ we also have a lower bound with the form of  the upper bound in \eqref{eq: stretched exponential decay for T_x intersecting with L_-n(x)} (in this case it is just $\mathbb{P}[|T_x\cap L_{-n}(x)|\geq k]\geq c_5$) and we omit the proof.  
 \end{remark}
\begin{question}
	Are there lower bounds for $\mathbb{P}[|\mathfrak{T}_x\cap L_{-n}(x)|\geq k]$  and $\mathbb{P}[|T_x\cap L_{-n}(x)|\geq k]$  with the same form as the upper bound in \eqref{eq: stretched exponential tail for x-wusf intersect with L_(-n) with many vertices} and \eqref{eq: stretched exponential decay for T_x intersecting with L_-n(x)} respectively for all $n,k\geq 1$? 
\end{question}

\begin{proof}[Proof of Proposition \ref{prop: sharp high moments for T_x intersect L_-n(x) on the toy model}]
	The proof of  \eqref{eq: sharp high moments for T_x intersect L_-n(x) on the toy model} is quite similar to \eqref{eq: c18 for high moments of T_x cap L_n}. By Observation \ref{obs: T_x under B_k,k'}, for $t\geq2$ one has 
	\begin{eqnarray}\label{eq: c22 for high moments of T_x cap L_(-n)}
	\mathbb{E}[|T_x\cap L_{-n}(x)|^t]
	&=&\sum_{k,k'=0}^{\infty} \mathbb{P}[B_{k,k'}]\mathbb{E}\bigg[ \big(X_{n}+\sum_{j=1}^{k-1}X_{j+n}'+\sum_{j=k+1}^{k+k'}X_{j+n}'+X_{k+n}''+\sum_{j=-n}^{k-1}\widetilde{X}_{j+n}'\big)^t\bigg]\nonumber\\
	&\leq&\sum_{k,k'=0}^{\infty}\frac{c_1}{b^{k+k'}}c_2^t\mathbb{E}\big[\big(\sum_{j=0}^{k+k'+n}X_j\big)^t\big]\nonumber\\
	&\stackrel{\eqref{eq: k-th moment of independent Z_0 plus dot to Z_m}}{\leq}&
	\sum_{k,k'=0}^{\infty}\frac{c_3^t}{b^{k+k'}}t!(k+k'+n)^t=\sum_{j=0}^{\infty}c_3^tt! \frac{j+1}{b^j}(j+n)^t\nonumber\\
	&\leq&c_3^tt!n^t+\sum_{j=1}^{\infty}c_3^tt!\frac{2j}{b^j}(j+n)^t\leq
	c_3^t t!n^t+ 2bc_3^tt!\int_{0}^{\infty}\frac{x(x+n)^t}{b^x}dx\nonumber\\
	&=&	c_3^t t!n^t+ 2bc_3^tt!\sum_{j=0}^{t}\binom{t}{j}\int_{0}^{\infty}x^{j+1}n^{t-j}b^{-x}dx\nonumber\\
	&=&c_3^t t!n^t+ 2bc_3^tt!\sum_{j=0}^{t}\binom{t}{j} n^{t-j}\frac{(j+1)!}{(\log b)^{j+2}}\nonumber\\
	&\leq&c_3^t t!n^t+ 2bc_4^tt!\sum_{j=0}^{t}t!(j+1) n^{t}\leq c_5^t(t!)^2n^t
	\end{eqnarray}
	
	Similar to the proof of Corollary \ref{cor: tail probability for x-wusf intersect with L_n with many vertices}, one has for $c_6<\frac{1}{2\sqrt{c_5}}$, 
	\begin{eqnarray}
	\mathbb{E}\bigg[\exp\bigg(c_6\sqrt{\frac{|T_x\cap L_{-n}(x)|}{n}} \bigg)\bigg]
	&\leq&2\sum_{k=0}^{\infty}\frac{c_6^{2k}}{(2k)!}\frac{1}{n^k}\mathbb{E}[|T_x\cap L_{-n}(x)|^k]\nonumber\\
	&\stackrel{\eqref{eq: c22 for high moments of T_x cap L_(-n)}}{\leq}&
	2\sum_{k=0}^{\infty}\frac{c_6^{2k}}{(2k)!} c_5^k(k!)^2<3
	\end{eqnarray}
	Hence by Markov's inequality one obtains \eqref{eq: stretched exponential decay for T_x intersecting with L_-n(x)}.
	
\end{proof}

For the toy model, the results in Proposition \ref{prop: stretched exponential decay for x-wusf intersect with L_n(x)} and Corollary \ref{cor: sharp high moments for x-wusf intersect L_n(x) on the toy model} also hold for the past $\mathfrak{P}(x)$ and we omit the details. For the future, it is not interesting since $|\mathfrak{F}(x,\infty)\cap L_n(x)|\in\{0,1,2\}$ almost surely.

\bibliography{revision}
\bibliographystyle{plain}

\end{document}